\documentclass[a4paper,12pt]{article}
\usepackage[margin=1.1in]{geometry}  % set the margins to 1in on all sides
\usepackage{graphicx}  
\usepackage{longtable}% to include figures
\usepackage[all,pdf]{xy}
\usepackage{amsmath}
\usepackage{amscd}               % great math stuff
\usepackage{amsfonts} 
\usepackage{amssymb}             % for blackboard bold, etc
\usepackage{amsthm}                % better theorem environment
\usepackage{mathrsfs}
\usepackage{amsbsy}
\usepackage{tikz}
\usepackage{booktabs}
\usetikzlibrary{decorations.pathmorphing}
\usepackage{hyperref}
\usepackage{siunitx}
\usepackage{float}

\hypersetup{pdfstartview=}

\numberwithin{equation}{section}

\newtheorem{thm}{Theorem}[section]
\newtheorem{defn}[thm]{Definition}
\newtheorem{lem}[thm]{Lemma}

\newtheorem{prop}[thm]{Proposition}
\newtheorem{cor}[thm]{Corollary}

\newtheorem{assumption}[thm]{Assumption}
\DeclareMathOperator{\id}{id}

\def\XXint#1#2#3{{\setbox0=\hbox{$#1{#2#3}{\int}$ }
\vcenter{\hbox{$#2#3$ }}\kern-.6\wd0}}

\makeatletter
\def\th@newremark{\th@remark\thm@headfont{\bfseries}}
\makeatletter

\theoremstyle{newremark}
\newtheorem{rmk}[thm]{Remark}

\newcommand{\CC}{\mathbb{C}}

\newcommand{\bB}{\mathcal{B}}
\newcommand{\cC}{\mathcal{C}}

\newcommand{\eE}{\mathcal{E}}
\newcommand{\fF}{\mathcal{F}}
\newcommand{\gG}{\mathcal{G}}
\newcommand{\hH}{\mathcal{H}}
\newcommand{\iI}{\mathcal{I}}

\newcommand{\kK}{\mathcal{K}}
\newcommand{\lL}{\mathcal{L}}
\newcommand{\mM}{\mathcal{M}}

\newcommand{\pP}{\mathcal{P}}

\newcommand{\rR}{\mathcal{R}}
\newcommand{\sS}{\mathcal{S}}
\newcommand{\tT}{\mathcal{T}}

\newcommand{\xX}{\mathcal{X}}
\newcommand{\yY}{\mathcal{Y}}
\newcommand{\zZ}{\mathcal{Z}}

\renewcommand{\b}{\mathbf{b}}
\newcommand{\e}{\mathbf{e}}
\newcommand{\p}{\mathbf{p}}
\newcommand{\s}{\textbf{s}}
\renewcommand{\v}{\mathbf{v}}
\newcommand{\w}{\mathbf{w}}

\newcommand{\E}{\mathbf{E}}
\newcommand{\I}{\mathbf{I}}

\newcommand{\N}{\mathbf{N}}
\renewcommand{\P}{\mathbf{P}}
\newcommand{\R}{\mathbf{R}}
\newcommand{\T}{\mathbf{T}}
\newcommand{\Z}{\mathbf{Z}}

\renewcommand{\id}{\mathrm{id}}
\newcommand{\eps}{\varepsilon}
\renewcommand{\d}{\partial}
\renewcommand{\div}{\mathrm{div }}
\newcommand{\supp}{\mathrm{supp }}

\newcommand{\Com}{{\operator@font Com}}
\newcommand{\Comad}{{\operator@font Comad}}
\newcommand{\Comi}{{\operator@font Comi}}
\newcommand{\diam}{{\operator@font diam}}
\newcommand{\Err}{{\operator@font Err}}

\newcommand{\Tr}{\text{Tr }}     %trace
\newcommand{\cov}{{\operator@font cov}}
\newcommand{\var}{{\operator@font var}}
\newcommand{\corr}{{\operator@font corr}}
\newcommand{\flux}{\mathsf{flux}}
\newcommand{\fs}{\mathfrak{s}}
\newcommand{\fC}{\mathfrak{C}}

\newcommand{\sT}{\mathsf{T}}

\renewcommand{\1}{\textbf{1}}

\newcommand{\Av}{{\operator@font Av}}
\newcommand{\trig}{{\operator@font trig}}
\newcommand{\PL}{{\operator@font PL}}

\definecolor{darkgreen}{rgb}{0.1,0.7,0.1}
\definecolor{darkred}{rgb}{0.7,0.1,0.1}
\definecolor{darkblue}{rgb}{0,0,0.7}
\def\scal#1{\langle#1\rangle}

\colorlet{symbols}{blue!90!black}
\colorlet{testcolor}{green!60!black}

\def\symbol#1{\textcolor{symbols}{#1}}
\def\1{\mathbf{\symbol{1}}}

\def\${|\!|\!|}

\definecolor{darkgreen}{rgb}{0.1,0.6,0.1}
\definecolor{darkblue}{rgb}{0.1,0,0.7}
\makeatletter
\pgfdeclareshape{crosscircle}
{
  \inheritsavedanchors[from=circle] % this is nearly a circle
  \inheritanchorborder[from=circle]
  \inheritanchor[from=circle]{north}
  \inheritanchor[from=circle]{north west}
  \inheritanchor[from=circle]{north east}
  \inheritanchor[from=circle]{center}
  \inheritanchor[from=circle]{west}
  \inheritanchor[from=circle]{east}
  \inheritanchor[from=circle]{mid}
  \inheritanchor[from=circle]{mid west}
  \inheritanchor[from=circle]{mid east}
  \inheritanchor[from=circle]{base}
  \inheritanchor[from=circle]{base west}
  \inheritanchor[from=circle]{base east}
  \inheritanchor[from=circle]{south}
  \inheritanchor[from=circle]{south west}
  \inheritanchor[from=circle]{south east}
  \inheritbackgroundpath[from=circle]
  \foregroundpath{
    \centerpoint%
    \pgf@xc=\pgf@x%
    \pgf@yc=\pgf@y%
    \pgfutil@tempdima=\radius%
    \pgfmathsetlength{\pgf@xb}{\pgfkeysvalueof{/pgf/outer xsep}}%  
    \pgfmathsetlength{\pgf@yb}{\pgfkeysvalueof{/pgf/outer ysep}}%  
    \ifdim\pgf@xb<\pgf@yb%
      \advance\pgfutil@tempdima by-\pgf@yb%
    \else%
      \advance\pgfutil@tempdima by-\pgf@xb%
    \fi%
    \pgfpathmoveto{\pgfpointadd{\pgfqpoint{\pgf@xc}{\pgf@yc}}{\pgfqpoint{-0.707107\pgfutil@tempdima}{-0.707107\pgfutil@tempdima}}}
    \pgfpathlineto{\pgfpointadd{\pgfqpoint{\pgf@xc}{\pgf@yc}}{\pgfqpoint{0.707107\pgfutil@tempdima}{0.707107\pgfutil@tempdima}}}
    \pgfpathmoveto{\pgfpointadd{\pgfqpoint{\pgf@xc}{\pgf@yc}}{\pgfqpoint{-0.707107\pgfutil@tempdima}{0.707107\pgfutil@tempdima}}}
    \pgfpathlineto{\pgfpointadd{\pgfqpoint{\pgf@xc}{\pgf@yc}}{\pgfqpoint{0.707107\pgfutil@tempdima}{-0.707107\pgfutil@tempdima}}}
  }
}
\makeatother

\colorlet{symbols}{blue}     %colorlet: define a new color based on an old one.
\colorlet{testcolor}{green!60!black}

\def\symbol#1{\textcolor{symbols}{#1}}
\def\1{\mathbf{\symbol{1}}}

\usetikzlibrary{calc}
\usetikzlibrary{shapes.misc}
\usetikzlibrary{shapes.symbols}
\usetikzlibrary{snakes}
\usetikzlibrary{decorations}
\usetikzlibrary{decorations.markings}

\tikzset{
	root/.style={circle,fill=testcolor,inner sep=0pt, minimum size=2mm},
	dot/.style={circle,fill=black,inner sep=0pt, minimum size=1mm},
	var/.style={circle,fill=black!10,draw=black,inner sep=0pt, minimum size=2mm},
	dotred/.style={circle,fill=black!50,inner sep=0pt, minimum size=2mm},
	generic/.style={semithick,shorten >=1pt,shorten <=1pt},
	dist/.style={ultra thick,testcolor,shorten >=1pt,shorten <=1pt},
	testfcn/.style={ultra thick,testcolor,shorten >=1pt,shorten <=1pt,<-},
	testfcnx/.style={ultra thick,testcolor,shorten >=1pt,shorten <=1pt,<-,
		postaction={decorate,decoration={markings,mark=at position 0.6 with {\drawx}}}},
	kernelprime/.style={semithick,shorten >=1pt,shorten <=1pt,densely dashed,->},
	kernelprimex/.style={semithick,shorten >=1pt,shorten <=1pt,densely dashed,->,
		postaction={decorate,decoration={markings,mark=at position 0.4 with {\drawx}}}},
	kernel/.style={semithick,shorten >=1pt,shorten <=1pt,->},
	kernel1/.style={->,semithick,shorten >=1pt,shorten <=1pt,postaction={decorate,decoration={markings,mark=at position 0.45 with {\draw[-] (0,-0.1) -- (0,0.1);}}}},
    akernel/.style={semithick,->},
	keps/.style={semithick,densely dashed,shorten >=1pt,shorten <=1pt,->},
	dots/.style={semithick,dotted,shorten >=1pt,shorten <=1pt},
	Deps/.style={semithick,draw=black!25,fill=black!25,shorten >=1pt,shorten <=1pt,->},
	kbase/.style={semithick,dotted,shorten >=1pt,shorten <=1pt,->},
	multx/.style={shorten >=1pt,shorten <=1pt,
		postaction={decorate,decoration={markings,mark=at position 0.5 with {\drawx}}}},
	kernelx/.style={semithick,shorten >=1pt,shorten <=1pt,->,
		postaction={decorate,decoration={markings,mark=at position 0.4 with {\drawx}}}},
	kernel1/.style={->,semithick,shorten >=1pt,shorten <=1pt,postaction={decorate,decoration={markings,mark=at position 0.45 with {\draw[-] (0,-0.1) -- (0,0.1);}}}},
	kernel2/.style={->,semithick,shorten >=1pt,shorten <=1pt,postaction={decorate,decoration={markings,mark=at position 0.45 with {\draw[-] (0.05,-0.1) -- (0.05,0.1);\draw[-] (-0.05,-0.1) -- (-0.05,0.1);}}}},
	kernelBig/.style={semithick,shorten >=1pt,shorten <=1pt,decorate, decoration={zigzag,amplitude=1.5pt,segment length = 3pt,pre length=2pt,post length=2pt}},
	rho/.style={dotted,semithick,shorten >=1pt,shorten <=1pt},
	renorm/.style={shape=circle,fill=white,inner sep=1pt},
	labl/.style={shape=rectangle,fill=white,inner sep=1pt},
	xi/.style={circle,fill=symbols!10,draw=symbols,inner sep=0pt,minimum size=1.0mm},
	xix/.style={crosscircle,fill=symbols!10,draw=symbols,inner sep=0pt,minimum size=1.2mm},
	xib/.style={circle,fill=symbols!10,draw=symbols,inner sep=0pt,minimum size=1.6mm},
	xibx/.style={crosscircle,fill=symbols!10,draw=symbols,inner sep=0pt,minimum size=1.6mm},
	not/.style={circle,fill=symbols,draw=symbols,inner sep=0pt,minimum size=0.5mm},
	>=stealth,
	graydot/.style={circle,fill=gray,inner sep=0pt, minimum size=1mm},
	zero/.style={circle,inner sep=0pt, minimum size=1mm, draw},
	kernelprimeeps/.style={densely dashed, semithick,shorten >=1pt,shorten <=1pt},
	smalldot/.style={circle,fill=black,draw=black, solid,inner sep=0pt,minimum size=0.5mm},
    tinydot/.style={circle,fill=black,draw=black, solid,inner sep=0pt,minimum size=0.2mm},
	}
 \makeatletter
 \def\DeclareSymbol#1#2#3{\expandafter\gdef\csname MH@symb@#1\endcsname{\tikz[baseline=#2,scale=0.15,draw=black]{#3}}
 \expandafter\gdef\csname MH@symb@#1s\endcsname{\scalebox{0.7}{\tikz[baseline=#2,scale=0.15,draw=black]{#3}}}}
 \def\<#1>{\csname MH@symb@#1\endcsname}
 \makeatother
\DeclareSymbol{1}{0}{\draw[white] (-.4,0) -- (.4,0); \draw (0,0)  -- (0,1.2) node[dot] {};}
\DeclareSymbol{2}{0}{\draw (-0.5,1.2) node[dot] {} -- (0,0) -- (0.5,1.2) node[dot] {};}
\DeclareSymbol{3}{0}{\draw (0,0) -- (0,1.2) node[smalldot] {}; \draw (-.7,1) node[smalldot] {} -- (0,0) -- (.7,1) node[smalldot] {};}
\DeclareSymbol{30}{-3}{\draw (0,0) -- (0,-1); \draw (0,0) -- (0,1.2) node[smalldot] {}; \draw (-.7,1) node[smalldot] {} -- (0,0) -- (.7,1) node[smalldot] {};}
\DeclareSymbol{31}{-3}{\draw (0,0) -- (0,-1) -- (1,0) node[smalldot] {}; \draw (0,0) -- (0,1.2) node[smalldot] {}; \draw (-.7,1) node[smalldot] {} -- (0,0) -- (.7,1) node[smalldot] {};}

\DeclareSymbol{32}{-3}{\draw (0,0) -- (0,-1) -- (1,0) node[smalldot] {}; \draw (0,0) -- (0,-1) -- (-1,0) node[smalldot] {}; \draw (0,0) -- (0,1.2) node[smalldot] {}; \draw (-.7,1) node[smalldot] {} -- (0,0) -- (.7,1) node[smalldot] {};}
\DeclareSymbol{20}{-3}{\draw (0,0) -- (0,-1);\draw (-.7,1) node[smalldot] {} -- (0,0) -- (.7,1) node[smalldot] {};}
\DeclareSymbol{22}{-3}{\draw (0,0.3) -- (0,-1) -- (1,0) node[smalldot] {}; \draw (0,0.3) -- (0,-1) -- (-1,0) node[smalldot] {};\draw (-.7,1) node[smalldot] {} -- (0,0.3) -- (.7,1) node[smalldot] {};}
\DeclareSymbol{31p}{-3}{\draw (0,0) -- (0,-1) -- (1,0) node[smalldot] {}; \draw (0,0) -- (0,1.2) node[smalldot] {}; \draw (-.7,1) node[smalldot] {} -- (0,0) -- (.7,1) node[smalldot] {};\fill[white] (0,-1) circle (8pt);\draw (0,-1) circle (8pt);}% \draw (0,-1) node{\scaleobj{0.5}{\pe}}; 
\DeclareSymbol{32p}{-3}{\draw (0,0) -- (0,-1) -- (1,0) node[smalldot] {}; \draw (0,0) -- (0,-1) -- (-1,0) node[smalldot] {}; \draw (0,0) -- (0,1.2) node[smalldot] {}; \draw (-.7,1) node[smalldot] {} -- (0,0) -- (.7,1) node[smalldot] {}; \fill[white] (0,-1) circle (8pt);\draw (0,-1) circle (8pt);}
\DeclareSymbol{22p}{-3}{\draw (0,0.3) -- (0,-1) -- (1,0) node[smalldot] {}; \draw (0,0.3) -- (0,-1) -- (-1,0) node[smalldot] {};\draw (-.7,1) node[smalldot] {} -- (0,0.3) -- (.7,1) node[smalldot] {}; \fill[white] (0,-1) circle (8pt); \draw (0,-1) circle (8pt);}
\DeclareSymbol{1}{0}{\draw[white] (-.4,0) -- (.4,0); \draw (0,0)  -- (0,1.2) node[smalldot] {};}

\DeclareSymbol{2}{0}{\draw (-0.5,1.2) node[smalldot] {} -- (0,0) -- (0.5,1.2) node[smalldot] {};}

\DeclareSymbol{11}{0}{\draw (0,1.8) node[smalldot] {} -- (-0.7,0.9) -- (0,0) -- (0.7,1) node[smalldot] {};}
\DeclareSymbol{10}{0}{\draw (0,1.8) node[smalldot] {} -- (-0.8,0.9) -- (0,0);}
\DeclareSymbol{21}{0.7}{\draw (-1,1.8) node[smalldot] {} -- (-0.5,0.9); \draw (0,1.8) node[smalldot] {} -- (-0.5,0.9) -- (0,0) -- (0.5,0.9) node[smalldot] {};}

\DeclareSymbol{20}{0.7}{\draw (-1,1.8) node[smalldot] {} -- (-0.5,0.9); \draw (0,1.8) node[smalldot] {} -- (-0.5,0.9) -- (-0.5,0);}

\DeclareSymbol{210}{1.1}{\draw (-0.5,2.4) node[smalldot] {} -- (-1,1.6);\draw (-1.5,2.4) node[smalldot] {} -- (-0.5,0.8); \draw (0,1.6) node[smalldot] {} -- (-0.5,0.8) -- (-0.5,0);}

\DeclareSymbol{211}{1.1}{\draw (-0.5,2.4) node[smalldot] {} -- (-1,1.6);\draw (-1.5,2.4) node[smalldot] {} -- (-0.5,0.8); \draw (0,1.6) node[smalldot] {} -- (-0.5,0.8) -- (0,0) -- (0.5,0.8) node[smalldot] {}; \fill[white] (0,0) circle (8pt); \draw (0,0) circle (8pt);}

\DeclareSymbol{22}{0.7}{\draw (-1.5,1.8) node[smalldot] {} -- (-1,0.9) -- (0,0) -- (1,0.9) -- (1.5,1.8) node[smalldot] {};
\draw (-0.5,1.8) node[smalldot] {} -- (-1,0.9);\draw (0.5,1.8) node[smalldot] {} -- (1,0.9);}
\DeclareSymbol{j1}{0.5}{\draw (-1, 2.4) node[smalldot] {} -- (-0.5,1.6) node[smalldot] {} -- (0, 2.4) node[smalldot] {}; \draw (-0.5,1.6) -- (0,0.8)  node[smalldot] {} -- (0.5, 1.6) node[smalldot] {}; \draw (0,0.8)
-- (0.5,0) node[smalldot] {};}
\DeclareSymbol{j2}{0.5}{\draw (0,0) node[smalldot] {} -- (0.5,0.8) node[smalldot] {};}
\DeclareSymbol{j3}{1.3}{\draw (-1, 2.4) node[smalldot] {} -- (-0.5,1.6) node[smalldot] {} -- (0, 2.4) node[smalldot] {}; \draw (-0.5,1.6) -- (0,0.8)  node[smalldot] {} -- (0.5, 1.6) node[smalldot] {}; }
\DeclareSymbol{j5}{0.5}{ \draw (0,0.8)  node[smalldot] {} -- (0.5, 1.6) node[smalldot] {}; \draw (0,0.8)
-- (0.5,0) node[smalldot] {};}

\DeclareSymbol{Xi22}{0.5}{\draw (0,0) node[xi] {} -- (-1,1) node[not] {} -- (0,2) node[xi] {};}

\DeclareSymbol{Xi2}{-2}{\draw (0,-0.25) node[xi] {} -- (-1,1) node[xi] {};}
\DeclareSymbol{Xi3}{0}{\draw (0,0) node[xi] {} -- (-1,1) node[xi] {} -- (0,2) node[xi] {};}
\DeclareSymbol{Xi4}{2}{\draw (0,0) node[xi] {} -- (-1,1) node[xi] {} -- (0,2) node[xi] {} -- (-1,3) node[xi] {};}
\DeclareSymbol{Xi2X}{-2}{\draw (0,-0.25) node[xi] {} -- (-1,1) node[xix] {};}
\DeclareSymbol{XXi2}{-2}{\draw (0,-0.25) node[xix] {} -- (-1,1) node[xi] {};}

\DeclareSymbol{IXi2}{0}{\draw (0,-0.25) node[not] {} -- (-1,1) node[xi] {} -- (0,2) node[xi] {};}
\DeclareSymbol{IXi^2}{-1}{\draw (-1,1) node[xi] {} -- (0,0) node[not] {} -- (1,1) node[xi] {};}

\DeclareSymbol{XiX}{-2.8}{\node[xibx] {};}
\DeclareSymbol{Xi}{-2.8}{\node[xib] {};}
\DeclareSymbol{IXiX}{-1}{\draw (0,-0.25) node[not] {} -- (-1,1) node[xix] {};}

\DeclareSymbol{Xi3b}{-1}{\draw (-1,1) node[xi] {} -- (0,0) node[xi] {} -- (1,1) node[xi] {};}

\DeclareSymbol{IXi3}{2}{\draw (0,-0.25) node[not] {} -- (-1,1) node[xi] {} -- (0,2) node[xi] {} -- (-1,3) node[xi] {};}
\DeclareSymbol{1p}{-2}{\draw (0,-0.7) node[not] {} -- (0,1.2) node[xi] {};}
\DeclareSymbol{1pd}{-2}{\draw (0,-0.7) node[not] {} -- (0,1.2) node[xi] {}; \draw (0.4,0.2)  {} -- (-0.4,0.2) {};}
\DeclareSymbol{1pdF}{-2}{\draw[red] (0,-0.7) node[not] {} -- (0,1.2) node[xi] {}; \draw (0.4,0.2)  {} -- (-0.4,0.2) {};}
\DeclareSymbol{1pdFP}{-2}{\draw[green] (0,-0.7) node[not] {} -- (0,1.2) node[xi] {}; \draw (0.4,0.2)  {} -- (-0.4,0.2) {};}
\DeclareSymbol{1pdF1pd}{-2}{\draw (0,-0.7) node[not] {} -- (-1.0,0.7) node[xi] {}; \draw[red](0,-0.7) node[not] {} -- (1.0,0.7) node[xi] {};\draw (-0.5,0.0)  {} -- (-0.22,0.2) {};\draw (-0.5,0.0)  {} -- (-0.78,-0.2) {};\draw (0.5,0.0)  {} -- (0.22,0.2) {};\draw (0.5,0.0)  {} -- (0.78,-0.2) {};}
\DeclareSymbol{1pdF1pd0}{-2}{\draw (0,-0.7) node[not] {} -- (-1.0,0.7) node[xi] {}; \draw[black](0,-0.7) node[not] {} -- (1.0,0.7) node[xi] {};\draw (-0.5,0.0)  {} -- (-0.22,0.2) {};\draw (-0.5,0.0)  {} -- (-0.78,-0.2) {};\draw (0.5,0.0)  {} -- (0.22,0.2) {};\draw (0.5,0.0)  {} -- (0.78,-0.2) {};}

\DeclareSymbol{11p}{-2}{\draw (0,-0.5) node [xi] {} -- (-1.2,1.2) node [xi] {};}
\DeclareSymbol{XiI}{-2}{\draw (0,-0.25) node[xi] {} -- (-1,1) node[not] {};}

\DeclareSymbol{Xi4b}{-1}{\draw(0,1.5) node[xi] {} -- (0,0); \draw (-1,1) node[xi] {} -- (0,0) node[xi] {} -- (1,1) node[xi] {};}
\DeclareSymbol{Xi4c}{0}{\draw (0,1) -- (0.8,2.2) node[xi] {};\draw (0,-0.25) node[xi] {} -- (0,1) node[xi] {} -- (-0.8,2.2) node[xi] {};}
\DeclareSymbol{Xi4d}{-4.5}{\draw (0,-1.5) node[not] {} -- (0,0); \draw (-1,1) node[xi] {} -- (0,0) node[xi] {} -- (1,1) node[xi] {};}
\DeclareSymbol{Xi4e}{0}{\draw (0,2) node[xi] {} -- (-1,1) node[xi] {} -- (0,0) node[xi] {} -- (1,1) node[xi] {};}

\DeclareSymbol{1gk}{0}{\draw[thick,green] (0,0)  -- (0,1.2) node[tinydot] {};}
\DeclareSymbol{1bk}{0}{\draw[thick,blue] (0,0)  -- (0,1.2) node[tinydot] {};} 
\DeclareSymbol{1rk}{0}{\draw[thick,red] (0,0)  -- (0,1.2) node[tinydot] {};}
\DeclareSymbol{2k}{0}{\draw[thick,red] (-0.5,1.2) node[tinydot] {} -- (0,0);\draw[thick,blue] (0,0)-- (0.5,1.2) node[tinydot] {};}
\DeclareSymbol{20gk}{0.7}{\draw[red,thick] (-1,1.8) node[tinydot] {} -- (-0.5,0.9); \draw[thick,green] (-0.5,0.9)-- (-0.5,0) ; \draw[thick,blue] (0,1.8) node[tinydot] {} -- (-0.5,0.9) ;}
\DeclareSymbol{20bk}{0.7}{\draw[red,thick] (-1,1.8) node[tinydot] {} -- (-0.5,0.9); \draw[thick,blue] (-0.5,0.9)-- (-0.5,0) ; \draw[thick,blue] (0,1.8) node[tinydot] {} -- (-0.5,0.9) ;}
\DeclareSymbol{20rk}{0.7}{\draw[red,thick] (-1,1.8) node[tinydot] {} -- (-0.5,0.9); \draw[thick,red] (-0.5,0.9)-- (-0.5,0) ; \draw[thick,blue] (0,1.8) node[tinydot] {} -- (-0.5,0.9) ;}
\DeclareSymbol{21k}{0.7}{\draw[red,thick] (-1,1.8) node[tinydot] {} -- (0,0); \draw[thick,blue] (0,1.8) node[tinydot] {} -- (-0.5,0.9); \draw[thick,blue] (0,0) -- (0.5,0.9) node[tinydot] {};}
\DeclareSymbol{210gk}{1.1}{\draw[thick,blue] (-0.5,2.4) node[tinydot] {} -- (-1,1.6);\draw[thick,red] (-1.5,2.4) node[tinydot] {} -- (-0.5,0.8); \draw[thick,blue] (0,1.6) node[tinydot] {} -- (-0.5,0.8) ;\draw[thick,green] (-0.5,0.8) -- (-0.5,0);}
\DeclareSymbol{210bk}{1.1}{\draw[thick,blue] (-0.5,2.4) node[tinydot] {} -- (-1,1.6);\draw[thick,red] (-1.5,2.4) node[tinydot] {} -- (-0.5,0.8); \draw[thick,blue] (0,1.6) node[tinydot] {} -- (-0.5,0.8);\draw[thick,blue] (-0.5,0.8) -- (-0.5,0);}
\DeclareSymbol{210rk}{1.1}{\draw[thick,blue] (-0.5,2.4) node[tinydot] {} -- (-1,1.6);\draw[thick,red] (-1.5,2.4) node[tinydot] {} -- (-0.5,0.8); \draw[thick,blue] (0,1.6) node[tinydot] {} -- (-0.5,0.8);\draw[thick,red] (-0.5,0.8) -- (-0.5,0);}
\DeclareSymbol{211k}{1.1}{\draw[thick,blue] (-0.5,2.4) node[tinydot] {} -- (-1,1.6);\draw[red,thick] (-1.5,2.4) node[tinydot] {} -- (0,0); \draw (0,1.6)[thick,blue] node[tinydot] {} -- (-0.5,0.8);\draw[thick,blue] (0.5,0.8)node[tinydot] {} -- (0,0) ; \fill[white] (0,0) circle (8pt); \draw (0,0) circle (8pt);}
\DeclareSymbol{2110gk}{1.0}{\draw[thick,blue] (-0.5,2.4) node[tinydot] {} -- (-1,1.6);\draw[red,thick] (-1.5,2.4) node[tinydot] {} -- (0,0); \draw (0,1.6)[thick,blue] node[tinydot] {} -- (-0.5,0.8);\draw[thick,blue] (0.5,0.8)node[tinydot] {} -- (0,0) ;\draw[green,thick] (0,0) -- (0,-1.2); \fill[white] (0,0) circle (8pt); \draw (0,0) circle (8pt);}
\DeclareSymbol{2110rk}{1.0}{\draw[thick,blue] (-0.5,2.4) node[tinydot] {} -- (-1,1.6);\draw[red,thick] (-1.5,2.4) node[tinydot] {} -- (0,0); \draw (0,1.6)[thick,blue] node[tinydot] {} -- (-0.5,0.8);\draw[thick,blue] (0.5,0.8)node[tinydot] {} -- (0,0) ;\draw[red,thick] (0,0) -- (0,-1.2); \fill[white] (0,0) circle (8pt); \draw (0,0) circle (8pt);}
\DeclareSymbol{22k}{0.7}{\draw[thick,red] (-1.5,1.8) node[tinydot] {} -- (-1,0.9)-- (0,0); \draw[thick,blue] (0,0) -- (1,0.9) -- (1.5,1.8) node[tinydot] {};
\draw[thick,blue] (-0.5,1.8) node[tinydot] {} -- (-1,0.9);\draw[thick,red] (0.5,1.8) node[tinydot] {} -- (1,0.9);}
\DeclareSymbol{220gk}{0.7}{\draw[thick,red] (-1.5,1.8) node[tinydot] {} -- (-1,0.9)-- (0,0); \draw[thick,blue] (0,0) -- (1,0.9) -- (1.5,1.8) node[tinydot] {};
\draw[thick,blue] (-0.5,1.8) node[tinydot] {} -- (-1,0.9);\draw[thick,red] (0.5,1.8) node[tinydot] {} -- (1,0.9);\draw[thick,green] (0,0) -- (0,-0.9);}
\DeclareSymbol{220rk}{0.7}{\draw[thick,red] (-1.5,1.8) node[tinydot] {} -- (-1,0.9)-- (0,0); \draw[thick,blue] (0,0) -- (1,0.9) -- (1.5,1.8) node[tinydot] {};
\draw[thick,blue] (-0.5,1.8) node[tinydot] {} -- (-1,0.9);\draw[thick,red] (0.5,1.8) node[tinydot] {} -- (1,0.9);\draw[thick,red] (0,0) -- (0,-0.9);}
\DeclareSymbol{101k}{1.0}{\draw[red,thick] (0,0) node[tinydot] {} -- (0,1.6);\draw[thick,blue] (0,1.6) node[tinydot] {} -- (1.6,2.4);\draw[thick,blue] (0,0) node[tinydot] {} -- (1.6,0.8); \fill[white] (0,0) circle (8pt); \draw (0,0) circle (8pt);}
\DeclareSymbol{10rk}{1.0}{\draw[red,thick] (0,0)  -- (0,1.6);\draw[thick,blue] (0,1.6) node[tinydot] {} -- (1.6,2.4);}
\DeclareSymbol{10bk}{1.0}{\draw[thick,blue] (0,0) node[tinydot] {} -- (0,1.6);\draw[thick,blue] (0,1.6) node[tinydot] {} -- (1.6,2.4);}
\DeclareSymbol{10gk}{1.0}{\draw[green,thick] (0,0)  -- (0,1.6);\draw[thick,blue] (0,1.6) node[tinydot] {} -- (1.6,2.4);}

\DeclareSymbol{1gk0}{0}{\draw[thick,green] (0,0)  -- (0,1.2) node[tinydot] {};}
\DeclareSymbol{1bk0}{0}{\draw[thick,black] (0,0)  -- (0,1.2) node[tinydot] {};} 
\DeclareSymbol{1rk0}{0}{\draw[thick,black] (0,0)  -- (0,1.2) node[tinydot] {};}
\DeclareSymbol{2k0}{0}{\draw[thick,black] (-0.5,1.2) node[tinydot] {} -- (0,0);\draw[thick] (0,0)-- (0.5,1.2) node[tinydot] {};}
\DeclareSymbol{20gk0}{0.7}{\draw[black,thick] (-1,1.8) node[tinydot] {} -- (-0.5,0.9); \draw[thick,green] (-0.5,0.9)-- (-0.5,0) ; \draw[thick] (0,1.8) node[tinydot] {} -- (-0.5,0.9) ;}
\DeclareSymbol{20bk0}{0.7}{\draw[black,thick] (-1,1.8) node[tinydot] {} -- (-0.5,0.9); \draw[thick,black] (-0.5,0.9)-- (-0.5,0) ; \draw[thick] (0,1.8) node[tinydot] {} -- (-0.5,0.9) ;}
\DeclareSymbol{20rk0}{0.7}{\draw[black,thick] (-1,1.8) node[tinydot] {} -- (-0.5,0.9); \draw[thick,black] (-0.5,0.9)-- (-0.5,0) ; \draw[thick] (0,1.8) node[tinydot] {} -- (-0.5,0.9) ;}
\DeclareSymbol{21k0}{0.7}{\draw[black,thick] (-1,1.8) node[tinydot] {} -- (0,0); \draw[thick] (0,1.8) node[tinydot] {} -- (-0.5,0.9); \draw[thick] (0,0) -- (0.5,0.9) node[tinydot] {};}
\DeclareSymbol{210gk0}{1.1}{\draw[thick] (-0.5,2.4) node[tinydot] {} -- (-1,1.6);\draw[thick,black] (-1.5,2.4) node[tinydot] {} -- (-0.5,0.8); \draw[thick] (0,1.6) node[tinydot] {} -- (-0.5,0.8) ;\draw[thick,green] (-0.5,0.8) -- (-0.5,0);}
\DeclareSymbol{210bk0}{1.1}{\draw[thick] (-0.5,2.4) node[tinydot] {} -- (-1,1.6);\draw[thick,black] (-1.5,2.4) node[tinydot] {} -- (-0.5,0.8); \draw[thick] (0,1.6) node[tinydot] {} -- (-0.5,0.8);\draw[thick] (-0.5,0.8) -- (-0.5,0);}
\DeclareSymbol{210rk0}{1.1}{\draw[thick] (-0.5,2.4) node[tinydot] {} -- (-1,1.6);\draw[thick,black] (-1.5,2.4) node[tinydot] {} -- (-0.5,0.8); \draw[thick] (0,1.6) node[tinydot] {} -- (-0.5,0.8);\draw[thick,black] (-0.5,0.8) -- (-0.5,0);}
\DeclareSymbol{211k0}{1.1}{\draw[thick] (-0.5,2.4) node[tinydot] {} -- (-1,1.6);\draw[black,thick] (-1.5,2.4) node[tinydot] {} -- (0,0); \draw (0,1.6)[thick] node[tinydot] {} -- (-0.5,0.8);\draw[thick] (0.5,0.8)node[tinydot] {} -- (0,0) ; \fill[white] (0,0) circle (8pt); \draw (0,0) circle (8pt);}
\DeclareSymbol{2110gk0}{1.0}{\draw[thick] (-0.5,2.4) node[tinydot] {} -- (-1,1.6);\draw[black,thick] (-1.5,2.4) node[tinydot] {} -- (0,0); \draw (0,1.6)[thick] node[tinydot] {} -- (-0.5,0.8);\draw[thick] (0.5,0.8)node[tinydot] {} -- (0,0) ;\draw[green,thick] (0,0) -- (0,-1.2); \fill[white] (0,0) circle (8pt); \draw (0,0) circle (8pt);}
\DeclareSymbol{2110rk0}{1.0}{\draw[thick] (-0.5,2.4) node[tinydot] {} -- (-1,1.6);\draw[black,thick] (-1.5,2.4) node[tinydot] {} -- (0,0); \draw (0,1.6)[thick] node[tinydot] {} -- (-0.5,0.8);\draw[thick] (0.5,0.8)node[tinydot] {} -- (0,0) ;\draw[black,thick] (0,0) -- (0,-1.2); \fill[white] (0,0) circle (8pt); \draw (0,0) circle (8pt);}
\DeclareSymbol{22k0}{0.7}{\draw[thick,black] (-1.5,1.8) node[tinydot] {} -- (-1,0.9)-- (0,0); \draw[thick] (0,0) -- (1,0.9) -- (1.5,1.8) node[tinydot] {};
\draw[thick] (-0.5,1.8) node[tinydot] {} -- (-1,0.9);\draw[thick,black] (0.5,1.8) node[tinydot] {} -- (1,0.9);}
\DeclareSymbol{220gk0}{0.7}{\draw[thick,black] (-1.5,1.8) node[tinydot] {} -- (-1,0.9)-- (0,0); \draw[thick] (0,0) -- (1,0.9) -- (1.5,1.8) node[tinydot] {};
\draw[thick] (-0.5,1.8) node[tinydot] {} -- (-1,0.9);\draw[thick,black] (0.5,1.8) node[tinydot] {} -- (1,0.9);\draw[thick,green] (0,0) -- (0,-0.9);}
\DeclareSymbol{220rk0}{0.7}{\draw[thick,black] (-1.5,1.8) node[tinydot] {} -- (-1,0.9)-- (0,0); \draw[thick] (0,0) -- (1,0.9) -- (1.5,1.8) node[tinydot] {};
\draw[thick] (-0.5,1.8) node[tinydot] {} -- (-1,0.9);\draw[thick,black] (0.5,1.8) node[tinydot] {} -- (1,0.9);\draw[thick,black] (0,0) -- (0,-0.9);}
\DeclareSymbol{101k0}{1.0}{\draw[black,thick] (0,0) node[tinydot] {} -- (0,1.6);\draw[thick] (0,1.6) node[tinydot] {} -- (1.6,2.4);\draw[thick] (0,0) node[tinydot] {} -- (1.6,0.8); \fill[white] (0,0) circle (8pt); \draw (0,0) circle (8pt);}
\DeclareSymbol{10rk0}{1.0}{\draw[black,thick] (0,0)  -- (0,1.6);\draw[thick] (0,1.6) node[tinydot] {} -- (1.6,2.4);}
\DeclareSymbol{10bk0}{1.0}{\draw[thick] (0,0) node[tinydot] {} -- (0,1.6);\draw[thick] (0,1.6) node[tinydot] {} -- (1.6,2.4);}
\DeclareSymbol{10gk0}{1.0}{\draw[green,thick] (0,0)  -- (0,1.6);\draw[thick] (0,1.6) node[tinydot] {} -- (1.6,2.4);}

\DeclareSymbol{1gks}{0}{\draw (0,0)  -- (0,1.2) node[tinydot] {};}
\DeclareSymbol{1bks}{0}{\draw (0,0)  -- (0,1.2) node[tinydot] {};} 
\DeclareSymbol{1rks}{0}{\draw (0,0)  -- (0,1.2) node[tinydot] {};}
\DeclareSymbol{2ks}{0}{\draw (-0.5,1.2) node[tinydot] {} -- (0,0);\draw  (0,0)-- (0.5,1.2) node[tinydot] {};}
\DeclareSymbol{20gks}{0.7}{\draw (-1,1.8) node[tinydot] {} -- (-0.5,0.9); \draw (-0.5,0.9)-- (-0.5,0) ; \draw  (0,1.8) node[tinydot] {} -- (-0.5,0.9) ;}
\DeclareSymbol{20iks}{0.7}{\draw (-1,1.8) node[tinydot] {} -- (-0.5,0.9); \draw (-0.5,0.9)-- (-0.5,0) ; \draw  (0,1.8) node[tinydot] {} -- (-0.5,0.9) ;}
\DeclareSymbol{20bks}{0.7}{\draw (-1,1.8) node[tinydot] {} -- (-0.5,0.9); \draw (-0.5,0.9)-- (-0.5,0) ; \draw  (0,1.8) node[tinydot] {} -- (-0.5,0.9) ;}
\DeclareSymbol{20rks}{0.7}{\draw (-1,1.8) node[tinydot] {} -- (-0.5,0.9); \draw (-0.5,0.9)-- (-0.5,0) ; \draw  (0,1.8) node[tinydot] {} -- (-0.5,0.9) ;}
\DeclareSymbol{21ks}{0.7}{\draw (-1,1.8) node[tinydot] {} -- (0,0); \draw  (0,1.8) node[tinydot] {} -- (-0.5,0.9); \draw  (0,0) -- (0.5,0.9) node[tinydot] {};}
\DeclareSymbol{210gks}{1.1}{\draw  (-0.5,2.4) node[tinydot] {} -- (-1,1.6);\draw (-1.5,2.4) node[tinydot] {} -- (-0.5,0.8); \draw  (0,1.6) node[tinydot] {} -- (-0.5,0.8) ;\draw (-0.5,0.8) -- (-0.5,0);}
\DeclareSymbol{210iks}{1.1}{\draw  (-0.5,2.4) node[tinydot] {} -- (-1,1.6);\draw (-1.5,2.4) node[tinydot] {} -- (-0.5,0.8); \draw  (0,1.6) node[tinydot] {} -- (-0.5,0.8) ;\draw (-0.5,0.8) -- (-0.5,0);}
\DeclareSymbol{210bks}{1.1}{\draw  (-0.5,2.4) node[tinydot] {} -- (-1,1.6);\draw (-1.5,2.4) node[tinydot] {} -- (-0.5,0.8); \draw  (0,1.6) node[tinydot] {} -- (-0.5,0.8);\draw  (-0.5,0.8) -- (-0.5,0);}
\DeclareSymbol{210rks}{1.1}{\draw  (-0.5,2.4) node[tinydot] {} -- (-1,1.6);\draw (-1.5,2.4) node[tinydot] {} -- (-0.5,0.8); \draw  (0,1.6) node[tinydot] {} -- (-0.5,0.8);\draw[ red] (-0.5,0.8) -- (-0.5,0);}
\DeclareSymbol{211ks}{1.1}{\draw  (-0.5,2.4) node[tinydot] {} -- (-1,1.6);\draw (-1.5,2.4) node[tinydot] {} -- (0,0); \draw (0,1.6)  node[tinydot] {} -- (-0.5,0.8);\draw  (0.5,0.8)node[tinydot] {} -- (0,0) ; \fill[white] (0,0) circle (8pt); \draw (0,0) circle (8pt);}
\DeclareSymbol{2110gks}{1.0}{\draw  (-0.5,2.4) node[tinydot] {} -- (-1,1.6);\draw (-1.5,2.4) node[tinydot] {} -- (0,0); \draw (0,1.6)  node[tinydot] {} -- (-0.5,0.8);\draw  (0.5,0.8)node[tinydot] {} -- (0,0) ;\draw (0,0) -- (0,-1.2); \fill[white] (0,0) circle (8pt); \draw (0,0) circle (8pt);}
\DeclareSymbol{2110rks}{1.0}{\draw  (-0.5,2.4) node[tinydot] {} -- (-1,1.6);\draw (-1.5,2.4) node[tinydot] {} -- (0,0); \draw (0,1.6)  node[tinydot] {} -- (-0.5,0.8);\draw  (0.5,0.8)node[tinydot] {} -- (0,0) ;\draw (0,0) -- (0,-1.2); \fill[white] (0,0) circle (8pt); \draw (0,0) circle (8pt);}
\DeclareSymbol{22ks}{0.7}{\draw (-1.5,1.8) node[tinydot] {} -- (-1,0.9)-- (0,0); \draw  (0,0) -- (1,0.9) -- (1.5,1.8) node[tinydot] {};
\draw  (-0.5,1.8) node[tinydot] {} -- (-1,0.9);\draw (0.5,1.8) node[tinydot] {} -- (1,0.9);}
\DeclareSymbol{220gks}{0.7}{\draw (-1.5,1.8) node[tinydot] {} -- (-1,0.9)-- (0,0); \draw  (0,0) -- (1,0.9) -- (1.5,1.8) node[tinydot] {};
\draw  (-0.5,1.8) node[tinydot] {} -- (-1,0.9);\draw (0.5,1.8) node[tinydot] {} -- (1,0.9);\draw (0,0) -- (0,-0.9);}
\DeclareSymbol{220rks}{0.7}{\draw (-1.5,1.8) node[tinydot] {} -- (-1,0.9)-- (0,0); \draw  (0,0) -- (1,0.9) -- (1.5,1.8) node[tinydot] {};
\draw  (-0.5,1.8) node[tinydot] {} -- (-1,0.9);\draw (0.5,1.8) node[tinydot] {} -- (1,0.9);\draw (0,0) -- (0,-0.9);}
\DeclareSymbol{101ks}{1.0}{\draw (0,0) node[tinydot] {} -- (0,1.6);\draw  (0,1.6) node[tinydot] {} -- (1.6,2.4);\draw  (0,0) node[tinydot] {} -- (1.6,0.8); \fill[white] (0,0) circle (8pt); \draw (0,0) circle (8pt);}
\DeclareSymbol{10rks}{1.0}{\draw (0,0)  -- (0,1.6);\draw  (0,1.6) node[tinydot] {} -- (1.6,2.4);}
\DeclareSymbol{10bks}{1.0}{\draw  (0,0) node[tinydot] {} -- (0,1.6);\draw  (0,1.6) node[tinydot] {} -- (1.6,2.4);}
\DeclareSymbol{10gks}{1.0}{\draw (0,0)  -- (0,1.6);\draw  (0,1.6) node[tinydot] {} -- (1.6,2.4);}
\begin{document}

\title{Operator dependent para-controlled calculus, periodic homogenisation and singular PDEs}

\author{
Yilin Chen
\and
Weijun Xu
}

\maketitle

\begin{abstract}
    We develop a variant of the para-controlled distributions framework based on operator dependent (generalised) Besov spaces. This in particular covers the class of oscillatory operators from periodic homogenisation. We then apply the framework to study periodic homogenisation problems for a class of singular stochastic PDEs. One novel viewpoint is to measure regularity via the geometry induced by the operators at hand rather than classical Fourier basis. 
    
    The generalised Besov spaces were introduced by Kerkyacharian-Petrushev (\cite{KP_general_Besov}). In contrast to the Fourier decomposition in the classical setting, they are built from spectral decomposition of general self-adjoint elliptic operators. A major difference is that products of two spectrally localised functions in general have ``frequencies" spread over the whole spectrum. We obtain a decay estimate on ``high frequency" component of the product of two spectrally-localised functions. This enables us to obtain uniform bounds for the naturally associated para-product and commutator operations over a class of operators. This in particular includes the family of periodic homogenisation operators uniform in the oscillation parameter. 

    Next, we show convergence properties of these operations in the periodic homogenisation setting. A key ingredient is a sharp quantitative convergence of the generalised Littlewood-Paley block operators. The proof utilises the Helffer-Sj\"ostrand's formula for self-adjoint operators. 

    Finally, we apply this framework to periodic homogenisation problems for dynamical $\phi^4_3$, $2$D gPAM, and $1$D KPZ. Assuming the convergence of stochastic objects to their homogenised limits, we show that the solutions and fluxes also converge to the corresponding limits. The convergence of the explicit stochastic terms with suitable renormalisations will be treated in a separate note. 
\end{abstract}

\tableofcontents

\bigskip
\bigskip

\section{Introduction}

\subsection{The problem and motivation}

The aim of this article is to develop a variant of the para-controlled distributions framework that are built from the self-adjoint elliptic operator for the problem at hand. This includes obtaining uniform bounds on operations such as generalised para-products, commutators and heat kernel convolutions over a class of operators, as well as obtaining convergences of these operations in the periodic homogenisation setting. 

The main motivation for developing such a variant is to study periodic homogenisation problem for singular parabolic stochastic PDEs. The linear operators in all these problems are given by
\begin{equation*}
    \lL_\eps := \div \big( A_\eps \nabla \big)\;,
\end{equation*}
where $A_\eps = A(\cdot / \eps)$ for some coefficient matrix $A$ satisfying Assumption~\ref{as:A} below, and $\eps \in \N^{-1}$ is an oscillation parameter taking values in inverse of integers to keep the problem on the torus. 

\begin{assumption} \label{as:A}
Let $A: \T^d \rightarrow \R^{d \times d}$ be symmetric, uniformly elliptic, $1$-periodic and H\"older-$\alpha$ continuous for some $\alpha>0$. 
\end{assumption}

The problems that can be dealt with the analytic framework in this article include the followings. 
\begin{itemize}
\item The dynamical $\phi^4_3$ model, formally given by
\begin{equation} \label{e:phi4_formal_intro}
    \d_t \phi_\eps = \lL_\eps \phi_\eps - \phi_\eps^3 + \xi + \infty_\eps \cdot \phi_\eps\;, \qquad (t,x) \in \R^{+} \times \T^3\;.
\end{equation}
Here, $\xi$ is the spacetime white noise on the three dimensional torus, and $\infty_\eps$ is a formal renormalisation function that depends both on the spatial location $x \in \T^3$ and the homogenisation parameter $\eps$, but independent of time. 

\item The two dimensional generalised parabolic Anderson model, formally given by
\begin{equation} \label{e:gPAM_formal_intro}
    \begin{split}
    \d_t u_\eps = \, &\lL_\eps u_\eps + \big(   \nabla^{\sT} u_\eps \, F_\eps \, \nabla u_\eps - \, \widetilde{\infty}_\eps g^2 (u_\eps) \big)\\
    &\qquad \qquad + g(u_\eps) \big( \xi - \infty_\eps g'(u_\eps) \big)\;, \qquad (t,x) \in \R^+ \times \T^2\;.
    \end{split}
\end{equation}
Here, $\xi$ is the \textit{spatial white noise} on two dimensional torus, and $F_\eps = F(\cdot / \eps)$ for some $F:\T^2\rightarrow \R^{2\times 2}$ bounded periodic function taking values in symmetric $2 \times 2$ matrices, and $g: \R \rightarrow \R$ is a smooth function with bounded derivatives. Again, $\infty_\eps$ and $\widetilde{\infty}_\eps$ are formal renormalisation functions independent of time. 

\item The one-dimensional KPZ equation, formally given by
\begin{equation} \label{e:kpz_formal_intro}
    \d_t h_\eps = \lL_\eps h_\eps + g_\eps \, (\d_x h_\eps)^2 + \xi - \infty_\eps - \widetilde{\infty}_\eps\cdot A_\eps\cdot \d_x h_\eps\;, \qquad (t,x) \in \R^+ \times \T\;.
\end{equation}
Here, $\xi$ is the one-dimensional spacetime white noise, $g_\eps := g(\cdot / \eps)$ is a bounded periodic function, and $\infty_\eps$ and $\widetilde{\infty}_\eps$ are formal renormalisation functions independent of time. 
\end{itemize}

All of the above are examples of periodic homogenisation problem for singular stochastic PDEs. At least in the constant coefficient case (that is, $A$ is a constant coefficient matrix, so that $\lL_\eps \equiv \lL = \div(A \nabla)$ is a constant change-of-coordinate of the Laplacian), this falls within the achievements of singular SPDEs in the past decade. With the developments of the theory of regularity structures (\cite{rs_theory}), para-controlled distributions (\cite{GIP_para_control}), and renormalisation group approaches (\cite{Kupiainen_RG, Duch_RG_cont}), one can now give rigorous constructions of solutions to a very large class of subcritical singular stochastic PDEs, including in particular all the examples above. 

The solutions obtained in this way are limits of solutions to classical PDEs. One insight from these developments is that although the noise-to-solution map is discontinuous in the topology where the smooth approximation $\xi^{(\delta)}$ converges to the original noise $\xi$, the solution does depend continuously on a collection of enhanced noises $\Upsilon^{(\delta)}$. These enhanced objects are built from the regularised version $\xi^{(\delta)}$ in a nonlinear but explicit way through stochastic analysis methods, possibly with renormalisations by diverging constants $C^{(\delta)}$. The possible choices of counter-terms (indexed by renormalisation constants) as well as enhanced stochastic objects $\Upsilon^{(\delta)}$ have also been systematically understood (\cite{rs_algebraic, rs_analytic, bphz_spectral_gap, rs_diagram_free}). 

The above approaches could in general also treat equations with uniformly elliptic operators whose coefficient matrix has sufficient regularity to ensure Schauder estimates. The analytic part of the programme has been achieved in \cite{SPDE_variable_nontrans} and more generally in \cite{rs_manifolds} with regularity structures framework. Very recently, a general BPHZ theorem for convergence of stochastic objects has been established in \cite{bphz_variable}. 

With this respect, for each fixed $\eps>0$, the solution $v_\eps$ to the above equations ($v_\eps$ can denote $\phi_\eps$, $u_\eps$ or $h_\eps$) should be understood as the limit of the regularised solutions $v_\eps^{(\delta)}$ after suitable renormalisations. This requires Schauder estimate of the linear operator, which in turn relied on at least some H\"older continuity of the coefficient matrix. This is simply not true for the coefficient $A_\eps$ if one wants uniformity in $\eps$ unless $A$ is constant. Hence, it is not a priori clear whether the solution $v_\eps$ obtained here as the $\delta \rightarrow 0$ limit of $v_\eps^{(\delta)}$ has uniform bounds in $\eps$. 

Understanding the asymptotic behaviour of the solutions as $\eps \rightarrow 0$ concerns the problem of homogenisation. This is another class of problems involving taking singular limits. Suppose $A$ satisfies Assumption~\ref{as:A}. Let $\lL_\eps = \div \big( A(\cdot/\eps) \nabla \big)$ as before. In a simple setup, periodic homogenisation is to consider the $\eps \rightarrow 0$ limit of the equation
\begin{equation} \label{e:homo_linear}
    \d_t v_\eps = \lL_\eps v_\eps + f\;, \qquad (t,x) \in \R^+ \times \T^d\;,
\end{equation}
with nice input $f$ and initial condition, and again with $\eps = \eps_N = \frac{1}{N}$ to keep the problem on the torus. It is well known that as $\eps \rightarrow 0$, $v_\eps$ converges to a limit $v_0$ in suitable sense, where $v_0$ solves the equation
\begin{equation*}
    \d_t v_0 = \div (\bar{A} \nabla v_0) + f\;, \qquad (t,x) \in \R^+ \times \T^d
\end{equation*}
for some \textit{constant elliptic matrix $A$}. The interesting point is that although $A_\eps = A(\cdot / \eps)$ converges weakly to the average of $A$ on the torus, the constant elliptic matrix $\bar{A}$ is not the average of $A$. This suggests that the coefficient-to-solution map $A_\eps \mapsto v_\eps$ is again not continuous in the topology where $A_\eps$ converges. 

It is then natural to ask how these two discontinuous solution maps / singular limiting procedures interact with each other. This leads to the consideration of periodic homogenisation problems for singular stochastic PDEs. 

Let $v_\eps^{(\delta)}$ denote the solution of a (general) periodic homogenisation problem for singular SPDEs with homogenisation parameter $\eps$ and regularisation parameter $\delta$, possibly after suitable renormalisation. The following questions are natural:
\begin{enumerate}
\item What are natural choices of renormalisation functions?

\item What is the asymptotic behaviour of the solution after suitable renormalisations?
\end{enumerate}

The first question is the probabilistic part while the second is more of a PDE / deterministic analytic nature. In this article, we mainly focus on the second question. 

We first remark that even if the enhanced stochastic objects $\Upsilon_\eps^{(\delta)}$ behave nicely, the analysis of the asymptotic behaviour of $v_\eps^{(\delta)}$ exhibits different levels of difficulties when sending $(\eps,\delta) \rightarrow (0,0)$ in different ways. Let us illustrate this point by looking at the two extreme ways of taking limits. 

For fixed $\delta>0$, the input on the right hand side is smooth, and hence it is standard that $v_\eps^{(\delta)}$ converges to the homogenised limit $v^{(\delta)}$ as $\eps \rightarrow 0$ (slightly trickier when the right hand side involves derivatives of the solution but still standard). Now, $v^{(\delta)}$ satisfies a standard singular SPDE with constant coefficient operator $\lL_0 = \div (\bar{A} \nabla)$, and hence is also well known that it converges to a limit $v$ as $\delta \rightarrow 0$. 

On the other hand, for fixed $\eps > 0$, a series of recent developments (\cite{rs_manifolds, SPDE_variable_nontrans, bphz_variable}) shows that for a large class of variable coefficient singular SPDEs (assuming sufficient regularity on the coefficient), there exist natural choices of counter-terms such that $v_\eps^{(\delta)}$ converges to a limit $v_\eps$ as $\delta \rightarrow 0$. The main question now concerns the asymptotic behaviour of $v_\eps$ as $\eps \rightarrow 0$. This is illustrated in the diagram below. 
\begin{center}
	\begin{tikzpicture}[scale=1.2,baseline=0cm]
	\node at (-3,1.5) [] (ul) {$v_\eps^{(\delta)}$};
	\node at (3,1.5) [] (ur) {$v^{(\delta)}$}; 
	\draw[akernel] (ul) to node[anchor = south]{\scriptsize $\eps \rightarrow 0$ \text{(homo.)}} (ur);
	\node at (3,-1.5) [] (br) {$v$}; 
	\draw[akernel] (ur) to node[labl]{\tiny $\delta \downarrow 0$ (SPDE)} (br); 
	\node at (-3,-1.5) [] (bl) {$v_\eps$}; 
	\draw[akernel] (ul) to node[labl]{\tiny $\delta \downarrow 0$ (SPDE)} (bl); 
	\draw[akernel] (bl) to node[anchor = north]{\scriptsize $\eps \rightarrow 0$ ?} (br);
	\end{tikzpicture}
\end{center}
The solution theory to variable coefficient singular SPDEs so far relies on sufficient regularity of the coefficient matrix, and hence does not directly give uniform-in-$\eps$ bounds in homogenisation problem. The bottom arrow has been a priori out of reach from the usual framework of singular stochastic PDEs or homogenisations. In general, with $(\eps,\delta) \rightarrow (0,0)$, the closer the limiting procedure is to the bottom arrow, the harder the analytic nature of the problem is, with the bottom arrow itself being the hardest. These are mainly because of the following two reasons:
\begin{enumerate}
    \item \textbf{Lack of uniform-in-$\eps$ Schauder estimate for $\d_t - \lL_\eps$.} Existing frameworks for singular SPDEs rely on sufficient smoothing properties of the linear operator in the equation. This improves regularities over iterations. But this is not possible for $\d_t - \lL_\eps$. The solution $v_\eps$ to \eqref{e:homo_linear} only converges weakly in $H^1$. In fact, even if $f$ is smooth, $v_\eps$ cannot have any regularity strictly above $1$ uniformly in $\eps$. 

    \item \textbf{Incompatible solution ansatzes.} Solutions to both singular SPDEs and periodic homogenisation problems are based on suitable ansatzes. But the ansatzes from singular SPDEs (Taylor-like expansions from regularity structures or para-controlled distributions) and the two-scale expansion from homogenisations have different formats. It is a priori not clear how to combine them to a suitable ansatz keeping features both from periodic homogenisation and from singular SPDEs. 
\end{enumerate}

To the best of our knowledge, the first works of periodic homogenisation for singular stochastic PDEs are \cite{st_homo_phi42, homo_pphi2} in the context of dynamical $\phi^4_2$. In that case, one only needs a Schauder estimate up to the level $\cC^{0+}$, which is available for the operator $\lL_\eps$ at hand. Furthermore, features from both singular SPDEs and homogenisations can be incorporated into a simple solution ansatz (Da Prato-Debussche decomposition). 

The next natural example is the two dimensional generalised parabolic Anderson model \eqref{e:gPAM_formal_intro}, which encounters both of the above difficulties (even if $F=0$). In the constant coefficient case, it requires a nontrivial ansatz from regularity structures or para-controlled distributions, and part of this ansatz needs to be of regularity $\cC^{1+}$ in order to multiply with the noise $\xi$. In \cite{homo_gPAM}, we identified a refined para-controlled ansatz for \eqref{e:gPAM_formal_intro} (with $F=0$) adapted to the homogenisation kernels. This enabled us to establish uniform bounds and convergence of the periodic homogenisation problem. The key is to identify the exact ``rough" part of the solution caused by oscillations from homogenisation, and treat its multiplication with the noise by hand. 

While similar techniques should be applicable to \eqref{e:gPAM_formal_intro} with general $F$ and to dynamical $\phi^4_3$, one expects the resulting calculations to be much more involved, and it is not clear to what extent these arguments can be applied to more singular and general examples. Hence, we hope to seek more unified framework to treat a class of problems. At the same time, Hairer-Singh (\cite{st_homo_phi42, homo_spde_rs}) also developed a systematic framework based on regularity structures. We give more detailed discussions in Section~\ref{sec:perspectives} below.

\subsection{Description of the framework}

The approach taken in this article is based on generalisations of the para-controlled distribution framework developed in \cite{GIP_para_control}. Since trigonometric functions form an eigen-basis of the Laplacian, the standard para-controlled calculus, which is based on these trigonometric functions, enjoys nice properties with the heat kernel and has been successful in treating singular SPDEs with the linear operator being Laplacian (or general constant coefficient elliptic operators). 

On the other hand, classical Fourier based methods are not so compatible with problems with variable coefficients. The generalised Besov spaces introduced in \cite{KP_general_Besov} are built from spectral decompositions of a general elliptic operator $\lL = \div (a \nabla)$ whose coefficient matrix $a$ and parabolic Green's function satisfy requirements in Definition~\ref{def:operator_class} below (namely that $a$ belongs to the class $\sS_d(\Lambda, M)$ for some parameters $\Lambda, M>1$). The spectral decomposition of $\lL_a = \div (a \nabla)$ gives natural definition of generalised Littlewood-Paley operators $\Delta_{j,a}$ and hence the spaces $\cC^{\alpha;a}$ for $\alpha \in \R$, all depending on the coefficient matrix $a$. One property for these spaces is that for every $\alpha \in (-1,1)$, the spaces $\cC^{\alpha;a}$ are equivalent among the coefficient matrices $a \in \sS_{d}(\Lambda, M)$, with uniform upper and lower bounds depending on $\Lambda$, $M$ and $d$ only. 

Back to our problems, it turns out that for every coefficient matrix $A$ satisfying Assumption~\ref{as:A}, there exist parameters $\Lambda$ and $M$ such that $A_\eps := A(\cdot/\eps)$ belongs to the class $\sS_d(\Lambda, M)$ for all $\eps$ (see Corollary~\ref{cor:coefficient_uniform}). Hence, it is natural to seek a suitable para-controlled calculus based on these generalised Besov spaces to tackle periodic homogenisation problems for singular SPDEs.  

Let $a \in \sS_d(\Lambda, M)$ for some fixed $\Lambda$ and $M$. We should think of $a = A(\cdot/\eps)$ for some $A$ satisfying Assumption~\ref{as:A} and $\eps \in \N^{-1}$. In the uniform boundedness part, we only seek bounds that depend on $\Lambda$ and $M$ but not on the actual element in $\sS_d(\Lambda, M)$. There are natural notions of Littlewood-Paley blocks ($\Delta_{j,a}$), para-products ($\prec_a$, $\circ_a$ and $\succ_a$) and commutators, all depending on $a \in \sS_d(\Lambda, M)$. They interact nicely with the heat kernel $\iI_a = (\d_t - \lL_a)^{-1}$. One would try to obtain uniform bounds on these operations over $a \in \sS_{d}(\Lambda, M)$ with regularity measured by the $a$-dependent norms $\cC^{\alpha;a}$. The main difficulty is that in general, the product of two spectrally-localised functions, $\Delta_{i,a} f \cdot \Delta_{j,a} g$, has frequencies spread over the whole spectrum. This is in complete contrast with the classical Fourier setting, where the frequency is still localised within a constant multiple of the support of the high-frequency part. A key ingredient here is Proposition~\ref{pr:decoupling_spectrum} below, which captures quantitatively the decay of the product $\Delta_{i,a} f \cdot \Delta_{j,a} g$ in its frequencies. This decay estimate allows us to obtain bounds for para-products and commutator operations, uniformly in $a \in \sS_d(\Lambda, M)$. Hence, these bounds also apply to the class of operators $\lL_\eps = \div \big( A_\eps \nabla \big)$, uniformly in $\eps\in \N^{-1}$. This enables us to obtain uniform-in-$\eps$ solution theory to these homogenisation problem for singular SPDEs. At this stage, assuming uniformity in enhanced stochastic objects after proper renormalisations, we can already establish a well-defined solution theory for each $\eps$, and show that they have uniform-in-$\eps$ bounds (with $\eps$-dependent operations and spaces), at least for $2$D gPAM with $F=0$ and for dynamical $\phi^4_3$. 

The cases for $2$D g-PAM with a general $F$ and KPZ are still a bit different since they involve the spatial derivatives in the nonlinearity, and that differentiation does not commute with $\lL_\eps$ (although it does commute with the Laplacian). Our treatments of the nonlinearities involving derivatives in these two examples are different. For $2$D gPAM, we employ integration by parts and ``completion of products" (similar to the techniques from \cite{homo_gPAM}), and arrive at an ansatz involving finer information on $\nabla u_\eps$. This relies on the fact that $\nabla u_\eps$ has a very mild divergence. For the KPZ equation, we are able to obtain uniform bounds involving derivatives in dimension $d=1$. This is due to the fact that the one dimensional derivative is the same as divergence, while this is not true in higher dimensions. Relevant bounds involving $1$D derivatives are treated in Section~\ref{sec:bounds_derivatives_1d}. They allow us to establish uniform-in-$\eps$ solution theory for KPZ as well. 

One nice feature with these generalised Besov spaces is that with the $\eps$-dependent spaces and operations, the solution ansatz and spaces are in complete analogy with the standard para-controlled calculus in the constant coefficient case (except the quadratic nonlinearity part in $2$D gPAM). 

Once we obtain the uniform bounds for the solutions, the next step is to investigate convergences of the solution and flux to the homogenised ones as $\eps \rightarrow 0$. In standard homogenisation problems, convergences are usually harder than merely getting uniform bounds. As purely analytic framework, we need to control the quantities
\begin{equation} \label{e:difference_operations}
    \bB_\eps (f_\eps, g_\eps) - \bB_0 (f_0, g_0) \quad \text{and} \quad \tT_\eps (f_\eps, g_\eps, h_\eps) - \tT_0 (f_0, g_0, h_0)
\end{equation}
for bilinear operations $\bB_\eps$ (such as para-products $\prec_\eps$, $\circ_\eps$ and $\succ_\eps$) and tri-linear operations $\tT_\eps$ (such as commutators), in terms of differences between the functions/distributions involved, as well as differences between the operations themselves. 

Here, the key ingredient is a quantitative convergence of the frequency-localised operators to the limiting one (including the convergence of generalised Littlewood-Paley blocks $\Delta_{j,\eps}$ to $\Delta_{j,0}$). This is stated in Theorem~\ref{th:convergence_single_block} below. Its proof uses the Helffer-Sj\"ostrand’s formula for self-adjoint operators\footnote{In an earlier version, we prove the convergence of these operators with a completely different way based on a clustering argument for eigenvalues, together with the convergence of individual eigenvalues by Kenig-LinShen (\cite{KLS2013}) and Weyl's asymptotic formula. The use of Helffer-Sj\"ostrand’s representation in the current version greatly simplifies the proof.}. Based on convergence of these block operators, one can establish effective control for quantities in \eqref{e:difference_operations} when all the inputs $(f_\eps, g_\eps, h_\eps) \in \cC^{\alpha;\eps} \times \cC^{\beta;\eps} \times \cC^{\gamma;\eps}$ for $|\alpha|, |\beta|, |\gamma| <1$. This restriction on exponents is due to that the spaces $\cC^{\alpha;\eps}$ and $\cC^\alpha$ are uniformly (in $\eps$) equivalent only when $\alpha \in (-1,1)$. In this case, one can compare $f_\eps$ and $f_0$ directly. 

On the other hand, the solution theories/ansatzes for all our examples involve components in $\cC^{\alpha;\eps}$ for $|\alpha| > 1$. Hence, one still needs to control \eqref{e:difference_operations} for some $\alpha, \beta, \gamma$ outside $(-1,1)$. It turns out that for $\alpha \in (1,2)$, one can compare $f_\eps \in \cC^{\alpha;\eps}$ with $f_0 \in \cC^{\alpha;0}$ by measuring the $\cC^{\alpha-2}$-norm of $\lL_\eps f_\eps - \lL_0 f_0$, where now $\alpha-2 \in (-1,0)$. Furthermore, the operations $\bB_\eps$ and $\tT_\eps$ are continuous in those norms. 
This allows us to establish very general continuity criterion for para-control related operations, which are further applicable to all the three equations we treat here. 

Finally, assuming the convergence of all the enhanced stochastic objects, we establish quantitative convergence of the oscillatory $\phi^4_3$, $2$D gPAM and KPZ solutions to their homogenised limits. Convergences of the fluxes are also established. The convergence of the stochastic objects is treated in the companion paper \cite{homo_stochastic}.

\subsection{Applications to periodic homogenisation for singular PDEs}

The framework described above allows us to establish periodic homogenisation for a class of singular stochastic PDEs. We give details for three models: dynamical $\phi^4_3$, $2$D gPAM and KPZ. They are respectively given by
\begin{equation} \label{e:phi4_approx_intro}
    \d_t \phi_\eps^{(\delta)} = \lL_\eps \phi_\eps^{(\delta)} - \big( \phi_\eps^{(\delta)} \big)^3 + \xi^{(\delta)} + C_\eps^{(\delta)} \phi_\eps^{(\delta)}\;, \quad (t,x) \in \R^+ \times \T^3\;,
\end{equation}
and
\begin{equation} \label{e:gpam_approx_intro}
    \begin{split}
    \d_t u_\eps^{(\delta)} = \, &\lL_\eps u_\eps^{(\delta)} + \Big( \nabla^\sT u_\eps^{(\delta)} \, F_\eps \, \nabla u_\eps^{(\delta)} - \Tr \big( F_\eps \widetilde{\fC}_\eps^{(\delta)} \big) \, g^2 \big( u_\eps^{(\delta)} \big) \Big)\\
    &+ g \big (u_\eps^{(\delta)} \big) \cdot \big( \xi^{(\delta)} - C_\eps^{(\delta)} g'\big(u_\eps^{(\delta)} \big) \big)\;, \qquad \quad (t,x) \in \R^+ \times \T^2\;,
    \end{split}
\end{equation}
and
\begin{equation} \label{e:kpz_approx_intro}
    \d_t h_\eps^{(\delta)} = \lL_\eps h_\eps^{(\delta)} + g_\eps \big( \d_x h_\eps^{(\delta)} \big)^2 + \xi^{(\delta)} - C_\eps^{(\delta)} - 4\widetilde{C}_\eps^{(\delta)} A_\eps \d_x h_\eps^{(\delta)}\;, \quad (t,x) \in \R^+ \times \T\;.
\end{equation}
In all of above, $\xi^{(\delta)}$ denotes smooth approximation to the spacetime white noise (or $2$D spatial white noise in \eqref{e:gpam_approx_intro}) $\xi$ at scale $\delta$, and $C_\eps^{(\delta)}$ and $\widetilde{C}_\eps^{(\delta)}$ are renormalisation functions that depend on the spatial variable $x$ but independent of time. 

In the KPZ equation \eqref{e:kpz_approx_intro}, we allow an oscillatory coupling coefficient $g_\eps = g(\cdot / \eps)$ for some bounded periodic function $g$. In $2$D gPAM \eqref{e:gpam_approx_intro}, $F_\eps = \big( f_{ij} (\cdot / \eps) \big)_{i,j=1}^{2}$ denotes an oscillatory symmetric $2 \times 2$ matrix with bounded periodic functions $f_{ij} = f_{ji}$, and
\begin{equation*}
    \nabla^\sT u_\eps^{(\delta)} \, F_\eps \nabla u_\eps^{(\delta)} = \sum_{i,j=1}^{2} f_{ij}(\cdot / \eps) \, \d_i u_\eps^{(\delta)} \, \d_j u_\eps^{(\delta)}\;.
\end{equation*}
Also, $\widetilde{\fC}_\eps^{(\delta)}$ is a $2 \times 2$ matrix-valued renormalisation function (independent of time and $F$), and $\Tr \big( F_\eps \widetilde{\fC}_\eps^{(\delta)} \big)$ is the trace of the matrix $F_\eps \widetilde{\fC}_\eps^{(\delta)}$. These are well defined versions of the formal equations \eqref{e:phi4_formal_intro}, \eqref{e:gPAM_formal_intro} and \eqref{e:kpz_formal_intro}. 

Let $\bar{A}$ be the homogenised constant matrix of $A$, and $\lL_0 = \div (\bar{A} \nabla)$. We have the following loose statements regarding the above three equations. Precise statements are given in Theorems~\ref{th:phi4_overall}, \ref{thm:PAM_overall} and ~\ref{thm:kpz_overall} respectively. 

\begin{thm} [Dynamical $\phi^4_3$]
\label{thm:phi43_informal}
Suppose the initial data $\phi_\eps^{(\delta)}(0)$ of \eqref{e:phi4_approx_intro} converges in $\cC^{-\frac{1}{2}-\eta}$ for a sufficiently small positive $\eta$ to a limit $\phi_0(0)$. There exists a choice of $C_\eps^{(\delta)}$ such that the solution $\phi_\eps^{(\delta)}$ converges to a limit $\phi_0$ as $(\eps,\delta) \rightarrow (0,0)$, where $\phi_0$ solves the standard dynamical $\phi^4_3$ model formally given by
\begin{equation*}
    \d_t \phi_0 = \div (\bar{A} \nabla) \phi_0 - \phi_0^3 + \xi + \infty \, \phi_0
\end{equation*}
with initial data $\phi_0(0)$. Furthermore, the flux $A_\eps \nabla \phi_\eps^{(\delta)}$ also converges to $\bar{A} \nabla \phi_0$. 
\end{thm}

We have a similar statement for the $2$D gPAM. 

\begin{thm} [$2$D gPAM] 
\label{thm:gpam_informal}
Suppose the initial data $u_\eps^{(\delta)}(0)$ of \eqref{e:gpam_approx_intro} converges in $\cC^{1-\eta}$ for a sufficiently small positive $\eta$ to a limit $u_0(0)$. Then there exist choices of $C_\eps^{(\delta)}$ and $\widetilde{\fC}_\eps^{(\delta)}$ such that the solution $u_\eps^{(\delta)}$ converges to a limit $u_0$ as $(\eps,\delta) \rightarrow (0,0)$, where $u_0$ solves the standard two dimensional gPAM formally given by
\begin{equation*}
    \d_t u_0 = \div (\bar{A} \nabla u_0) + \big( \nabla^\sT u_0 \, \bar{F} \, \nabla u_0 - \widetilde{\infty} \cdot g^2(u_0) \big) + g(u_0) \big( \xi - \infty \cdot g'(u_0) \big)
\end{equation*}
with initial data $u_0(0)$. Here, $\bar{F}$ is a $2 \times 2$ constant matrix given by
\begin{equation*}
    \bar{F} := \int_{\T^2} \big( \id + (\nabla \chi)(x) \big)^{ \sT } \cdot F(x) \cdot \big( \id + (\nabla \chi)(x) \big) {\rm d}x\;,
\end{equation*}
where $\chi$ is the homogenisation corrector associated to the matrix $A$. Furthermore, the flux $A_\eps \nabla u_\eps^{(\delta)}$ converges to $\bar{A} \nabla u_0$ as well. 
\end{thm}

Finally, we give the statement for the KPZ equation. 

\begin{thm} [The KPZ equation]
\label{thm:kpz_informal}

Suppose the initial data $h_\eps(0)$ of \eqref{e:kpz_approx_intro} converges in $\cC^{\frac{1}{2}-\eta}$ for a sufficiently small positive $\eta$ to a limit $h_0(0)$. Then there exist choices of $C_\eps^{(\delta)}$ and $\widetilde{C}_\eps^{(\delta)}$ such that the solution $h_\eps^{(\delta)}$ converges to a limit $h_0$ as $(\eps,\delta) \rightarrow (0,0)$, where $h_0$ solves the standard KPZ equation formally given by
\begin{equation*}
    \d_t h_0 = \bar{A} \d_x^2 h_0 + \lambda (\d_x h_0)^2 + \xi - \infty
\end{equation*}
with initial data $h_0(0)$ and coupling constant
\begin{equation*}
    \lambda = \bar{A}^2 \int_\T \frac{g(x)}{A^2(x)} {\rm d}x\;.
\end{equation*}
Furthermore, the flux $A_\eps \d_x h_\eps^{(\delta)}$ also converges to $\bar{A} \d_x h_0$. 
\end{thm}

The above statements are purely analytic. In fact, we obtain bounds of the type
\begin{equation*}
    \|v_\eps^{(\delta)} - v_0\| \lesssim \eps^\theta + \|v_\eps(0) - v_0(0)\| + \|\Upsilon_\eps^{(\delta)} - \Upsilon_0\|
\end{equation*}
for proper norms, showing that the solution depends continuously on the initial data, the enhanced noise as well as the homogenisation parameter. In particular, it allows the noise $\Upsilon_\eps^{(\delta)}$ to converge to different limits along different subsequences $(\eps,\delta) \rightarrow (0,0)$. In this case, the main theorems imply that the solution $v_\eps^{(\delta)}$ will also converge to different limits along different subsequences, while each limit is determined by the corresponding limiting enhanced noise. 

In the forthcoming work \cite{homo_stochastic}, we give natural choices of $C_\eps^{(\delta)}$ and $\widetilde{C}_\eps^{(\delta)}$ in the above three examples under which $\Upsilon_\eps^{(\delta)}$ converges to the same limit $\Upsilon_0$ as $(\eps,\delta) \rightarrow (0,0)$ in arbitrary ways. However, we do not claim that these are the only natural choices for renormalisation functions which gives uniformity in $(\eps,\delta)$. In principle, there is an infinite degree of freedom of choices which can possibly lead to different limits of the stochastic objects. See \cite{homo_stoch_rp} for an interesting example in rough paths where the Brownian motion is correlated with the (random) coefficient field and the L\'evy area does not commute with homogenisation.

\subsection{Related works and further perspectives}
\label{sec:perspectives}

Recently, Hairer-Singh (\cite{st_homo_phi42, homo_spde_rs}) developed a general framework for homogenisation for singular SPDEs based on regularity structures. They utilised Green's function estimates in periodic homogenisation problem for general coefficient matrices (see for example \cite{GengShen2020, Shen2018}), and implemented the oscillatory kernels into regularity structures. This allowed them to obtain a general continuity statement of the solution maps to a class of singular SPDEs with respect to enhanced stochastic objects, initial data and oscillatory integration kernels. This approach leads to a systematic framework for general spacetime periodic coefficients $A$ without assuming symmetry, and is also convenient in taking derivatives. The examples covered in \cite{homo_spde_rs} include $2$D g-PAM and $\phi^4_3$, which require suitable ansatz with expansion up to $\cC^{1+}$. 

The approach presented here heavily relies on spectral theory for self-adjoint operators. Hence, compared to \cite{homo_spde_rs}, our approach is restricted to problems with symmetric and time-independent coefficient matrices. At a more philosophical level, the main novelty of this article is the viewpoint to measure regularities adapted to the background geometry of the problem (see also the recent work \cite{2d_YMH} on the development of a covariant para-controlled calculus to study stochastic quantisation for two dimensional Yang-Mills-Higgs). In addition, we also give a complete treatment of the KPZ equation, which requires an ansatz describing local expansion up to $\cC^{\frac{3}{2}+}$. 

The para-product and commutator (for heat kernel) estimates obtained in this article require the restriction $\alpha + \beta > -1$ on the regularities of the two inputs (see Proposition~\ref{pr:paraproducts_gBesov} and Lemma~\ref{le:heat_comm_spacetime} below). They allow to cover dynamical $\phi^4_3$, $2$D g-PAM and KPZ, but not $3$D g-PAM yet. 

The works \cite{Bailleul2016, Bailleul2018, Bailleul2019} are close in spirit to the generalised para-controlled calculus developed in this article. \cite{Bailleul2016} used heat semi-groups to define a generalised version of Besov spaces. The requirements there for the coefficient matrix are similar to Definition~\ref{def:operator_class} below, and hence also includes $A_\eps = A(\cdot/\eps)$ from periodic homogenisation setting. Unfortunately, the para-product estimates there have restrictions on the regularity exponents (namely all exponents should be strictly less than $1$), and hence not applicable to singular SPDE problems whose solution ansatz has a component with regularity strictly above $\cC^{1}$ (including standard examples like $2$D g-PAM, $\phi^4_3$ and KPZ). The approaches in \cite{Bailleul2018, Bailleul2019} have removed this restriction on the regularity exponent, but require higher regularity of the coefficient matrix, and hence are not applicable to homogenisation problems.

\subsection{Organisation of this article}

In Section~\ref{sec:para_general}, we develop uniform bounds for para-product, commutator and heat kernel operations over a class of operators. This class includes in particular the family $\lL_\eps$ of oscillatory operators generated from the same coefficient matrix $A$. In Section~\ref{sec:para_convergence}, we establish convergence properties of these operations in the homogenisation setting. Section~\ref{sec:bounds_derivatives_1d} is devoted to uniform bounds and convergences for operations involving derivatives specific in dimension $d=1$. 

Finally in Sections~\ref{sec:phi4}, ~\ref{sec:gpam} and~\ref{sec:kpz}, we apply the framework to periodic homogenisation problems for dynamical $\phi^4_3$, $2$D gPAM and KPZ. In particular, we show that the solutions are uniform in the oscillation parameter, and that both the solutions and fluxes converge to the homogenised ones if the initial data and enhanced noise converge.

\subsection{Notations and definitions}

\subsection*{General notations}

For functions $f$ on $\T^d$, we let
\begin{equation*}
    \Pi_0 f := \int_{\T^d} f(x) {\rm d}x\;, \qquad \Pi_0^\perp := \id - \Pi_0\;.
\end{equation*}
For $T>0$, function space $\zZ$ and $p \in [1,+\infty]$, we write $L_T^p \zZ := L^p \big( [0,T]; \zZ \big)$ as the spacetime norm for $L^p$-in-time evolution in $\zZ$. For $\sigma \geq 0$, we define the weighted norm 
$L_{\sigma,T}^{\infty} \zZ$ by
\begin{equation} \label{eq:spacetime_norm_weight}
    \|f\|_{L_{\sigma,T}^{\infty} \zZ} := \sup_{t \in [0,T]} \Big( t^{\frac{\sigma}{2}} \|f(t)\|_{\zZ} \Big)\;.
\end{equation}
For $\beta \in (0,1)$, we use $\cC_T^{\beta,\sigma} \zZ$ to denote the norm
\begin{equation} \label{eq:spacetime_Holder_norm_weight}
    \|f\|_{\cC_{\sigma,T}^{\beta} \zZ} := \|f\|_{L_{\sigma,T}^{\infty} \zZ} + \sup_{0<s<t<T} \bigg( s^{\frac{\sigma}{2}} \cdot \frac{\|f(t) - f(s)\|_{\zZ}}{|t-s|^\beta} \bigg)\;.
\end{equation}

\subsection*{Homogenisation quantities}

In all places except in Section~\ref{sec:para_general}, we fix $A: \T^d \rightarrow \R^{d\times d}$ be a symmetric, H\"older continuous and uniformly elliptic coefficient matrix satisfying Assumption~\ref{as:A}. Let $\chi: \T^d \rightarrow \R^d$ be the corrector associated to its homogenisation problem in the sense that $\chi$ is the unique solution to
\begin{equation}\label{def:chi_homo_corrector}
    \div \big( A (\id + \nabla \chi) \big) = 0\;, \qquad \int_{\T^d} \chi = 0\;.
\end{equation}
Let
\begin{equation}\label{def:bar_a_homogenised_matrix}
    \bar{A} := \int_{\T^d} A(y) \big( \id + \nabla \chi (y) \big) {\rm d}y
\end{equation}
be the homogenised matrix. Throughout, $\eps\in \N^{-1}$ always denotes an inverse integer (that is, $\eps = \frac{1}{n}$ for some $n \in \N^{+}$). We write $A_\eps := A(\cdot/\eps)$, and
\begin{equation*}
    \lL_\eps := \div ( A_\eps \nabla )\;, \qquad \lL_0 := \div ( \bar{A} \nabla )\;.
\end{equation*}
Also write
\begin{equation*}
    \Phi := \id + \nabla \chi\;, \qquad \Phi_\eps := \Phi (\cdot / \eps)\;.
\end{equation*}
The Littlewood-Paley operators $\Delta_{j,\eps}$ and $\Delta_{j,0}$ are defined the same way as $\Delta_{j,a}$ in \eqref{eq:operator_LP_block} with $a = A_\eps$ and $a = \bar{A}$ respectively. Similar as \eqref{eq:operator_S}, define the operators $S_{j,\eps}$ and $S_{j,\eps}^{\perp}$ as
\begin{equation} \label{eq:operator_S_oscillatory}
    S_{j,\eps} := \sum_{i \leq j-2} \Delta_{j,\eps}\;, \qquad S_{j,\eps}^{\perp} := \id - S_{j,\eps}\;.
\end{equation}
The para-product operations $\prec_\eps$, $\circ_\eps$ and $\succ_\eps$ are given by
The space $\cC^{\alpha;\eps}$ denotes the $\eps$-dependent Besov space associated with the operator $\lL_\eps$. We will simply write $\|\cdot\|_{\alpha;\eps}$ for $\|\cdot\|_{\cC^{\alpha;\eps}}$. 

We will show in Section~\ref{sec:para_general} below that spaces $\{\cC^{\alpha;\eps}\}_{\eps\in \N^{-1}}$ are equivalent to each other if $\alpha \in (-1,1)$. For $\alpha \in (1,2)$, we compare a function $f_\eps \in \cC^{\alpha;\eps}$ with $f_0 \in \cC^{\alpha;0}$ by
\begin{equation} \label{e:norm_comparison_high_reg}
    |\!|\!| f_\eps; f_0 |\!|\!|_{\alpha;\eps} := |\Pi_0 (f_\eps - f_0)| + \|\lL_\eps f_\eps - \lL_0 f_0\|_{\alpha-2}\;.
\end{equation}
For $\alpha \in (-2,-1)$, we compare a function $g_\eps \in \cC^{\alpha;\eps}$ with $g_0 \in \cC^{\alpha;0}$ by
\begin{equation} \label{e:norm_comparison_low_reg}
    |\!|\!| g_\eps; g_0 |\!|\!|_{\alpha;\eps} := |\Pi_0 (g_\eps - g_0)| + \|\lL_\eps^{-1} \Pi_0^\perp g_\eps - \lL_0^{-1} \Pi_0^\perp g_0\|_{\alpha+2}\;.
\end{equation}

\subsection*{Green functions}

For $\eps \in \N^{-1}$ (including $\eps=0$), write $\iI_\eps = (\d_t - \lL_\eps)^{-1}$ and $Q_\eps$ be its kernel in the sense that
\begin{equation*}
    (\iI_\eps f)(t) := \int_{0}^{t} e^{(t-r) \lL_\eps} f(r) {\rm d}r = \int_{0}^{t} \int_{\T^d} Q_\eps (t-r;x,y) f(r,y) {\rm d}r {\rm dy}\;.
\end{equation*}
Let $G_\eps$ (including $\eps=0$) be the elliptic Green function in the sense that
\begin{equation*}
    (-\lL_\eps^{-1} f)(x) = \int_{\T^d} G_\eps (x,y) f(y) {\rm d}y
\end{equation*}
for $f$ with $\Pi_0 f = 0$. Note that we have the identities
\begin{equation*}
    \begin{split}
    - \lL_\eps^{-1} \Pi_0^\perp = \int_{\R^+} e^{r \lL_\eps} \Pi_0^\perp {\rm d}r\;, \qquad G_\eps = \int_{\R^+} \big( Q_\eps(r; \cdot, \cdot) - 1 \big) {\rm d}r\;.
    \end{split}
\end{equation*}
Write $\rR_\eps^{\p}$ and $(\rR_\eps^{\p})^*$ as
\begin{equation} \label{eq:Green_remainder_parabolic}
    \begin{split}
    \rR_\eps^{\p} (t;x,y) &:= (\nabla_x Q_\eps)(t;x,y) - \Phi_\eps (x) \cdot (\nabla_x Q_0)(t;x,y)\;,\\
    (\rR_\eps^{\p})^* (t;x,y) &:= (\nabla_y Q_\eps)(t;x,y) - \Phi_\eps (y) \cdot (\nabla_y Q_0)(t;x,y)\;.
    \end{split}
\end{equation}
Similarly, we use $\rR_\eps^{\e}$ and $(\rR_\eps^{\e})^*$ to denote the elliptic versions in the sense that
and $\rR_\eps^{\e}$ as
\begin{equation} \label{eq:Green_remainder_elliptic}
    \begin{split}
    \rR_\eps^{\e}(x,y) &:= (\nabla_x G_\eps) (x,y) - \Phi_\eps(x) \cdot (\nabla_x G_0) (x,y)\;,\\
    (\rR_\eps^{\e})^*(x,y) &:= (\nabla_y G_\eps) (x,y) - \Phi_\eps(y) \cdot (\nabla_y G_0) (x,y)\;.
    \end{split}
\end{equation}
Note that here we have $\Phi_\eps (y)$ instead of $\Phi_\eps^*
(y)$ since $A$ is symmetric (and hence $\Phi = \Phi^*$).

\subsection{Comparison with the previous version and statement on use of AI}

The mathematical contents as well as the write up of this article are developed without assistance of AI. A previous version of this article was submitted to and appeared on arXiv in August 2025. That version already contained essentially all the main ideas and statements. 

In this version, there are four main changes compared to the 2025 arXiv version: 
\begin{itemize}
\item we added Theorem~\ref{thm:gpam_informal} on the $2$D gPAM;
\item we extended the oscillatory coupling coefficient in the KPZ equation (in Theorem~\ref{thm:kpz_informal}) to general oscillatory function $g_\eps$; 
\item we simplified the proof of the uniform decay estimates (Proposition~\ref{pr:decoupling_spectrum}); 
\item we also significantly simplified the proof for the convergence of block operators (Theorem~\ref{th:convergence_single_block}). 
\end{itemize}
None of these improvements/extensions used AI tools. We give more detailed descriptions of these changes in Appendix~\ref{sec:comparison}.

\subsection*{Acknowledgements}

We thank Martin Hairer, Harprit Singh and Jinping Zhuge for helpful discussions. In particular, Martin Hairer pointed out several mathematical issues in a much earlier version of the draft. Fixing these issues results in a complete re-writing of the manuscript. 

We learned the Helffer-Sj\"ostrand's formula for self-adjoint operators from Siqi Fu's lectures on spectral theory of the complex Laplacian at Westlake during March-April 2026. This formula enables us to greatly simplify the proof of Theorem~\ref{th:convergence_single_block}. 

W. Xu was supported by Ministry of Science and Technology via the National Key R\&D Program of China (no.2023YFA1010100) and National Science Foundation China via Standard Project Grant (no.8200906145).

\section{Operator-dependent para-controlled calculus -- uniform bounds}
\label{sec:para_general}

\subsection{Set up}

We fix a space dimension $d$. Let $\T^d = (\R / \Z)^d$ be the $d$ dimensional torus of side length $1$. For $\Lambda > 1$, let $\sS_d(\Lambda)$ be the set of matrix-valued functions $a: \T^d \rightarrow \R^{d \times d}$ such that $a = a^\sT$ and
\begin{equation}\label{eq: sS d lambda}
    \Lambda^{-1} |v|^2 \leq v^\sT a v \leq \Lambda |v|^2
\end{equation}
for all $v \in \R^d$. For $a \in \sS_{d}(\Lambda)$, define $\lL_a :=\div( a \nabla )$. Let $Q_a = (\d_t - \lL_a)^{-1}$ be the Green's function. 

We aim to obtain uniform para-differential calculus associated with the operators $\lL_a$ over the following class of coefficient matrices $a$. 

\begin{defn} \label{def:operator_class}
For $M>0$ and $\Lambda>1$, let $\sS_d (\Lambda, M)$ be the 
subset of $\sS_d (\Lambda)$ such that the parabolic Green's function $Q_a$ for $a \in \sS_d(\Lambda, M)$ satisfies the following:
\begin{itemize}
\item Small time Gaussian bound:
    \begin{equation}
    |Q_{a}(t;x,y)|\leq M t^{-\frac{d}{2}}e^{-\frac{|x - y|^2}{M t}}\;; \quad t \leq 1\;.
    \end{equation}
    
\item $(\nabla_x Q_a)(t;x,y)$ and $(\nabla_y Q_a)(t;x,y)$ exist on $\R^+ \times \T^d \times \T^d$. Furthermore, $Q_a$ satisfies the bounds
\begin{equation}
    \begin{split}
    |Q_{a}(t;x,y) - Q_{a}(t;x',y)| &\leq M |x-x'| \, t^{-\frac{d+1}{2}} \, e^{-\frac{|x-y|^{2}}{M t}}\;,\\
    |Q_{a}(t;x,y) - Q_{a}(t;x,y')| &\leq M |y-y'| \, t^{-\frac{d+1}{2}} \, e^{-\frac{|x-y|^{2}}{M t}}\;
    \end{split}
\end{equation}
for $|x-x'| \leq \sqrt{t}$ (respectively $|y-y'| \leq \sqrt{t}$) and $t \leq 1$. 
\end{itemize} 
The above bounds hold uniformly in spacetime variables with $t \leq 1$ and $a \in \sS_{d}(\Lambda, M)$. 
\end{defn}

In the above definition, $M$ is the constant that controls the behaviour of Green functions $Q_a$ on the time interval $[0,1]$. By semi-group property, bounds on the time interval $[0,T]$ for any $T>0$ can be obtained from the bounds on $[0,1]$. The associated constant $M=M_T$ then depends on $T$. 

For $a \in \sS_d(\Lambda, M)$, the operator $\sqrt{-\lL_a}$ has a discrete spectrum
\begin{equation*}
    0 = \lambda_{0,a} < \lambda_{1,a} \leq \cdots \lambda_{n,a} \leq \cdots\;.
\end{equation*}
Each eigenvalue $\lambda_{n,a}$ has an associated eigenfunction $\psi_{n,a}$. The set of eigenfunctions $\{\psi_{n,a}\}_{n \in \Z^+}$ forms an orthonormal basis of $L^2(\T^d)$. 

For any function $F: \R^+ \rightarrow \R$, we define the operator $F(\sqrt{-\lL_a})$ by
\begin{equation} \label{e:functional_operator_form}
    F(\sqrt{-\lL_{a}})f := \sum_{n=0}^{\infty} F(\lambda_{n,a}) \, \scal{f,\psi_{n,a}} \, \psi_{n,a}\;.
\end{equation}
This definition does not depend on the actual choice of eigenfunctions $\psi_{n,a}$'s. The following proposition (\cite[Theorem~3.1]{KP_general_Besov}) is the basis for all our subsequent analysis. 

\begin{prop} \label{pr:kernel_pointwise_bound}
Fix $\Lambda > 1$ and $M>0$. Let $\sS_d(\Lambda, M)$ be the set of coefficient matrices given in Definition~\ref{def:operator_class}. Let $F$ be a smooth even function on $\R$ with compact support in $(-R,R)$ for some $R>1$. For every $k \geq d+1$, there exist $C_k$ and $C_k'$ such that for all $a \in \sS_d(\Lambda, M)$, we have the bounds
\begin{equation} \label{e:kernel_pointwise_bound}
    \left| F \big( \delta \sqrt{-\lL_a} \big)(x,y) \right| \leq C_k \delta^{-d} \Big( 1 + \frac{|x-y|}{\delta} \Big)^{-k}
\end{equation}
uniformly over $x,y \in \T^d$ and $\delta \in (0,1)$, and
\begin{equation} \label{e:kernel_difference_bound}
    \left| F \big( \delta \sqrt{-\lL_a} \big)(x,y) - F \big( \delta \sqrt{-\lL_a} \big)(x,y') \right| \leq C_k' \cdot \frac{|y-y'|}{\delta^{d+1}} \cdot \Big( 1 + \frac{|x-y|}{\delta} \Big)^{-k}
\end{equation}
uniformly over $x,y,y' \in \T^d$ and $\delta \in (0,1)$ with $|y-y'| \leq \delta$. 

As a consequence, $\nabla_x F\big( \delta \sqrt{-\lL_a} \big)(x,y)$ and $\nabla_y F\big( \delta \sqrt{-\lL_a} \big)(x,y)$ exist pointwise, and satisfy the bound
\begin{equation} \label{e:kernel_gradient_bound}
    \left| \nabla_x F \big( \delta \sqrt{-\lL_a} \big)(x,y) \right| + \left| \nabla_y F \big( \delta \sqrt{-\lL_a} \big)(x,y) \right| \leq \frac{C_k'}{\delta^{d+1}} \Big( 1 + \frac{|x-y|}{\delta} \Big)^{-k}\;.
\end{equation}
The constants $C_k$ and $C_k'$ are given by
\begin{equation*}
    C_k := R^d \big((c_1 \, k)^{k} \|F\|_{L^{\infty}} + (c_2 R)^{k} \|F^{(k)}\|_{L^{\infty}}\big)\;, \quad C_k' = c_3 C_k R\;,
\end{equation*}
where $c_1, c_2, c_3$ are constants that depend on the values of $\Lambda$ and $M$ only. 

Moreover, we have the identity
\begin{equation} \label{e:functional_operator_integral}
    \int_{\T^d} F \big( \sqrt{-\lL_{a}} \big)(x,y) {\rm d} y = F(0)\;.
\end{equation}
\end{prop}
\begin{proof}
The bounds \eqref{e:kernel_pointwise_bound} and \eqref{e:kernel_difference_bound} are proven in \cite[Theorem~3.1]{KP_general_Besov} (with $\alpha=1$ in their case). The pointwise existence of $\nabla_y F\big( \delta \sqrt{-\lL_a} \big)(x,y)$ follows from the representation
\begin{equation*}
    F \big( \delta \sqrt{-\lL_a} \big)(x,y) = \int_{\T^d} \big( F \big( \delta \sqrt{-\lL_a} \big) e^{-\delta^2 \lL_a} \big)(x,z) \cdot e^{ \delta^2 \lL_a} (z,y) \, {\rm d} z\;,
\end{equation*}
the assumption that $\nabla_y e^{\delta^2 \lL_a}(z,y)$ exists pointwise for $a \in \sS_d(\Lambda, M)$, the bound \eqref{e:kernel_difference_bound} together with Lebesgue's dominated convergence theorem. The existence of the gradient $\nabla_x F\big( \delta \sqrt{-\lL_a} \big)(x,y)$ can be obtained similarly. The bound \eqref{e:kernel_gradient_bound} for the gradients then follows from \eqref{e:kernel_difference_bound}. 

Finally, applying $f \equiv 1$ to the definition \eqref{e:functional_operator_form}, and using the fact that $\psi_{0,a} \equiv 1$ and the orthogonality of the $\psi_{n,a}$'s, we then deduce the identity \eqref{e:functional_operator_integral}. 
\end{proof}

\begin{rmk}
    In this section, we use $\T^d$ as the base manifold for convenience, although similar estimates for para-products and commutators can also be obtained on compact Riemannian manifolds without boundaries. More general settings can be found in \cite{KP_general_Besov}.
\end{rmk}

As in classical Littlewood-Paley theory, we need the existence of a dyadic partition of unity. Let $\rho$ be an even, smooth function on $\R$ taking values between $[0,1]$ with the further properties that $\supp (\rho) \subset (-1,1)$ and $\rho = 1$ on $[-\frac{3}{4}, \frac{3}{4}]$. 

Define the functions $\{\varphi_j\}_{j \geq -1}$ by
\begin{equation} \label{e:rho_rescaling}
    \varphi_{-1} := \rho\;, \qquad \varphi_0 :=\rho(\cdot/2) - \rho, \qquad \varphi_j := \varphi_0(\cdot / 2^j) \quad \text{for} \; j \geq 1\;.
\end{equation}
The properties of $\rho$ ensure that
\begin{equation} \label{e:dyadic_partition}
    \sum_{j=-1}^{+\infty} \varphi_j = 1\;.
\end{equation}
Furthermore, for each $j \geq 0$, we have $\supp(\varphi_j) \subset ( - 2^{j+1}, - \frac{3}{4} \times 2^j) \cup (\frac{3}{4} \times 2^j, 2^{j+1})$. 

We now introduce the operator dependent (generalised) Besov spaces, as defined in \cite[Theorem~3.1]{KP_general_Besov}. For $a \in \sS_{d}(\Lambda, M)$ and the functions $\varphi_j$ given in \eqref{e:rho_rescaling}, we write
\begin{equation} \label{eq:operator_LP_block}
    \Delta_{j,a} := \varphi_j \big( \sqrt{-\lL_a} \big)\;.
\end{equation}
For $\alpha\in \R$, $p,q\in [1,\infty]$, define the $\bB^{\alpha; a}_{p,q}$-norm by
\begin{equation}\label{eq: gBesov_defn}
    \|f\|_{\bB^{\alpha; a}_{p,q}} := \bigg [\underset{ j \geq -1 }{\sum} \Big(2^{\alpha j} \, \|\Delta_{j,a} f\|_{L^p}\Big)^{q} \bigg]^{1/q}\;. 
\end{equation}
For $\alpha \in \R$, define the $\cC^{\alpha;a}$-norm by
\begin{equation} \label{eq: Holder_Besov_defn}
    \|f\|_{\cC^{\alpha; a}} := \underset{ j \geq -1 }{ \sup } \Big(2^{\alpha j} \, \|\Delta_{j,a} f\|_{L^\infty}\Big)\;. 
\end{equation}
For convenience, we write $\|\cdot\|_{\alpha;a}$ for $\|\cdot\|_{\cC^{\alpha;a}}$. For $j \geq -1$, define the operators $S_{j,a}$ and $S_{j,a}^{\perp}$ by
\begin{equation} \label{eq:operator_S}
    S_{j,a} f:= \sum_{i=-1}^{j-2} \Delta_{i,a} f\;, \qquad S_{j,a}^{\perp} := \id - S_{j,a}\;.
\end{equation}
We then define the para-products $\prec_a$, $\circ_a$ and $\succ_a$ by
\begin{equation}\label{eq: paraproduct_a}
    f \prec_a g := \sum_{j \geq -1} S_{j,a}f \cdot \Delta_{j,a}g\;, \qquad f \circ_a g := \sum_{|i-j| \leq 1} \Delta_{i,a}f \cdot \Delta_{j,a}g\;,
\end{equation}
and $f \succ_a g := g \prec_a f$. Finally, we define
\begin{equation}\label{eq: mathbf Pi}
    \mathbf{\Pi}_a(f,g):=\sum_{j=-1}^{\infty}\sum_{|i-j|\leq 1}^{\infty} \big(a\nabla \Delta_{i,a}f\big)\cdot \nabla \Delta_{j,a} g\;.
\end{equation}
Classical Littlewood-Paley theory takes $a = \id$, and the blocks $\Delta_{j}$ are given by $\Delta_j = \varphi_j \big(\sqrt{-\Delta}\big)$. It has the nice property that for every two functions $f, g$ on $\T^d$ and $i \leq j - 2$, one has
\begin{equation*}
    \supp \big( \Delta_i f \cdot \Delta_j g \big) \subset B(0, 2^{j+2}) \setminus B(0, 2^{j-2})\;.
\end{equation*}
However, this property is no longer true for general $a \in \sS_{d}(\Lambda, M)$, and $\Delta_{i,a} f \cdot \Delta_{j,a}g$ in general spreads over the whole spectrum of $\sqrt{-\lL_a}$. 

Also, in classical settings, one has the identity $\nabla \Delta_j = \Delta_j \nabla$, and for every $\alpha \in \R$, 
\begin{equation*}
    \|\nabla f\|_{\cC^{\alpha-1}} \lesssim \|f\|_{\cC^\alpha}
\end{equation*}
uniformly over all $f \in \cC^\alpha$. But neither of this is true for general $a \in \sS_{d}(\Lambda, M)$. In particular, differentiation does not make the regularity one degree under the norm induced by $a$ uniformly over $\sS_{d}(\Lambda, M)$. These facts make the analysis of the KPZ equation particularly subtle, as it involves the square of a derivative. Section~\ref{sec:bounds_derivatives_1d} is devoted to the effect of differentiations in space dimension one, and Section~\ref{sec:kpz} is devoted to the full analysis of the KPZ equation.

\subsection{Paraproduct estimates}

Let $a \in \sS_d(\Lambda,M)$. The constants appearing in this subsection depend only on $d$, $M$ and $\Lambda$ unless otherwise specified. For $\alpha \in \R$, we denote the classical H\"older-Besov space in $\T^d$ by $\cC^{\alpha}(\T^d)$ as usual. We start with the following useful lemma, which is a consequence of Proposition~\ref{pr:kernel_pointwise_bound}.

\begin{lem}\label{lem:blocks_derivatives_pointwise}
    For $\ell \in \{0,1\}$, every $m \in \Z$ and every $k \in \N$, we have the bound
    \begin{equation*}
        \big| \big( \nabla_x^{\ell} (-\lL_a)^{m} \Delta_{ i , a } \big)( x , y ) \big| + \big| \big( \nabla_y^{\ell} (-\lL_a)^{m} \Delta_{ i , a } \big)( x , y ) \big| \lesssim_m 2^{(2m+\ell)i}\frac{2^{ d i }}{(1+2^{i}|x-y|)^{k}}\;,
    \end{equation*}
    uniformly over the pair $(i,m) \in \big( \N \times \Z \big) \, \cup \, \big( \{-1\} \times \N \big)$. As a consequence, by choosing $k \geq d+1$, we have
    \begin{equation*}
        \sup_{x \in \T^d} \Big( \big\| \big( \nabla_x^\ell (-\lL_a)^{m} \Delta_{i,a}\big)(x,\cdot) \big\|_{L^1} + \big\| \big( \nabla_y^\ell (-\lL_a)^{m} \Delta_{i,a}\big)(x,\cdot) \big\|_{L^1} \Big) \lesssim_m 2^{(2m+\ell)i}
    \end{equation*}
    for $\ell = 0,1$.
    \end{lem}
    \begin{proof}
        We first consider the case $i \geq 0$. Since $\supp(\varphi_0)$ is in the annulus $B(0,2) \setminus B(0,\frac{3}{4})$, for every $m \in \Z$ (including negative values of $m$), the function
        \begin{equation*}
            F_{m}(x) := x^{2m} \varphi_0(x)
        \end{equation*}
        is smooth and has compact support in the same annulus, and $\|F_m^{(k)}\|_{L^\infty} \lesssim 2^{\frac{5|m|}{2}}$, where the proportionality constant depends on $k$ but is independent of $m$. We have the precise relation
        \begin{equation*}
            (-\lL_a)^{m} \Delta_{i,a} = 2^{2m i} F_m \big( 2^{-i} \sqrt{-\lL_a} \big)\;.
        \end{equation*}
        The pointwise bounds for $\ell=0,1$ then follow from Proposition~\ref{pr:kernel_pointwise_bound} and the control on $F_m^{(k)}$. The $L_x^\infty L_y^1$ bound follows from the above pointwise bounds and choosing $k \geq d+1$. 

        For $m \geq 0$, the function $G_m(x) = x^{2m} \rho(x)$ is smooth with compact support, and
        \begin{equation*}
            (-\lL_a)^{m} \Delta_{-1,a} = G_m(\sqrt{-\lL_a})\;.
        \end{equation*}
        Hence the same conclusions hold uniformly over $(i,m) \in \{-1\} \times \N$. 
    \end{proof}

\begin{lem} \label{lem:Holder_equivalence_prelim}
For every $p \in [1,+\infty]$, we have the bounds
\begin{equation*}
    \begin{split}
    \|\Delta_{j,a} \Delta_{k} f\|_{L^p} &\lesssim 2^{-|j-k|} \sum_{i=k-1}^{k+1}  \|\Delta_i f\|_{L^p}\;,\\
    \|\Delta_j \Delta_{k,a} f\|_{L^p} &\lesssim 2^{-|j-k|} \sum_{i=k-1}^{k+1} \|\Delta_{i,a} f\|_{L^p}
    \end{split}
\end{equation*}
uniformly over all $a \in \sS_{d}(\Lambda, M)$, $k,j \geq -1$ and $f \in \sS'(\T^d)$. The proportionality constant is also independent of $p \in [1,+\infty]$. 
\end{lem}
\begin{proof}
Without loss of generality, it suffices to prove that
\begin{equation*}
    \begin{split}
    \|\Delta_{j,a} \Delta_{k} f\|_{L^p}+\|\Delta_j \Delta_{k,a} f\|_{L^p}  \lesssim 2^{-|j-k|} \|f\|_{L^p}
    \end{split}
\end{equation*}
since the following identities hold
\begin{equation*}
    \Delta_{k,a} = \sum_{i=k-1}^{k+1}\Delta_{k,a}\Delta_{i,a}\;,\quad \Delta_{k} = \sum_{i=k-1}^{k+1}\Delta_{k}\Delta_{i}\;.
\end{equation*}
We give details for the first bound. 
Define for $x\in \T^d$, $k\geq -1$,
\begin{equation*}
    V_{k}(x):=\frac{2^{dk}}{(1+2^{k} |x|)^{d+10}}\;.
\end{equation*}
If $|j-k| \leq 1$, applying Lemma~\ref{lem:blocks_derivatives_pointwise} (with $\ell = m = 0$) to the kernel $\Delta_{j,a}$, we have
\begin{equation} \label{e:equivalence_a0_0}
    |(\Delta_{j,a} \Delta_k f)(x)| \lesssim (V_{j}*V_{k}*|f|)(x)\;.
\end{equation}
Using Young's convolution inequality, we have
\begin{equation*}
    \|\Delta_{j,a} \Delta_k f\|_{L^{p}}\lesssim \|f\|_{L^{p}}\;.
\end{equation*}
Hence, it remains to consider the situation $|j-k| \geq 2$.
For $j \geq k+2$, we write
\begin{equation*}
    \Delta_{j,a} \Delta_{k}f(x) = \int_{\T^d}\Delta_{j,a}(x,y)(\Delta_{k} f(y)-\Delta_{k}f(x)) {\rm d}y\;.
\end{equation*}
By Proposition \ref{pr:kernel_pointwise_bound}, we get
\begin{equation*}
   |\Delta_{k} f(y)- \Delta_{k} f(x)|\lesssim  |x-y|\cdot 2^{k}\cdot \Big[(V_{k}*f)(y)+(V_{k}*f)(x)\Big]\;.
\end{equation*}
Plugging this back, we have
\begin{equation*}
    \begin{split}
    |\Delta_{j,a} \Delta_{k}f(x)| \lesssim  \int_{\T^d}V_{j}(x-y)|x-y|\cdot 2^{k}\cdot \Big[(V_{k}*f)(y)+(V_{k}*f)(x)\Big] {\rm d}y\;. 
    \end{split}
\end{equation*}
This together with Young's inequality deduces
\begin{equation} \label{e:equivalence_a0_2}
    \|\Delta_{j,a} \Delta_{k} f\|_{L^p} \lesssim 2^{-(j-k)} \|f\|_{L^{p}}\;, \qquad j \geq k+2\;.
\end{equation}
 If $k \geq j+2$, for every $h\in L^{p'}(\T^d)$ where $p'$ is the H\"older conjugate index of $p$,
\begin{equation*}
    \scal{h,\Delta_{j,a} \Delta_k f} = \scal{\Delta_k \Delta_{j,a}h,  f}\;.
\end{equation*}
Analogously, we have
\begin{equation*}
    \|\Delta_k \Delta_{j,a}h\|_{L^{p'}}\lesssim 2^{-(k-j)}\|h\|_{L^{p'}}\;.
\end{equation*}
Hence, we obtain
\begin{equation} \label{e:equivalence_a0_1}
    |\scal{h,\Delta_{j,a} \Delta_k f}| \lesssim 2^{-(k-j)} \|f\|_{L^p} \cdot \|h\|_{L^{p'}}\;, \qquad k \geq j+2\;.
\end{equation}
Combining \eqref{e:equivalence_a0_0}, \eqref{e:equivalence_a0_1} and \eqref{e:equivalence_a0_2} gives the desired claim for $\Delta_{j,a} \Delta_k f$. The bound for $\Delta_j \Delta_{k,a} f$ can be obtained in essentially the same way. 
\end{proof}

The following result illustrates that for every $\alpha \in (-1,1)$, the spaces $\bB^{\alpha}_{p,q}(\T^{d})$ are all equivalent uniformly over $a \in \sS_d(\Lambda, M)$.

\begin{lem} \label{lem:Holder_equivalence}
For $\alpha \in (-1,1)$ and $p,q \in [1,+\infty]$, we have the bound
    \begin{equation*}
        (1-|\alpha|) \|f\|_{\bB_{p,q}^{\alpha}} \lesssim \|f\|_{\bB_{p,q}^{\alpha;a}} \lesssim \frac{1}{1-|\alpha|} \cdot \|f\|_{\bB_{p,q}^\alpha}
    \end{equation*}
    uniformly over all $a \in \sS_{d}(\Lambda, M)$ and all $f \in \bB_{p,q}^{\alpha;a}$. Furthermore, the proportionality constants are also independent of $\alpha \in (-1,1)$ and $p,q \in [1,+\infty]$, and depend on the dimension $d$ and $\Lambda, M$ only. 
\end{lem}
\begin{proof}
Fix $\alpha \in (-1,1)$ and $p,q \in [1,+\infty]$. By definition of the norm $\bB_{p,q}^{\alpha;a}$ and duality, we have
\begin{equation*}
    \|f\|_{\bB_{p,q}^{\alpha;a}} = \sup \bigg\{  \sum_{j \geq -1} \lambda_j 2^{\alpha j} \|\Delta_{j,a} f\|_{L^p}\,: \; \|\vec{\lambda}\|_{\ell^{q'}} \leq 1 \bigg\}\;,
\end{equation*}
where $\vec{\lambda} = (\lambda_j)_{j \geq -1}$ and $q'$ is the conjugate of $q$. By Lemma~\ref{lem:Holder_equivalence_prelim}, we have
\begin{equation*}
    \|\Delta_{j,a} f\|_{L^p} \leq \sum_{k \geq -1} \|\Delta_{j,a} \Delta_k f\|_{L^p} \lesssim \sum_{k \geq -1} 2^{-|j-k|} \sum_{i=k-1}^{k+1} \|\Delta_i f\|_{L^p}\;.
\end{equation*}
Hence, we deduce
\begin{equation*}
    \sum_{j \geq -1} |\lambda_j| \, 2^{\alpha j} \|\Delta_{j,a} f\|_{L^p} \lesssim \sum_{k \geq -1} \bigg[ \Big( \sum_{i=k-1}^{k+1}  2^{\alpha i} \|\Delta_i f\|_{L^p} \Big) \Big( \sum_{j \geq -1} |\lambda_j| \cdot 2^{-\alpha(k-j) - |k-j|} \Big) \bigg]\;.
\end{equation*}
Applying H\"older inequality to the right hand side above, we get
\begin{equation} \label{e:equivalence_a0_bound}
    \sum_{j \geq -1} |\lambda_j| \, 2^{\alpha j} \|\Delta_{j,a} f\|_{L^p} \lesssim \bigg\| \sum_{i=k-1}^{k+1} 2^{\alpha i} \|\Delta_i f\|_{L^p} \bigg\|_{\ell_k^q} \cdot \bigg\| \sum_{j \geq -1} |\lambda_j| \cdot 2^{-\alpha(k-j) - |k-j|} \bigg\|_{\ell_k^{q'}}\;.
\end{equation}
The first factor on the right hand side above can be directly bounded by
\begin{equation} \label{e:equivalence_a0_factor1}
    \bigg\| \sum_{i=k-1}^{k+1} 2^{\alpha i} \|\Delta_i f\|_{L^p} \bigg\|_{\ell_k^q} \lesssim \bigg( \sum_{k \geq -1} 2^{\alpha k q} \|\Delta_k f\|_{L^p}^q \bigg)^{\frac{1}{q}} \lesssim \|f\|_{\bB_{p,q}^\alpha}\;.
\end{equation}
For the second factor, the quantity inside the norm is a discrete convolution in $j$. Hence, by Young's inequality, we control it by
\begin{equation} \label{e:equivalence_a0_factor2}
    \Big\| \vec{\lambda} * \big( 2^{-\alpha(k-\cdot) - |k-\cdot|} \big) \Big\|_{\ell_k^{q'}} \lesssim \|\vec{\lambda}\|_{\ell^{q'}} \|2^{-\alpha \cdot - |\cdot|}\|_{\ell^1} \lesssim \frac{1}{1-|\alpha|}\;,
\end{equation}
where we have used that $\|\vec{\lambda}\|_{\ell^{q'}} \leq 1$ in the constraint. Plugging \eqref{e:equivalence_a0_factor1} and \eqref{e:equivalence_a0_factor2} back into \eqref{e:equivalence_a0_bound}, we obtain the bound
\begin{equation*}
    \|f\|_{\bB_{p,q}^{\alpha;a}} \lesssim \frac{1}{1-|\alpha|} \cdot \|f\|_{\bB_{p,q}^{\alpha}}\;.
\end{equation*}
The other direction can be obtained in the same way, and hence we omit the details. 
\end{proof}

The following estimate is the key ingredient for most of the properties of the generalised Besov spaces. 

\begin{prop} \label{pr:decoupling_spectrum}
Let $\alpha, \beta \in \R$. We have the bound
\begin{equation*}
    \left\| \Delta_{k,a} \big( \Delta_{i,a} f \cdot \Delta_{j,a} g \big)  \right\|_{L^\infty} \lesssim 2^{-2k + (1-\alpha)i + (1-\beta)j} \, \|f\|_{\alpha;a} \|g\|_{\beta;a}\;, \quad k \geq \max\{i,j\} + 2\;,
\end{equation*}
and
\begin{equation*}
    \left\| \Delta_{k,a} \big( \Delta_{i,a} f \cdot \Delta_{j,a} g \big)  \right\|_{L^\infty} \lesssim 2^{k + (1-\alpha)i - (2+\beta)j} \, \|f\|_{\alpha;a} \|g\|_{\beta;a}\;, \quad j \geq \max\{i,k\} + 2\;.
\end{equation*}
Both proportionality constants are independent of $i,j,k \geq -1$ within the specified range. 
\end{prop}
\begin{proof}
For $k\geq \max\{i,j\}+2$, $\Delta_{k,a} \big( \Delta_{i,a} f \cdot \Delta_{j,a} g \big)(x)$ equals
    \begin{equation*}
         \int_{\T^d}\Delta_{k,a}(x,y)\left[(\Delta_{i,a} f (y)-\Delta_{i,a} f (x))\cdot (\Delta_{j,a} g(y)-\Delta_{j,a}g(x))\right] {\rm d}y\;.
    \end{equation*}
    Since Proposition \ref{pr:kernel_pointwise_bound} gives 
    \begin{equation*}
        |\Delta_{i,a} f (y)-\Delta_{i,a} f (x)|\lesssim |x-y| \cdot  2^{(1-\alpha)i} \cdot \|f\|_{{\alpha;a}}\;,
    \end{equation*}
    and
    \begin{equation*}
        |\Delta_{j,a} g (y)-\Delta_{j,a} g (x)|\lesssim  |x-y| \cdot  2^{(1-\beta)j} \cdot \|g\|_{{\beta;a}}\;,
    \end{equation*}
    then we have
    \begin{equation*}
        \|\Delta_{k,a} \big( \Delta_{i,a} f \cdot \Delta_{j,a} g \big)\|_{L^{\infty}}\lesssim  2^{-2k + (1-\alpha)i + (1-\beta)j} \, \|f\|_{\alpha;a} \|g\|_{\beta;a}\;.
    \end{equation*}
    For $j\geq \max\{i,k\}+2$, for every $h\in L^{1}(\T^d)$, by self-adjointness of $\Delta_{k,a}$ and $\Delta_{j,a} = \Delta_{j,a}\sum_{\ell=j-1}^{j+1}\Delta_{\ell,a}$, we have the following equality. 
    \begin{equation*}
        \scal{\; h, \, \Delta_{k,a} \big( \Delta_{i,a} f \cdot \Delta_{j,a} g \big)\;} = \sum_{\ell=j-1}^{j+1}\scal{\;\Delta_{j,a}\big(\Delta_{k,a} h\cdot \Delta_{i,a} f \big), \;\Delta_{\ell,a}g\;}\;.
    \end{equation*}
    $\Delta_{j,a}\big(\Delta_{k,a} h\cdot \Delta_{i,a} f \big)(x)$ equals
    \begin{equation*}
         \int_{\T^d}\Delta_{j,a}(x,y)\left[(\Delta_{i,a} f (y)-\Delta_{i,a} f (x))\cdot (\Delta_{k,a} h(y)-\Delta_{k,a}h(x))\right] {\rm d}y\;.
    \end{equation*}
    Proposition \ref{pr:kernel_pointwise_bound} shows that 
    \begin{equation*}
        |\Delta_{i,a} f (y)-\Delta_{i,a} f (x)|\lesssim |x-y| \cdot  2^{(1-\alpha)i} \cdot \|f\|_{{\alpha;a}}\;,
    \end{equation*}
    and 
    % \begin{equation*}
    %     |\Delta_{k,a} h (y)-\Delta_{k,a} h (x)|\leq |x-y| \cdot  \int_{\T^d}\frac{2^{(d+1)k}|h(z)|}{(1+2^{k}\min\{|z-y|,|z-x|\})^{d+10}} {\rm d}z \;.
    % \end{equation*}
    \begin{equation*}
        |\Delta_{k,a} h (y)-\Delta_{k,a} h (x)|\lesssim |x-y| \cdot 2^{k} \cdot \big[(V_{k}*|h|)(x)+(V_{k}*|h|)(y) \big]\;,
    \end{equation*}
    where define for $x\in \T^d$, $\ell\geq -1$,
\begin{equation*}
    V_{\ell}(x):=\frac{2^{d\ell}}{(1+2^{\ell} |x|)^{d+10}}\;.
\end{equation*}
    Plugging this back, $\Delta_{j,a}\big(\Delta_{k,a} h\cdot \Delta_{i,a} f \big)(x)$ has the upper bound
    % \begin{equation*}
    %      2^{(1-\alpha)i} \cdot \|f\|_{{\alpha;a}} \int_{\T^d} \frac{2^{dj}|x-y|^2}{(1+2^{j}|x-y|)^{d+10}}\cdot\left(\int_{\T^d}\frac{2^{(d+1)k}|h(z)|}{(1+2^{k}\min\{|z-y|,|z-x|\})^{d+10}} {\rm d}z\right){\rm d}y\;.
    % \end{equation*}
    \begin{equation*}
         2^{(1-\alpha)i+k} \cdot \|f\|_{{\alpha;a}} \int_{\T^d} |x-y|^2 \cdot V_{j}(x-y) \cdot \big[(V_{k}*|h|)(x)+(V_{k}*|h|)(y) \big]{\rm d}y\;.
    \end{equation*}
    This yields
    \begin{equation*}
        \|\Delta_{j,a}\big(\Delta_{k,a} h\cdot \Delta_{i,a} f \big)\|_{L^{1}}\lesssim 2^{-2 j+(1-\alpha) i+k} \|h\|_{L^1}\|f\|_{\alpha;a}\;.
    \end{equation*}
    This together with the bound for $g$ gives the second bound. This completes the proof.
\end{proof}

\begin{cor} \label{cor:decoupling_spectrum_S}
Let $\alpha < 1$ and $\beta \in \R$. Recall the operator $S_{j,a}$ from \eqref{eq:operator_S}. We have the bound
\begin{equation*}
    \left\| \Delta_{k,a} \big( S_{j,a} f \cdot \Delta_{j,a} g \big) \right\|_{L^\infty} \lesssim 2^{-2k + (2-\alpha-\beta)j} \|f\|_{\alpha;a} \|g\|_{\beta;a}\;, \quad k \geq j + 2\;,
\end{equation*}
and the bound
\begin{equation*}
    \left\| \Delta_{k,a} \big( S_{j,a} f \cdot \Delta_{j,a} g \big) \right\|_{L^\infty} \lesssim 2^{k-(1+\alpha+\beta)j} \|f\|_{\alpha;a} \|g\|_{\beta;a}\;, \quad j \geq k + 2\;.
\end{equation*}
Both proportionality constants are uniform over $j,k \geq -1$ within the specified range. 
\end{cor}
\begin{proof}
We have
\begin{equation*}
    \Delta_{k,a} \big( S_{j,a} f \cdot \Delta_{j,a} g \big) = \sum_{i \leq j - 2} \Delta_{k,a} \big( \Delta_{i,a} f \cdot \Delta_{j,a} g \big)\;.
\end{equation*}
The claims then follow directly from the two bounds in Propositions~\ref{pr:decoupling_spectrum} and that $\alpha<1$. 
\end{proof}

% By virtue of the method in Lemma \ref{pr:decoupling_spectrum}, with the results of  Corollary~\ref{cor: L p estimate for high order derivatives} in hand, we have $L^{p}$ estimates for coupled blocks below.
% \begin{cor}\label{le: decoupling by analyticity for different frequencies for L p estimates}
% $\alpha,\beta\in \R$, $r,p,q\geq 1$ and $\frac{1}{r}=\frac{1}{p}+\frac{1}{q}$, $\exists C>0$, $\forall j,k\leq i-4$, 
%     \begin{equation*}
%          \|\Delta_{ i , a }( \Delta_{ j , a } f \Delta_{ k , a } g )\|_{L^{r}(\T^d)}\leq C 2^{  -2i+j+k}\sum_{l=j-1}^{j+1}\|\Delta_{l,a} f\|_{L^{p}(\T^d)}\sum_{m=k-1}^{k+1}\|\Delta_{m,a} g\|_{L^{q}(\T^d)},
%     \end{equation*}
% and $\exists C>0$, $\forall i,j\leq k-4$,
%     \begin{equation*}
%          \|\Delta_{ i , a }( \Delta_{ j , a } f \Delta_{ k , a } g )\|_{L^{r}(\T^d)}\leq C 2^{ -2k+i+j}\sum_{l=j-1}^{j+1}\|\Delta_{l,a} f\|_{L^{p}(\T^d)}\sum_{m=k-1}^{k+1}\|\Delta_{m,a} g\|_{L^{q}(\T^d)}.
%     \end{equation*}
% \end{cor}

\begin{prop} \label{pr:paraproducts_gBesov}
    We have the following generalised paraproduct estimates:
    \begin{equation} \label{e:bony_generalised_low_high}
        \|f \prec_a g\|_{\beta;a} \lesssim \big( \|f\|_{\alpha;a} + \|f\|_{L^\infty} \big) \, \|g\|_{\beta;a}\;, \quad \alpha \in [0,1) \; \; \text{and} \; \; \alpha + \beta \in (-1,2)\;,
    \end{equation}
    and
    \begin{equation} \label{e:bony_generalised_resonance}
        \|f \circ_a g\|_{\alpha+\beta;a} \lesssim \|f\|_{\alpha;a} \, \|g\|_{\beta;a}\;, \quad \alpha+\beta \in (0,2)\;,
    \end{equation}
    and finally
    \begin{equation} \label{e:bony_generalised_high_low}
        \|f \succ_a g\|_{\alpha+\beta;a} \lesssim \|f\|_{\alpha;a} \, \|g\|_{\beta;a}\;, \quad \beta<0 \; \; \text{and} \; \; \alpha+\beta \in (-1,2)\;.
    \end{equation}
    All proportionality constants are uniform for $a \in \sS_{d}(\Lambda, M)$. 
\end{prop}
\begin{proof}
    We start with \eqref{e:bony_generalised_low_high}. By definition, we have
    \begin{equation*}
        \Delta_{k,a} (f \prec_a g) = \sum_{j \geq -1} \Delta_{k,a} \big( S_{j,a} f \cdot \Delta_{j,a} g \big)\;.
    \end{equation*}
    We split the sum into three disjoint parts: $j \leq k-2$, $k-1 \leq j \leq k+1$ and $j \geq k+2$. For the first one, by Corollary~\ref{cor:decoupling_spectrum_S} and that $\alpha<1$ and $\alpha+\beta<2$, we have
    \begin{equation*}
        \begin{split}
        \sum_{j \leq k-2} \|\Delta_{k,a} \big( S_{j,a} f \cdot\, \Delta_{j,a} g \big)\|_{L^\infty} \lesssim \sum_{j \leq k-2}2^{-2k +(2-\alpha-\beta)j}\|f\|_{\alpha;a}\|g\|_{\beta;a}\lesssim 2^{-(\alpha + \beta) k} \|f\|_{\alpha;a} \|g\|_{\beta;a}\;.
        \end{split}
    \end{equation*}
    For the second one, we have
    \begin{equation*}
        \sum_{j=k-1}^{k+1} \|\Delta_{k,a} \big( S_{j,a} f \cdot\, \Delta_{j,a} g \big)\|_{L^\infty} \lesssim \sum_{j=k-1}^{k+1} \|S_{j,a} f\|_{L^\infty} \|\Delta_{j,a} g\|_{L^\infty} \lesssim 2^{-\beta k} \|f\|_{L^\infty} \|g\|_{\beta;a}\;,
    \end{equation*}
    where we have used that $\|\Delta_{k,a}\|_{L^\infty \rightarrow L^\infty} \lesssim 1$ and $\|S_{j,a}\|_{L^\infty \rightarrow L^\infty} \lesssim 1$. Finally, by Corollary~\ref{cor:decoupling_spectrum_S} and the assumption $\alpha+\beta>-1$, we have
    \begin{equation*}
        \begin{split}
        \sum_{j \geq k+2} \|\Delta_{k,a} \big( S_{j,a} f \cdot\, \Delta_{j,a} g \big)\|_{L^\infty} \lesssim \sum_{j \geq k+2} 2^{k -(1+\alpha+\beta)j}\|f\|_{\alpha;a}\|g\|_{\beta;a}\lesssim 2^{-(\alpha + \beta) k} \|f\|_{\alpha;a} \|g\|_{\beta;a}\;.
        \end{split}
    \end{equation*}
    Combining the above three bounds, we deduce that
    \begin{equation} \label{e:bony_generalised_low_high_other}
        \|\Delta_{k,a} (f \prec_a g)\|_{L^\infty} \lesssim 2^{-\beta k} \big( 2^{-\alpha k} \|f\|_{\alpha;a} + \|f\|_{L^\infty} \big) \|g\|_{\beta;a}\;,
    \end{equation}
    which implies \eqref{e:bony_generalised_low_high}. 

    For the bound \eqref{e:bony_generalised_resonance} concerning the resonance product, we have the expression
    \begin{equation*}
        \Delta_{k,a} (f \circ_a g) = \sum_{|i-j| \leq 1} \Delta_{k,a} \big( \Delta_{i,a} f \cdot \Delta_{j,a} g \big)\;.
    \end{equation*}
    We further split the set $\{|i-j| \leq 1\}$ into intersections with the ranges $\{j \leq k-3\}$ and $\{j \geq k-2\}$. For the first one, applying the first bound in Proposition~\ref{pr:decoupling_spectrum} and using the assumption $\alpha+\beta<2$, we have
    \begin{equation*}
        \begin{split}
        \sum_{\stackrel{|i-j| \leq 1}{j \leq k-3}} \left\| \Delta_{k,a} \big( \Delta_{i,a} f \cdot \Delta_{j,a} g \big) \right\|_{L^\infty} &\lesssim \sum_{j \leq k-3} 2^{-2k + (2-\alpha-\beta)j} \|f\|_{\alpha;a} \|g\|_{\beta;a}\\
        &\lesssim 2^{-(\alpha+\beta)k} \|f\|_{\alpha;a} \|g\|_{\beta;a}\;.
        \end{split}
    \end{equation*}
    For the second one, we directly use $\|\Delta_{k,a}\|_{L^\infty \rightarrow L^\infty} \lesssim 1$, the definitions of $\|f\|_{\alpha;a}$ and $\|g\|_{\beta;a}$ as well as the assumption $\alpha+\beta>0$ to get
    \begin{equation*}
        \begin{split}
        \sum_{\stackrel{|i-j| \leq 1}{j \geq k-2}} \left\| \Delta_{k,a} \big( \Delta_{i,a} f \cdot \Delta_{j,a} g \big) \right\|_{L^\infty} &\lesssim \sum_{\stackrel{|i-j| \leq 1}{j \geq k-2}} 2^{-\alpha i - \beta j} \|f\|_{\alpha;a} \|g\|_{\beta;a}\\
        &\lesssim 2^{-(\alpha+\beta)k} \|f\|_{\alpha;a} \|g\|_{\beta;a}\;.
        \end{split}
    \end{equation*}
    The above two bounds together give
    \begin{equation*}
        \left\| \Delta_{k,a} (f \circ_a g) \right\|_{L^\infty} \lesssim 2^{-(\alpha+\beta) k} \|f\|_{\alpha;a} \|g\|_{\beta;a}\;,
    \end{equation*}
    which implies the bound \eqref{e:bony_generalised_resonance}. 

    For the final bound \eqref{e:bony_generalised_high_low}, we have the expression
    \begin{equation*}
        \Delta_{k,a} (f \succ_a g) = \sum_{j \geq -1} \Delta_{k,a} \big( \Delta_{j,a} f \cdot S_{j,a} g \big)\;.
    \end{equation*}
    Similar as the proof for \eqref{e:bony_generalised_low_high}, we decompose the sum into $j \leq k-2$, $k-1 \leq j \leq k + 1$ and $j \geq k+ 2$. Derivation of the bounds for the sums over $j \leq k-2$ and $j \geq k+2$ proceeds in exactly the same way as their counterparts in \eqref{e:bony_generalised_low_high}. For the range $|j-k| \leq 1$, we have
    \begin{equation*}
        \sum_{j=k-1}^{k+1} \left\| \Delta_{k,a} \big( \Delta_{j,a} f \cdot S_{j,a} g \big) \right\|_{L^\infty} \leq \sum_{j=k-1}^{k+1} \|\Delta_{j,a} f\|_{L^\infty} \|S_{j,a} g\|_{L^\infty} \lesssim 2^{-(\alpha+\beta) k} \|f\|_{\alpha;a} \|g\|_{\beta;a}\;,
    \end{equation*}
    where we have used the assumption $\beta<0$ to get $\|S_{j,a} g\|_{L^\infty} \lesssim 2^{-\beta j} \|g\|_{\beta;a}$. This gives the desired bound \eqref{e:bony_generalised_high_low} and hence completes the proof of the proposition. 
\end{proof}

\begin{rmk} \label{rm:para_bony}
Note that \eqref{e:bony_generalised_low_high} is an extension of the classical paraproduct estimate. Indeed, since $\|f\|_{0;a} \lesssim \|f\|_{L^\infty}$, taking $\alpha=0$ in \eqref{e:bony_generalised_low_high} then gives the same form of the classical Bony's estimate with the exponent range $\beta \in (-1,2)$. In the dynamical $\phi^4_3$ case (and also $2$D g-PAM), we need to deal with the situation where the ``high frequency" component of the para-product has regularity (measured by the operator $\lL_{\eps}$) below $-1$. In this case, we need the low frequency part to have higher regularity. This is encoded in $\|f\|_{\alpha;a}$ with $\alpha \in (0,1)$. 

On the other hand, the restriction $\alpha+\beta>-1$ in \eqref{e:bony_generalised_low_high} and \eqref{e:bony_generalised_high_low} excludes the application to $3$D gPAM for the moment. 
\end{rmk}

For $ f$, $g\in C^{\infty}(\T^d)$, define
\begin{equation} \label{e:operator_Pa}
    \textbf{P}_{a}(f,g) :=  \sum_{j=-1}^{\infty} a \nabla S_{j,a} f \cdot \nabla \Delta_{j,a} g\;.
\end{equation}
The following lemma says that for regularity exponents within a certain range, the two $\nabla$ operators in $\P_a (f,g)$ indeed give the desired effects. 

\begin{lem} \label{le:prec_alternative_bound}
Let $\alpha \in (0,1)$ and $\beta > 0$ be such that $\alpha + \beta \in (1,2)$. Then, we have
    \begin{equation*}
        \|\mathbf{P}_{a}( f , g )\|_{\alpha+\beta-2;a}\lesssim \|f\|_{\alpha;a}\|g\|_{ \beta ; a }\;.
    \end{equation*}
\end{lem}
\begin{proof}
It suffices to show that
\begin{equation*}
    \sum_{j \geq -1} \left\| \Delta_{k,a} \big( a \nabla S_{j,a} f \cdot \nabla \Delta_{j,a} g \big) \right\|_{L^\infty} \lesssim 2^{(2-\alpha-\beta) k} \, \|f\|_{\alpha;a} \|g\|_{\beta;a}\;.
\end{equation*}
We split the sum into $j \leq k + 2$ and $j \geq k+3$. For $j \leq k+2$, we use the bound
\begin{equation*}
    \left\| \Delta_{k,a} \big( a \nabla S_{j,a} f \cdot \nabla \Delta_{j,a} g \big) \right\|_{L^\infty} \lesssim \|\nabla S_{j,a} f\|_{L^\infty} \|\nabla \Delta_{j,a} g\|_{L^\infty} \lesssim 2^{(2-\alpha-\beta) j} \|f\|_{\alpha;a} \|g\|_{\beta;a}\;,
\end{equation*}
where in the last bound follows from Lemma~\ref{lem:blocks_derivatives_pointwise} on the behaviour of the kernel $\nabla \Delta_{j,a}$ (with $\ell=1$ and $m=0$ in that lemma) and that $\alpha < 1$. Noting that $\alpha+\beta<2$, summing over $j \leq k + 2$ gives
\begin{equation*}
    \sum_{j=-1}^{k+2} \left\| \Delta_{k,a} \big( a \nabla S_{j,a} f \cdot \nabla \Delta_{j,a} g \big) \right\|_{L^\infty} \lesssim 2^{(2-\alpha-\beta)k} \|f\|_{\alpha;a} \|g\|_{\beta;a}\;.
\end{equation*}
For $j \geq k+3$, we write
\begin{equation*}
    a \nabla S_{j,a} f \cdot \nabla \Delta_{j,a}g = \div \big(\, a \nabla S_{j,a} f\, \cdot\, \Delta_{j,a} g \,\big) - S_{j,a} (\lL_a f) \cdot \Delta_{j,a} g\;.
\end{equation*}
For the first part, we have
\begin{equation*}
    \begin{split}
    \left\| \Delta_{k,a} \div \big( a \nabla S_{j,a} f \cdot \Delta_{j,a} g \big) \right\|_{L^\infty} &= \left\| \int_{\T^d} \nabla_y \Delta_{k,a} (\cdot,y) \cdot \big( a \nabla S_{j,a} f \cdot \Delta_{j,a} g \big)(y) {\rm d}y \right\|_{L^\infty}\\
    &\lesssim 2^{k+(1-\alpha-\beta)j} \|f\|_{\alpha;a} \|g\|_{\beta;a}\;,
    \end{split}
\end{equation*}
where in the last bound we have again used Lemma~\ref{lem:blocks_derivatives_pointwise} and that $\alpha<1$. For the second part, using the second estimate in Proposition~\ref{pr:decoupling_spectrum} and $\alpha<1$, we have
\begin{equation*}
    \begin{split}
    \|\Delta_{k,a} \big( S_{j,a} (\lL_a f) \cdot \Delta_{j,a} g \big)\|_{L^\infty} &\lesssim 2^{k-(1+\alpha-2+\beta) j} \|\lL_a f\|_{\alpha-2;a} \|g\|_{\beta;a}\\
    &\lesssim 2^{k+(1-\alpha-\beta)j} \|f\|_{\alpha;a} \|g\|_{\beta;a}\;.
    \end{split}
\end{equation*}
Summing the above two bounds over $j \geq k+3$ and using the assumption $\alpha+\beta>1$, we get
\begin{equation*}
    \sum_{j=k+3}^{+\infty} \left\| \Delta_{k,a} \big( a \nabla S_{j,a} f \cdot \nabla \Delta_{j,a} g \big) \right\|_{L^\infty} \lesssim 2^{(2-\alpha-\beta) k} \|f\|_{\alpha;a} \|g\|_{\beta;a}\;.
\end{equation*}
This completes the proof. 
\end{proof}

% Similarly, we have the following estimates for para-products of gBesov. The norm of $\bB_{p,q}^{s}$ is sometimes used in SPDE. For example, $\bB_{p,\infty}^{\alpha}$ is used in \cite[Section 6]{GubinelliKPZreloaded}. 
% \begin{lem}\label{lem: paraproduct for gBesov s p q}
% $\forall p,q,p_1,q_1,p_2,q_2\in [1,+\infty]$, $\forall s,s_1,s_2\in \R$ satisfy
%     \begin{equation*}
%         \frac{1}{p}=\frac{1}{p_1}+\frac{1}{p_2}; \quad \frac{1}{q}=\frac{1}{q_1}+\frac{1}{q_2}; \quad s_{1}+s_{2}=s,
%     \end{equation*}
%     we have
%     \begin{equation*}
%         \|f\prec_{a}g\|_{\bB_{p,q}^{s,a}}\lesssim\|f\|_{\bB_{p_{1},q_1}^{s_1,a}}\|g\|_{\bB_{p_2,q_2}^{s_2,a}}
%     \end{equation*}
%     and
%     \begin{equation*}
%         \|f\prec_{a}g\|_{\bB_{p,q}^{s,a}}\lesssim\|f\|_{\bB_{p_{1},q_1}^{s_1,a}}\|g\|_{\bB_{p_2,q_2}^{s_2,a}}
%     \end{equation*}
%     provided that $s_1<0$.
%     Moreover,
%     \begin{equation*}
%         \|f\circ_{a}g\|_{\bB_{p,q}^{s,a}}\lesssim\|f\|_{\bB_{p_{1},q_1}^{s_1,a}}\|g\|_{\bB_{p_2,q_2}^{s_2,a}}
%     \end{equation*}
%     provided that $s>0$.
% \end{lem}

\subsection{Commutator and para-linearisation}

The aim of this subsection is to establish uniform estimates for the operator $\Com_a$ and $\PL_{a}^F$ given respectively by
\begin{equation} \label{e:commutator_defn}
    \Com_a (f; g; h) := (f \prec_a g) \circ_a h - f (g \circ_a h)\;,
\end{equation}
and
\begin{equation} \label{e:para_linear_defn}
    \PL_a^F (u) := F(u) - F'(u) \prec_a u\;.
\end{equation}
An essential role is played by the operator
\begin{equation} \label{e:R_op_defn}
    R_{j,a}(f,g) := \Delta_{j,a} (f \prec_a g) - f \, \Delta_{j,a} g\;.
\end{equation}
We start with the following lemma for a commutator involving $\Delta_{j,a}$. 

\begin{lem}\label{le:block_multiplication_comm}
    For $\alpha\in (0,1)$, we have 
    \begin{equation*}
        \big\| [\Delta_{j,a},f] \, g \big\|_{L^\infty}\lesssim 2^{-\alpha j}\|f\|_{\dot{\cC}^{\alpha}} \|g\|_{L^{\infty}}\;.
    \end{equation*}
    Moreover,
    \begin{equation*}
        \big\| [\Delta_{j,a},f] \, g \big\|_{L^\infty}\lesssim 2^{- j}\|\nabla f\|_{L^{\infty}} \|g\|_{L^{\infty}}\;.
    \end{equation*}
    where $[\Delta_{j,a},f] \, g := \Delta_{j,a}(f g) - f\Delta_{j,a}g$.
\end{lem}
\begin{proof} By definition, we have
    \begin{equation*}
        [\Delta_{j,a},f] \, g = \int_{\T^d}\Delta_{j,a}(x,y) \, \big(f(y)-f(x) \big) g(y) {\rm d}y\;.
    \end{equation*}
    The desired bounds then follow from the pointwise bound for $\Delta_{j,a}(x,y)$ in Lemma~\ref{lem:blocks_derivatives_pointwise}. 
\end{proof}

We have the following lemma. 

\begin{lem} \label{le:R_Linfty_bound}
    Let $\alpha \in (0,1)$ and $\beta \in \R$ be such that $\alpha + \beta \in (-1,2)$. We have the bound
    \begin{equation*}
        \|R_{j,a}(f,g)\|_{L^\infty} \lesssim 2^{-(\alpha+\beta) j} \|f\|_{\alpha;a} \|g\|_{\beta;a}\;.
    \end{equation*}
\end{lem}
\begin{proof}
We have the expression
\begin{equation*}
    R_{j,a} ( f , g ) = \sum_{|p-j|\leq 1}\big([\Delta_{j,a}, S_{p,a}f]\Delta_{p,a}g -S_{p,a}^{\perp}f\cdot\Delta_{j,a}\Delta_{p,a}g\big)+ \sum_{|p-j| \geq 2} \Delta_{j,a}( S_{p,a}f\,\cdot\,\Delta_{p,a}g\,) \;.
\end{equation*}
We first treat the two terms in the first sum above. For the first one, by Lemma~\ref{le:block_multiplication_comm} and the assumption $\alpha <1$, we have
\begin{equation*}
    \|[\Delta_{j,a}, S_{p,a}f]\Delta_{p,a}g\|_{L^{\infty}}\lesssim 2^{-j}\|\nabla S_{p,a}f\|_{L^{\infty}}\|\Delta_{p,a}g\|_{L^{\infty}}\lesssim 2^{-j+(1-\alpha-\beta)p}\|f\|_{\alpha;a}\|g\|_{\beta;a}\;.
\end{equation*}
For the second one, since $\alpha>0$, we have
\begin{equation*}
    \|S_{p,a}^{\perp}f\cdot\Delta_{j,a}\Delta_{p,a}g\|_{L^{\infty}} \lesssim \|S_{p,a}^{\perp} f\|_{L^\infty} \|\Delta_{p,a} g\|_{L^\infty}  \lesssim 2^{-(\alpha +\beta) p}\|f\|_{\alpha;a}\|g\|_{\beta;a}\;.
\end{equation*}
Summing the above two bounds over $j-1 \leq p \leq j+1$ gives
\begin{equation}\label{eq: first term in le:R_Linfty_bound}
   \sum_{|p-j| \leq 1} \left\| [\Delta_{j,a}, S_{p,a}f] \Delta_{p,a} g - S_{p,a}^{\perp}f\cdot\Delta_{j,a}\Delta_{p,a}g \right\|_{L^{\infty}}\lesssim 2^{-(\alpha+\beta)j}\|f\|_{\alpha;a}\|g\|_{\beta;a}\;.
\end{equation}
For the other term on the right hand side of the expression of $R_{j,a}(f,g)$, by Corollary~\ref{cor:decoupling_spectrum_S} and the assumption $\alpha<1$, we have the bounds
\begin{equation*}
    \begin{split}
    \left\| \Delta_{j,a} \big( S_{p,a}f \cdot \Delta_{p,a} g \big) \right\|_{L^{\infty}} &\lesssim 2^{j-(1+\alpha+\beta)p} \, \|f\|_{\alpha;a} \|g\|_{\beta;a}\;, \qquad \quad p \geq j+2\;,\\
    \left\| \Delta_{j,a} \big( S_{p,a}f \cdot \Delta_{p,a}g \big) \right\|_{L^{\infty}} &\lesssim 2^{-2j + (2-\alpha-\beta)p} \, \|f\|_{\alpha;a} \|g\|_{\beta;a}\;, \qquad p \leq j-2\;.
    \end{split}
\end{equation*}
Summing over $p>j+2$ and $p<j-2$ respectively and using the assumption $\alpha + \beta \in (-1,2)$, we get
\begin{equation}\label{eq: second term in le:R_Linfty_bound}
    \sum_{|p-j|>2} \left\| \Delta_{j,a} \big( S_{p,a}f \cdot \Delta_{p,a}g \big) \right\|_{L^{\infty}} \lesssim 2^{-(\alpha+\beta)j} \|f\|_{\alpha;a} \|g\|_{\beta;a}\;.
\end{equation}
The desired result follows from \eqref{eq: first term in le:R_Linfty_bound} and \eqref{eq: second term in le:R_Linfty_bound}.
\end{proof}

\begin{lem} \label{le: estimate for R}
For $\alpha \in (0,1)$ and $\beta \in \R$ such that $\alpha+\beta\in(-1,2)$, we have 
    \begin{equation*}
        \|\Delta_{k,a} \big(R_{j,a} ( f , g )\big)\|_{ L^{ \infty } }\lesssim 2^{ - \alpha \max\{ j, k\}-\beta j } \|f\|_{{\alpha;a}}\|g\|_{{\beta;a}}\;.
    \end{equation*}
\end{lem}
\begin{proof}
For $k\geq j+ 2$, we have $\Delta_{k,a} \Delta_{j,a} = 0$. This, together with Lemma~\ref{le:block_multiplication_comm}, gives
     \begin{equation*}
         \|\Delta_{k,a} \big(R_{j,a} ( f , g )\big)\|_{L^{\infty}} = \|[f,\Delta_{k,a}]\Delta_{j,a}g\|_{L^{\infty}}\lesssim 2^{-k\alpha-j\beta}\|f\|_{{\alpha;a}}\|g\|_{{\beta;a}}\;, \quad k \geq j+2\;.
     \end{equation*}
     For $k \leq j+1$, the bound follows directly from Lemma~\ref{le:R_Linfty_bound}. This completes the proof.
\end{proof}

\begin{lem} \label{le:com_block}
For every $\alpha \in (0,1)$, $\beta, \gamma \in \R$ such that $\alpha+\beta\in(-1,2)$, we have
    \begin{equation*}
        \|\Delta_{k,a} \big( R_{j,a}(f,g) \cdot \Delta_{i,a}h \big)\|_{L^\infty} \lesssim  2^{ - \alpha \max\{ j, k\}-(\beta+\gamma)j} \, \|f\|_{{\alpha;a}}\|g\|_{{\beta;a}}\|h\|_{\gamma;a}
    \end{equation*}
    uniformly over $i,j,k \geq -1$ with the constraint $|i-j| \leq 1$. 
\end{lem}
\begin{proof}
    For $k\leq j+ 2$, it follows from Lemma~\ref{le:R_Linfty_bound} directly that
    \begin{equation*}
        \|\Delta_{k,a} \big( R_{j,a}(f,g) \cdot \Delta_{i,a}h \big)\|_{L^\infty} \lesssim \|R_{j,a}(f,g)\|_{L^\infty} \|\Delta_{i,a} h\|_{L^\infty} \lesssim 2^{-(\alpha+\beta+\gamma) j} \|f\|_{\alpha;a} \|g\|_{\beta;a} \|h\|_{\gamma;a}\;.
    \end{equation*}
    It remains to consider the situation $k \geq j + 3$. We have
    \begin{equation*}
        \Delta_{k,a} \big( R_{j,a}(f,g) \cdot \Delta_{i,a}h \big) = \sum_{p \geq -1} \Delta_{k,a} \Big( \Delta_{p,a} \big(R_{j,a}(f,g)\big) \cdot \Delta_{i,a}h \Big)\;.
    \end{equation*}    
    If $p \geq k - 1$, then $p \geq j + 2$, and hence one has the bound
    \begin{equation*}
    \begin{split}
    \left\| \Delta_{k,a} \big( \Delta_{p,a} \big( R_{j,a}(f,g) \big) \cdot \Delta_{i,a} h \big) \right\|_{L^\infty} &\lesssim \|\Delta_{p,a} \big( R_{j,a}(f,g) \big)\|_{L^\infty} \|\Delta_{i,a} h\|_{L^\infty}\\
    &\lesssim 2^{-\alpha p - (\beta+\gamma)j} \|f\|_{\alpha;a} \|g\|_{\beta;a} \|h\|_{\gamma;a}\;,
    \end{split}
    \end{equation*}
    where we used Lemma~\ref{le: estimate for R} in the second inequality. If $p \leq k-2$, then by the first bound in Proposition~\ref{pr:decoupling_spectrum} and Lemma~\ref{le:R_Linfty_bound}, we have
    \begin{equation*}
        \begin{split}
        \left\| \Delta_{k,a} \big( \Delta_{p,a} \big( R_{j,a}(f,g) \big) \cdot \Delta_{i,a} h \big) \right\|_{L^\infty} &\lesssim 2^{-2k + p + (1-\gamma)j} \| R_{j,a}(f,g) \|_{L^\infty} \|h\|_{\gamma;a}\\
        &\lesssim 2^{-2k+p} \cdot 2^{(1-\alpha-\beta-\gamma) j} \|f\|_{\alpha;a} \|g\|_{\beta;a} \|h\|_{\gamma;a}\;.
        \end{split}
    \end{equation*}
    Summing the above two bounds over $p \geq k-1$ and $p \leq k-2$ respectively, we get
    \begin{equation*}
        \|\Delta_{k,a} \big( R_{j,a}(f,g) \cdot \Delta_{i,a}h \big)\|_{L^\infty} \lesssim 2^{-\alpha k - (\beta+\gamma)j} \|f\|_{\alpha;a} \|g\|_{\beta;a} \|h\|_{\gamma;a}\;,
    \end{equation*}
    where we used $k \geq j + 3$ and $\alpha<1$ in performing the sum. This completes the proof of the lemma.
\end{proof}

We then have the following direct consequence for the commutator $\Com_a$.

\begin{cor} \label{cor:com_prec_circ}
    Recall the commutator $\Com_a$ given in \eqref{e:commutator_defn}. For $\alpha\in (0,1)$ and $\beta,\gamma\in \R$ such that $\alpha+\beta+\gamma>0$, $\alpha+\beta\in(-1,2)$ and $\beta+\gamma<0$, we have
    \begin{equation*}
        \|\Com_a(f;\, g;\, h)\|_{\alpha+\beta+\gamma;a}\lesssim \|f\|_{\alpha;a} \|g\|_{\beta;a} \|h\|_{\gamma;a}\;.
    \end{equation*}
\end{cor}
\begin{proof}
It suffices to note that one has the expression
\begin{equation*}
    \Com_a (f;g;h) = \sum_{|i-j| \leq 1} R_{j,a}(f,g) \cdot \Delta_{i,a} h\;.
\end{equation*}
The bound then follows directly from Lemma~\ref{le:com_block}. 
\end{proof}

\begin{lem}
For $\alpha\in(0,1)$ and $\beta<0$ such that $\alpha+\beta\in (-1,1)$, 
    \begin{equation*}
        \left\| \, [ \Delta_{ j , a } , f ] g - \Delta_{j,a} (f \circ_a g) \right\|_{L^{\infty}} \lesssim 2^{-(\alpha+\beta)j} \|f\|_{\alpha;a} \|g\|_{\beta;a}
    \end{equation*}
    uniformly over $j\geq -1$.
\end{lem}
\begin{proof}
By explicit computation, we have
\begin{equation*}
    [\Delta_{j,a},f] g-\Delta_{j,a}(f\circ_a g)= \Delta_{j,a}(f\succ_a g) + R_{j,a}(f,g)\;.
\end{equation*}
Then the assertion follows from the bound \eqref{e:bony_generalised_high_low} and Lemma~\ref{le:R_Linfty_bound}.
\end{proof}

We now turn to controlling the quantity $\PL_a^F$ defined in \eqref{e:para_linear_defn}. We first give a preliminary lemma. 

\begin{lem} \label{lem:Holder_composition}
Let $\alpha \in (0,1)$. We have
\begin{equation*}
    \|G(u)\|_\alpha \lesssim \|G\|_{\cC^1} \big( 1 + \|u\|_\alpha \big)\;,
\end{equation*}
and
\begin{equation*}
    \|G(u_1) - G(u_2)\|_\alpha \lesssim \|G\|_{\cC^2} \big( 1 + \|u_1\|_\alpha + \|u_2\|_\alpha \big) \|u_1 - u_2\|_\alpha\;.
\end{equation*}
\end{lem}
\begin{proof}
Both bounds follow directly from the expression
\begin{equation*}
    G \big( u_j(x) \big) - G \big( u_j(y) \big) = \big( u_j(x) - u_j(y) \big) \int_0^1 G' \big( s \, u_j(x) + (1-s) u_j(y) \big) {\rm d}s\;.
\end{equation*}
We omit the details. 
\end{proof}

We are now ready to state the para-linearisation lemma. 

\begin{lem} \label{lem:para_linearisation}
Recall the operator $\PL_a^F$ in \eqref{e:para_linear_defn}. Let $\alpha \in (0,1)$. We have
\begin{equation} \label{eq:para_linear}
    \|\PL_a^F(u)\|_{2\alpha;a} \lesssim \|F\|_{\cC^2} \big(1 + \|u\|_\alpha^2 \big)\;,
\end{equation}
and
\begin{equation} \label{eq:para_linear_diff}
    \left\| \PL_a^F (u_1) - \PL_a^F(u_2) \right\|_{2\alpha;a} \lesssim \|F\|_{\cC^3} \big( 1 + \|u_1\|_\alpha^2 + \|u_2\|_\alpha^2 \big) \|u_1-u_2\|_\alpha\;.
\end{equation}
\end{lem}
\begin{proof}
We start with \eqref{eq:para_linear}. We need to show that
\begin{equation} \label{eq:para_linear_block}
    \big\| \Delta_{k,a} \big( F(u) - F'(u) \prec_a u \big) \big\|_{L^\infty} \lesssim 2^{-2\alpha k} \|F\|_{\cC^1} \big( 1 + \|u\|_\alpha^2 \big)\;.
\end{equation}
According to definition of the operator $R_{k,a}$ in \eqref{e:R_op_defn}, we have the identity
\begin{equation*}
    \Delta_{k,a} \big( F(u) - F'(u) \prec_a u \big) = \big( \Delta_{k,a} F(u) - F'(u) \, \Delta_{k,a} u \big) - R_{k,a} \big( F'(u), u \big)\;.
\end{equation*}
For the term with $R_{k,a}$, by Lemmas~\ref{le:R_Linfty_bound} and~\ref{lem:Holder_composition}, we have
\begin{equation*}
    \big\|R_{k,a} \big( F'(u), u \big) \big\|_{L^\infty} \lesssim 2^{-2\alpha k} \|F'(u)\|_\alpha \|u\|_\alpha \lesssim 2^{-2\alpha k} \|F\|_{\cC^1} \big( 1 + \|u\|_\alpha^2 \big)\;.
\end{equation*}
For the first term, noting that $\Delta_{k,a}(x,\cdot)$ integrates to $0$ for $k \geq 0$, we have
\begin{equation*}
    \big( \Delta_{k,a} F(u) \big)(x) - F'\big( u(x) \big) \cdot (\Delta_{k,a}u)(x) = \int_{\T^d} \Delta_{k,a}(x,y) \Theta^F_{u}(x,y) {\rm d}y + \1_{k=-1} \, \widetilde{F} \big( u(x) \big)\;,
\end{equation*}
where $\widetilde{F}(u) = F(u) - F'(u) \cdot u$, and $\Theta^F_{u}(x,y)$ is the Taylor remainder
\begin{equation*}
    \begin{split}
    \Theta^F_{u}(x,y):=  F\big( u(y) \big) - F \big( u(x) \big) - F'\big( u(x) \big) \, \big( u(y) - u(x) \big)\;,
    \end{split}
\end{equation*}
which satisfies the pointwise bound
\begin{equation*}
    \big| \Theta^F_{u}(x,y) \big| \leq \|F''\|_{L^\infty} \, \|u\|_\alpha^2 \, |x-y|^{2\alpha}\;.
\end{equation*}
Combining with the pointwise bound for the kernel $\Delta_{k,a}$ in Lemma~\ref{lem:blocks_derivatives_pointwise} (with $\ell=m=0$), we get
\begin{equation*}
    \big\| \Delta_{k,a} \big( F(u) \big) - F'(u) \, \Delta_{k,a} u \big\|_{L^\infty} \lesssim 2^{-2\alpha k} \, \|F\|_{\cC^2} \big( 1 + \|u\|_\alpha^2 \big)\;.
\end{equation*}
This proves the bound \eqref{eq:para_linear_block} and hence also \eqref{eq:para_linear}. 

We now turn to the difference bound \eqref{eq:para_linear_diff}. According to the above decomposition of $\Delta_{k,a} \big( \PL_a^F (u_j) \big)$, we have
\begin{equation*}
    \begin{split}
    \big\| \Delta_{k,a} \big( &\PL_a^F (u_1) - \PL_a^F (u_2) \big) \big\|_{L^\infty} \leq \sup_{x \in \T^d} \int_{\T^d} \big| \Delta_{k,a}(x,y) \cdot (\Theta^F_{u_1}-\Theta^F_{u_2})(x,y) \big| \, {\rm dy}\\
    &+ \1_{k=-1} \big\| \widetilde{F}(u_1) - \widetilde{F}(u_2) \big\|_{L^\infty} + \big\| R_{k,a} \big( F'(u_1), u_1 \big) - R_{k,a} \big( F'(u_2), u_2 \big) \|_{L^\infty}\;.
    \end{split}
\end{equation*}
For the term with $\widetilde{F}$, it is straightforward that
\begin{equation} \label{eq:para_linear_difference_constant}
    \1_{k=-1} \|\widetilde{F}(u_1) - \widetilde{F}(u_2)\|_{L^\infty} \lesssim \1_{k=-1} \|F\|_{\cC^2} \big( 1 + \|u_1\|_{L^\infty} + \|u_2\|_{L^\infty} \big) \|u_1 - u_2\|_{L^\infty}\;.
\end{equation}
As for the term involving $R_{k,a}$, by the bi-linearity of $R_{k,a}$, Lemmas~\ref{le:R_Linfty_bound} and~\ref{lem:Holder_composition} (with $G=F'$), we have
\begin{equation} \label{eq:para_linear_difference_R}
    \begin{split}
    &\phantom{111}\big\| R_{k,a} \big( F'(u_1), u_1 \big) - R_{k,a} \big( F'(u_2), u_2 \big) \|_{L^\infty}\\
    &\leq \big\| R_{k,a} \big( F'(u_1) - F'(u_2), \, u_1 \big) \big\|_{L^\infty} + \big\| R_{k,a} \big( F'(u_2), \, u_1 - u_2 \big) \big\|_{L^\infty}\\
    &\lesssim 2^{-2\alpha k} \|F\|_{\cC^3} \big( 1 + \|u_1\|_\alpha^2 + \|u_2\|_\alpha^2 \big) \cdot \|u_1 - u_2\|_\alpha\;.
    \end{split}
\end{equation}
For the term involving the difference between $\Theta^F_{u_j}$, we first note that by the expression
\begin{equation*}
    \Theta^F_{u_j}(x,y) = \big( u_j(y) - u_j(x) \big)^2 \iint_{0<s<t<1} F'' \big( (1-s) u_j(x) + s u_j(y) \big) \, {\rm d}s \, {\rm d}t\;,
\end{equation*}
we have the pointwise bound
\begin{equation*}
    \big| \Theta^F_{u_1}(x,y) - \Theta^F_{u_2}(x,y) \big| \lesssim \|F\|_{\cC^3} \big( 1 + \|u_1\|_\alpha^2 + \|u_2\|_\alpha^2 \big) \, \|u_1 - u_2\|_\alpha \, |x-y|^{2\alpha}\;.
\end{equation*}
Combining with the pointwise bound for $\Delta_{k,a}$ from Lemma~\ref{lem:blocks_derivatives_pointwise}, we deduce
\begin{equation} \label{eq:para_linear_difference_comm}
    \begin{split}
    &\phantom{111}\sup_{x \in \T^d} \int_{\T^d} \big| \Delta_{k,a}(x,y) \cdot (\Theta^F_{u_1}-\Theta^F_{u_2})(x,y) \big|\, {\rm d}y \\
    &\lesssim 2^{-2\alpha k} \|F\|_{\cC^3} \big( 1 + \|u_1\|_\alpha^2 + \|u_2\|_\alpha^2 \big) \, \|u_1 - u_2\|_\alpha\;.
    \end{split}
\end{equation}
The desired bound \eqref{eq:para_linear_diff} then follows from combining \eqref{eq:para_linear_difference_constant}, \eqref{eq:para_linear_difference_R} and \eqref{eq:para_linear_difference_comm}. This completes the proof of the lemma. 
\end{proof}

\subsection{Heat kernel bounds}

We define
\begin{equation} \label{e:I_a_defn}
    \iI_a(f) := \int_{0}^{t}e^{(t-\tau)\lL_a} f(\tau,\cdot) {\rm d}\tau\;.
\end{equation}
Recall from \eqref{eq:spacetime_norm_weight} and \eqref{eq:spacetime_Holder_norm_weight} the weighted spacetime norms $L_{\sigma,T}^\infty \zZ$ and $\cC_{\sigma,T}^\beta \zZ$. We start with the following lemma.

\begin{lem} \label{le:heat_smoothing_space}
    For every $\alpha \in \R$ and $\sigma \geq 0$, we have
    \begin{equation*}
        \|e^{t \lL_a} f\|_{\alpha+\sigma;a} \lesssim t^{-\frac{\sigma}{2}} \|f\|_{\alpha;a}\;.
    \end{equation*}
    Furthermore, if $\theta \in [0,2]$, we also have the bound
    \begin{equation*}
        \| (e^{t \lL_a} - e^{s \lL_a}) f\|_{\alpha+\sigma-\theta;a} \lesssim (t-s)^{\frac{\theta}{2}} s^{-\frac{\sigma}{2}} \|f\|_{\alpha;a}\;.
    \end{equation*}
\end{lem}
\begin{proof}
By Lemma~\ref{lem:heat_convolution_dyadic_bound}, there exists $c>0$ such that
\begin{equation*}
    \|\Delta_{j,a} e^{t \lL_a}\|_{L^\infty \rightarrow L^\infty} \lesssim e^{- c t 2^{2j}} \lesssim t^{-\frac{\sigma}{2}} 2^{-\sigma j}\;.
\end{equation*}
The proportionality constant depends on $\sigma \geq 0$ but is uniform in $t$ and $j$. Hence, we have
\begin{equation*}
    \|\Delta_{j,a} e^{t \lL_a} f\|_{L^\infty} \lesssim t^{-\frac{\sigma}{2}} 2^{-\sigma j} \sum_{\ell=j-1}^{j+1} \|\Delta_{\ell,a} f\|_{L^\infty} \lesssim t^{-\frac{\sigma}{2}} 2^{-(\alpha+\sigma) j} \|f\|_{\alpha;a}\;.
\end{equation*}
This proves the first bound. As for second one, we first consider the case $s=\sigma=0$, where we interpret $0^0=1$. Since $\|e^{t \lL_a}\|_{L^\infty \rightarrow L^\infty} \lesssim 1$, we have
\begin{equation*}
    \|\Delta_{j,a} (e^{t \lL_a} - \id) f\|_{L^\infty} =  \|(e^{t \lL_a} - \id) \Delta_{j,a} f\|_{L^\infty} \lesssim \|\Delta_{j,a} f\|_{L^\infty} \lesssim 2^{-\alpha j} \|f\|_{\alpha;a}\;.
\end{equation*}
On the other hand, using the first claim with $\sigma=0$, we also have
\begin{equation*}
    \|\Delta_{j,a} (e^{t \lL_a} - \id) f\|_{L^\infty} \leq \int_{0}^{t} \|e^{r \lL_a} \Delta_{j,a} \lL_a f \|_{L^\infty} {\rm d}r \lesssim t 2^{(2-\alpha) j} \|f\|_{\alpha;a}\;.
\end{equation*}
Combining the two bounds together, we deduce
\begin{equation*}
    \|\Delta_{j,a} (e^{t \lL_a} - \id) f\|_{L^\infty} \lesssim (1 \wedge t 2^{2j}) 2^{-\alpha j} \|f\|_{\alpha;a} \lesssim t^{\frac{\theta}{2}} 2^{(\theta-\alpha) j} \|f\|_{\alpha;a}\;,
\end{equation*}
where the last bound holds if $\theta \in [0,2]$. For general $s<t$ and $\sigma \geq 0$, replacing $f$ by $e^{s \lL_a} f$ and $\alpha$ by $\alpha+\sigma$, we then have
\begin{equation*}
    \begin{split}
    \left\|\Delta_{j,a} \big( e^{(t-s) \lL_a} - \id \big) e^{s \lL_a} f \right\|_{L^\infty} &\lesssim (t-s)^{\frac{\theta}{2}} 2^{-(\alpha + \sigma -\theta) j} \|e^{s \lL_a} f\|_{\alpha+\sigma;a}\\&\lesssim (t-s)^{\frac{\theta}{2}} 2^{-(\alpha + \sigma - \theta) j} s^{-\frac{\sigma}{2}} \|f\|_{\alpha;a}\;.
    \end{split}
\end{equation*}
This completes the proof of the second bound and hence the lemma. 
\end{proof}
With the same method in Lemma \ref{le:heat_smoothing_space}, we have
\begin{equation*}
    \|e^{\cdot\lL_{a}}f\|_{\cC^{(\alpha+\sigma)/2}_{\sigma,T}L^{\infty}}\lesssim \|f\|_{\alpha;a}\;.
\end{equation*}

\begin{lem} \label{le:heat_smoothing_spacetime}
Let $\alpha, \beta \in \R$ with $\alpha \in [\beta, \beta+2]$, and $\sigma \in [0,2)$. Then, we have
\begin{equation} \label{e:heat_smoothing_spacetime_fixed}
    \|(\iI_a f)(t)\|_{\alpha;a} \lesssim t^{1-\frac{\sigma+\alpha-\beta}{2}} \|f\|_{L_{\sigma,t}^\infty \cC_{x}^{\beta;a}}\;.
\end{equation}
Furthermore, for $\alpha \in [\beta, \beta+2)$, every $\eta \in (0,2-\alpha+\beta]$ and $T>0$, we also have
\begin{equation} \label{e:heat_smoothing_spacetime_time}
    \left\| (\iI_a f)(t) - (\iI_a f)(s) \right\|_{\alpha;a} \lesssim (t-s)^{\frac{\eta}{2}} \big( s^{1-\frac{\sigma+\alpha+\eta-\beta}{2}} + t^{1-\frac{\sigma+\alpha+\eta-\beta}{2}} \big) \|f\|_{L_{\sigma,t}^{\infty} \cC_{x}^{\beta;a}}\;.
\end{equation}
The proportionality constant depends on $T$, but is independent of $0<s<t<T$. 
\end{lem}
\begin{proof}
By Lemma~\ref{lem:heat_convolution_dyadic_bound}, we have
\begin{equation*}
    \begin{split}
    \| \Delta_{j,a} (\iI_a f)(t) \|_{L^\infty} &\leq \sum_{\ell=j-1}^{j+1} \int_{0}^{t} \left\| \Delta_{j,a} e^{(t-r) \lL_a} \Delta_{\ell,a} f(r) \right\|_{L^\infty} {\rm d}r\\
    &\lesssim 2^{- \beta j} \|f\|_{L_{\sigma,t}^{\infty} \cC_{x}^{\beta;a}} \int_{0}^{t} e^{-c (t-r) \cdot 2^{2j} } r^{-\frac{\sigma}{2}} {\rm d}r\;.
    \end{split}
\end{equation*}
Splitting the domain of the above time integral into $[0,\frac{t}{2}]$ and $[\frac{t}{2},t]$, we get
\begin{equation*}
    \|\Delta_{j,a} (\iI_a f)(t)\|_{L^\infty} \lesssim 2^{-\beta j}  \cdot \Big( t^{1-\frac{\sigma}{2}} e^{-c_1 t \cdot 2^{2j}} + 2^{-2j} t^{-\frac{\sigma}{2}} \big( 1 - e^{-c_2 t \cdot 2^{2j}} \big) \Big) \cdot \|f\|_{L_{\sigma,t}^{\infty} \cC_x^{\beta;a}}
\end{equation*}
for some $c_1, c_2 > 0$. For the two terms in the parenthesis, we bound them separately by
\begin{equation*}
    \begin{split}
    t^{1-\frac{\sigma}{2}} e^{-c_1 t \cdot 2^{2j}} \lesssim t^{1-\frac{\sigma}{2}} (t \cdot 2^{2j})^{-\frac{\alpha-\beta}{2}} \lesssim 2^{-(\alpha-\beta)j} \, t^{1-\frac{\sigma+\alpha-\beta}{2}}\;,\\
    2^{-2j} t^{-\frac{\sigma}{2}} \big( 1 - e^{-c_2 t \cdot 2^{2j}} \big) \lesssim 2^{-2j} t^{-\frac{\sigma}{2}} (t \cdot 2^{2j})^{1-\frac{\alpha-\beta}{2}} \lesssim 2^{-(\alpha-\beta)j} \, t^{1-\frac{\sigma+\alpha-\beta}{2}}\;,
    \end{split}
\end{equation*}
where we have used the assumption $\alpha \geq \beta$ in the first one, and $\alpha-\beta \leq 2$ in the second one. Together with the factor $2^{-\beta j}$, this shows
\begin{equation*}
    \|\Delta_{j,a} (\iI_a f)(t)\|_{L^\infty} \lesssim 2^{-\alpha j} t^{1-\frac{\sigma+\alpha-\beta}{2}} \|f\|_{L_{\sigma,t}^{\infty} \cC_{x}^{\beta;a}}\;,
\end{equation*}
which proves the first claim. 

As for continuity in time, for $t \geq s$, we have
\begin{equation*}
    (\iI_a f)(t) - (\iI_a f)(s) = \big( e^{(t-s) \lL_a} - \id \big) \big( (\iI_a f)(s) \big) + \int_{s}^{t} e^{(t-r) \lL_a} f(r) {\rm d}r\;.
\end{equation*}
For the first term above, by Lemma~\ref{le:heat_smoothing_space} and \eqref{e:heat_smoothing_spacetime_fixed}, we have
\begin{equation*}
    \begin{split}
    \left\| \big( e^{(t-s) \lL_a} - \id \big) \big( (\iI_a f)(s) \big) \right\|_{\alpha;a} 
    &\lesssim (t-s)^{\frac{\eta}{2}} \|(\iI_a f)(s)\|_{\alpha+\eta;a}\\
    &\lesssim (t-s)^{\frac{\eta}{2}} s^{1-\frac{\sigma+\alpha+\eta-\beta}{2}}  \|f\|_{L_{\sigma,t}^{\infty} \cC_{x}^{\beta;a}}\;.
    \end{split}
\end{equation*}
For the second term, again by Lemma~\ref{le:heat_smoothing_space} and the assumption $\eta < 2-(\alpha-\beta)$, we have
\begin{equation*}
    \begin{split}
    \int_{s}^{t} \left\| e^{(t-r) \lL_a} f(r)  \right\|_{\alpha;a} {\rm d}r &\lesssim \|f\|_{L_{\sigma,t}^{\infty} \cC_{x}^{\beta;a}} \int_{s}^{t} (t-r)^{-\frac{\alpha-\beta}{2}} r^{-\frac{\sigma}{2}} {\rm d}r\\
    &\lesssim (t-s)^{1-\frac{\alpha-\beta}{2}} s^{-\frac{\sigma}{2}} \|f\|_{L_{\sigma,t}^{\infty} \cC_{x}^{\beta;a}} \lesssim (t-s)^{\frac{\eta}{2}} t^{1-\frac{\alpha-\beta+\eta}{2}} s^{-\frac{\sigma}{2}} \|f\|_{L_{\sigma,t}^{\infty} \cC_{x}^{\beta;a}}\;.
    \end{split}
\end{equation*}
Combining the above two bounds finishes the proof of \eqref{e:heat_smoothing_spacetime_time} and hence also the lemma. 
\end{proof}

For for every $T>0$ and $\sigma,\alpha \in \R$, define 
\begin{equation}\label{eq: weighted_mathfrakL}
    \|f\|_{\mathfrak{L}^{\alpha;a}_{\sigma}((0,T]\times \T^d)}:=\|f\|_{L_{\sigma,T}^{\infty}\cC^{\alpha;a}_{x}}+\|f\|_{ \cC^{{\alpha}/{2}}_{\sigma,T}L^{\infty}_{x}}\;.
\end{equation}
For convenience, we denote $\| \cdot \|_{ \mathfrak{ L }^{ \alpha ; a }_{ \sigma }(( 0 , T ] \times \T^d ) }$ by $\|\cdot\|_{\mathfrak{L}^{\alpha;a}_{\sigma,T}}$. This estimate below is in spirit of the proof in \cite[Lemma 18]{Gubinelli2017lectures}.

\begin{lem} \label{le:heat_comm_space}
For $\alpha\in (0,1)$ and $\beta<0$ such that $\alpha+\beta \in (-1,0)$, there exist $C$, $c>0$ such that
\begin{equation*}
    \left\|\Delta_{k,a}([e^{t\lL_{a}},f\prec_{a}]g) \right\|_{L^{\infty}}\leq C B_{k}(t)\|f\|_{\alpha;a}\|g\|_{\beta;a}
\end{equation*}
where
\begin{equation} \label{eq:heat_comm_space_Akt}
    B_{k}(t):=2^{-2k}\sum_{j=-1}^{k}2^{(2-\alpha-\beta)j}e^{-ct2^{2j}}\;,
\end{equation}
and the proportionality constants are uniform in $k$ and $t>0$. 
\end{lem}
\begin{proof}
By definition of the paraproduct $\prec_a$, we have the expression
\begin{equation*}
    \Delta_{k,a} \big(\, [\, e^{t \lL_a}\, ,\, f \prec_a]\, g\, \big) = \sum_{j \geq -1} \Delta_{k,a} \big( [\, e^{t \lL_a},\, S_{j,a} f\,] \Delta_{j,a} g \,\big)\;.
\end{equation*}
We split the sum into three regions: $j \leq k-2$, $k-1 \leq j \leq k + 1$ and $j \geq k+2$. For both $j \leq k-2$ and $j \geq k+2$, we write
\begin{equation*}
    \begin{split}
    \Delta_{k,a} \big( [e^{t \lL_a}, S_{j,a} f] \Delta_{j,a} g \big) = &\sum_{\ell=k-1}^{k+1} \Delta_{\ell,a} e^{t \lL_a} \Delta_{k,a} \big( S_{j,a}f \, \Delta_{j,a} g \big)\\
    &- \sum_{\ell=j-1}^{j+1} \Delta_{k,a} \Big( S_{j,a} f \cdot \Delta_{j,a} \big( e^{t \lL_a} \Delta_{\ell,a} g \big) \Big)\;.
    \end{split}
\end{equation*}
Applying Corollary~\ref{cor:decoupling_spectrum_S} to the right hand side above and combining with the operator bounds (implied by Lemma~\ref{lem:heat_convolution_dyadic_bound})
\begin{equation} \label{e:heat_block_composition}
    \|e^{t \lL_a} \Delta_{\ell,a}\|_{L^\infty \rightarrow L^\infty} \lesssim e^{-ct 2^{2\ell}}\;, \qquad  \|e^{t \lL_a} \Delta_{j,a} \|_{\cC^{\beta;a} \rightarrow \cC^{\beta;a}} \lesssim e^{-ct 2^{2j}}\;,
\end{equation}
we get (with a factor $\|f\|_{\alpha;a} \|g\|_{\beta;a}$)
\begin{equation} \label{e:heat_comm_space_ends}
   \left\| \Delta_{k,a} \big( [e^{t \lL_a}, S_{j,a} f] \Delta_{j,a} g \big)  \right\|_{L^\infty} \lesssim
   \begin{cases}
2^{-2 k + (2-\alpha-\beta) j}  e^{-ct \cdot 2^{2j}}, & j \leq k-2\;,\\
2^{k-(1+\alpha+\beta)j} e^{-ct \cdot 2^{2 k}}, & j \geq k+2\;,
\end{cases} 
\end{equation}
As for $k - 1 \leq j \leq k + 1$, we write
\begin{equation*}
    \Delta_{k,a} \big( [e^{t \lL_a}, S_{j,a} f] \Delta_{j,a} g \big) = [\Delta_{k,a} e^{t\lL_a}, S_{j,a} f] \Delta_{j,a} g + [S_{j,a} f, \Delta_{k,a}] \big( e^{t \lL_a} \Delta_{j,a} g \big)\;.
\end{equation*}
By Lemma~\ref{le:block_multiplication_comm} (and the same argument for its proof) and Lemma~\ref{lem:heat_convolution_dyadic_bound} on the kernel of the operator $\Delta_{k,a} e^{t \lL_a}$ (and $\Delta_{j,a} e^{t\lL_a}$), we have the bounds
\begin{equation*}
    \begin{split}
    \left\| [\Delta_{k,a} e^{t\lL_a}, S_{j,a} f] \Delta_{j,a} g \right\|_{L^\infty} &\lesssim 2^{-k} e^{-c t \cdot 2^{2k}} \|\nabla S_{j,a} f\|_{L^\infty} \|\Delta_{j,a} g\|_{L^\infty}\\
    &\lesssim 2^{-k + (1-\alpha-\beta) j} e^{-c t \cdot 2^{2k}} \|f\|_{\alpha;a} \|g\|_{\beta;a}\;,
    \end{split}
\end{equation*}
and
\begin{equation*}
    \begin{split}
    \left\| [S_{j,a} f, \Delta_{k,a}] \big( e^{t \lL_a} \Delta_{j,a} g \big) \right\|_{L^\infty} &\lesssim  2^{-k} \|\nabla S_{j,a} f\|_{L^\infty} \sum_{\ell=j-1}^{j+1} \left\| \Delta_{j,a} e^{t \lL_a} \big( \Delta_{\ell,a} g\big)  \right\|_{L^\infty}\\
    &\lesssim 2^{-k + (1-\alpha-\beta) j} e^{-ct \cdot 2^{2j}} \|f\|_{\alpha;a} \|g\|_{\beta;a}\;.
    \end{split}
\end{equation*}
Hence, we have
\begin{equation*}
    \left\| \Delta_{k,a} \big( [e^{t \lL_a}, S_{j,a} f] \Delta_{j,a} g \big) \right\|_{L^\infty} \lesssim 2^{-(\alpha+\beta)k} e^{-ct \cdot 2^{2k}} \|f\|_{\alpha;a} \|g\|_{\beta;a}\;, \quad \text{for} \; k-1\leq j \leq k + 1\;.
\end{equation*}
Combining the above bound and \eqref{e:heat_comm_space_ends} and using the assumption $1+\alpha+\beta>0$, we arrive at the conclusion. 
\end{proof}

\begin{lem}
\label{le:heat_comm_spacetime}
Let $\alpha \in (0,1)$ and $\beta<0$ be such that $\alpha + \beta \in (-1,0)$. Let $\sigma \in (0,2)$ and $T>0$. Then we have
    \begin{equation*}
        \| \big([\iI_{a},f\prec_{a}]g \big)(t)\|_{\alpha+\beta+2;a} \lesssim t^{-\frac{\sigma}{2}} \big( \|f\|_{\cC_{\sigma,t}^{\alpha/2} L_x^\infty} + \|f\|_{L_{\sigma,t}^{\infty} \cC_{x}^{\alpha;a}} \big)  \|g\|_{L^{\infty}_{t}\cC^{\beta;a}_{x}}
    \end{equation*}
    uniformly over $t \in [0,T]$. 
\end{lem}
\begin{proof}
    We write
    \begin{equation} \label{e:heat_prec_com_decomposition}
        \begin{split}
        \big([\iI_a, f \prec_a] g\big)(t) = &\int_{0}^{t} e^{(t-r) \lL_a} \Big( \big( f(r) - f(t) \big) \prec_a g(r) \Big) {\rm d}r\\
        &+ \int_{0}^{t} [e^{(t-r) \lL_a}, f(t) \prec_a] g(r) {\rm d}r\;.
        \end{split}
    \end{equation}
    For the first term on the right hand side above, by \eqref{e:heat_block_composition} and \eqref{e:bony_generalised_low_high_other}, the $L^\infty$-norm of the $k$-th block of its integrand can be controlled by
    \begin{equation*}
        \begin{split}
        &\phantom{111}\left\| \Delta_{k,a} e^{(t-r) \lL_a} \Big( \big( f(r) - f(t) \big) \prec_a g(r) \Big) \right\|_{L^\infty}\\
        &\lesssim e^{-c (t-r) \cdot 2^{2 k}} \sum_{\ell=k-1}^{k+1} \left\| \Delta_{\ell,a} \Big( \big( f(r) - f(t) \big) \prec_a  g(r) \Big) \right\|_{L^\infty}\\
        &\lesssim 2^{- \beta k} e^{-c(t-r) \cdot 2^{2 k}} r^{-\frac{\sigma}{2}} \big((t-r)^{\frac{\alpha}{2}}  \|f\|_{\cC_{\sigma,t}^{\alpha/2}L^{\infty}_{x}}+ 2^{-\alpha k}\|f\|_{L^{\infty}_{\sigma,t}\cC_x^{\alpha;a}} \big) \|g\|_{L_t^\infty \cC_x^{\beta;a}}\;.
        \end{split}
    \end{equation*}
    Integrating $r \in [0,\frac{t}{2}]$ and $r \in [\frac{t}{2},t]$ separately and summing them up, we get
    \begin{equation*}
        \begin{split}
        &\phantom{11}\int_{0}^{t} \left\| \Delta_{k,a} e^{(t-r) \lL_a} \Big( \big( f(r) - f(t) \big) \prec_a g(r) \Big) \right\|_{L^\infty} {\rm d}r\\
        &\lesssim 2^{-(\alpha+\beta+2)k} t^{-\frac{\sigma}{2}} (\|f\|_{\cC_{\sigma,t}^{\alpha/2} L_x^\infty}+\|f\|_{L^{\infty}_{\sigma,t}\cC_x^{\alpha;a}}) \|g\|_{L_t^\infty \cC_x^{\beta;a}}\;.
        \end{split}
    \end{equation*}
    For the second term on the right hand side of \eqref{e:heat_prec_com_decomposition}, by Lemma~\ref{le:heat_comm_space}, we can control the $L^\infty$-norm of the $k$-th block of its integrand by
    \begin{equation*}
        \left\| \Delta_{k,a} \big( [e^{(t-r) \lL_a}, f(t) \prec_a] g(r)  \big) \right\|_{L^\infty} \lesssim B_{k}(t-r) \cdot t^{-\frac{\sigma}{2}} \|f\|_{L_{\sigma,t}^{\infty} \cC_{x}^{\alpha;a}} \|g\|_{L_t^\infty \cC_x^{\beta;a}}\;,
    \end{equation*}
    where $B_{k}$ is given in \eqref{eq:heat_comm_space_Akt}. This time, integrating $r \in [0,t]$ directly gives
    \begin{equation*}
        \int_{0}^{t} \left\| \Delta_{k,a} \big( [e^{(t-r) \lL_a}, f(t) \prec_a] g(r)  \big) \right\|_{L^\infty} {\rm d}r \lesssim 2^{-(\alpha+\beta+2) k}\cdot t^{-\frac{\sigma}{2}} \|f\|_{L_{\sigma,t}^{\infty} \cC_{x}^{\alpha;a}} \|g\|_{L_t^\infty \cC_x^{\beta;a}}\;,
    \end{equation*}
    where we have used the assumption that $\alpha+\beta<0$. This completes the proof. 
\end{proof}

In dynamical $\Phi^{4}_{3}$, we need use the operator
\begin{equation*}
    (\widetilde{\iI}_{a}f)(t) := \int_{0}^{t} e^{(t-r)(\lL_{a}-1)} f(r){\rm d} r
\end{equation*}
instead of $\iI_a$. But all the bounds above for $\iI_a$ are still valid for $\widetilde{\iI}_a$.

\section{Convergences}
\label{sec:para_convergence}

In this section, we fix a coefficient matrix $A$ satisfying Assumption~\ref{as:A}, and let $\bar{A}$ be its homogenised matrix as given in \eqref{def:bar_a_homogenised_matrix}. By Corollary~\ref{cor:coefficient_uniform}, there exist constants $\Lambda, M>0$ such that $A_\eps := A(\cdot / \eps) \in \sS_d(\Lambda, M)$ for all $\eps \in \N^{-1}$. This guarantees that all the bounds in Section~\ref{sec:para_general} are true for the operator $\lL_\eps$ uniformly in $\eps \in \N^{-1}$ (including $\lL_0$). 

We use the subscript $\eps$ in $\Delta_{j,\eps}$, $\prec_{\eps}$, $\circ_{\eps}$, $\Com_\eps$ and so on to indicate that these operations are associated to $-\lL_{\eps} = - \div \big( A_\eps \nabla \big)$. For operations associated to the homogenised operator $-\lL_0 = - \div (\bar{A} \nabla)$, we drop the subscript $0$ and simply write $\Delta_{j}$, $\prec$, $\circ$, $\Com$, etc.

\subsection{Convergence of blocks}

The main convergence theorem in this section is the following. 

\begin{thm} \label{th:convergence_single_block}
For every $m \in \Z$ and every $\kappa>0$, there exists $\theta>0$ such that
\begin{equation*}
    \| (\lL_\eps^m \Delta_{j,\eps} - \lL_0^m \Delta_{j,0}) \Pi_0^\perp f \|_{L^\infty} \lesssim_m \eps^\theta 2^{(2m+\kappa)j} \|\Pi_0^\perp f\|_{L^\infty}
\end{equation*}
uniformly over all $j \geq -1$, $\eps \in \N^{-1}$ and all $f \in L^\infty (\T^d)$. 
\end{thm}

The proof of this theorem is postponed after Corollary~\ref{cor:convergence_L2_rate}. We first give an almost immediate consequence.

\begin{cor} \label{cor:convergence_single_block}
Let $\alpha \in (-1,1)$. For every $\kappa>0$, there exists $\theta>0$ such that for every $m \in \Z$, we have the bounds
\begin{equation} \label{eq:convergence_single_block_positive}
    \|(\lL_\eps^m \Delta_{j,\eps} - \lL_0^m \Delta_{j,0}) f\|_{L^\infty} \lesssim \eps^\theta 2^{(2m-\alpha+\kappa) j} \|f\|_{\alpha}\;, \qquad m \geq 0\;,
\end{equation}
and
\begin{equation} \label{eq:convergence_single_block_negative}
    \|(\lL_\eps^m \Delta_{j,\eps} - \lL_0^m \Delta_{j,0}) \Pi_0^\perp f\|_{L^\infty} \lesssim \eps^\theta 2^{(2m-\alpha+\kappa)j} \|\Pi_0^\perp f\|_{\alpha}\;, \qquad m \leq -1\;.
\end{equation}
The proportionality constants depend on $\alpha, \kappa, \theta$ and $m$, but are uniform over $j \geq -1$, $\eps \in \N^{-1}$ and $f \in \cC^\alpha (\T^d)$. 
\end{cor}
\begin{proof}
We first note it suffices to show that
\begin{equation} \label{eq:convergence_single_block_inclusive}
    \|(\lL_\eps^m \Delta_{j,\eps} - \lL_0^m \Delta_{j,0}) f\|_{L^\infty} \lesssim \eps^\theta 2^{(2m-\alpha+\kappa) j} \|f\|_{\alpha}
\end{equation}
uniformly over all $j \geq -1$, $\eps \in \N^{-1}$ and $f \in \cC^{\alpha}(\T^d)$ such that $\Pi_0 f = 0$. For $m \geq -1$, this is exactly the same as \eqref{eq:convergence_single_block_negative}. For $m \geq 0$, we note that
\begin{equation*}
    (\lL_\eps^m \Delta_{j,\eps} - \lL_0 \Delta_{j,0}) \Pi_0 \equiv 0\;,
\end{equation*}
and hence it also implies \eqref{eq:convergence_single_block_positive}. We now focus on \eqref{eq:convergence_single_block_inclusive} with $f$ such that $\Pi_0 f = 0$. Also, since $\|f\|_0 \lesssim \|f\|_{L^\infty}$, the case $\alpha=0$ is included in Theorem~\ref{th:convergence_single_block}. Hence, we only consider the cases where $\alpha \in (0,1)$ or $\alpha \in (-1,0)$. 

For $\alpha \in (0,1)$, by Lemma~\ref{lem:blocks_derivatives_pointwise} (with $\ell=0$), we have the uniform bound
\begin{equation*}
    \begin{split}
    \|\lL_\eps^m \Delta_{j,\eps} f\|_{L^\infty} &\leq \sum_{k=j-1}^{j+1} \|\lL_\eps^m \Delta_{j,\eps} \Delta_{k,\eps} f\|_{L^\infty}\\
    &\lesssim 2^{2mj} \sum_{k=j-1}^{j+1} \|\Delta_{k,\eps} f\|_{L^\infty} \lesssim 2^{(2m-\alpha) j} \|f\|_\alpha\;.
    \end{split}
\end{equation*}
Interpolating the above bound with Theorem~\ref{th:convergence_single_block} (and noting that $\|f\|_{L^\infty} \lesssim \|f\|_\alpha$ for $\alpha>0$) gives the desired bound \eqref{eq:convergence_single_block_inclusive} for $\alpha \in (0,1)$. 

For $\alpha \in (-1,0)$, write $\gamma = - \alpha \in (0,1)$. Also write $\Gamma_\eps^{(m,j)}:= \lL_\eps^m \Delta_{j,\eps} - \lL_0^m \Delta_{j,0}$. We then have
\begin{equation*}
    \left| \big( (\lL_\eps^m \Delta_{j,\eps} - \lL_0^m \Delta_{j,0}) f \big)(x) \right| \lesssim \|\Gamma_\eps^{(m,j)}(x,\cdot)\|_{L^1}^{1-\gamma} \, \|(\nabla_y \Gamma_\eps^{(m,j)})(x,\cdot)\|_{L^1}^\gamma \, \|f\|_{\cC^{-\gamma}}
\end{equation*}
By Lemma~\ref{lem:blocks_derivatives_pointwise} (with $\ell=1$), we have
\begin{equation*}
    \|(\nabla_y \Gamma_\eps^{(m,j)})(x,\cdot)\|_{L^1} \lesssim 2^{(2m+1)j}\;.
\end{equation*}
For the other term, by duality and Theorem~\ref{th:convergence_single_block}, for every $\kappa'>0$, there exists $\theta'>0$ such that
\begin{equation*}
    \|\Gamma_\eps^{(m,j)}(x,\cdot)\|_{L^1} = \sup_{\|g\|_{L^\infty}=1} \big| \big\langle \Gamma_{\eps}^{(m,j)}(x,\cdot), \, g \big\rangle \big| \lesssim \eps^{\theta'} 2^{(2m+\kappa')j}\;.
\end{equation*}
Together, they imply
\begin{equation*}
    \|\Gamma_\eps^{(m,j)} f\|_{L^\infty} \lesssim \big( \eps^{\theta'} 2^{(2m+\kappa')j} \big)^{1-\gamma} 2^{(2m+1)\gamma j} \|f\|_{\cC^\gamma} \lesssim \eps^{(1-\gamma) \theta'} 2^{(2m-\alpha+(1-\gamma) \kappa')} \|f\|_\alpha\;.
\end{equation*}
Since $\kappa'$ is arbitrary, this then implies \eqref{eq:convergence_single_block_inclusive} for $\alpha \in (-1,0)$. The proof of the corollary is thus complete. 
\end{proof}

% The resolvent set of a linear operator $T: \hH \rightarrow \hH$ on a Hilbert space $\hH$, denoted by $\rR(T)$, is the set of all complex numbers $\lambda \in \mathbf{C}$ such that $\lambda - T$ is a bijection on $\hH$, and that $\|(\lambda - T)^{-1}\|_{\hH \rightarrow \hH}$ is finite. The spectrum of the linear operator $T$, denoted by $\spec(T)$, is $\mathbf{C}\backslash \rR(T)$.

% Let $T$ be a linear, self-adjoint and compact operator from a separable Hilbert space $\hH$ to itself. Then $T$ has a discrete set of real eigenvalues, each with finite multiplicity. In this case, $\spec(T)$ is the set of all eigenvalues of $T$ together with their accumulation points. For every subset $\xX \subset \R$, let $P_{\xX,T}$ denote the projection onto the direct sum of invariant subspaces associated to eigenvalues of $T$ in $\xX$. In our situation, we will use $T = (\id - \lL_\eps)^{-1}$ or $(\id - \lL_0)^{-1}$. 

We now start to prove Theorem~\ref{th:convergence_single_block}. The main ingredient is the following Helffer-Sj\"ostrand's formula (see for example \cite{Helffer-Sjostrand}). For readers' convenience, we include a proof for it. 

\begin{lem} [Helffer-Sj\"ostrand's formula]
\label{lem:Helffer-Sjostrand}

For every smooth function $f\in \cC^{\infty}_0(\R)$, every self-adjoint operator $T$ on a Hilbert space $\hH$, and every $N\in \N^{+}$, we have the identity
    \begin{equation*}
        f(T) = -\frac{1}{\pi} \iint_{\CC} \, (\d_{\bar{z}} \widetilde{f}_N) (z) \, (z-T)^{-1} \, {\rm d}x \wedge {\rm d}y\;,
    \end{equation*}
    where $\rho\in \cC^{\infty}_0(-1,1)$ satisfies that $\rho|_{(-\frac{1}{2},\frac{1}{2})}\equiv 1$, $z= x + iy$, and $\widetilde{f}_{N}: \CC \to \CC$ is given by
    \begin{equation}\label{eq:almost_analytic_extension}
        \widetilde{f}_{N}(x+i y) := \sum_{n=0}^{N} \, \frac{f^{(n)}(x)}{n !} \, (i y)^{n} \, \rho\left(\frac{y}{\sqrt{1+x^2}}\right).
    \end{equation}
\end{lem}
\begin{proof}
    By Cauchy-Pompeiu formula, for any $\lambda\in \R$, we have
    \begin{equation*}
        f(\lambda) = -\frac{1}{\pi} \iint_{\R^2} \, \frac{\partial \tilde{f}_{N}}{\partial \bar{z}}(z) \, (z-\lambda)^{-1} \, {\rm d}x \wedge {\rm d}y\;,
    \end{equation*}
    where $z=  x + iy$. By spectral decomposition theorem, we have for $u,v\in \hH$,
    \begin{equation*}
        \scal{f(T)u,v} = \int_{\R}f(\lambda) \mu_{u,v}({\rm d} \lambda)\;,
    \end{equation*}
    where $\mu_{u,v}$ is the spectral measure with respect to the self-adjoint operator $T$. Combining these estimates, we have
    \begin{equation*}
        \scal{f(T)u,v} = -\frac{1}{\pi}\int_{\R} \; \Big[ \iint_{\R^2} \, \frac{(\d_{\bar{z}} \widetilde{f}_{N})(z)}{z-\lambda}  \, {\rm d}x \wedge {\rm d}y \Big]\; \mu_{u,v}({\rm d} \lambda)\;.
    \end{equation*}
    By the expression \eqref{eq:almost_analytic_extension}, $\d_{\bar{z}} \widetilde{f}_N$ satisfies the bound
    \begin{equation*}
            \left| (\d_{\bar{z}} \widetilde{f}_N)(z)  \right| \lesssim |y|^{N}\;.
    \end{equation*}
    Take $N = 1$. By Fubini theorem\footnote{Although $\mu_{u,v}$ can be a signed measure for generic $u,v \in \hH$, the identity holds
    \begin{equation*}
        \mu_{u,v} = \frac{1}{4}( \, \mu_{u+v,u+v} \, - \, \mu_{u-v,u-v} \, )\;,
    \end{equation*}
    where $\mu_{w,w}$ is a finite positive measure for every $w\in \hH$ and $\mu_{w,w}(\R)=\|w\|^2_{\hH}$.}, we can interchange the order of integrations to obtain
    \begin{equation*}
        \scal{f(T)u,v} = -\frac{1}{\pi} \iint_{\CC} \, (\d_{\bar{z}} \widetilde{f}_N)(z) \;\left[ \, \int_{\R}(z-\lambda)^{-1}\mu_{u,v}({\rm d} \lambda) \right]\, {\rm d}x \wedge {\rm d}y\; \;.
    \end{equation*}
    For every $z\in \CC\backslash \R$, by spectral decomposition theorem again, we have
    \begin{equation*}
        \scal{(z-T)^{-1}u,v} = \int_{\R}(z-\lambda)^{-1}\mu_{u,v}({\rm d} \lambda)\;.
    \end{equation*}
    This completes the proof.
\end{proof}

\begin{lem}\label{lem:HS}
    Let $f\in \cC^{\infty}_0(\R)$. Then we have
    \begin{equation*}
        \|f(T_1)-f(T_2)\|_{\hH\to \hH}\leq C \|T_1-T_2\|_{\hH\to \hH}
    \end{equation*}
    for all bounded self-adjoint operators $T_1, T_2$ on $\hH$. The proportinality constant depends on $f$, but is independent of the operators $T_1$ and $T_2$. 
\end{lem}
\begin{proof}
    All operator norms $\|\cdot\|$ appearing here are $\|\cdot\|_{\hH\to \hH}$. By Helffer-Sj\"ostrand's formula (Lemma~\ref{lem:Helffer-Sjostrand}), we have
    \begin{equation*}
        f(T_1) - f(T_2) = -\frac{1}{\pi} \iint_{\CC} \, \d_{\bar{z}} (\widetilde{f}_N)(z) \cdot \big( (z-T_1)^{-1} - (z-T_2)^{-1} \big) \, {\rm d}x \wedge {\rm d}y\;,
    \end{equation*}
    where $\tilde{f}_{N}$ is given in \eqref{eq:almost_analytic_extension}. Note that we have the bound
    \begin{equation*}
        \big\| (z-T_1)^{-1} - (z-T_2)^{-1}  \big\| \leq \|z-T_2\|^{-1} \|T_1 - T_2\| \|z-T_1^{-1}\| \lesssim \frac{\|T_1 - T_2\|}{|\text{Im} (z)|^2}\;,
    \end{equation*}
    where the second inequality follows from self-adjointness of $T_j$ so that $\|z-T_{j}\|^{-1} \lesssim |\text{Im} (z)|^{-1}$. Also note that for every $N \in \N^+$, we have
    \begin{equation*}
            \big|  (\d_{\bar{z}} \widetilde{f}_N)(z) \big| \lesssim |\text{Im} (z)|^{N}\;.
    \end{equation*}
    The proof is complete by combining the above two bounds (taking $N=3$) and plugging back to the above expression of $f(T_1) - f(T_2)$. 
\end{proof}

\begin{cor} \label{cor:convergence_L2_rate}
For every $m\in \Z$, there exists $C>0$ such that for every $\eps\in \N^{-1}$ and $j\geq -1$,
    \begin{equation*}
        \|\lL_{\eps}^{m}\Delta_{j,\eps}\Pi_0^{\perp}-\lL_{0}^{m}\Delta_{j,0}\Pi_0^{\perp}\|_{L^{2}\to L^2}\leq C \eps 2^{2j(1+m)}\;.
    \end{equation*}
\end{cor}
\begin{proof}
For $j=-1$, the desired result follows from \cite[Lemma~4.2]{Zhuge2020} and contour integrations. We now focus on $j \geq 0$. Let
\begin{equation*}
    f(x) := x^{-m}\, \varphi_0 \big( 1 / \sqrt{x} \big) 
\end{equation*}
where $\varphi_0$ is given in \eqref{e:rho_rescaling}. Then we have
    \begin{equation*}
        \lL_{\eps}^{m}\Delta_{j,\eps}\Pi_0^{\perp} = 2^{2mj} \, f(-2^{2j}\lL^{-1}_{\eps})\Pi_0^{\perp}\;.
    \end{equation*}
    Together with the fact that
    \begin{equation*}
        \|\lL_{\eps}^{-1}-\lL_{0}^{-1}\|\lesssim \eps\;,
    \end{equation*}
    applying Lemma \ref{lem:HS} to $-2^{2j}\lL^{-1}_{\eps}$ in $\Pi_0^{\perp}L^{2}(\T^d)$ gives the desired result.
\end{proof}

We are now ready to complete the proof of Theorem~\ref{th:convergence_single_block}. 

\begin{proof} [Proof of Theorem~\ref{th:convergence_single_block}]
For $\alpha \in (0,1)$, we have
\begin{equation*}
    \big\| \big(\lL_\eps^m \Delta_{j,\eps} - \lL_0^m \Delta_{j,0} \big) \big\|_{L^\infty \rightarrow \cC^\alpha} \lesssim 2^{(2m+\alpha) j}\;.
\end{equation*}
The desired claim then follows immediately by interpolating the above bound and the $L^2 \rightarrow L^2$ bound in Corollary~\ref{cor:convergence_L2_rate}. 
\end{proof}

\begin{lem} \label{le:convergence_gradient_block}
For every $\gamma>0$, there exist $\theta,\nu>0$ such that
\begin{equation} \label{eq:convergence_gradient_block}
    \big\| \nabla \Delta_{j,\eps} - \Phi_{\eps}\nabla \Delta_{j,0} \big\|_{L^\infty \rightarrow \cC^{\nu}} \lesssim \eps^\theta 2^{(1+\gamma)j}\;.
\end{equation}
As a consequence, for every $\kappa>0$ and $\gamma>0$, there exists $\theta>0$ such that
\begin{equation} \label{eq:convergence_flux_block}
    \|A_{\eps} \nabla \Delta_{j,\eps}- \bar{A} \nabla \Delta_{j,0}\|_{L^{\infty}\to \cC^{-\kappa}} \lesssim \eps^{\theta}2^{(1+\gamma)j}\;.
\end{equation}
The proportionality constants above depend on $\kappa$ and $\gamma$ but are independent of $\eps\in \N^{-1}$ and $j \geq -1$. 
\end{lem}
\begin{proof}
We first show the bound \eqref{eq:convergence_gradient_block}. By Lemma~\ref{lem:blocks_derivatives_pointwise}, we have the uniform bound
\begin{equation} \label{eq:convergence_gradient_block_uniform}
    \big\| \nabla \Delta_{j,\eps} \big\|_{L^\infty \rightarrow L^\infty} + \big\|\nabla \Delta_{j,0} \big\|_{L^\infty \rightarrow L^\infty} \lesssim 2^j\;.
\end{equation}
To get the convergence rate, we decompose the difference into
    \begin{equation*}
        \begin{split}
        \nabla \Delta_{j,\eps} - \Phi_{\eps} \nabla \Delta_{j,0}= \nabla \lL_\eps^{-1} \big( \lL_\eps \Delta_{j,\eps} - \lL_0 \Delta_{j,0} \big) + \Big( \nabla \lL_\eps^{-1} - \Phi_{\eps} \nabla \lL_0^{-1}  \Big) \lL_0 \Delta_{j,0}\;.
        \end{split}
    \end{equation*}
    We estimate the two terms on the right hand side above separately. For the first one, using the Green's function estimates for $\nabla \lL_\eps^{-1}$ and then applying Corollary~\ref{cor:convergence_single_block} with $m=1$, we get
    \begin{equation*}
        \big\| \nabla \lL_\eps^{-1} \big( \lL_\eps \Delta_{j,\eps} - \lL_0 \Delta_{j,0} \big) f \big\|_{L^\infty} \lesssim \big\| \big( \lL_\eps \Delta_{j,\eps} - \lL_0 \Delta_{j,0} \big) f \big\|_{L^\infty} \lesssim \eps^\theta 2^{3j} \|f\|_{L^\infty}\;.
    \end{equation*}
    For the second term, by Lemma~\ref{lem:elliptic_homo_convergence_gradient}, we have
    \begin{equation*}
        \Big\| \Big( \nabla \lL_\eps^{-1} - \Phi_{\eps} \nabla \lL_0^{-1}  \Big)\lL_{0}\Delta_{j,0} f \Big\|_{L^\infty} \lesssim \eps^\theta \|\lL_0 \Delta_{j,0} f\|_{\cC^1} \lesssim \eps^\theta 2^{3j} \|f\|_{L^\infty}\;.
    \end{equation*}
    Combining the above two bounds and then interpolating with \eqref{eq:convergence_gradient_block_uniform}, we deduce the bound
    \begin{equation}
        \big\| \nabla \Delta_{j,\eps} - \Phi_{\eps}\nabla \Delta_{j,0} \big\|_{L^\infty \rightarrow L^{\infty}} \lesssim \eps^\theta 2^{(1+\gamma)j}\;.
    \end{equation}
    By Schauder estimates, there exists $M>0$ such that for every $\eps\in \N^{-1}$, $j\geq -1$,
    \begin{equation*}
         \big\| \nabla \lL_{\eps}^{-1} (\lL_{\eps}\Delta_{j,\eps}) \big\|_{L^\infty \rightarrow \cC^{\frac{1}{2}}} + \big\| \Phi_{\eps} \nabla \Delta_{j,0} \big\|_{L^\infty \rightarrow \cC^{\frac{1}{2}}} \lesssim \eps^{-M}\cdot 2^{M j}\;.
    \end{equation*}
    These estimates give \eqref{eq:convergence_gradient_block}. 

    We now turn to \eqref{eq:convergence_flux_block}. We write the left hand side as
    \begin{equation*}
        \begin{split}
        A_\eps \nabla \Delta_{j,\eps} - \bar{A} \nabla \Delta_{j,0} = A_\eps \Big( \nabla \Delta_{j,\eps} - \Phi_{\eps} \nabla \Delta_{j,0} \Big)+ \Pi_0^{\perp}(A_{\eps}\Phi_{\eps})\nabla \Delta_{j,0}\;.
        \end{split}
    \end{equation*}
    For the first term on the right hand side above, it follows directly from \eqref{eq:convergence_gradient_block} that
    \begin{equation*}
        \Big\| A_\eps \left( \nabla \Delta_{j,\eps} - \Phi_{\eps}\nabla \Delta_{j,0}\right) \Big\|_{L^\infty \rightarrow L^\infty} \lesssim \eps^\theta 2^{(1+\gamma)j}\;.
    \end{equation*}
    For the second term, we have
    \begin{equation*}
        \begin{split}
        \Big\| \Pi_0^{\perp}\left( A_\eps \Phi_{\eps} \right) \nabla \Delta_{j,0} \Big\|_{L^\infty \rightarrow \cC^{-\kappa}}\lesssim \big\|\Pi_0^{\perp}\left( A_\eps \Phi_{\eps} \right) \big\|_{\cC^{-\kappa'}} \|\nabla \Delta_{j,0}\|_{L^\infty \rightarrow \cC^{2\kappa'}} \lesssim \eps^{\theta'} 2^{2j}\;.
        \end{split}
    \end{equation*}
    The desired claim \eqref{eq:convergence_flux_block} then follows from interpolation. This completes the proof. 
    \end{proof}

\subsection{Para-products}

\begin{lem} \label{le:convergence_circ_good}
Let $\alpha,\beta\in (-1,1)$ be such that $\alpha+\beta\in (0,1)$. Then for every $\kappa>0$, there exists $\theta>0$ such that
\begin{equation*}
    \|f_\eps \circ_\eps g_\eps - f_0 \circ g_0\|_{\alpha+\beta-\kappa} \lesssim \|f_0\|_{\alpha} \big( \eps^\theta \|g_0\|_\beta + \|g_\eps - g_0\|_\beta \big) + \|f_\eps - f_0\|_{\alpha} \|g_\eps\|_\beta\;,
\end{equation*}
uniformly over functions in relevant spaces and $\eps\in \N^{-1}$. 
\end{lem}
\begin{proof}
We split the difference into
\begin{equation*}
    f_\eps \circ_\eps g_\eps - f_0 \circ g_0 = (f_\eps - f_0) \circ_\eps g_\eps + f_0 \circ_\eps (g_\eps - g_0) + \big( f_0 \circ_\eps g_0 - f_0 \circ g_0 \big)\;.
\end{equation*}
For the first two terms on the right hand side, it follows directly from Proposition~\ref{pr:paraproducts_gBesov} and Lemma~\ref{lem:Holder_equivalence} that
\begin{equation*}
    \|(f_\eps - f_0) \circ_\eps g_\eps\|_{\alpha+\beta} \lesssim \|f_\eps - f_0\|_{\alpha} \|g_\eps\|_{\beta}\;, \quad \|f_0 \circ_\eps (g_\eps - g_0)\|_{\alpha+\beta} \lesssim \|f_0\|_{\alpha} \|g_\eps - g_0\|_{\beta}\;.
\end{equation*}
For the last term, on one hand, we have the uniform bound
\begin{equation} \label{e:convergence_circ_good_uniform}
    \| f_0 \circ_\eps g_0 - f_0 \circ g_0 \|_{\alpha+\beta} \lesssim \|f_0\|_{\alpha} \|g_0\|_{\beta}\;.
\end{equation}
For smallness in $\eps$, we expand the difference as
\begin{equation*}
    f_0 \circ_\eps g_0 - f_0 \circ g_0 = \sum_{|i-j| \leq 1} \Big( \big( \Delta_{i,\eps} - \Delta_{i,0} \big) f_0 \cdot \Delta_{j,\eps} g_0 + \Delta_{i,0} f_0 \cdot \big( \Delta_{j,\eps} - \Delta_{j,0} \big) g_0 \Big)\;.
\end{equation*}
By Lemma~\ref{lem:Holder_equivalence}, Corollary~\ref{cor:convergence_single_block} (with $m=0$) and the assumption $\alpha+\beta>0$, we can control its $L^\infty$-norm by
\begin{equation*}
    \| f_0 \circ_\eps g_0 - f_0 \circ g_0 \|_{L^\infty} \lesssim \eps^{\theta'} \|f_0\|_\alpha \|g_0\|_\beta
\end{equation*}
for some $\theta'>0$. Interpolating it with the uniform bound \eqref{e:convergence_circ_good_uniform} gives
\begin{equation*}
    \| f_0 \circ_\eps g_0 - f_0 \circ g_0 \|_{\alpha+\beta-\kappa} \lesssim \eps^{\theta} \|f_0\|_{\alpha} \|g_0\|_{\beta}\;.
\end{equation*}
This completes the proof. 
\end{proof}

We now give a situation where the functions involved in the resonance products have regularities outside the range of equivalence of norms as in Lemma~\ref{lem:Holder_equivalence}. 

\begin{lem}
\label{le:convergence_circ_outside1}
Let $\alpha \in (-1,0)$ and $\beta \in (0,1)$ be such that $\alpha + \beta > 0$. Then for every $\kappa>0$, there exists $\theta>0$ such that
\begin{equation} \label{e:convergence_circ_outside1}
    \left\| \big( \lL_{\eps}^{-1} \Pi_0^\perp f\big) \circ_{\eps} (\lL_{\eps}g) - \big( \lL_{0}^{-1} \Pi_0^\perp f \big)\circ (\lL_{0}g) \right\|_{\alpha+\beta-\kappa} \lesssim \eps^{\theta}\|\Pi_0^\perp f\|_{\alpha}\|g\|_{\beta}\;,
\end{equation}
where the proportionality constant is independent of $\eps\in \N^{-1}$, $j \geq -1$ and $f,g$ in respective function classes. 
\end{lem}
\begin{proof}
We assume without loss of generality that $\Pi_0 f = 0$, so that $\Pi_0^\perp f = f$. The assumption on $\alpha$ and $\beta$ ensures $\alpha+\beta \in (0,1)$. Hence, by Proposition~\ref{pr:paraproducts_gBesov} and Lemma~\ref{lem:Holder_equivalence}, we have the bound
\begin{equation*}
\left\| \big( \lL_{\eps}^{-1} f\big) \circ_{\eps} (\lL_{\eps}g) \right\|_{\alpha+\beta} \lesssim  \|\lL_\eps^{-1} f \|_{\alpha+2;\eps} \|\lL_\eps g\|_{\beta-2;\eps} \lesssim \|f\|_\alpha \|g\|_\beta
\end{equation*}
uniformly over $\eps \in \N^{-1}$, where we have used the equivalence of norms for $\alpha$, $\beta$ and $\alpha+\beta$. 

On the other hand, by definition of $\circ_\eps$ and $\circ$, we have the expression
\begin{equation*}
    \begin{split}
    &\phantom{1111111}\big( \lL_{\eps}^{-1} f\big) \circ_{\eps} (\lL_{\eps}g) - \big( \lL_0^{-1} f\big) \circ (\lL_0 g)\\
    &= \sum_{|i-j| \leq 1} \Big( \big( \Delta_{j,\eps} \lL_\eps^{-1} - \Delta_{j,0} \lL_0^{-1} \big) f \cdot \Delta_{i,\eps} \lL_\eps g + \Delta_{j,0} \lL_0^{-1} f \cdot (\Delta_{i,\eps} \lL_\eps - \Delta_{i,0} \lL_0) g \Big)\;.
    \end{split}
\end{equation*}
Applying Corollary~\ref{cor:convergence_single_block} with $m = \pm 1$, we get
\begin{equation*}
    \| \big( \lL_{\eps}^{-1} f\big) \circ_{\eps} (\lL_{\eps}g) - \big( \lL_0^{-1} f\big) \circ (\lL_0 g) \|_{L^\infty} \lesssim \eps^{\theta'} \|f\|_\alpha \|g\|_\beta\
\end{equation*}
for some $\theta'>0$. The desired bound \eqref{e:convergence_circ_outside1} then follows from interpolation. 
\end{proof}  

% For $f_\eps \in \cC^{\alpha;\eps}$, we introduce the notation
% \begin{equation*}
%     F_\eps := \Pi_0 f_\eps - (-\lL_\eps)^{-1} \Pi_0^\perp f_\eps\;,
% \end{equation*}
% so that we have
% \begin{equation*}
%     f_\eps = \Pi_0 f_\eps - \lL_\eps F_\eps\;.
% \end{equation*}
% With the above notation, we have the following lemma. 

\begin{lem} \label{lem:homo_diff_circ_estimate}
Let $\alpha \in (1,2)$ and $\beta \in (-2,-1)$ be such that $\alpha+\beta>0$. Then for every $\kappa>0$, there exists $\theta>0$ such that
\begin{equation*}
    \begin{split}
    \| f_{\eps} \circ_{\eps} g_{\eps} - f_{0} \circ g_0\|_{\alpha+\beta-\kappa} \lesssim (\eps^\theta \| f_\eps\|_{\alpha;\eps}+|\!|\!|f_\eps;f_0|\!|\!|_{\alpha;\eps}) \|g_\eps\|_{\beta;\eps}+ 
 \| f_0\|_{\alpha;0} |\!|\!|g_\eps;g_0|\!|\!|_{\beta;\eps}\;,
    \end{split}
\end{equation*}
where the norms $|\!|\!|f_\eps;f_0|\!|\!|_{\alpha;\eps})$ and $|\!|\!|g_\eps;g_0|\!|\!|_{\beta;\eps}$ are given in \eqref{e:norm_comparison_high_reg} and \eqref{e:norm_comparison_low_reg}. 
\end{lem}
\begin{proof}
Writing $F_\eps = \lL_\eps f_\eps$ and $G_\eps = \lL_\eps^{-1} \Pi_0^\perp g_\eps$, we have
\begin{equation*}
    f_{\eps} = \Pi_{0} f_{\eps} +\lL_{\eps}^{-1} F_\eps\;, \quad g_{\eps} = \Pi_{0} g_{\eps} + \lL_{\eps} G_\eps\;.
\end{equation*}
We have the expression
\begin{equation} \label{e:homo_diff_circ_decomposition}
    \begin{split}
    &\phantom{111}f_\eps \circ_\eps g_\eps - f_0 \circ g_0\\
    &= \Big( \Pi_0 f_\eps \circ_\eps \Pi_0 g_\eps - \Pi_0 f_0 \circ \Pi_0 g_0 \Big) + \Big( \Pi_0 f_\eps \circ_\eps \lL_\eps G_\eps - \Pi_0 f_0 \circ \lL_0 G_0 \Big)\\
    &\phantom{111}+ \Big( \lL_\eps^{-1} F_\eps \circ_\eps \Pi_0 g_\eps - \lL_0^{-1} F_0 \circ \Pi_0 g_0 \Big) + \Big( \lL_\eps^{-1} F_\eps \circ_\eps \lL_\eps G_\eps - \lL_0^{-1} F_0 \circ \lL_0 G_0 \Big)\;.
    \end{split}
\end{equation}
We control the four terms on the right hand side above one by one. The first one is just multiplication of two constants, so it is straightforward that
\begin{equation*}
    \left\| \Pi_0 f_\eps \circ_\eps \Pi_0 g_\eps - \Pi_0 f_0 \circ \Pi_0 g_0 \right\|_{\alpha+\beta-\kappa} \leq \left| \Pi_0 (f_\eps - f_0) \right| \cdot |\Pi_0 g_\eps| + |\Pi_0 f_0| \cdot \left| \Pi_0 (g_\eps - g_0) \right|\;.
\end{equation*}
The second one equals
\begin{equation*}
    \Pi_0 (f_\eps - f_0) \cdot S_{3,\eps} \lL_\eps G_\eps + \Pi_0 f_0 \cdot (S_{3,\eps} \lL_\eps - S_{3,0} \lL_0) G_\eps + \Pi_0 f_0 \cdot S_{3,0} \lL_0 (G_\eps - G_0)\;.
\end{equation*}
The $\|\cdot\|_{\alpha+\beta-\kappa}$-norm for the first and third term above can be controlled by
\begin{equation*}
    |\Pi_0 (f_\eps - f_0)| \cdot \|G_\eps\|_{\beta+2} + |\Pi_0 f_0| \cdot \|G_\eps - G_0\|_{\beta+2}\;,
\end{equation*}
while by Corollary~\ref{cor:convergence_single_block} (with $m=0$) and Lemma~\ref{lem:blocks_derivatives_pointwise}, the middle term can be controlled by (for some $\gamma \in (0,1)$)
\begin{equation*}
    |\Pi_0 f_0| \cdot  \left\| (S_{3,\eps} \lL_\eps - S_{3,0} \lL_0) G_\eps \right\|_{L^\infty}^\gamma \left\| \nabla (S_{3,\eps} \lL_\eps - S_{3,0} \lL_0) G_\eps \right\|_{L^\infty}^{1-\gamma} \lesssim \eps^\theta |\Pi_0 f_0| \cdot \|G_\eps\|_{\beta+2}\;.
\end{equation*}
Similarly, for the third term on the right hand side of \eqref{e:homo_diff_circ_decomposition}, we can write it as
\begin{equation*}
    \begin{split}
        \Pi_{0}g_{\eps}\big((S_{3,\eps}\lL_{\eps}^{-1}-S_{3,0}\lL_{0}^{-1})\lL_{\eps}f_{\eps} + S_{3,0}\lL_{0}^{-1}(\lL_{\eps}f_{\eps}-\lL_{0}f_{0})\big)+\Pi_{0}(g_{\eps}-g_0)S_{3,0}\Pi_0^{\perp}f_{0}\;.
    \end{split}
\end{equation*}
Similar as before, the $\|\cdot\|_{\alpha+\beta-\kappa}$-norms of the latter two terms can be controlled by
\begin{equation*}
    |\Pi_0 g_\eps| \cdot \|\lL_{\eps}f_\eps - \lL_{0}f_0\|_{\alpha-2} + |\Pi_0 (g_\eps - g_0)| \cdot \|\Pi_0^\perp f_0\|_{\alpha}\;,
\end{equation*}
while the first term can be bounded by (for some $\gamma \in (0,1)$)
\begin{equation*}
    \begin{split}
    &\phantom{111}|\Pi_0 g_\eps| \cdot \| \big( \Delta_{j,\eps}\lL_\eps^{-1} - \Delta_{j,0} \lL_0^{-1} \big) \lL_{\eps} f_\eps \|_{L^\infty}^{\gamma} \| \nabla \big( \Delta_{j,\eps}\lL_\eps^{-1} - \Delta_{j,0} \lL_0^{-1} \big)\lL_{\eps} f_\eps \|_{L^\infty}^{1-\gamma}\\
    &\lesssim \eps^\theta |\Pi_0 g_\eps| \cdot \|\lL_{\eps} f_\eps\|_{\alpha-2}\;.
    \end{split}
\end{equation*}
Finally we turn to the last term on the right hand side of \eqref{e:homo_diff_circ_decomposition}. We decompose it into the sum of two difference terms as
\begin{equation*}
    \Big( \lL_\eps^{-1} F_\eps \circ_\eps \lL_\eps G_\eps - \lL_0^{-1} F_\eps \circ \lL_0 G_\eps \Big) + \Big( \lL_0^{-1} F_\eps \circ \lL_0 G_\eps -  \lL_0^{-1} F_0 \circ \lL_0 G_0 \Big)\;.
\end{equation*}
By Lemma~\ref{le:convergence_circ_outside1} and Proposition~\ref{pr:paraproducts_gBesov}, the $\|\cdot\|_{\alpha+\beta-\kappa}$-norms of the two terms can be together bounded by
\begin{equation*}
    \eps^\theta \|\lL_{\eps} f_\eps\|_{\alpha-2} \|G_\eps\|_{\beta+2} + \| \lL_{\eps}f_\eps - \lL_{0} f_0\|_{\alpha-2} \|G_\eps\|_{\beta+2} + \|\lL_{0} f_0\|_{\alpha-2} \|G_\eps - G_0\|_{\beta+2}\;.
\end{equation*}
The proof of the lemma is then complete by combining all the four bounds above. 
\end{proof}

\begin{lem} \label{le:convergence_prec_good}
    Let $\alpha \in (-1,0)$ and $\beta \in (-1,1)$ be such that $\alpha + \beta \in (-1,1)$. Then for every $\kappa>0$, there exists $\theta>0$ such that
    \begin{equation} \label{e:convergence_prec_good}
    	\|f_\eps \prec_{\eps} g_\eps - f_0 \prec g_0\|_{\alpha+\beta-\kappa} \lesssim \big( \eps^{\theta} \|g_0\|_{\beta} + \|g_\eps - g_0\|_{\beta} \big) \|f_0\|_\alpha + \|f_\eps - f_0\|_{\alpha} \, \|g_\eps\|_\beta\;.
    \end{equation}
\end{lem}
\begin{proof}
We first deal with the special situation when $f_\eps = f_0 = f$ and $g_\eps = g_0 = g$. By definition, we need to show that for every $\kappa>0$, there exists $\theta>0$ such that
\begin{equation*}
    \left\| \Delta_{k,0} \big( f \prec_\eps g - f \prec g \big) \right\|_{L^\infty} \lesssim \eps^\theta \cdot 2^{-(\alpha+\beta-\kappa) k} \|f\|_{\alpha} \|g\|_\beta\;.
\end{equation*}
We write
\begin{equation} \label{e:convergence_prec_good_decomposition}
    \begin{split}
    \Delta_{k,0}(f\prec_{\eps}g-f \prec g ) = &(\Delta_{k,0} - \Delta_{k,\eps})( f \prec_\eps g)\\
    &+ \Big( \Delta_{k,\eps }( f \prec_{ \eps } g )-\Delta_{k,0} ( f \prec g) \Big)\;.
    \end{split}
\end{equation}
For the first part, it follows from Corollary~\ref{cor:convergence_single_block} with $m=0$ and Proposition~\ref{pr:paraproducts_gBesov} that
\begin{equation*}
    \left\| (\Delta_{k,0}-\Delta_{k,\eps})( f \prec_\eps g) \right\|_{L^\infty} \lesssim \eps^\theta 2^{-(\alpha+\beta-\kappa) k} \|f \prec_\eps g\|_{\alpha+\beta} \lesssim \eps^\theta 2^{-(\alpha+\beta-\kappa) k} \|f\|_\alpha \|g\|_\beta\;.
\end{equation*}
Here, we have used the equivalence of norms for $\alpha+\beta \in (-1,1)$. We now turn to the second part in \eqref{e:convergence_prec_good_decomposition}. By definition, we have
\begin{equation} \label{e:converegence_prec_good_expansion}
    \Delta_{k,\eps }( f \prec_{ \eps } g ) -\Delta_{k,0} ( f \prec g)  = \sum_{j \geq -1} \Big( \Delta_{k,\eps} \big( S_{j,\eps} f \cdot \Delta_{j,\eps} g \big) - \Delta_{k,0} \big( S_{j,0} f \cdot \Delta_{j,0} g \big) \Big)\;.
\end{equation}
By the assumption $\alpha<1$, Corollaries~\ref{cor:decoupling_spectrum_S} and~\ref{cor:convergence_single_block} (with $m=0$) and equivalence of norms, we have
\begin{equation*}\begin{split}
    &\left\| \Delta_{k,\eps }( S_{j,\eps} f \cdot \Delta_{j,\eps} g )  - \Delta_{k,0} \big( S_{j,0} f \cdot \Delta_{j,0} g \big) \right\|_{L^\infty}\\
    \lesssim &
    \begin{cases}
\eps^\theta \, 2^{(\kappa-2)k} \, 2^{(2-\alpha-\beta)j} \, \|f\|_\alpha \|g\|_\beta, & \text{if} \; j \leq k-2\;;\\
\eps^\theta \, 2^k \, 2^{-(\alpha+\beta+1-\kappa)j} \, \|f\|_\alpha \|g\|_\beta, & \text{if} \; j \geq k + 2\;.
\end{cases}
\end{split}\end{equation*}
 For $j$ such that $|j-k| \leq 1$, by Corollary~\ref{cor:convergence_single_block} (with $m=0$), the assumption $\alpha<0$ and that $\|\Delta_{k,0}\|_{L^\infty \rightarrow L^\infty} \lesssim 1$, we have
\begin{equation*}
    \begin{split}
    \| \Delta_{ k , \eps }( S_{j,\eps} f \cdot &\Delta_{j,\eps} g)- \Delta_{ k , 0 }( S_{j,0} f\cdot\Delta_{j,0} g)\|_{L^{\infty}} \leq \left\| (\Delta_{k,\eps} - \Delta_{k,0}) ( S_{j,\eps} f \cdot \Delta_{j,\eps} g) \right\|_{L^\infty}\\
    &+ \left\| \Delta_{k,0} \big( (S_{j,\eps} - S_{j,0}) f \cdot \Delta_{j,\eps} g \big) \right\|_{L^\infty} + \left\| \Delta_{k,0} \big( S_{j,0} f \cdot (\Delta_{k,\eps} - \Delta_{k,0}) g \big) \right\|_{L^\infty}\\
    &\lesssim \eps^\theta 2^{-(\alpha+\beta-\kappa)k} \|f\|_\alpha \|g\|_\beta\;.
    \end{split}
\end{equation*}
Substituting the above bounds back into \eqref{e:converegence_prec_good_expansion} gives
\begin{equation*}
    \left\| \Delta_{k,\eps }( f \prec_{ \eps } g ) -\Delta_{k,0} ( f \prec g) \right\|_{L^\infty} \lesssim \eps^\theta 2^{-(\alpha+\beta-\kappa) k} \|f\|_\alpha \|g\|_\beta
\end{equation*}
for arbitrary $\kappa>0$ and $\theta>0$ depending on $\kappa$. This proves \eqref{e:convergence_prec_good} when $f_\eps = f_0$ and $g_\eps = g_0$. To treat the general situation, we just write
\begin{equation*}
    f_\eps \prec_\eps g_\eps - f_0 \prec g_0 = (f_\eps - f_0) \prec_\eps g_\eps + f_0 \prec_\eps (g_\eps - g_0) + \big( f_0 \prec_\eps g_0 - f_0 \prec g_0 \big)\;,
\end{equation*}
and apply Proposition~\ref{pr:paraproducts_gBesov} for the first two terms and the previous bound for the special case to the third term. This completes the proof of the lemma. 
\end{proof}

\begin{lem} \label{le:convergence_prec_reg-2}
Let $\alpha, \beta \in (0,1)$ be such that $\alpha+\beta>1$. Then for every $\kappa>0$, there exists $\theta>0$ such that
\begin{equation}\label{e:convergence_prec_reg-2}
    \|(\lL_\eps f) \prec_{\eps} g-(\lL_0 f) \prec g\|_{\alpha+\beta-2-\kappa} \lesssim \eps^{\theta}\|f\|_{\alpha}\|g\|_{\beta}\;,
    \end{equation}
uniformly over $\eps\in \N^{-1}$, $f \in \cC^\alpha$ and $g \in \cC^\beta$. 
\end{lem}
\begin{proof}
By definition, we need to show for every $\kappa>0$, there exists $\theta>0$ such that
\begin{equation} \label{e:convergence_prec_reg-2_block}
    \left\| \Delta_{k,0} \big( (\lL_\eps f) \prec_{\eps} g-(\lL_0 f) \prec g \big) \right\|_{L^\infty} \lesssim \eps^\theta 2^{(2-\alpha-\beta+\kappa) k} \|f\|_\alpha \|g\|_\beta\;.
\end{equation}
We write the left hand side of \eqref{e:convergence_prec_reg-2_block} as
\begin{equation} \label{e:convergence_prec_reg-2_decomposition}
    \begin{split}
    \Delta_{k,0} \big( (\lL_\eps f) \prec_{\eps} g- &(\lL_0 f) \prec g \big) = (\Delta_{k,0} - \Delta_{k,\eps}) \big( (\lL_\eps f) \prec_\eps g \big)\\
    &+ \Big( \Delta_{k,\eps} \big( (\lL_\eps f) \prec_{\eps} g \big) - \Delta_{k,0} \big( (\lL_0 f) \prec g \big)  \Big)\;.
    \end{split}
\end{equation}
By Corollary~\ref{cor:convergence_single_block} (with $m=0$) and Proposition~\ref{pr:paraproducts_gBesov}, the first term above can be directly bounded by
\begin{equation*}
    \begin{split}
    \left\| (\Delta_{k,0} - \Delta_{k,\eps}) \big( (\lL_\eps f) \prec_\eps g \big) \right\|_{L^\infty} &\lesssim \eps^\theta 2^{(2-\alpha-\beta+\kappa)k} \|(\lL_\eps f) \prec_\eps g\|_{\alpha+\beta-2}\\
    &\lesssim \eps^\theta 2^{(2-\alpha-\beta+\kappa)k} \|f\|_\alpha \|g\|_\beta\;.
    \end{split}
\end{equation*}
For the second term on the right hand side of \eqref{e:convergence_prec_reg-2_decomposition}, by definition, we have the expansion
\begin{equation} \label{e:convergence_prec_reg-2_expansion}
    \begin{split}
    &\phantom{111}\Delta_{k,\eps} \big( (\lL_\eps f) \prec_{\eps} g \big) - \Delta_{k,0} \big( (\lL_0 f) \prec g \big)\\
    &= \sum_{j \geq -1} \Big( \Delta_{k,\eps} \big( ( S_{j,\eps} \lL_\eps f ) \cdot (\Delta_{j,\eps} g) \big) - \Delta_{k,0} \big( ( S_{j,0} \lL_0 f ) \cdot (\Delta_{j,0} g) \big) \Big)\;.
    \end{split}
\end{equation}
By Corollaries~\ref{cor:decoupling_spectrum_S} and~\ref{cor:convergence_single_block} (with $m=0$) and that $\alpha-2<0$, we have
\begin{equation}\begin{split} \label{e:convergence_prec_reg-2_uniform}
    &\left\| \Delta_{k,\eps } \big( S_{j,\eps} \lL_\eps f \cdot \Delta_{j,\eps} g \big) - \Delta_{k,0} \big( ( S_{j,0} \lL_0 f ) \cdot (\Delta_{j,0} g) \big) \right\|_{L^\infty}\\
    \lesssim &
    \begin{cases}
\eps^\theta \, 2^{(\kappa-2)k} \, 2^{(4-\alpha-\beta)j} \, \|f\|_{\alpha} \|g\|_\beta, & \text{if} \; j \leq k-2\;;\\
\eps^\theta \, 2^k \, 2^{(1-\alpha-\beta+\kappa)j} \, \|f\|_{\alpha} \|g\|_\beta, & \text{if} \; j \geq k + 2\;,
\end{cases}
\end{split}\end{equation}
where we also used $\|\lL_\eps f\|_{\alpha-2;\eps} \lesssim \|f\|_\alpha$ and equivalence of norms for $\alpha, \beta \in (0,1)$.

On the other hand, similar to the proof in Lemma~\ref{le:convergence_prec_good}, we can also use triangle inequality and apply Corollary~\ref{cor:convergence_single_block} (with $m=0$) to control
\begin{equation*}
    \left\| \Delta_{k,\eps} \big( ( S_{j,\eps} \lL_\eps f ) \cdot (\Delta_{j,\eps} g) \big) - \Delta_{k,0} \big( ( S_{j,0} \lL_0 f ) \cdot (\Delta_{j,0} g) \big) \right\|_{L^\infty} \lesssim \eps^{\theta} 2^{(2-\alpha-\beta+\kappa)k} \|f\|_\alpha \|g\|_\beta
\end{equation*}
for $|j-k| \leq 1$. Substituting the above bound as well as \eqref{e:convergence_prec_reg-2_uniform} back into the sum \eqref{e:convergence_prec_reg-2_expansion} gives the desired bound for the second term on the right hand side of \eqref{e:convergence_prec_reg-2_decomposition}. We have thus proved \eqref{e:convergence_prec_reg-2_block} and hence the lemma. 
\end{proof}

\begin{lem} \label{le:convergence_prec_below}
    Let $\alpha \in (-2,-1)$ and $\beta \in (0,1)$ be such that $\alpha+\beta>-1$. Then for every $\kappa>0$, there exists $\theta>0$ such that
    \begin{equation} \label{e:homo_diff_prec_less_minus_1}
        \| f_\eps \prec_\eps g_\eps - f_0 \prec g_0 \|_{\alpha+\beta-\kappa} \lesssim \big( \eps^\theta \|f_0\|_{\alpha}+ |\!|\!| f_\eps;f_0 |\!|\!|_{\alpha} \big) \|g_0\|_\beta + \|f_\eps\|_{\alpha;\eps} \, \|g_\eps-g_0\|_\beta\;.
    \end{equation}
\end{lem}
\begin{proof}
Similar as before, we decompose
\begin{equation*}
    f_\eps \prec_\eps g_\eps - f_0 \prec g_0 = f_\eps \prec_\eps (g_\eps - g_0) + \big( f_\eps \prec_\eps g_0 - f_0 \prec g_0 \big)\;.
\end{equation*}
The first term can be controlled directly with Proposition~\ref{pr:paraproducts_gBesov}. Hence, it suffices to control the second term. We write $g_0 = g$ for simplicity. 

Writing $F_\eps := \lL_\eps^{-1} \Pi_0^\perp f_\eps$ and $F_0 := \lL_0^{-1} \Pi_0^\perp f_0$, we have
\begin{equation*}
    \begin{split}
    f_\eps \prec_\eps g - f_0 \prec g = &\big( \Pi_0 (f_\eps - f_0) \big) \prec_\eps g + \big( \Pi_0 f_0 \prec_\eps g - \Pi_0 f_0 \prec g \big)\\
    &+ \big( \lL_\eps (F_\eps - F_0) \big) \prec_\eps g + \big( (\lL_\eps F_0) \prec_\eps g - (\lL_0 F_0) \prec g \big)\;.
    \end{split}
\end{equation*}
For the first term, by Proposition~\ref{pr:paraproducts_gBesov}, we have
\begin{equation*}
    \|\Pi_0 (f_\eps - f_0) \prec_\eps g\|_{\alpha+\beta} \leq |\Pi_0 (f_\eps - f_0)| \cdot \|1 \prec_\eps g\|_{\alpha+\beta} \lesssim |\Pi_0 (f_\eps - f_0)| \cdot \|g\|_\beta\;.
\end{equation*}
For the second term, it follows from Lemma~\ref{le:convergence_prec_good} that
\begin{equation*}
    \left\| \Pi_0 f_0 \prec_\eps g - \Pi_0 f_0 \prec g \right\|_{\alpha+\beta-\kappa} \lesssim \eps^\theta |\Pi_0 f_0| \cdot \|g\|_\beta\;,
\end{equation*}
where we have used that $\alpha+\beta-\kappa>-1$. For the third term, by Proposition~\ref{pr:paraproducts_gBesov} and that $\alpha+\beta>-1$, we have
\begin{equation*}
    \left\| \big( \lL_\eps (F_\eps - F_0) \big) \prec_\eps g \right\|_{\alpha+\beta} \lesssim \|\lL_\eps (F_\eps - F_0)\|_{\alpha;\eps} \|g\|_\beta \lesssim \|F_\eps - F_0\|_{\alpha+2} \|g\|_\beta\;.
\end{equation*}
Finally, by Lemma~\ref{le:convergence_prec_reg-2}, we have
\begin{equation*}
    \left\|  (\lL_\eps F_0) \prec_\eps g - (\lL_0 F_0) \prec g \right\|_{\alpha+\beta-\kappa} \lesssim \eps^\theta \|F_0\|_{\alpha+2} \|g\|_\beta\;.
\end{equation*}
This completes the proof of the lemma. 
\end{proof}

Recall from \eqref{e:operator_Pa} that
\begin{equation*}
    \P_\eps (f,g) := \sum_{j=-1}^{+\infty} A_\eps \nabla S_{j,\eps} f \cdot \nabla \Delta_{j,\eps} g\;, \quad \P_0 (f,g) := \sum_{j=-1}^{+\infty} \bar{A} \nabla S_{j,0} f \cdot \nabla \Delta_{j,0} g
\end{equation*}
We have the following bound. 

\begin{lem} \label{le:convergence_prec_alternative}
Let $\alpha, \beta \in (0,1)$ be such that $\alpha + \beta > 1$. Then for every $\kappa>0$, there exists $\theta>0$ such that
\begin{equation} \label{e:convergence_prec_alternative}
    \begin{split}
	\| \P_{\eps}(f_\eps , g_\eps) - \P_{0}(f_0 , g_0)\|_{\alpha+\beta-2-\kappa } \lesssim &\big( \eps^{\theta} \|g_0\|_{\beta} + \|g_\eps - g_0\|_\beta \big) \|f_0\|_\alpha\\
    &+ \|f_\eps - f_0\|_\alpha \, \|g_\eps\|_\beta\;,
    \end{split}
\end{equation}
where the proportionality constant is independent of $\eps\in \N^{-1}$ and functions in their relevant classes. 
\end{lem}
\begin{proof}
Similar as before, it suffices to consider the case $f_\eps = f_0 = f$ and $g_\eps = g_0 = g$. The general case then follows from the special one, the decomposition
\begin{equation*}
    \P_\eps (f_\eps , g_\eps) - \P_0 (f_0 , g_0) = \P_\eps (f_\eps - f_0 , g_\eps) + \P_\eps (f_0 , g_\eps - g_0) + \big( \P_\eps (f_0 , g_0) - \P_0 (f_0 , g_0) \big)
\end{equation*}
as well as Lemma~\ref{le:prec_alternative_bound}. For the special case $f_\eps = f_0 = f$ and $g_\eps = g_0 = g$, we need to show that
\begin{equation*}
    \left\| \Delta_{k,0} \big( \P_{\eps}(f,g)- \P_{0}(f,g) \big) \right\|_{L^\infty} \lesssim \eps^\theta 2^{(2-\alpha-\beta+\kappa) k} \|f\|_\alpha \|g\|_\beta
\end{equation*}
for arbitrary $\kappa > 0$ and some $\theta>0$ depending on $\kappa$. We first decompose it as
\begin{equation*}
    \Delta_{k,0} \big( \P_{\eps}(f,g)- \P_{0}(f,g) \big) = (\Delta_{k,0} - \Delta_{k,\eps}) \P_\eps(f,g) + \big( \Delta_{k,\eps} \P_\eps (f,g) - \Delta_{k,0} \P_0(f,g) \big)\;.
\end{equation*}
For the first term on the left hand side above, by Corollary~\ref{cor:convergence_single_block} (with $m=0$) and Lemma~\ref{le:prec_alternative_bound}, we have
\begin{equation*}
    \left\| (\Delta_{k,0} - \Delta_{k,\eps}) \P_\eps(f,g) \right\|_{L^\infty} \lesssim \eps^{\theta} 2^{-(\alpha+\beta-2-\kappa)k} \|\P_\eps (f,g)\|_{\alpha+\beta-2} \lesssim \eps^{\theta} 2^{-(\alpha+\beta-2-\kappa)} \|f\|_\alpha \|g\|_\beta\;.
\end{equation*}
For the second term, we write it as
\begin{equation*}
    \begin{split}
    &\phantom{111}\Delta_{k,\eps} \P_\eps (f,g) - \Delta_{k,0} \P_0(f,g)\\
    &= \sum_{j \geq -1} \Big( \Delta_{k,\eps} \big( A_\eps \nabla S_{j,\eps} f \cdot \nabla \Delta_{j,\eps} g \big) - \Delta_{k,0} \big( \bar{A} \nabla S_{j,0} f \cdot \nabla \Delta_{j,0} g \big) \Big)\;.
    \end{split}
\end{equation*}
By the proof in Lemma~\ref{le:prec_alternative_bound}, we have
\begin{equation} \label{e:convergence_prec_alternative_uniform}
    \left\| \Delta_{k,\eps} \big( A_\eps \nabla S_{j,\eps} f \cdot \nabla \Delta_{j,\eps} g \big) \right\|_{L^\infty} \lesssim
    \begin{cases}
2^{(2-\alpha-\beta)j} \|f\|_\alpha \|g\|_\beta, & \text{if} \; j \leq k\;;\\
2^{k +(1-\alpha-\beta)j} \|f\|_\alpha \|g\|_\beta, & \text{if} \; j \geq k+1\;.\\
\end{cases}
\end{equation}
On the other hand, using symmetry of $A_\eps$ and $\bar{A}$, we have
\begin{equation*}
    A_\eps \nabla S_{j,\eps} f \cdot \nabla \Delta_{j,\eps} g = \frac{1}{2} \Big( \lL_\eps \big( S_{j,\eps} f \cdot \Delta_{j,\eps} g \big) - (S_{j,\eps} \lL_\eps f) \cdot (\Delta_{j,\eps} g) - (S_{j,\eps} f) \cdot (\Delta_{j,\eps} \lL_\eps g) \Big)\;.
\end{equation*}
We can also taking the difference term-wise to get the very loose bound
\begin{equation*}
    \left\| \Delta_{k,\eps} \big( A_\eps \nabla S_{j,\eps} f \cdot \nabla \Delta_{j,\eps} g \big) - \Delta_{k,0} \big( \bar{A} \nabla S_{j,0} f \cdot \nabla \Delta_{j,0} g \big)  \right\|_{L^\infty} \lesssim \eps^{\theta'} 2^{M(j+k)} \|f\|_\alpha \|g\|_\beta
\end{equation*}
for some $\theta'$ and $M>0$. Interpolating this bound with \eqref{e:convergence_prec_alternative_uniform} and then summing over $j \geq -1$ gives
\begin{equation*}
    \left\| \Delta_{k,\eps} \P_\eps (f,g) - \Delta_{k,0} \P_0(f,g)  \right\|_{L^\infty} \lesssim \eps^\theta 2^{(2-\alpha-\beta+\kappa) k} \|f\|_\alpha \|g\|_\beta\;.
\end{equation*}
This completes the proof. 
\end{proof}

\begin{lem} \label{lem:prec_normal_high}
Let $\alpha \in (0,1)$ and $\beta \in (1,2)$. Then for every $\kappa>0$, there exists $\theta>0$ such that
\begin{equation*}
    \begin{split}
    &\phantom{111}\|\lL_\eps (f_\eps \prec_\eps g_\eps) - \lL_0 (f_0 \prec g_0)\|_{(\alpha-1) \wedge (\beta-2) - 2\kappa}\\
    &\lesssim \big( \eps^\theta \|\lL_0 g_0\|_{\beta-2} + \|\lL_\eps g_\eps - \lL_0 g_0\|_{\beta-2} \big) \|f_0\|_\alpha + \|f_\eps - f_0\|_\alpha \|\lL_\eps g_\eps\|_{\beta-2}\;.
    \end{split}
\end{equation*}
\end{lem}
\begin{proof}
By definition, we have
\begin{equation*}
    \lL_\eps (f_\eps \prec_\eps g_\eps) = (\lL_\eps f_\eps) \prec_\eps g_\eps + f_\eps \prec_\eps (\lL_\eps g_\eps) + 2 \P_\eps(f_\eps; \, g_\eps)\;.
\end{equation*}
We control the differences one by one. For the second one, we have
\begin{equation*}
    \begin{split}
    &\phantom{111}\|f_\eps \prec_\eps \lL_\eps g_\eps - f_0 \prec \lL_0 g_0\|_{\beta-2-\kappa}\\
    &\lesssim \big( \eps^\theta \|\lL_0 g_0\|_{\beta-2} + \|\lL_\eps g_\eps - \lL_0 g_0\|_{\beta-2} \big) \|f_0\|_\alpha + \|f_\eps - f_0\|_\alpha \|\lL_\eps g_\eps\|_{\beta-2}\;.
    \end{split}
\end{equation*}
For the first term, by Lemma~\ref{le:convergence_prec_below}, we have
\begin{equation*}
    \begin{split}
    &\phantom{111}\|(\lL_\eps f_\eps) \prec_\eps g_\eps - (\lL_0 f_0) \prec g_0\|_{\alpha-1-\kappa}\\
    &\lesssim \big( \eps^\theta \|\lL_0 f_0\|_{\alpha-2; 0} + |\!|\!| \lL_\eps f_\eps; \, \lL_0 f_0 |\!|\!|_{\alpha-2,\, \eps} \big) \|g_0\|_{1-\kappa} + \|\lL_\eps f_\eps\|_{\alpha-2;\eps} \|g_\eps - g_0\|_{1-\kappa}\\
    &\lesssim \big( \eps^\theta \|f_0\|_{\alpha; 0} + \| f_\eps - f_0 \|_{\alpha} \big) \|g_0\|_{1-\kappa} + \|f_\eps\|_{\alpha} \|g_\eps - g_0\|_{1-\kappa}\;.
    \end{split}
\end{equation*}
For the third term, by Lemma~\ref{le:convergence_prec_alternative}, we have
\begin{equation*}
    \begin{split}
    \|\P_\eps(f_\eps; g_\eps) - \P_0(f_0; g_0)\|_{\alpha-1-2\kappa} &\lesssim \big( \eps^\theta \|g_0\|_{1-\kappa} + \|g_\eps - g_0\|_{1-\kappa} \big) \|f_0\|_\alpha\\
    &\phantom{111}+ \|f_\eps - f_0\|_\alpha \|g_\eps\|_{1-\kappa}\;.
    \end{split}
\end{equation*}
Note that in order to apply Lemmas~\ref{le:convergence_prec_below} and~\ref{le:convergence_prec_alternative}, we used $\|\cdot\|_{1-\kappa}$-norm on $g_0$, $g_\eps$ and $g_\eps - g_0$. But we have
\begin{equation*}
    \begin{split}
    \|\Pi_0^\perp (g_\eps - g_0)\|_{1-\kappa} &= \|(\lL_\eps^{-1} - \lL_0^{-1}) \lL_\eps g_\eps\|_{1-\kappa} + \|\lL_0^{-1} \big( \lL_\eps g_\eps - \lL_0 g_0 \big) \|_{1-\kappa}\\
    &\lesssim \eps^\theta \|\lL_\eps g_\eps\|_{\beta-2} + \|\lL_\eps g_\eps - \lL_0 g_0\|_{\beta-2}\;.
    \end{split}
\end{equation*}
The claim then follows. 
\end{proof}

\subsection{The commutator}

In this subsection, we provide error bounds for $\Com_\eps (f_\eps; g_\eps; h_\eps) - \Com_0 (f_0; g_0; h_0)$ under different regularity assumptions on the inputs. The other commutator involving the heat kernel will be dealt with in the next subsection. 

\begin{lem} \label{le:convergence_comm_good}
    Let $\alpha \in (0,1)$ and $\beta, \gamma \in (-1,1)$ be such that $\alpha + \beta \in (-1,2)$, $\beta + \gamma < 0$ and $\alpha+\beta+\gamma>0$. Then for every $\kappa>0$, there exists $\theta>0$ such that
    \begin{equation*}
        \begin{split}
        \big\| \Com_\eps (f_\eps; g_\eps; h_\eps) - &\Com_0 (f_0; g_0; h_0) \big\|_{\alpha+\beta+\gamma-\kappa} \lesssim \big( \eps^\theta \|h_0\|_{\gamma} + \|h_\eps - h_0\|_{\gamma} \big) \, \|f_0\|_\alpha \|g_0\|_\beta\\
        &+ \big( \|f_\eps - f_0\|_\alpha \|g_\eps\|_\beta + \|f_0\|_\alpha \|g_\eps - g_0\|_\beta \big) \, \|h_\eps\|_\gamma\;.
        \end{split}
    \end{equation*}
\end{lem}
\begin{proof}
We split the difference into
\begin{equation*}
    \begin{split}
    \Com_\eps (f_\eps; g_\eps; h_\eps) - &\Com_0 (f_0; g_0; h_0) = \Com_\eps (f_\eps - f_0; g_\eps; h_\eps) + \Com_\eps (f_\eps; g_\eps - g_0; h_\eps)\\
    &+ \Com_\eps (f_0; g_0; h_\eps - h_0) + \big( \Com_\eps (f_0; g_0; h_0) - \Com_0 (f_0; g_0; h_0) \big)\;.
    \end{split}
\end{equation*}
The control for $\|\cdot\|_{\alpha+\beta+\gamma}$-norm of the first three terms on the right hand side follow directly from Corollary~\ref{cor:com_prec_circ}, and hence we focus on the last one. 

Writing $(f,g,h) = (f_0, g_0, h_0)$ for simplicity, we need to show
\begin{equation} \label{e:convergence_comm_good_block}
    \big\| \Delta_{k,0} \big( \Com_\eps (f; g; h) - \Com_0 (f; g; h) \big) \big\|_{L^\infty} \lesssim \eps^\theta 2^{-(\alpha+\beta+\gamma-\kappa)k} \|f\|_\alpha \|g\|_\beta \|h\|_\gamma\;.
\end{equation}
We further split the quantity on the left hand side above into two terms as
\begin{equation*}
    (\Delta_{k,0} - \Delta_{k,\eps}) \Com_\eps (f;g;h) + \Big( \Delta_{k,\eps} \big( \Com_\eps (f;g;h) \big) - \Delta_{k,0} \big( \Com_0(f;g;h) \big) \Big)\;.
\end{equation*}
The assumption on $\alpha, \beta$ and $\gamma$ guarantees that $\alpha + \beta + \gamma \in (0,1)$. Hence, by Corollaries~\ref{cor:convergence_single_block} (with $m=0$) and~\ref{cor:com_prec_circ}, we have
\begin{equation} \label{e:convergence_comm_good_block_1}
    \begin{split}
    \big\| (\Delta_{k,0} - \Delta_{k,\eps}) \Com_\eps (f;g;h) \big\|_{L^\infty} &\lesssim \eps^\theta 2^{-(\alpha+\beta+\gamma-\kappa)} \|\Com_\eps (f;g;h)\|_{\alpha+\beta+\gamma}\\
    &\lesssim \eps^\theta 2^{-(\alpha+\beta+\gamma-\kappa)} \|f\|_\alpha \|g\|_\beta \|h\|_\gamma\;.
    \end{split}
\end{equation}
For the other term, we expand it as
\begin{equation} \label{e:convergence_comm_good_expansion}
    \begin{split}
    &\phantom{111}\Delta_{k,\eps} \big( \Com_\eps (f;g;h) \big) - \Delta_{k,0} \big( \Com_0(f;g;h) \big)\\
    &= \sum_{j \geq -1} \sum_{i: |i-j| \leq 1} \Big( \Delta_{k,\eps} \big( R_{j,\eps}(f,g) \Delta_{i,\eps} h \big) - \Delta_{k,0} \big( R_{j,0}(f,g) \Delta_{i,0} h \big) \Big)\;,
    \end{split}
\end{equation}
where $R_{j,\eps}$ is defined in \eqref{e:R_op_defn}. By Lemma~\ref{le:com_block}, for each term in the sum on the right hand side above (noting that $|i-j| \leq 1$), we have the uniform bound
\begin{equation*}
    \big\| \Delta_{k,\eps} \big( R_{j,\eps}(f,g) \Delta_{i,\eps} h \big) - \Delta_{k,0} \big( R_{j,0}(f,g) \Delta_{i,0} h \big) \big\|_{L^\infty} \lesssim 2^{-\alpha (j \vee k) - (\beta+\gamma) j} \|f\|_\alpha \|g\|_\beta \|h\|_\gamma\;.
\end{equation*}
On the other hand, we also have the bound (again noting that $|i-j| \leq 1$)
\begin{equation*}
    \big\| \Delta_{k,\eps} \big( R_{j,\eps}(f,g) \Delta_{i,\eps} h \big) - \Delta_{k,0} \big( R_{j,0}(f,g) \Delta_{i,0} h \big) \big\|_{L^\infty} \lesssim \eps^{\theta'} 2^{M(j+k)} \|f\|_\alpha \|g\|_\beta \|h\|_\gamma
\end{equation*}
for some $\theta', M>0$. Interpolating it with the above uniform bound and substituting it into the sum \eqref{e:convergence_comm_good_expansion}, we deduce that for every $\kappa>0$, there exists $\theta>0$ such that
\begin{equation*}
    \big\| \Delta_{k,\eps} \big( \Com_\eps (f;g;h) \big) - \Delta_{k,0} \big( \Com_0(f;g;h) \big) \big\|_{L^\infty} \lesssim \eps^\theta 2^{-(\alpha+\beta+\gamma-\kappa) k} \|f\|_\alpha \|g\|_\beta \|h\|_\gamma\;.
\end{equation*}
Together with \eqref{e:convergence_comm_good_block_1}, it implies the bound \eqref{e:convergence_comm_good_block}. We have thus completed the proof of the lemma. 
\end{proof}

\begin{lem} \label{le:convergence_comm_outside1}
Let $\alpha, \beta \in (0,1)$ and $\gamma \in (-2,-1)$ be such that $\alpha+\beta+\gamma>0$. Then for every $\kappa>0$, there exists $\theta>0$ such that
\begin{equation} \label{e:convergence_comm_outside1}
    \|\Com_\eps(f; g; \lL_{\eps} H)-\Com_{0}(f;g;\lL_{0}H)\|_{\alpha+\beta+\gamma-\kappa} \lesssim \eps^{\theta} \|f\|_{\alpha} \|g\|_{\beta} \|H\|_{\gamma+2}\;,
\end{equation}
uniformly over $\eps\in \N^{-1}$, $f \in \cC^\alpha$, $g \in \cC^\beta$ and $H \in \cC^{\gamma+2}$. 
\end{lem}
\begin{proof} 
By definition, we need to show
\begin{equation} \label{e:convergence_comm_outside1_bound}
    \begin{split}
    &\phantom{111}\left\| \Delta_{k,0} \big( \Com_\eps(f; g; \lL_{\eps} H)-\Com_0(f;g;\lL_{0}H) \big) \right\|_{L^\infty}\\
    &\lesssim \eps^\theta 2^{-(\alpha+\beta+\gamma-\kappa) k} \|f\|_{\alpha} \|g\|_{\beta} \|H\|_{\gamma+2}\;.
    \end{split}
\end{equation}
The proof is similar to that of the second half of Lemma~\ref{le:convergence_comm_good}. We again split the left hand side above into two terms as
\begin{equation*}
    (\Delta_{k,0} - \Delta_{k,\eps}) \Com_\eps(f, g, \lL_{\eps} H) + \Big( \Delta_{k,\eps} \big( \Com_\eps (f,g,\lL_\eps H) \big) - \Delta_{k,0} \big( \Com_0 (f,g,\lL_0 H) \big) \Big)\;.
\end{equation*}
The bound for the first term follows from Corollaries~\ref{cor:convergence_single_block} (with $m=0$) and~\ref{cor:com_prec_circ} that
\begin{equation*}
    \begin{split}
    \big\| (\Delta_{k,0} - \Delta_{k,\eps}) \Com_\eps(f; g; \lL_{\eps} H) \big\|_{L^\infty} &\lesssim \eps^\theta 2^{-(\alpha+\beta+\gamma-\kappa)k} \big\|  \Com_\eps(f; g; \lL_{\eps} H) \big\|_{\alpha+\beta+\gamma}\\
    &\lesssim \eps^\theta 2^{-(\alpha+\beta+\gamma-\kappa)k} \|f\|_\alpha \|g\|_\beta \|H\|_{\gamma+2}\;.
    \end{split}
\end{equation*}
For the second term, we have a slightly different expansion than \eqref{e:convergence_comm_good_expansion} that
\begin{equation} \label{e:convergence_comm_outside1_expansion}
    \begin{split}
    &\phantom{111}\Delta_{k,\eps} \big( \Com_\eps (f;g;\lL_\eps H) \big) - \Delta_{k,0} \big( \Com_0 (f;g;\lL_0 H) \big)\\
    &= \sum_{j \geq -1} \sum_{i:|i-j| \leq 1} \Big( \Delta_{k,\eps} \big( R_{j,\eps}(f,g) \cdot \Delta_{i,\eps} \lL_\eps H \big) - \Delta_{k,0} \big( R_{j,0}(f,g) \cdot \Delta_{i,0} \lL_0 H \big) \Big)\;,
    \end{split}
\end{equation}
where the input $h$ in $\Com_\eps$ and $\Com_0$ are replaced by $\lL_\eps H$ and $\lL_0 H$ respectively, whose regularities are measured by the norms $\|\cdot\|_{\gamma;\eps}$ and $\|\cdot\|_{\gamma}$ for $\gamma \in (-2,1)$. As in the previous lemma, for each term in the sum above, we again interpolate between a sharp uniform bound and a loose bound with a positive power in $\eps$. Since the assumption on $\alpha, \beta$ and $\gamma$ ensures that $\alpha + \beta + \gamma \in (0,1)$, by Lemma~\ref{le:com_block}, we have the uniform bound
\begin{equation} \label{e:convergence_comm_outside1_uniform}
    \begin{split}
    &\phantom{111}\left\| \Delta_{k,\eps} \big( R_{j,\eps}(f,g) \cdot \Delta_{i,\eps} \lL_\eps H \big) - \Delta_{k,0} \big( R_{j,0}(f,g) \cdot \Delta_{i,0} \lL_0 H \big) \right\|_{L^\infty}\\
    &\lesssim 2^{-\alpha (j \vee k) - (\beta+\gamma)j} \|f\|_\alpha \|g\|_\beta \big( \|\lL_\eps H\|_{\gamma;\eps} + \|\lL_0 H\|_{\gamma} \big)\\
    &\lesssim 2^{-\alpha (j \vee k) - (\beta+\gamma)j} \|f\|_\alpha \|g\|_\beta \|H\|_{\gamma+2}\;,
    \end{split}
\end{equation}
where we have used the equivalence of norms for $\gamma+2 \in (0,1)$. On the other hand, using additionally (than before) Corollary~\ref{cor:convergence_single_block} (with $m=0, 1$), we also have the bound
\begin{equation} \label{e:convergence_comm_outside1_small_loose}
    \begin{split}
    &\phantom{111}\left\| \Delta_{k,\eps} \big( R_{j,\eps}(f,g) \cdot \Delta_{i,\eps} \lL_\eps H \big) - \Delta_{k,0} \big( R_{j,0}(f,g) \cdot \Delta_{i,0} \lL_0 H \big) \right\|_{L^\infty}\\
    &\lesssim \eps^{\theta'} 2^{M(j+k)} \|f\|_\alpha \|g\|_\beta \|H\|_{\gamma+2}\;.
    \end{split}
\end{equation}
for some $\theta', M>0$. Interpolating the bounds \eqref{e:convergence_comm_outside1_uniform} and \eqref{e:convergence_comm_outside1_small_loose}, and then substituting it back to the sum in \eqref{e:convergence_comm_outside1_expansion} gives the desired bound for the $L^\infty$-norm of the quantity in \eqref{e:convergence_comm_outside1_expansion}. 
\end{proof}

\begin{lem} \label{le:convergence_comm_below}
Let $\alpha, \beta \in (0,1)$ and $\gamma \in (-2,-1)$ be such that $\alpha+\beta+\gamma>0$. Then for every $\kappa>0$, there exists $\theta>0$ such that
\begin{equation*}
    \begin{split}
    \big\| \Com_\eps (f_\eps ;g_\eps; h_\eps) - \Com_0 &(f_0; g_0; h_0) \big\|_{\alpha+\beta+\gamma-\kappa} \lesssim \big( \eps^\theta  \|h_0\|_{\gamma} + |\!|\!|h_{\eps};h_0|\!|\!|_{\gamma; \eps} \big) \, \|f_0\|_\alpha \|g_0\|_\beta\\
    &+ \big( \|f_\eps - f_0\|_\alpha \, \|g_\eps\|_\beta + \|f_0\|_\alpha \, \|g_\eps - g_0\|_\beta \big) \, \|h_\eps\|_{\gamma;\eps}
    \end{split}
\end{equation*}
uniformly over $\eps\in \N^{-1}$ and functions in respective spaces. 
\end{lem}
\begin{proof}
We first treat the situation when $f_\eps \equiv f_0 = f$ and $g_\eps \equiv g_0 = g$. Define $H_\eps := \lL_\eps^{-1} \Pi_0^\perp h_\eps$ for $\eps\in \N^{-1}$, so that
\begin{equation*}
    h_\eps = \Pi_0 h_\eps + \lL_\eps H_\eps\;.
\end{equation*}
We then decompose the object into several differences as
\begin{equation*}
    \begin{split}
    \Com_\eps (f,g,h_\eps) - &\Com_0 (f,g,h_0) = \Pi_0 h_0 \cdot \big( \Com_\eps (f,g,1) - \Com_0 (f,g,1) \big)\\
    &+ \big( \Pi_0 (h_\eps - h_0) \big) \cdot \Com_\eps (f,g,1) + \Com_\eps \big(f,g, \lL_\eps (H_\eps - H_0) \big)\\
    &+ \big( \Com_\eps (f,g, \lL_\eps H_0) - \Com_0 (f,g, \lL_0 H_0)  \big)\;.
    \end{split}
\end{equation*}
We need to control the $\|\cdot\|_{\alpha+\beta+\gamma-\kappa}$-norm for each of the four terms on the right hand side. By Corollary~\ref{cor:com_prec_circ}, the $\|\cdot\|_{\alpha+\beta+\gamma}$-norm of the second and third terms are bounded by
\begin{equation*}
    \|f\|_\alpha \|g\|_\beta  \; |\!|\!|h_{\eps};h_0|\!|\!|_{\gamma;\eps}\;.
\end{equation*}
By Lemma~\ref{le:convergence_comm_outside1}, the $\|\cdot\|_{\alpha+\beta+\gamma-\kappa}$-norm of the fourth term is bounded by
\begin{equation*}
    \eps^\theta \|f\|_\alpha \|g\|_\beta \|H_0\|_{\gamma+2}\;.
\end{equation*}
For the first one, let $\gamma' \in (-1,0)$ be such that $\beta + \gamma' < 0$. Then the triple $(\alpha, \beta, \gamma')$ satisfies the assumption of Lemma~\ref{le:convergence_comm_good}. Hence, we have
\begin{equation*}
    \begin{split}
    &\phantom{111}|\Pi_0 h_0| \cdot \big\| \big( \Com_\eps (f,g,1) - \Com_0 (f,g,1) \big) \big\|_{\alpha+\beta+\gamma-\kappa}\\
    &\lesssim |\Pi_0 h_0| \cdot \big\| \big( \Com_\eps (f,g,1) - \Com_0 (f,g,1) \big) \big\|_{\alpha+\beta+\gamma'-\kappa} \lesssim \eps^\theta |\Pi_0 h_0| \cdot \|f\|_\alpha \|g\|_\beta\;.
    \end{split}
\end{equation*}
Combining the bounds for the above four terms then gives the bound
\begin{equation} \label{e:convergence_comm_below_special}
    \left\| \Com_\eps (f,g,h_\eps) - \Com_0 (f,g,h_0) \right\|_{\alpha+\beta+\gamma-\kappa}\lesssim \big( \eps^\theta  \|h_0\|_{\gamma} + |\!|\!|h_{\eps};h_0|\!|\!|_{\gamma; \eps} \big) \|f\|_\alpha \|g\|_\beta\;.
\end{equation}
For the general case, we write it as
\begin{equation} \label{e:convergence_comm_below_split}
    \begin{split}
    \Com_\eps (f_\eps; &g_\eps; h_\eps) - \Com_0 (f_0; g_0; h_0) = \Com_\eps (f_\eps - f_0; g_\eps; h_\eps)\\
    &+ \Com_\eps (f_0; g_\eps - g_0; h_\eps) + \big( \Com_\eps (f_0; g_0; h_\eps) - \Com_0 (f_0; g_0; h_0) \big)\;.
    \end{split}
\end{equation}
Applying Corollary~\ref{cor:com_prec_circ} to the first two terms and \eqref{e:convergence_comm_below_special} to the last term on the right hand side above gives the desired bound. 
\end{proof}

\begin{lem} \label{lem:para_linearisation_diff}
Let $\alpha \in (\frac{1}{2},1)$. For every $\kappa>0$, there exists $\theta>0$ such that
    \begin{equation*}
        \|\lL_{\eps}\PL_{\eps}^{F}(f_{\eps}) - \lL_0 \PL_{0}^{F}(f_0)\|_{{2\alpha-2-\kappa}} \lesssim (1+\|f\|_{\alpha}+\|f_{\eps}\|)^4\|F\|_{\cC^3}(\eps^{\theta}  +\|f_{\eps}-f_0\|_{\alpha})
    \end{equation*}
    where $\PL_{\eps}^{F} := \PL_{A_{\eps}}^{F}$ is given in \eqref{e:para_linear_defn}.
\end{lem}
\begin{proof}
By definition,
\begin{equation*}
    \begin{split}
    \lL_{\eps}\PL^{F}_{\eps}(f_{\eps}) - \lL_0 \PL^{F}_{0}(f_0)=\lL_{\eps}(\PL^{F}_{\eps}(f_{\eps})-\PL^{F}_{\eps}(f_{0}))+ \big(\lL_{\eps}\PL^{F}_{\eps}(f_{0})-\lL_0 \PL^{F}_{0}(f_0)\big)\;.
    \end{split}
\end{equation*}
By Lemma \ref{lem:para_linearisation}, we have
\begin{equation*}
    \|\lL_{\eps}(\PL^{F}_{\eps}(f_{\eps})-\PL^{F}_{\eps}(f_{0}))\|_{2\alpha-2}\lesssim \|F\|_{\cC^3}(1+\|f_{\eps}\|_{\alpha}+\|f_0\|_{\alpha})^4\|f_{\eps}-f_0\|_{\alpha}\;.
\end{equation*}
By definition,
\begin{equation*}
    \begin{split}
    &\Delta_{k,0}\big(\lL_{\eps}\PL_{\eps}^{F}(f_{0})-\lL_0 \PL^{F}_{0}(f_0)\big)= (\Delta_{k,0}-\Delta_{k,\eps})\lL_{\eps}\PL^{F}_{\eps}(f_0)\\
    +&(\lL_{\eps}\Delta_{k,\eps}-\lL_0\Delta_{k,0})\PL^{F}_{\eps}(f_0)+\lL_0\Delta_{k,0}(\PL^{F}_{\eps}(f_0)-\PL^{F}_{0}(f_0))\;.
    \end{split}
\end{equation*}
By Lemma~\ref{lem:para_linearisation} and Corollary~\ref{cor:convergence_single_block},
\begin{equation*}
    \|\Delta_{k,0}\big(\lL_{\eps}\PL^{F}_{\eps}(f_{0})-\lL_0 \PL^{F}_{0}(f_0)\big)\|_{L^{\infty}}\lesssim  \eps^{\theta} (1+\|f_0\|_{\alpha}+\|f_{\eps}\|)^4\|F\|_{\cC^3}\;.
\end{equation*}
This completes the proof.
\end{proof}

\subsection{Heat kernel}

\begin{lem} \label{le:convergence_heat_normal}
Let $\beta \in (-1,1)$ and $\eta \in [0,2)$. Let $\alpha \in (\beta,1)$. For every $\kappa>0$, there exists $\theta>0$ such that
\begin{equation*}
    \|e^{t \lL_\eps} f_\eps - e^{t \lL_0} f_0\|_{\cC_{\sigma,T}^{\eta/2} \cC_{x}^{\alpha-\kappa}} \lesssim \eps^\theta \|f_0\|_{\beta} + \|f_\eps - f_0\|_{\beta}\;,
\end{equation*}
where $\sigma = (\eta+\alpha-\beta) \vee 0$. 

On the other hand, if $\gamma \in (-1,1)$, then for every $\kappa>0$, there exists $\theta>0$ such that
\begin{equation*}
    \|\lL_\eps e^{t \lL_\eps} f_\eps - \lL_0 e^{t \lL_0} f_0\|_{\cC_{\sigma,T}^{\eta/2} \cC_{x}^{\gamma-\kappa}} \lesssim \eps^\theta \|f_0\|_{\beta} + \|f_\eps - f_0\|_{\beta}\;,
\end{equation*}
where now $\sigma = (\eta+\gamma+2-\beta) \vee 0$. 
\end{lem}
\begin{proof}
For the first claim, we write
\begin{equation*}
    e^{t \lL_\eps} f_\eps - e^{t \lL_0} f_0 = e^{t \lL_\eps} (f_\eps - f_0) + (e^{t \lL_\eps} - e^{t \lL_0}) f_0\;.
\end{equation*}
The bound for the first term follows from Lemma~\ref{le:heat_smoothing_space}, while the bound for the second term follows the $L_{\sigma,T}^{\infty} L_x^\infty$ bound from standard homogenisation and then interpolation with Lemma~\ref{le:heat_smoothing_space}. 

For the second claim, we write
\begin{equation*}
    \lL_\eps e^{t \lL_\eps} f_\eps - \lL_0 e^{t \lL_0} f_0 = \lL_\eps e^{t \lL_\eps} (f_\eps - f_0) + (\lL_\eps e^{t \lL_\eps} - \lL_0 e^{t \lL_0}) f_0\;.
\end{equation*}
For the first term on the right hand side above, by Lemma~\ref{le:heat_smoothing_space}, we have
\begin{equation*}
    \|\lL_\eps e^{t \lL_\eps} (f_\eps - f_0)\|_{\cC_{\sigma,T}^{\eta/2} \cC_{x}^{\gamma}} \lesssim \|e^{t \lL_\eps} (f_\eps - f_0)\|_{\cC_{\sigma,T}^{\eta/2} \cC_{x}^{\gamma+2;\eps}} \lesssim \|f_\eps - f_0\|_{\beta}\;.
\end{equation*}
As for the second term, again by Lemma~\ref{le:heat_smoothing_space}, one has the uniform bound
\begin{equation*}
    \begin{split}
    \| \Delta_{k,0} \big( \lL_\eps e^{t \lL_\eps} - \lL_0 e^{t \lL_0} \big) f_0 \|_{\cC_{\sigma,T}^{\eta/2} L_x^\infty} &\lesssim 2^{-\gamma k} \big( \|\lL_\eps e^{t \lL_\eps} f_0\|_{\cC_{\sigma,T}^{\eta/2} \cC_x^\gamma} + \|\lL_0 e^{t \lL_0} f_0\|_{\cC_{\sigma,T}^{\eta/2} \cC_x^\gamma} \big)\\
    &\lesssim 2^{-\gamma k} \|f_0\|_{\beta}\;.
    \end{split}
\end{equation*}
To get smallness in $\eps$, we decompose it as (omitting the input $f_0$)
\begin{equation*}
    \begin{split}
    &\phantom{111}\Delta_{k,0} \big( \lL_\eps e^{t \lL_\eps} - \lL_0 e^{t \lL_0} \big)\\
    &= (\Delta_{k,0} - \Delta_{k,\eps}) \lL_\eps e^{t \lL_\eps} + (\Delta_{k,\eps} \lL_\eps - \Delta_{k,0} \lL_0) e^{t \lL_\eps} + \Delta_{k,0} \lL_0 (e^{t \lL_\eps} - e^{t \lL_0})\;.
    \end{split}
\end{equation*}
The $\cC_x^\beta \rightarrow \cC_{\sigma,T}^{\eta/2} L_x^\infty$ operator norms of the above three terms are all bounded by $\eps^{\theta'} 2^{Mk}$ for some $\theta', M>0$. Interpolating it with the above uniform bound gives
\begin{equation*}
    \left\| \big( \lL_\eps e^{t \lL_\eps} - \lL_0 e^{t \lL_0} \big) f_0 \right\|_{\cC_{\sigma,T}^{\eta/2} \cC_x^{\gamma-\kappa}} \lesssim \eps^\theta \|f_0\|_{\beta}
\end{equation*}
for some $\theta>0$. This completes the proof. 
\end{proof}

\begin{lem} \label{le:convergence_heat_spacetime_normal}
Fix $\beta \in (-1,1)$. Let $\alpha \in (\beta, 1)$ and $\eta \in [0,2)$ be such that $\alpha + \eta \leq 2 + \beta$. Then for every $\kappa>0$ sufficiently small, there exists $\theta>0$ such that
\begin{equation*}
    \left\| \iI_\eps (\kK_\eps) - \iI_0 (\kK_0) \right\|_{\cC_{\widetilde{\sigma}, T}^{\eta/2} \cC_x^{\alpha-\kappa}} \lesssim \eps^{\theta} \|\kK_0\|_{L_{\sigma, T}^{\infty} \cC_x^\beta} + \|\kK_\eps - \kK_0\|_{L_{\sigma,T}^{\infty} \cC_x^\beta}\;,
\end{equation*}
where $\widetilde{\sigma} = (\sigma+\eta+\alpha-\beta-2) \vee 0$. 

On the other hand, if $\gamma \in (-1,1)$ and $\eta \in [0,2)$ are such that $\eta + \gamma - \beta \in [-2,0]$, then for every $\kappa>0$, there exists $\theta>0$ such that
\begin{equation*}
    \|\lL_\eps \iI_\eps \kK_\eps - \lL_0 \iI_0 \kK_0\|_{\cC_{\widetilde{\sigma},T}^{\eta/2} \cC_x^{\gamma-\kappa}} \lesssim \eps^\theta \|\kK_0\|_{L_{\sigma,T}^{\infty} \cC_x^\beta} + \|\kK_\eps - \kK_0\|_{L_{\sigma,T}^{\infty} \cC_x^\beta}\;,
\end{equation*}
where this time $\widetilde{\sigma} = (\sigma + \eta + \gamma - \beta) \vee 0$. 
\end{lem}
\begin{proof}
For the first claim, we write
\begin{equation*}
    \iI_\eps (\kK_\eps) - \iI_0 (\kK_0) = \iI_\eps (\kK_\eps - \kK_0) + (\iI_\eps - \iI_0) \kK_0\;.
\end{equation*}
The bound for the first term follows from Lemma~\ref{le:heat_smoothing_spacetime}, while the bound for the second term follows from the $L_{\widetilde{\sigma},T}^{\infty} L_x^\infty$ bound from standard homogenisation and then interpolation with Lemma~\ref{le:heat_smoothing_spacetime}. 

For the second claim, we write
\begin{equation*}
    \lL_\eps \iI_\eps \kK_\eps - \lL_0 \iI_0 \kK_0 = \lL_\eps \iI_\eps (\kK_\eps - \kK_0) + (\lL_\eps \iI_\eps - \lL_0 \iI_0) \kK_0\;.
\end{equation*}
By Lemma~\ref{le:heat_smoothing_spacetime}, we have
\begin{equation*}
    \|\lL_\eps \iI_\eps (\kK_\eps - \kK_0)\|_{\cC_{\widetilde{\sigma},T}^{\eta/2} \cC_{x}^{\gamma-\kappa}} \lesssim \|\iI_\eps (\kK_\eps - \kK_0)\|_{\cC_{\widetilde{\sigma},T}^{\eta/2} \cC_{x}^{\gamma+2-\kappa; \eps}} \lesssim \|\kK_\eps - \kK_0\|_{L_{\sigma,T}^{\infty} \cC_{x}^{\beta}}
\end{equation*}
Finally, we estimate the $k$-th block of $(\lL_\eps \iI_\eps - \lL_0 \iI_0) \kK_0$ (omitting the input function) by
\begin{equation*}
    \Delta_{k,0} (\lL_\eps \iI_\eps - \lL_0 \iI_0) = (\Delta_{k,0} - \Delta_{k,\eps}) \lL_\eps \iI_\eps + (\Delta_{k,\eps} \lL_\eps - \Delta_{k,0} \lL_0) \iI_\eps + \Delta_{k,0} \lL_0 (\iI_\eps - \iI_0)\;.
\end{equation*}
Their $L_{\sigma, T}^{\infty} \cC_{x}^{\beta} \rightarrow \cC_{\widetilde{\sigma},T}^{\eta/2} L_x^\infty$ operator norms are bounded by $\eps^{\theta'} 2^{Mk}$ for some $\theta', M > 0$. On the other hand, we have the uniform bound
\begin{equation*}
    \begin{split}
    \|\Delta_{k,0} \big( \iI_\eps \lL_\eps \kK_0 - \iI_0 \lL_0 \kK_0 \big)\|_{\cC_{\widetilde{\sigma},T}^{\eta/2} L_x^\infty} &\lesssim  2^{-\gamma k} \big( \|(\iI_\eps \lL_\eps \kK_0)\|_{\cC_{\widetilde{\sigma},T}^{\eta/2} \cC_{x}^{\gamma}} + \|\iI_0 \lL_0 \kK_0\|_{\cC_{\widetilde{\sigma},T}^{\eta/2} \cC_{x}^{\gamma}} \big)\\
    &\lesssim 2^{-\gamma k} \|\kK_0\|_{L_{\sigma, T}^\infty \cC_x^{\beta}}\;.
    \end{split}
\end{equation*}
Interpolating the two bounds then gives
\begin{equation*}
    \| (\iI_\eps \lL_\eps - \iI_0 \lL_0) \kK_0 \|_{\cC_{\widetilde{\sigma},T}^{\eta/2} \cC_{x}^{\gamma-\kappa}} \lesssim \eps^{\theta} \|\kK_0\|_{L_{\sigma,T}^{\infty} \cC_x^{\beta+2}}\;.
\end{equation*}
This completes the proof. 
\end{proof}

\begin{lem} \label{le:convergence_heat_below}
Let $\beta \in (-2,-1)$. For every $\kappa > 0$, there exists $\theta>0$ such that
\begin{equation*}
    \| \iI_\eps (g_\eps) - \iI_0 (g_0) \|_{L_T^\infty \cC_{x}^{\beta+2-\kappa}} \lesssim T \, \|\Pi_0 (g_\eps - g_0)\|_{L_T^\infty} + \|G_\eps - G_0\|_{L_T^\infty \cC_x^{\beta+2}} + \eps^\theta \|G_0\|_{L_{T}^{\infty} \cC_{x}^{\beta+2}}\;,
\end{equation*}
where $G_\eps := \lL_\eps^{-1} \Pi_0^\perp g_\eps$ (including $\eps=0$). The proportionality constant is independent of $T>0$ and $\eps\in \N^{-1}$. 
\end{lem}
\begin{proof}
We decompose $\iI_\eps (g_\eps) - \iI_0 (g_0)$ as
\begin{equation*}
    \iI_\eps (g_\eps) - \iI_0 (g_0) = \int_{0}^{t} \Pi_0 \big( g_\eps(r) - g_0(r) \big) {\rm d}r + \iI_\eps \lL_\eps G_\eps - \iI_0 \lL_0 G_0\;,
\end{equation*}
where we have omitted the time $t$ in most of the expressions. For the ``constant" term, one directly has
\begin{equation*}
    \left| \int_{0}^{t} \Pi_0 \big( g_\eps (r) - g_0 (r) \big) {\rm d}r \right| \leq t \cdot \sup_{r \in [0,t]} \left| \Pi_0 \big( g_\eps(r) - g_0(r) \big) \right|\;.
\end{equation*}
For the difference of the other two terms, note that $\beta + 2 \in (-1,1)$, then the desired bound follows from Lemma~\ref{le:convergence_heat_spacetime_normal} with $\beta \mapsto \beta + 2$, $\gamma = \beta+2$, $\eta=0$ and $\sigma = \tilde{\sigma}=0$. This completes the proof of the lemma. 
\end{proof}

\begin{lem} \label{le:convergence_heat_prec_below}
Let $\beta \in (-2,-1)$ and $\sigma \in [0,2]$. Then for every $\kappa>0$, there exists $\theta>0$ such that
\begin{equation*}
    \begin{split}
    \big\| f_\eps (t) \prec_\eps (\iI_\eps g_\eps)(t) &- f_0 (t) \prec (\iI_0 g_0)(t) \big\|_{\beta+2-\kappa} \lesssim t^{-\frac{\sigma}{2}} \Big( \eps^\theta \|f_0\|_{L_{\sigma,T}^{\infty} L_x^\infty} \|g_0\|_{L_t^\infty \cC_x^\beta}\\
    &+\|f_\eps - f_0\|_{L_{\sigma,T}^{\infty} L_x^\infty} \|g_\eps\|_{\beta;\eps} + \|f_0\|_{L_{\sigma,T}^{\infty} L_x^\infty} |\!|\!| g_\eps; g_0 |\!|\!|_{L_{T}^{\infty} \cC_{x}^{\beta;\eps}} \Big)\;,
    \end{split}
\end{equation*}
where the proportionality constant is uniform over $\eps\in \N^{-1}$, $t>0$. 
\end{lem}
\begin{proof}
Note that $\Pi_0 g_\eps$ and $\Pi_0 g_0$ play no rule in the quantities concerned, so we can assume without loss of generality that $\Pi_0 g_\eps = \Pi_0 g_0 = 0$, and hence $g_\eps = \Pi_0^\perp g_\eps$ and $g_0 = \Pi_0^\perp g_0$. Also, since the bound is pointwise in time and that the quantities involving $g$ all have $L^\infty$-in-time bounds, we will omit the time $t$ below and consider spatial norms only, and add weight in time in the end. 

We start with the decomposition
\begin{equation*}
    \begin{split}
    f_\eps \prec_\eps \iI_\eps (g_\eps) - f_0 \prec \iI_0 (g_0) = &(f_\eps - f_0) \prec_\eps \iI_\eps (g_\eps) + f_0 \prec_\eps \big( \iI_\eps(g_\eps) - \iI_0 (g_0) \big)\\
    &+ \Big( f_0 \prec_\eps \iI_0 (g_0) - f_0 \prec \iI_0 (g_0) \Big)\;.
    \end{split}
\end{equation*}
By Proposition~\ref{pr:paraproducts_gBesov} and Lemmas~\ref{le:heat_smoothing_spacetime} and~\ref{le:convergence_heat_below}, the first two terms on the right hand side above are bounded respectively by
\begin{equation*}
    \begin{split}
    \left\| (f_\eps - f_0) \prec_\eps \iI_\eps (g_\eps) \right\|_{\beta+2} &\lesssim \|f_\eps - f_0\|_{L^\infty} \cdot \|g_\eps\|_{\beta;\eps}\\
    \left\| f_0 \prec_\eps \big( \iI_\eps (g_\eps) - \iI_0 (g_0) \big) \right\|_{\beta+2-\kappa} &\lesssim  \big( \eps^\theta \|g_0\|_{\beta}+ |\!|\!| g_\eps; g_0 |\!|\!|_{\beta;\eps} \big) \cdot \|f_0\|_{L^\infty}\;.
    \end{split}
\end{equation*}
For the third one, by Lemmas~\ref{le:convergence_prec_good} and~\ref{le:heat_smoothing_spacetime}, we have
\begin{equation*}
    \left\| f_0 \prec_\eps \iI_0 (g_0) - f_0 \prec \iI_0 (g_0) \right\|_{\beta+2-\kappa} \lesssim \eps^\theta \|f_0\|_{L^\infty} \|g_0\|_{\beta}\;.
\end{equation*}
These together imply the desired bound without time. Adding the weight in time from $f_\eps$ and $f_0$ then gives the conclusion. 
\end{proof}

\begin{lem} \label{le:convergence_prec_heat_below}
Let $\alpha \in (0,1)$ and $\beta \in (-2,-1)$ be such that $\alpha+\beta>-1$. Let $\eta, \gamma \in (0,1)$ be such that $\eta + \gamma \in (0,\beta+2)$. Let $\sigma \in [0,2)$. Then for every $\kappa>0$, there exists $\theta>0$ such that
\begin{equation} \label{e:convergence_prec_heat_below}
    \begin{split}
    \big\| \iI_\eps (f_\eps \prec_\eps g_\eps) - &\iI_0 (f_0 \prec g_0) \big\|_{\cC_{\widetilde{\sigma},T}^{\eta/2} \cC_{x}^{\gamma-\kappa}} \lesssim \|f_\eps-f_0\|_{L_{\sigma,T}^\infty \cC_x^\alpha} \|\Pi_0^\perp g_\eps\|_{L_T^\infty \cC_x^{\beta;\eps}}\\
    &+ \big( \eps^\theta \|\Pi_0^\perp g_0\|_{L_{T}^{\infty} \cC_x^{\beta}}  + |\!|\!| \Pi_0^\perp g_\eps; \Pi_0^\perp g_0 |\!|\!|_{L_{T}^{\infty} \cC_x^{\beta;\eps}} \big) \|f_0\|_{L_{\sigma,T}^{\infty} \cC_{x}^{\alpha}}\;,
    \end{split}
\end{equation}
where $\widetilde{\sigma} = (\sigma + \eta + \gamma - \beta-2) \vee 0$. The proportionality constant is independent of $\eps\in \N^{-1}$ and functions in their relative classes. 
\end{lem}
\begin{proof}
Similar as before, we can assume without loss of generality that $\Pi_0 g_\eps = \Pi_0 g_0 = 0$, and we write $G_\eps = \lL_\eps^{-1} g_\eps$ (including $\eps=0$). We first treat the situation when $f_\eps = f_0 = f$. In that case, we have
\begin{equation*}
    \iI_\eps (f \prec_\eps \lL_\eps G_\eps) = \iI_\eps \lL_\eps (f \prec_\eps G_\eps) - \iI_\eps \big( (\lL_\eps f) \prec_\eps G_\eps  \big) - 2 \iI_\eps \big( \P_\eps (f, G_\eps) \big)\;.
\end{equation*}
The bound for the difference from first term above follows from Lemma~\ref{le:convergence_heat_spacetime_normal} with $\kK_\eps = f \prec_\eps G_\eps$ and $\kK_0 = f \prec G_0$ (replacing $\beta$ there by $\beta+2$). 

For the second term, by Lemma~\ref{le:convergence_heat_spacetime_normal}, Proposition~\ref{pr:paraproducts_gBesov} and Lemma~\ref{le:convergence_prec_reg-2}, we have
\begin{equation*}
    \begin{split}
    &\phantom{111}\left\| \iI_\eps \big( (\lL_\eps f) \prec_\eps G_\eps  \big) - \iI_0 \big( (\lL_0 f) \prec G_0  \big) \right\|_{\cC_{\widetilde{\sigma},T}^{\eta/2} \cC_x^{\gamma-\kappa}}\\
    &\lesssim \eps^\theta \| (\lL_0 f) \prec G_0 \|_{L_{\sigma,T}^{\infty} \cC_{x}^{\alpha+\beta}} + \|(\lL_\eps f) \prec_\eps G_\eps - (\lL_0 f) \prec G_0 \|_{L_{\sigma,T}^{\infty} \cC_{x}^{\alpha+\beta-\frac{\kappa}{2}}}\\
    &\lesssim \big( \eps^\theta \|G_0\|_{L_{T}^{\infty} \cC_x^{\beta+2}}  + \|G_\eps - G_0\|_{L_{T}^{\infty} \cC_x^{\beta+2}} \big) \|f\|_{L_{\sigma,T}^{\infty} \cC_{x}^{\alpha}}\;.
    \end{split}
\end{equation*}
For the third term, by Lemmas~\ref{le:convergence_heat_spacetime_normal} and~\ref{le:convergence_prec_alternative}, we have
\begin{equation*}
    \begin{split}
    &\phantom{111}\left\| \iI_\eps \big( \P_\eps (f, G_\eps) \big) - \iI_0 \big( \P_0 (f, G_0) \big) \right\|_{\cC_{\widetilde{\sigma},T}^{\eta/2} \cC_x^{\gamma-\kappa}}\\
    &\lesssim \eps^\theta \|\P_0 (f, G_0)\|_{L_{\sigma,T}^{\infty} \cC_{x}^{\alpha+\beta-\frac{\kappa}{2}}} + \|\P_\eps (f, G_\eps) - \P_0 (f, G_0)\|_{L_{\sigma,T}^{\infty} \cC_{x}^{\alpha+\beta-\frac{\kappa}{2}}}\\
    &\lesssim \big( \eps^\theta \|G_0\|_{L_{T}^{\infty} \cC_x^{\beta+2}} + \|G_\eps - G_0\|_{L_{T}^{\infty} \cC_{x}^{\beta+2}} \big) \|f\|_{L_{\sigma,T}^{\infty} \cC_{x}^{\alpha}}\;.
    \end{split}
\end{equation*}
Collecting the above bounds gives the bound \eqref{e:convergence_prec_heat_below} when $f_\eps = f_0 = f$. To treat the general situation, we just write
\begin{equation*}
    \iI_\eps (f_\eps \prec_\eps g_\eps) - \iI_0 (f_0 \prec g_0) = \iI_\eps \big( (f_\eps - f_0) \prec_\eps g_\eps \big) + \Big( \iI_\eps (f_0 \prec_\eps g_\eps) - \iI_0 (f_0 \prec g_0) \Big)\;,
\end{equation*}
and apply Lemma~\ref{le:heat_smoothing_spacetime} and Proposition~\ref{pr:paraproducts_gBesov} for the first term. The second term has been treated above. This completes the proof of the lemma. 
\end{proof}

\begin{lem} \label{le:convergence_heat_comm_below}
Let $\alpha \in (0,1)$, $\beta \in (-2,-1)$ be such that $\alpha+\beta>-1$. Then for every $\kappa>0$, there exists $\theta>0$ such that
\begin{equation*}\begin{split}
    &\left\| \lL_\eps \, \big( [\iI_\eps, f \prec_\eps] g_\eps \big) - \lL_0 \, \big( [\iI_0, f \prec] g_0 \big)  \right\|_{L^{\infty}_{\sigma, t}\cC^{\alpha+\beta-\kappa}_{x}}\\
    \lesssim &  \|f\|_{\mathfrak{L}_{\sigma, t}^{\alpha}} (\eps^{\theta}\|\Pi_{0}^{\perp}g_{0}\|_{L^{\infty}_{t}\cC^{\beta;0}_{x}}+  \|\lL_\eps^{-1}\Pi_{0}^{\perp} g_\eps-\lL_0^{-1}\Pi_{0}^{\perp} g_0\|_{L^{\infty}_{t}\cC^{\beta+2}_{x}})\;.
\end{split}\end{equation*}
\end{lem}
\begin{proof}
We need to control
\begin{equation*}\begin{split}
    &\left\| \Delta_{k,0} \Big( \lL_\eps \, \big( [\iI_\eps, f \prec_\eps] g_\eps \big)(r) - \lL_0 \, \big( [\iI_0, f \prec] g_0 \big)(r) \Big) \right\|_{L^\infty}\\
    \lesssim & 2^{-(\alpha+\beta-\kappa) k} r^{-\frac{\sigma}{2}} \|f\|_{\mathfrak{L}_{\sigma, t}^{\alpha}} (\eps^{\theta}\|g_{0}\|_{L^{\infty}_{t}\cC^{\beta;0}_{x}}+ \|\lL_\eps^{-1}\Pi_{0}^{\perp} g_\eps-\lL_0^{-1}\Pi_{0}^{\perp} g_0\|_{L^{\infty}_{t}\cC^{\beta+2}_{x}})\;.
\end{split}\end{equation*}
We decompose the object of study into
\begin{equation*}
    \begin{split}
     &\Delta_{k,0} \Big( \lL_\eps \, \big( [\iI_\eps, f \prec_\eps] g_\eps \big) - \lL_0 \, \big( [\iI_0, f \prec] g_0 \big) \Big)=\Delta_{k,\eps} \lL_\eps \, \big( [\iI_\eps, f \prec_\eps]\lL_{\eps}(\lL_{\eps}^{-1}\Pi_{0}^{\perp} g_\eps-\lL_{0}^{-1}\Pi_{0}^{\perp} g_0) \big)\\
    +&(\Delta_{k,0}-\Delta_{k,\eps}) \lL_\eps \, \big( [\iI_\eps, f \prec_\eps]g_{\eps}\big)+ \Big(\Delta_{k,\eps} \lL_\eps \, \big( [\iI_\eps, f \prec_\eps]\lL_{\eps}\lL_{0}^{-1}\Pi_0^{\perp} g_0 \big) - \Delta_{k,0}\lL_0 \, \big( [\iI_0, f \prec] g_0 \big) \Big) \;.
    \end{split}
\end{equation*}
Fix $\kappa'>0$ arbitrary but such that $\alpha+\beta-2\kappa'>-1$. By Lemma~\ref{le:heat_comm_spacetime} and Corollary~\ref{cor:convergence_single_block} (with $m=0$), there exists $\theta'>0$ such that the $L^\infty$-norm of the first two terms are bounded by
\begin{equation*}
    \begin{split}
    2^{-(\alpha+\beta-\kappa) k} r^{-\frac{\sigma}{2}} \|f\|_{\mathfrak{L}_{t,\sigma}^{\alpha}} (\eps^{\theta'}\|\Pi_{0}^{\perp}g_{0}\|_{L^{\infty}_{t}\cC^{\beta;0}_{x}}+ \|\lL_\eps^{-1}\Pi_{0}^{\perp} g_\eps-\lL_0^{-1}\Pi_{0}^{\perp} g_0\|_{L^{\infty}_{t}\cC^{\beta+2}_{x}})\;.
    \end{split}
\end{equation*}
We now turn to the third term. Since $\beta + 2 \in (0,1)$, we have
\begin{equation*}
    \|\lL_\eps \lL_0^{-1} \Pi_0^\perp g_0\|_{\beta;\eps} \lesssim \|\lL_0^{-1} \Pi_0^\perp g_0\|_{\beta+2} \lesssim \|g_0\|_{\beta;0}\;.
\end{equation*}
Hence, by Lemma~\ref{le:heat_comm_spacetime}, the third term is bounded by
    \begin{equation*}
         2^{-(\alpha+\beta) k} r^{-\frac{\sigma}{2}} \|f\|_{\mathfrak{L}_{t,\sigma}^{\alpha}} \|\Pi_0^{\perp}g_{0}\|_{L^{\infty}_{t}\cC^{\beta;0}_{x}}\;.
    \end{equation*}
On the other hand, we rewrite it as
\begin{equation*}\begin{split}
    &(\Delta_{k,\eps} \lL_\eps-\Delta_{k,0}\lL_0) \, \big( [\iI_\eps, f \prec_\eps]\lL_{\eps}\lL_{0}^{-1}\Pi_0^{\perp} g_0 \big)\\
    +&\Delta_{k,0}\lL_0\Big(\big( [\iI_\eps, f \prec_\eps]\lL_{\eps}\lL_{0}^{-1}\Pi_0^{\perp} g_0 \big) - \, \big( [\iI_0, f \prec] \lL_{0}\lL_{0}^{-1}\Pi_0^{\perp} g_0  \big) \Big)\;.
\end{split}\end{equation*}
Writing out the commutators $[\iI_\eps, f \prec_\eps]$ (including $\eps=0$) and applying Lemmas~\ref{le:convergence_heat_prec_below} and~\ref{le:convergence_prec_heat_below}, we see its $L^\infty$-norm is bounded by
\begin{equation*}
    \eps^{\theta'}2^{Mk} r^{-\frac{\sigma}{2}} \|f\|_{\mathfrak{L}_{t,\sigma}^{\alpha}} \|\Pi_{0}^{\perp}g_{0}\|_{L^{\infty}_{t}\cC^{\beta;0}_{x}}\;.
\end{equation*}
The proof is then complete with suitable interpolations. 
\end{proof}

Similar as above, we also have the following convergence with $\beta$ in a different range. 

\begin{lem} \label{le:convergence_heat_comm_good}
    Let $\alpha \in (0,1)$, $\beta \in (-1,0)$ and $\sigma \in [0,2)$. Then for every $\kappa>0$, there exists $\theta>0$ such that
    \begin{equation*}
        \begin{split}
        &\phantom{111}\|\lL_\eps \big( [\iI_\eps, f_\eps \prec_\eps] g_\eps \big) - \lL_0 \big( [\iI_0, f_0 \prec] g_0 \big)\|_{L_{\sigma,t}^{\infty} \cC_{x}^{\alpha+\beta-\kappa}}\\
        &\lesssim \|f_0\|_{\mathfrak{L}_{\sigma, t}^{\alpha}} \big( \eps^\theta \|g_0\|_{L_t^\infty \cC_x^\beta} + \|g_\eps - g_0\|_{L_t^\infty \cC_x^\beta} \big) + \|f_\eps - f_0\|_{\mathfrak{L}_{\sigma, t}^{\alpha}} \|g_\eps\|_{L_t^\infty \cC_x^\beta}\;.
        \end{split}
    \end{equation*}
\end{lem}
\begin{proof}
The proof is essentially the same as that for Lemma~\ref{le:convergence_heat_comm_below} but only simpler, since we have $\beta \in (-1,0)$ and hence can compare $g_\eps$ with $g_0$ directly. We omit the proof here. 
\end{proof}

\section{Bounds involving derivatives specific to $d=1$}
\label{sec:bounds_derivatives_1d}

In this section, we restrict to the one dimensional torus $\T$. As before, for uniform bounds, operations are associated with coefficient matrix $a \in \sS_1(\Lambda, M)$. For convergence, operations are associated with $A_\eps$ and $\bar{A}$.

\subsection{Uniform bounds}

All proportionality constants below are independent of the choice of $a \in \sS_1(\Lambda,M)$.

\begin{lem}\label{le:derivative_a_1d}
    For $\alpha\in (1,2)$, we have the bound
    \begin{equation*}
        \|a\d_{x}f\|_{\alpha-1} \lesssim \|f\|_{\alpha;a}\;.
    \end{equation*}
\end{lem}
\begin{proof}
We need to show that both $\|\d_x f\|_{L^\infty}$ and $\|a \d_x f\|_{\dot{\cC}^{\alpha-1}}$ are controlled by $\|f\|_{\alpha;a}$. Note that
\begin{equation*}
    f'(x) = \sum_{j \geq -1} \sum_{\ell=j-1}^{j+1} \int_{\T} (\d_x \Delta_{j,a})(x,y) \cdot \big( \Delta_{\ell,a} f \big)(y) {\rm d}y\;.
\end{equation*}
By Proposition~\ref{pr:kernel_pointwise_bound}, we have
\begin{equation*}
    \|a \d_x f\|_{L^\infty} \lesssim \|\d_x f\|_{L^\infty} \lesssim \sum_{j \geq -1} 2^{j} \sum_{\ell=j-1}^{j+1} \|\Delta_{\ell,a} f\|_{L^\infty} \lesssim \|f\|_{\alpha;a}\;.
\end{equation*}
We now turn to the homogeneous $\dot{\cC}^{\alpha-1}$-norm. It suffices to note that the derivative coincides with the divergence in dimension one, so that
\begin{equation*}
    \|a \d_x f\|_{\dot{\cC}^{\alpha-1}} \lesssim \|\d_x \big( a \d_x f \big)\|_{\alpha-2} = \|\lL_a f\|_{\alpha-2} \lesssim \|f\|_{\alpha;a}\;,
\end{equation*}
where in the last step we have used $\alpha \in (1,2)$ so that $\|\cdot\|_{\alpha-2} \sim \|\cdot\|_{\alpha-2;a}$. This completes the proof. 
\end{proof}

Recall definitions for $\mathbf{\Pi}_{a}$ and $\P_{a}$.
\begin{equation*}
    \mathbf{\Pi}_{a}(f,g) := \sum_{j=-1}^{\infty}\sum_{|i-j|\leq 1}a\d_{x}\Delta_{i,a}f\cdot\d_{x}\Delta_{j,a}g\;;\quad \P_{a}(f,g):=\sum_{j=-1}^{\infty}a\d_{x}S_{j,a}f\cdot\d_{x}\Delta_{j,a}g\;.
\end{equation*}
We have the following lemma. 

\begin{lem} \label{le:alternative_prec_circ_1d}
    Let $\alpha \in (1,2)$ and $\beta \in (0,1)$ be such that $\alpha+\beta \in (1,2)$. Then we have the bounds
    \begin{equation} \label{e:alternative_prec_1d}
        \| \P_{a} ( f , g )-( a \d_{ x } f ) \prec_{a}\d_{x}g\|_{\alpha+\beta-2;a}\lesssim \|f\|_{\alpha;a}\|g\|_{\beta;a}
    \end{equation}
    and
    \begin{equation} \label{e:alternative_circ_1d}
        \| \mathbf{\Pi}_{a}(f,g)-(a\d_{x}f)\circ_{a}\d_{x}g\|_{\alpha+\beta-2;a}\lesssim \|f\|_{\alpha;a}\|g\|_{\beta;a}.
    \end{equation}
\end{lem}
\begin{proof}
We need to show that the $L^\infty$-norms of the $k$-th block of both quantities are controlled by $2^{(2-\alpha-\beta)k}$. We start with \eqref{e:alternative_prec_1d}. 

For every $k \geq -1$, we have the identity
\begin{equation} \label{e:alternative_prec_1d_expression}
    \begin{split}
    \Delta_{k,a} \big( (a \, \d_x f) &\prec_a \d_x g - \P_a (f,g) \big) = R_{k,a} \big(a \, \d_x f, \d_x g\big) + \Delta_{k,a} \big( (\lL_a f) \prec_a g \big)\\
    &+ \sum_{j \geq -1} \Big(  [a \d_x S_{j,a}f\,, \Delta_{k,a}\nabla ]  \Delta_{j,a} g + a \d_x S_{j,a}^{\perp} f \cdot \Delta_{k,a} \big(\d_x \Delta_{j,a} g\big) \Big)\;,
    \end{split}
\end{equation}
where $R_{k,a}$ is given in \eqref{e:R_op_defn}. We treat the four terms on the right hand side above one by one. For the first one, by Lemmas~\ref{le:R_Linfty_bound} and~\ref{le:derivative_a_1d}, we have
\begin{equation*}
    \left\| R_{k,a}(a \, \d_x f, \d_x g)  \right\|_{L^\infty} \lesssim 2^{(2-\alpha-\beta) k} \|a \, \d_x f\|_{\alpha-1;a} \|\d_x g\|_{\beta-1; a} \lesssim 2^{(2-\alpha-\beta) k} \|f\|_{\alpha;a} \|g\|_{\beta;a}\;.
\end{equation*}
For the second term, it follows from Proposition~\ref{pr:paraproducts_gBesov} and $\alpha<2$ and $\alpha-2+\beta>1$ that
\begin{equation*}
    \big\| \Delta_{k,a} \big( (\lL_a f) \prec_a g \big) \big\|_{L^\infty} \lesssim 2^{-(\alpha-2+\beta)k} \big\| (\lL_a f) \prec_a g \big\|_{\alpha-2+\beta;a} \lesssim 2^{(2-\alpha-\beta)k} \|f\|_{\alpha;a} \|g\|_{\beta;a}\;.
\end{equation*}
For the third term on the right hand side of \eqref{e:alternative_prec_1d_expression}, integrating by parts, we have
\begin{equation*}
    \begin{split}
    &\phantom{111}\big([a \nabla S_{j,a} f, \Delta_{k,a} \nabla] \Delta_{j,a} g\big)(x)\\
    &= \int_\T (\d_y \Delta_{k,a})(x,y) \big( (a \nabla S_{j,a} f)(y) - (a \nabla S_{j,a} f)(x) \big) (\Delta_{j,a} g)(y) {\rm d}y\;,
    \end{split}
\end{equation*}
where we use $\nabla$ to replace $\d_x$ or $\d_y$ to avoid confusions. By the kernel estimate in Lemma~\ref{lem:blocks_derivatives_pointwise}, for $\gamma \in (0,1]$, we have
\begin{equation*}
    \begin{split}
    \big\| [a \nabla S_{j,a} f, \Delta_{k,a} \nabla] \Delta_{j,a} g \big\|_{L^\infty} &\lesssim 2^{(1-\gamma)k} \|a \d_x S_{j,a} f\|_{\gamma} \|\Delta_{j,a} g\|_{L^\infty}\\
    &\lesssim 2^{(1-\gamma) k} 2^{-\beta j} \|a \d_x S_{j,a} f\|_{\dot{\cC}^\gamma} \|g\|_{\beta;a}\;.
    \end{split}
\end{equation*}
For $\gamma=1$ and $\gamma=\alpha-1$, we respectively have
\begin{equation*}
    \begin{split}
    \|a \d_x S_{j,a} f\|_{\dot{\cC}^1} &\lesssim \|S_{j,a} \lL_a f\|_{L^\infty} \lesssim 2^{(2-\alpha) j} \|f\|_{\alpha;a}\;\\
    \|a \d_x S_{j,a} f\|_{\alpha-1} &\lesssim \|S_{j,a} f\|_{\alpha;a} \lesssim \|f\|_{\alpha;a}\;.
    \end{split}
\end{equation*}
Hence, for $j \leq k$ and $j > k$, we take $\gamma=1$ and $\gamma= \alpha-1$ respectively so that
\begin{equation*}
    \begin{split}
    \big\| [a \nabla S_{j,a} f, \Delta_{k,a} \nabla] \Delta_{j,a} g \big\|_{L^\infty} &\lesssim 2^{(2-\alpha-\beta)j} \|f\|_{\alpha;a} \|g\|_{\beta;a}\\
    \big\| [a \nabla S_{j,a} f, \Delta_{k,a} \nabla] \Delta_{j,a} g \big\|_{L^\infty} &\lesssim 2^{(2-\alpha) k} 2^{-\beta j} \|f\|_{\alpha;a} \|g\|_{\beta;a}\;.
    \end{split}
\end{equation*}
Summing the above bounds over $j \geq -1$, we get
\begin{equation*}
    \sum_{j \geq -1} \big\| [a \nabla S_{j,a} f, \Delta_{k,a} \nabla] \Delta_{j,a} g \big\|_{L^\infty} \lesssim 2^{(2-\alpha-\beta)k} \|f\|_{\alpha;a} \|g\|_{\beta;a}\;.
\end{equation*}
Now we turn to the last term on the right hand side of \eqref{e:alternative_prec_1d_expression}. By Proposition~\ref{pr:kernel_pointwise_bound} and that $\alpha>1$, we have
\begin{equation*}
    \|a \d_x S_{j,a}^{\perp} f\|_{L^\infty} \lesssim \sum_{\ell \geq j-1} \sum_{q=\ell-1}^{\ell+1} \|\nabla \Delta_{\ell,a} \Delta_{q,a} f\|_{L^\infty}  \lesssim 2^{(1-\alpha) j} \|f\|_{\alpha;a}\;,
\end{equation*}
and
\begin{equation*}
    \|\Delta_{k,a} \big( \d_x \Delta_{j,a} g \big)\|_{L^\infty} \lesssim 2^{j \wedge k} \cdot 2^{-\beta j} \|g\|_{\beta;a}\;.
\end{equation*}
Summing over $j \geq -1$ then gives the desired bound on the last term. This completes the proof for \eqref{e:alternative_prec_1d}. 

As for \eqref{e:alternative_circ_1d}, note that we have the identity
\begin{equation*}
    \begin{split}
    &\phantom{111}\mathbf{\Pi}_a(f,g) - (a \d_x f) \circ_a \d_x g\\
    &= (a \d_x f) \succ_a \d_x g + \Big( (a \d_x f) \prec_a (\d_x g) - \P_a(f,g) \Big) - \P_a (g,f)\;.
    \end{split}
\end{equation*}
The desired bounds for the three terms follow from Lemma~\ref{le:derivative_a_1d}, the bound \eqref{e:alternative_prec_1d} and Lemma~\ref{le:prec_alternative_bound}. The proof is then complete. 
\end{proof}

\begin{lem} \label{le:prec_kpz}
Let $\alpha \in (0,1)$ and $\beta \in (1,2)$ be such that $\alpha + \beta < 2$. Then, we have
    \begin{equation} \label{e:prec_kpz}
    \sum_{j=-1}^{\infty} \left\| a \, \d_{x}(S_{j,a}f) \cdot \Delta_{j,a}g \right\|_{\alpha+\beta-1;a}\lesssim \|f\|_{\alpha;a}\|g\|_{\beta;a}\;.
    \end{equation}
\end{lem}
\begin{proof}
We split the sum into $j \leq k-2$ and $j \geq k-1$. For $j \leq k-2$, we have
\begin{equation} \label{e:prec_kpz_low}
    \begin{split}
    \left\| \Delta_{k,a} \big( a \, \d_{x}(S_{j,a}f) \cdot \Delta_{j,a}g \big)  \right\|_{L^\infty} &= \left\| \, \big[ \Delta_{k,a}\,, a \d_x (S_{j,a} f) \, \big] \Delta_{j,a} g \right\|_{L^\infty}\\
    &\lesssim 2^{-k} \|\lL_a S_{j,a} f\|_{L^\infty} \|\Delta_{j,a} g\|_{L^\infty}\\
    &\lesssim 2^{-k + (2-\alpha-\beta)j} \|f\|_{\alpha;a} \|g\|_{\beta;a}\;,
    \end{split}
\end{equation}
where the identity on the first line holds since $\Delta_{k,a} \Delta_{j,a} = 0$ for $j \leq k-2$, and the inequality on the second line follows from Lemma~\ref{le:block_multiplication_comm}. 

For $j \geq k-1$, we have
\begin{equation} \label{e:prec_kpz_high}
    \begin{split}
    \left\| \Delta_{k,a} \big( a \, \d_{x}(S_{j,a}f) \cdot \Delta_{j,a}g \big)  \right\|_{L^\infty} &\lesssim \|a \, \d_{x}(S_{j,a}f)\|_{L^\infty} \|\Delta_{j,a} g\|_{L^\infty}\\
    &\lesssim 2^{-(\alpha+\beta-1) j} \|f\|_{\alpha;a} \|g\|_{\beta;a}\;.
    \end{split}
\end{equation}
Summing over $j \leq k-2$ with the bound \eqref{e:prec_kpz_low} and $j \geq k-1$ with the bound \eqref{e:prec_kpz_high} gives the desired bound and hence completes the proof of the lemma. 
\end{proof}

\begin{lem} \label{le: R_operator_KPZ}
  Let $\alpha \in (0,1)$ and $\beta>1$ be such that $\alpha + \beta \in (1,2)$. Then we have the bound
    \begin{equation*}\begin{split}
         \big \|\Delta_{k,a}\big(S_{j,a}f\cdot a\d_{x}\Delta_{j,a}g\big)- f \Delta_{k,a} \big(a\d_{x} \Delta_{j,a}g \big)\big\|_{L^{\infty}}\lesssim 2^{(1-\alpha) (k \wedge j) -k+(1-\beta)j}\|f\|_{\alpha;a}\|g\|_{\beta;a}
    \end{split}
    \end{equation*}
    for all $k \geq -1$. 
\end{lem}
\begin{proof}
We have the expression
\begin{equation} \label{e:R_operator_kpz_decomposition}
    \begin{split}
    &\phantom{111}\Delta_{k,a} \big(S_{j,a} f \cdot a \, \d_{x} \Delta_{j,a}g\big)- f \Delta_{k,a} \big( a \, \d_{x} \Delta_{j,a} g \big)\\
    &= [\Delta_{k,a}\,, S_{j,a} f] \, \big( a \, \d_x \Delta_{j,a} g \big) - S_{j,a}^{\perp}f \cdot \Delta_{k,a} \big( a \, \d_x \Delta_{j,a} g \big)\;.
    \end{split}
\end{equation}
We treat the two terms on the right hand side above separately. For the first one, by the first estimate in Lemma~\ref{le:block_multiplication_comm}, we have
\begin{equation*}
    \left\| [\Delta_{k,a}\,, S_{j,a} f] \, \big( a \, \d_x \Delta_{j,a} g \big)  \right\|_{L^\infty} \lesssim 2^{-\alpha k} \|S_{j,a} f\|_{\alpha;a} \|a\d_{x}\Delta_{j,a} g\|_{L^\infty} \lesssim 2^{-\alpha k + (1-\beta) j} \|f\|_{\alpha;a} \|g\|_{\beta;a}\;.
\end{equation*}
On the other hand, we can also use the second estimate in Lemma~\ref{le:block_multiplication_comm} to control it by
\begin{equation*}
    \left\| [\Delta_{k,a}\,, S_{j,a} f] \, \big( a \, \d_x \Delta_{j,a} g \big)  \right\|_{L^\infty} \lesssim 2^{-k} \|\d_x S_{j,a} f\|_{L^\infty} \|a\d_{x}\Delta_{j,a} g\|_{L^\infty} \lesssim 2^{-k + (2-\alpha-\beta) j} \|f\|_{\alpha;a} \|g\|_{\beta;a}\;.
\end{equation*}
They together give the desired bound for the first term on the right hand side of \eqref{e:R_operator_kpz_decomposition}. We now turn to the second term. It is bounded by
\begin{equation*}
\begin{split}
    \left\| (S_{j,a}^{\perp}f) \cdot \Delta_{k,a} \big( a \, \d_x \Delta_{j,a} g \big) \right\|_{L^\infty} \lesssim & \|S_{j,a}^{\perp} f\|_{L^\infty} \|\Delta_{k,a} \big( a \, \d_x \Delta_{j,a} g \big) \|_{L^\infty}\\
    \lesssim & 2^{-\alpha j}\|f\|_{\alpha;a} \|\Delta_{k,a} \big( a \, \d_x \Delta_{j,a} g \big) \|_{L^\infty}\;.\end{split}
\end{equation*}
There are two ways to control the norm on the right hand side above. First, it can be directly controlled by
\begin{equation*}
    \|\Delta_{k,a} \big( a \, \d_x \Delta_{j,a} g \big) \|_{L^\infty} \lesssim \|\d_x \Delta_{j,a} g\|_{L^\infty} \lesssim 2^{(1-\beta) j} \|g\|_{\beta;a}\;.
\end{equation*}
Also, noting that
\begin{equation*}
    \Delta_{k,a} \big( a \, \d_x \Delta_{j,a} g \big) = [\Delta_{k,a}\,, a \d_x \big( \Delta_{j,a} g \big)] \1_{\T} + \1_{k=-1} \cdot a \d_x \Delta_{j,a} g\;,
\end{equation*}
we can also bound it by Lemma~\ref{le:block_multiplication_comm} as
\begin{equation*}
    \begin{split}
    \|\Delta_{k,a} \big( a \, \d_x \Delta_{j,a} g \big) \|_{L^\infty} &\lesssim 2^{-k} \|\lL_a \Delta_{j,a} g\|_{L^\infty} + \1_{k=-1} \|a \d_x \Delta_{j,a} g\|_{L^\infty}\\
    &\lesssim \big( 2^{-(k-j)} + \1_{k=-1} \big) \cdot 2^{(1-\beta)j} \|g\|_{\beta;a}\;.
    \end{split}
\end{equation*}
Combing the above and using $\alpha \in (0,1)$, we see that the $L^\infty$-norm of the second term on the right hand side of \eqref{e:R_operator_kpz_decomposition} also satisfies the desired bound. This completes the proof. 
\end{proof}

\begin{prop}\label{pr:a_derivative_prec_comm}
Let $\alpha \in (0,1)$ and $\beta \in (1,2)$ be such that $\alpha + \beta \in (1,2)$. Then we have
    \begin{equation}\label{eq:a_derivative_prec_comm}
        \left\| a \, \d_{x}(f\prec_{a}g)-f\prec_{a}(a \, {\d_{x}}g) \right\|_{\alpha+\beta-1;a}\lesssim \|f\|_{\alpha;a}\|g\|_{\beta;a}\;.
    \end{equation}
\end{prop}
\begin{proof}
We have the identity
\begin{equation} \label{e:a_derivative_prec_comm_decomposition}
    \begin{split}
    a \, \d_{x}(f\prec_{a}g) - &f\prec_{a}(a \, {\d_{x}}g) = \sum_{j \geq -1} a \, (\d_x S_{j,a} f) \cdot (\Delta_{j,a} g) \\
    &+ \sum_{j \geq -1} (S_{j,a} f) \cdot a \d_x ( \Delta_{j,a} g) - f \prec_a (a \, \d_x g)\;.
    \end{split}
\end{equation}
The desired bound for the first term above is precisely Lemma~\ref{le:prec_kpz}. For the terms on the second line above, its $k$-th block can be written as
\begin{equation*}
    \sum_{j \geq -1} \Big( \Delta_{k,a} \big( (S_{j,a} f) \cdot a \d_x ( \Delta_{j,a} g) \big) - f \Delta_{k,a} \big(a \d_x \Delta_{j,a} g\big) \Big) - R_{k,a} (f, \d_x g)\;.
\end{equation*}
It follows from Lemma~\ref{le: R_operator_KPZ} and Lemmas~\ref{le:block_multiplication_comm} and~\ref{le:derivative_a_1d} that the $L^\infty$-norms of the above two terms are both bounded by
\begin{equation*}
    2^{{(1-\alpha-\beta)k}} \|f\|_{\alpha;a} \|g\|_{\beta;a}\;.
\end{equation*}
The proof of the proposition is thus complete. 
\end{proof}

\subsection{Convergences}

In this subsection, we fix the coefficient matrix $A$ satisfying Assumption~\ref{as:A} (on one dimensional torus $\T$), and let $A_\eps = A (\cdot / \eps)$. As before, we write $\lL_\eps = \div(A_\eps \nabla)$ and $\lL_0 = \div (\bar{A} \nabla)$ for the homogenised matrix $\bar{A}$. All proportionality constants below are independent of $\eps \in \N^{-1}$.

\begin{lem} \label{le:convergence_derivative_a_1d}
    Let $\alpha \in (-1,0)$. Then for every $\kappa>0$, there exists $\theta>0$ such that
    \begin{equation} \label{e:convergence_derivative_a_1d}
        \|A_{\eps}\d_{x}\lL_{\eps}^{-1}g_{\eps}- \bar{A} \d_{x}\lL_{0}^{-1} g_0\|_{\alpha+1-\kappa} \lesssim \eps^{\theta} \|g_0\|_{\alpha}+\|g_{\eps}-g_{0}\|_{\alpha}
    \end{equation}
    uniformly over $\eps\in \N^{-1}$ and $g_{\eps}, g_0 \in \cC^{\alpha}$ such that $\Pi_0 g_{\eps} = \Pi_0 g_0 = 0$. 
\end{lem}
\begin{proof}
We write the left hand side of \eqref{e:convergence_derivative_a_1d} as
\begin{equation*}
    A_{\eps}\d_{x}\lL_{\eps}^{-1}g_{\eps}- \bar{A} \d_{x}\lL_{0}^{-1} g_0 = A_\eps \d_x \lL_\eps^{-1} (g_\eps - g_0) + \big( A_\eps \d_x \lL_\eps^{-1} - \bar{A} \d_x \lL_0^{-1} \big) g_0\;.
\end{equation*}
For the first term, it follows from Lemma~\ref{le:derivative_a_1d} directly that
\begin{equation*}
    \big\| A_\eps \d_x \lL_\eps^{-1} (g_\eps - g_0) \big\|_{\alpha+1} \lesssim \|\lL_\eps^{-1} (g_\eps - g_0)\|_{\alpha+2; \eps} \lesssim \|g_\eps - g_0\|_\alpha\;.
\end{equation*}
For the second one, similarly, we have the uniform bound
\begin{equation} \label{e:convergence_derivative_a_1d_uniform}
    \big\| \big( A_\eps \d_x \lL_\eps^{-1} - \bar{A} \d_x \lL_0^{-1} \big) g_0 \big\|_{\alpha+1} \lesssim \|g_0\|_\alpha\;.
\end{equation}
To obtain smallness in $\eps$, we rewrite it as
\begin{equation*}
    \big( A_\eps \d_x \lL_\eps^{-1} - \bar{A} \d_x \lL_0^{-1} \big) g_0 = \sum_{j \geq -1} \big( A_\eps \d_x \lL_\eps^{-1} - \bar{A} \d_x \lL_0^{-1} \big) \Delta_{j,0} g_0\;.
\end{equation*}
Let $\kappa'>0$ be arbitrarily small. For each term in the sum, again by Lemma~\ref{le:derivative_a_1d}, we have the uniform-in-$\eps$ bound
\begin{equation*}
    \|A_\eps \d_x \lL_\eps^{-1} \Delta_{j,0} g_0\|_{\kappa'} \lesssim \|\lL_\eps^{-1} \Delta_{j,0} g_0\|_{1+\kappa'; \eps} \lesssim \|\Delta_{j,0} g_0\|_{-1+\kappa'} \lesssim 2^{-(1+\alpha-\kappa')j} \|g_0\|_\alpha\;.
\end{equation*}
On the other hand, by Lemma~\ref{lem:elliptic_homo_Green_pointwise_expansion}, we have
\begin{equation*}
    \| \big( A_\eps \d_x \lL_\eps^{-1} - \bar{A} \d_x \lL_0^{-1} \big) \Delta_{j,0} g_0 \|_{-\kappa''} \lesssim \sqrt{\eps} 2^{Mj} \|g_0\|_\alpha
\end{equation*}
for some $M>0$. Interpolating the above two bounds and summing over $j \geq -1$, we obtain
\begin{equation} \label{e:convergence_derivative_a_1d_small}
    \big\| \big( A_\eps \d_x \lL_\eps^{-1} - \bar{A} \d_x \lL_0^{-1} \big) g_0 \big\|_{L^\infty} \lesssim \eps^\theta \|g_0\|_\alpha\;.
\end{equation}
Interpolating \eqref{e:convergence_derivative_a_1d_uniform} and \eqref{e:convergence_derivative_a_1d_small} gives the desired bound and thus completes the proof of the lemma. 
\end{proof}

\begin{prop} \label{pr:homo_diff_a_derivative_prec_comm}
Let $\alpha \in (0,1)$ and $\beta \in (1,2)$ be such that $\alpha + \beta \in (1,2)$. Then for every $\kappa>0$, there exists $\theta>0$ such that
    \begin{equation}\label{eq:homo_diff_a_derivative_prec_comm}
    \begin{split}
        &\left\| [A_{\eps} \, \d_{x},f_{\eps}\prec_{\eps}]G_{\eps} - [\bar{A} \, \d_{x},f_0\prec]G_{0}\right\|_{\alpha+\beta-1-\kappa}\\
        \lesssim &\|f_0\|_{\alpha}(\eps^{\theta}\|G_0\|_{\beta;0}+\|\lL_{\eps}G_{\eps}-\lL_{0}G_0\|_{\beta-2})+\|f_{\eps}-f\|_{\alpha}\|G_{\eps}\|_{\beta;\eps}\;.
    \end{split}
    \end{equation}
\end{prop}
\begin{proof}
Since $\Pi_0 G_\eps$ and $\Pi_0 G_0$ play no role in the quantity of study, we can assume without loss of generality that $\Pi_0 G_\eps = \Pi_0 G_0 = 0$. Let $g_\eps = \lL_\eps G_\eps$ (including $\eps=0$). Then the left hand side of \eqref{eq:homo_diff_a_derivative_prec_comm} can be expressed by
\begin{equation*}
    \begin{split}
    [A_\eps \, \d_x, (f_\eps - &f_0) \prec_\eps ] \lL_\eps^{-1} g_\eps + [A_\eps \, \d_x, \, f_0 \prec_\eps] \big( \lL_\eps^{-1} (g_\eps - g_0) \big)\\
    &+ \Big( [A_\eps \d_x, \, f_0 \prec_\eps] \lL_\eps^{-1} g_0 - [\bar{A} \, \d_x, \, f_0 \prec] \lL_0^{-1} g_0 \Big)\;.
    \end{split}
\end{equation*}
By Proposition~\ref{pr:a_derivative_prec_comm}, the first two terms above can be bounded by
\begin{equation*}
    \big\| \, [A_\eps \, \d_x, (f_\eps - f_0) \prec_\eps ] \lL_\eps^{-1} g_\eps \big\|_{\alpha+\beta-1} \lesssim \|f_\eps-f_0\|_\alpha \|G_\eps\|_{\beta;\eps}\;,
\end{equation*}
and
\begin{equation*}
    \begin{split}
    \big\| \, [A_\eps \, \d_x, \, f_0 \prec_\eps] \big( \lL_\eps^{-1} (g_\eps - g_0) \big) \big\|_{\alpha+\beta-1} &\lesssim \|f_0\|_\alpha \, \|\lL_\eps^{-1} (g_\eps - g_0)\|_{\beta;\eps}\\
    &\lesssim \|f_0\|_\alpha \, \|g_\eps - g_0\|_{\beta-2}\;.
    \end{split}
\end{equation*}
So it remains to consider the third term. To simplify notations, we drop the subscript $0$ from $f_0$ and $g_0$, and write $f=f_0$ and $g=g_0$. Again by Proposition~\ref{pr:a_derivative_prec_comm}, we have the uniform-in-$\eps$ bound (including $\eps=0$)
\begin{equation} \label{e:convergence_1d_comm_uniform}
    \big\| \, [A_\eps \d_x, \, f \prec_\eps] \lL_\eps^{-1} g  \big\|_{\alpha+\beta-1} \lesssim \|f\|_\alpha \, \|\lL_\eps^{-1} g\|_{\beta;\eps} \lesssim \|f\|_\alpha \, \|g\|_{\beta-2}\;.
\end{equation}
So it suffices to control the difference in a weaker norm (we choose $L^\infty$ in this case). Writing out the commutator and taking the difference term by term, we get
\begin{equation} \label{e:convergence_1d_comm_expand}
    \begin{split}
    [&A_\eps \d_x, \, f \prec_\eps] \lL_\eps^{-1} g - [\bar{A} \, \d_x, \, f \prec] \lL_0^{-1} g = \Big( A_\eps \d_x \big( f \prec_\eps \lL_\eps^{-1} g \big) - \bar{A} \d_x \big( f \prec \lL_0^{-1} g \big) \Big)\\
    &-  f \prec_\eps \big( ( A_\eps \d_x \lL_\eps^{-1} - \bar{A} \d_x \lL_0^{-1} ) g \big) - \Big( f \prec_\eps \big( \bar{A} \d_x \lL_0^{-1} g \big) - f \prec \big( \bar{A} \d_x \lL_0^{-1} g \big) \Big)\;.
    \end{split}
\end{equation}
The two terms on the second line above can be treated directly. For example, by Proposition~\ref{pr:paraproducts_gBesov} and Lemma~\ref{le:convergence_derivative_a_1d}, we have
\begin{equation*}
    \begin{split}
    \left\| f \prec_\eps \big( ( A_\eps \d_x \lL_\eps^{-1} - \bar{A} \d_x \lL_0^{-1} ) g \big) \right\|_{L^\infty} &\lesssim \|f\|_\alpha \, \big\| ( A_\eps \d_x \lL_\eps^{-1} - \bar{A} \d_x \lL_0^{-1} ) g \big\|_{\beta-1-\kappa'}\\
    &\lesssim \eps^\theta \|f\|_\alpha \, \|g\|_{\beta-2}\;,
    \end{split}
\end{equation*}
where $\kappa'>0$ is such that $\beta-1-\kappa' \in (0,1)$. Also, since $\beta-1 \in (0,1)$, by Lemma~\ref{le:convergence_prec_good}, we have
\begin{equation*}
    \begin{split}
    \left\| f \prec_\eps \big( \bar{A} \d_x \lL_0^{-1} g \big) - f \prec \big( \bar{A} \d_x \lL_0^{-1} g \big) \right\|_{L^\infty} &\lesssim \eps^\theta \|f\|_\alpha \, \|\bar{A} \d_x \lL_0^{-1} g\|_{\beta-1}\\
    &\lesssim \eps^\theta \|f\|_\alpha \, \|g\|_{\beta-2}\;.
    \end{split}
\end{equation*}
For the first term on the right hand side of \eqref{e:convergence_1d_comm_expand}, we write it as
\begin{equation*}
    \begin{split}
    &\phantom{111}A_\eps \d_x \big( f \prec_\eps \lL_\eps^{-1} g \big) - \bar{A} \d_x \big( f \prec \lL_0^{-1} g \big)\\
    &= \sum_{j \geq -1} \bigg[ A_\eps \d_x S_{j,\eps} f \cdot \big( \Delta_{j,\eps} \lL_\eps^{-1} g - \Delta_{j,0} \lL_0^{-1} g \big) + \big( A_\eps \d_x S_{j,\eps} f - \bar{A} \d_x S_{j,0} f \big) \cdot \Delta_{j,0} \lL_0^{-1} g\\
    &+ (S_{j,\eps} f - S_{j,0} f) \cdot A_\eps \d_x \lL_\eps^{-1} \Delta_{j,\eps} g + S_{j,0} f \big( A_\eps \d_x \lL_\eps^{-1} \Delta_{j,\eps} g - \bar{A} \d_x \lL_0^{-1} \Delta_{j,0} g \big) \bigg]\;.
    \end{split}
\end{equation*}
For each $j$, the $L^\infty$-norm for each term in the sum above can be controlled by $\eps^\theta 2^{-\eta j} \|f\|_\alpha \|g\|_{\beta-2}$ for some $\theta, \eta>0$. For example, by Corollary~\ref{cor:convergence_single_block} (with $m=-1$) and Lemma~\ref{le:derivative_a_1d}, we have
\begin{equation*}
    \begin{split}
    \|A_\eps \d_x S_{j.\eps} f\|_{L^\infty} \cdot \|(\Delta_{j,\eps} \lL_\eps^{-1} - \Delta_{j,0} \lL_0^{-1}) g\|_{L^\infty} &\lesssim \|S_{j,\eps} f\|_{1+\kappa;\eps} \cdot \eps^\theta 2^{-(\beta-\kappa)j} \|g\|_{\beta-2}\\
    &\lesssim \eps^\theta 2^{(1-\alpha-\beta+2\kappa)j} \|f\|_\alpha \, \|g\|_{\beta-2}\;.
    \end{split}
\end{equation*}
Also, writing $A_\eps \d_x S_{j,\eps} f = A_\eps \d_x \lL_\eps^{-1} \lL_\eps S_{j,\eps} f$ and using Lemma~\ref{le:convergence_derivative_a_1d}, we have
\begin{equation*}
    \begin{split}
    \|(A_\eps \d_x S_{j,\eps} - \bar{A} \d_x S_{j,0}) f\|_{L^\infty} &\lesssim \eps^\theta \|\lL_0 S_{j,0} f\|_{-\gamma} + \|(\lL_\eps S_{j,\eps} - \lL_0 S_{j,0)} f\|_{-\gamma}\\
    &\lesssim \eps^\theta 2^{(2-\alpha-\gamma+\kappa)j} \|f\|_\alpha\;,
    \end{split}
\end{equation*}
and
\begin{equation*}
    \|\Delta_{j,0} \lL_0^{-1} g\|_{L^\infty} \lesssim 2^{-\beta j} \|g\|_{\beta-2}\;.
\end{equation*}
Hence, taking $\gamma \in (0,1)$ such that $\alpha+\beta+\gamma>1$ and $\kappa$ sufficiently small, we get
\begin{equation*}
    \left\| \big( A_\eps \d_x S_{j,\eps} f - \bar{A} \d_x S_{j,0} f \big) \cdot \Delta_{j,0} \lL_0^{-1} g \right\|_{L^\infty} \lesssim \eps^\theta 2^{-\kappa' j} \|f\|_\alpha \, \|g\|_{\beta-2}
\end{equation*}
for some $\theta, \kappa' > 0$. The other two terms can be dealt with similarly. Summing over $j \geq -1$, we get
\begin{equation*}
    \left\| A_\eps \d_x \big( f \prec_\eps \lL_\eps^{-1} g \big) - \bar{A} \d_x \big( f \prec \lL_0^{-1} g \big) \right\|_{L^\infty} \lesssim \eps^\theta \|f\|_\alpha \, \|g\|_{\beta-2}\;.
\end{equation*}
Plugging all the above bounds back into \eqref{e:convergence_1d_comm_expand} and interpolating with \eqref{e:convergence_1d_comm_uniform}, we get
\begin{equation*}
    \left\| [A_\eps \d_x, \, f \prec_\eps] \lL_\eps^{-1} g - [\bar{A} \, \d_x, \, f \prec] \lL_0^{-1} g \right\|_{\alpha+\beta-1-\kappa} \lesssim \eps^\theta \|f\|_\alpha \, \|g\|_{\beta-2}\;.
\end{equation*}
We have thus completed the proof of the proposition. 
\end{proof}

We are now ready to prove the following main proposition. 

\begin{prop} \label{prop:derivative_heat_prec_comm}
Let $\alpha \in (0,1)$, $\beta \in (-1,0)$ and $\sigma \in (0,2)$. We have
\begin{equation*}
    \begin{split}
    &\phantom{111}\big\| [A_\eps \iI_\eps', f_\eps \prec_\eps] g_\eps - [\bar{A} \iI_0', f_0 \prec] g_0 \big\|_{L_{\sigma,T}^\infty \cC_x^{\alpha+\beta+1-2\kappa}}\\
    &\lesssim \|f_0\|_{\mathfrak{L}_{\sigma,T}^\alpha} \big( \eps^\theta \|g_0\|_{L_T^\infty \cC_x^\beta} + \|g_\eps - g_0\|_{L_T^\infty \cC_x^\beta} \big) + \|f_\eps - f\|_{\mathfrak{L}_{\sigma,T}^\alpha} \|g_\eps\|_{L_{T}^\infty \cC_x^\beta}\;.
    \end{split}
\end{equation*}
\end{prop}
\begin{proof}
For simplicity, we write $\gamma=\alpha+\beta$. Decompose the commutators $[A_\eps \iI_\eps', f_\eps \prec_\eps] g_\eps$ and $[\bar{A} \iI_0', f_0 \prec] g_0$ by
\begin{equation*}
    \begin{split}
    [A_\eps \iI_\eps', f_\eps \prec_\eps] g_\eps &= A_\eps \d_x \big( [\iI_\eps, f_\eps \prec_\eps] g_\eps \big) + [A_\eps \d_x, f_\eps \prec_\eps] \iI_\eps (g_\eps)\;,\\
    [\bar{A} \iI_0', f_0 \prec] g_0 &= \bar{A} \d_x \big( [\iI_0, f_0 \prec] g_0 \big) + [\bar{A} \d_x, f_0 \prec] \iI_0 (g_0)\;.
    \end{split}
\end{equation*}
We compare them term-wise. For the first one, by Lemmas~\ref{le:convergence_derivative_a_1d}, ~\ref{le:heat_comm_spacetime} and~\ref{le:convergence_heat_comm_good}, we have
\begin{equation*}
    \begin{split}
    &\phantom{111}\big\| A_\eps \d_x \big([ \iI_\eps, f_\eps \prec_\eps  ] g_\eps \big) - \bar{A} \d_x \big( [\iI_0, f_0 \prec] g_0 \big) \big\|_{L_{\sigma,T}^\infty \cC_x^{\gamma+1-2\kappa}}\\
    &\lesssim \eps^\theta \big\| \lL_0 \big( [\iI_0, f_0 \prec] g_0 \big) \big\|_{L_{\sigma,T}^\infty \cC_x^{\gamma-\kappa}} + \big\| \lL_\eps \big( [\iI_\eps, f_\eps \prec_\eps] g_\eps \big) - \lL_0 \big( [\iI_0, f_0 \prec] g_0 \big) \big\|_{L_{\sigma,T}^\infty \cC_x^{\gamma-\kappa}}\\
    &\lesssim \|f_0\|_{\mathfrak{L}_{\sigma,T}^\alpha} \big( \eps^\theta \|g_0\|_{L_T^\infty \cC_x^\beta} + \|g_\eps - g_0\|_{L_T^\infty \cC_x^\beta} \big) + \|f_\eps - f\|_{\mathfrak{L}_{\sigma,T}^\alpha} \|g_\eps\|_{L_{T}^\infty \cC_x^\beta}\;.
    \end{split}
\end{equation*}
For the second one, by Proposition~\ref{pr:homo_diff_a_derivative_prec_comm} and Lemmas~\ref{le:heat_comm_spacetime} and~\ref{le:convergence_heat_spacetime_normal}, we have
\begin{equation*}
    \begin{split}
    &\phantom{111}\big\| [A_\eps \d_x, f_\eps \prec_\eps] \iI_\eps (g_\eps) - [\bar{A} \d_x, f_0 \prec] \iI_0(g_0)  \big\|_{L_{\sigma,T}^\infty \cC_x^{\gamma+1-2\kappa}}\\
    &\lesssim \|f_0\|_{L_{\sigma,T}^\infty \cC_x^\alpha} \Big( \eps^\theta \|\iI_0 (g_0)\|_{L_T^\infty \cC_x^{\beta+2-\kappa}} + \|\lL_\eps \iI_\eps (g_\eps) - \lL_0 \iI_0 (g_0)\|_{L_T^\infty \cC_x^{\beta-\kappa}} \Big)\\
    &\phantom{11111}+ \|f_\eps - f\|_{L_{\sigma,T}^\infty \cC_x^\alpha} \|\iI_\eps (g_\eps)\|_{L_T^\infty \cC_x^{\beta+2-\kappa}}\\
    &\lesssim \|f_0\|_{L_{\sigma,T}^\infty \cC_x^\alpha} \big( \eps^\theta \|g_0\|_{L_T^\infty \cC_x^\beta} + \|g_\eps - g_0\|_{L_{T}^\infty \cC_x^\beta} \big) + \|f_\eps - f_0\|_{L_{\sigma,T}^\infty \cC_x^\alpha} \|g_\eps\|_{L_T^\infty \cC_x^\beta}\;.
    \end{split}
\end{equation*}
Combining the above gives the claim of the proposition. 
\end{proof}

\section{The dynamical $\phi^4_3$ model}
\label{sec:phi4}

In this section, we study periodic homogenisation for dynamical $\phi^4_3$ model, given by
\begin{equation} \label{e:phi4_hom_full}
    \d_t \phi_\eps^{(\delta)} = (\lL_\eps-1) \phi_\eps^{(\delta)} - (\phi_\eps^{(\delta)})^3 + \xi^{(\delta)} + C_\eps^{(\delta)} \phi_\eps^{(\delta)}\;, \quad (t,x) \in \R^+ \times \T^3\;.
\end{equation}
Here, $\xi^{(\delta)}$ is a regularised version of the spacetime white noise $\xi$, $\lL_\eps = \div (A_\eps \nabla)$ with $A_\eps = A(\cdot / \eps)$ for a fixed coefficient matrix $A$ satisfying Assumption~\ref{as:A}, and $C_\eps^{(\delta)}$ is a renormalisation function (depending on $x$ but independent of time) specified in \eqref{eq:phi4_renormalisation} below. We use $\lL_\eps -1$ instead of $\lL_\eps$ to make the stochastic objects stationary-in-time. This additional linear term can be added back in the renormalisation $C_\eps^{(\delta)}$. 

The aim of this section is to show that assuming the convergence of proper stochastic objects (which is shown in \cite{homo_stochastic} and stated in Proposition~\ref{prop:phi43_stochastic} below), the solution $\phi_\eps^{(\delta)}$ to \eqref{e:phi4_hom_full} converges to the solution $\phi_0$ to the standard dynamical $\phi^4_3$ model, formally given by
\begin{equation} \label{e:phi4_limiting_equation}
    \d_t \phi = (\lL_0 - 1) \phi - \phi^3 + \xi + \infty \phi\;.
\end{equation}
Here, $\lL_0 = \div(\bar{A} \nabla)$ with $\bar{A}$ being the homogenised matrix of $A$. 

We introduce some notations before proceeding to details. For $\eps\in \N^{-1}$, we use $\iI_\eps$ and $\textbf{I}_\eps$ to denote the heat kernel operations
\begin{equation}\label{eq: definition for heat kernel Italian I}
    (\iI_{\eps}f)(t) := \int_{0}^{t}e^{(t-\tau)(\lL_{\eps}-1)}f(\tau){\rm d}\tau
\end{equation}
and
\begin{equation} \label{eq: definition for heat kernel textbf I}
       (\textbf{I}_{\eps}f)(t):=\int_{-\infty}^{t}e^{(t-\tau)(\lL_{\eps}-1)}f(\tau){\rm d}\tau\;.
\end{equation}
In the rest of this section, we fix the parameters
\begin{equation} \label{eq:parameters_phi43}
    \eta \in \left(0, \frac{1}{1000}\right)\;, \quad \sigma = \frac{1}{2} + \eta\;, \quad \kappa \in \left(0, \frac{\eta}{1000}\right)\;.
\end{equation}

\subsection{Definition of stochastic objects and main result}
\label{sec:phi4_overall}

For $\delta \in (0,1)$, let $\xi^{(\delta)}$ be a regularisation in space of $\xi$ at scale $\delta$. Define the stationary-in-time process $\<1>_\eps^{(\delta)}$ by
\begin{equation*}
    \<1>_\eps^{(\delta)} (t) := (\I_\eps \xi^{(\delta)})(t) = \int_{-\infty}^{t} e^{(t-r) (\lL_\eps -1)} \xi^{(\delta)}(r) {\rm d}r\;.
\end{equation*}
We define two other processes $\<2>_\eps^{(\delta)}$ and $\<3>_\eps^{(\delta)}$ by
\begin{equation*}
    \<2>_\eps^{(\delta)} := (\<1>_\eps^{(\delta)})^2 - C_{\eps,1}^{(\delta)}\;, \qquad \<3>_\eps := (\<1>_\eps^{(\delta)})^3 - 3\, C_{\eps,1}^{(\delta)}\,\cdot \, \<1>_\eps^{(\delta)}\;.
\end{equation*}
We write $\<20>_\eps^{(\delta)}$ and $\<30>_\eps^{(\delta)}$ for the objects
\begin{equation*}
    \begin{split}
    \<20>_\eps^{(\delta)} (t) := (\I_\eps \<2>_\eps^{(\delta)})(t) = \int_{-\infty}^{t} e^{(t-r)(\lL_\eps -1)} \<2>_\eps^{(\delta)} (r) {\rm d}r\;,\\
    \<30>_\eps^{(\delta)} (t) := (\I_\eps \<3>_\eps^{(\delta)})(t) = \int_{-\infty}^{t} e^{(t-r)(\lL_\eps -1)} \<3>_\eps^{(\delta)} (r) {\rm d}r\;.
    \end{split}
\end{equation*}
We further define three more processes $\<31p>_\eps^{(\delta)}$, $\<22p>_\eps^{(\delta)}$ and $\<32p>_\eps^{(\delta)}$ by
\begin{equation*}
    \begin{split}
    &\<31p>_\eps^{(\delta)} := \<30>_\eps^{(\delta)} \circ_\eps \<1>_\eps^{(\delta)}\;, \qquad \<22p>_\eps^{(\delta)} := \<20>_\eps^{(\delta)} \circ_\eps \<2>_\eps^{(\delta)} - C_{\eps,2}^{(\delta)}\;,\\
    &\<32p>_\eps^{(\delta)} := \<30>_\eps^{(\delta)} \circ_\eps \<2>_\eps^{(\delta)} - 3 \, C_{\eps,2}^{(\delta)} \cdot \<1>_\eps^{(\delta)}\;.
    \end{split}
\end{equation*}
Here, $C_{\eps,1}^{(\delta)}$ and $C_{\eps,2}^{(\delta)}$ are renormalisation functions with some flexibilities in choices. For each fixed choice of the pair $(C_{\eps,1}^{(\delta)}, C_{\eps,2}^{(\delta)})$, 
\begin{equation} \label{eq:phi4_renormalisation}
    C_\eps^{(\delta)} := 3 C_{\eps,1}^{(\delta)} - 9 C_{\eps,2}^{(\delta)}
\end{equation}
is the choice for the renormalisation function in \eqref{e:phi4_hom_full}. This is in complete analogy with the standard constant coefficient case. 

For every $\eps\in \N^{-1}$ and $\delta \in (0,1)$, define
\begin{equation} \label{e:phi4_stochastic_collection}
    \Upsilon_\eps^{(\delta)} := \big( \<1>_\eps^{(\delta)}\,, \; \Pi_0 \<2>_\eps^{(\delta)}\,, \; \lL_\eps^{-1} \Pi_0^\perp \<2>_\eps^{(\delta)}\,, \;\<30>_\eps^{(\delta)}\,, \; \<31p>_\eps^{(\delta)}\,, \;\<22p>_\eps^{(\delta)}\,, \; \<32p>_\eps^{(\delta)} \big)\;.
\end{equation}
As in the standard constant coefficient dynamical $\phi^4_3$ case, we expect that for each fixed $\eps\in \N^{-1}$, the components to $\Upsilon_\eps^{(\delta)}$ converge as $\delta \rightarrow 0$ in their corresponding regularity spaces. More precisely, for each $T>0$, define the space $\yY_{T}$ by
\begin{equation*}
    \begin{split}
    \yY_{T} := &L_{T}^{\infty} \cC_x^{-\frac{1}{2}-\kappa} \, \times \, L_{T}^{\infty} \, \times L_T^\infty \cC_x^{1-\kappa} \times \, (L_{T}^{\infty} \cC_{x}^{\frac{1}{2}-\kappa} \cap \cC_{T}^{\eta/2} \cC_{x}^{\eta})\\
    &\times L_T^\infty \cC_x^{-\kappa} \, \times \, L_T^\infty \cC_x^{-\kappa} \, \times \, L_T^\infty \cC_x^{-\frac{1}{2}-\kappa}\;.
    \end{split}
\end{equation*}
A difference between the collection \eqref{e:phi4_stochastic_collection} and the standard constant coefficient case is that we replace the object $\<2>_\eps^{(\delta)}$ by the pair $(\Pi_0 \<2>_\eps^{(\delta)}, \lL_\eps^{-1} \Pi_0^\perp \<2>_\eps^{(\delta)})$. This is because the process $\<2>_\eps^{(\delta)}$ is expected to live in $\cC_x^{-1-\kappa;\eps}$ uniformly in $\delta$, and we need to transform it into a space ($L_T^\infty \times L_T^\infty \cC_x^{1-\kappa}$ in this case) where one can compare with its limit $(\Pi_0 \<2>_0, \lL_0^{-1} \Pi_0^\perp \<2>_0)$ directly. Also note that the original process $\<2>_\eps^{(\delta)}$ can be recovered from $(\Pi_0 \<2>_\eps^{(\delta)}, \lL_\eps^{-1} \Pi_0^\perp \<2>_\eps^{(\delta)})$ by
\begin{equation*}
    \<2>_\eps^{(\delta)} = \Pi_0 \<2>_\eps^{(\delta)} + \lL_\eps \lL_\eps^{-1} \Pi_0^\perp \<2>_\eps^{(\delta)}\;.
\end{equation*}
We also define the flux related stochastic objects $\fF_\eps^{(\delta)}$ by
\begin{equation*}
    \fF_\eps^{(\delta)} := \big( A_\eps \nabla \<1>_\eps^{(\delta)}, \, A_\eps \nabla \<20>_\eps^{(\delta)}, \, A_\eps \nabla \<30>_\eps^{(\delta)} \big)\;,
\end{equation*}
and the space $\yY_T^{\flux}$ by
\begin{equation*}
    \yY_T^\flux := \cC^{ - \kappa }_{ T } \cC^{ - \frac{ 3 }{ 2 } - \kappa }_{ x } \times L_T^\infty \cC_x^{-\kappa} \times L_T^\infty \cC_x^{-\frac{1}{2}-\kappa}\;.
\end{equation*}
The following proposition \cite{homo_stochastic} states the convergence of $\Upsilon_\eps^{(\delta)}$ to the corresponding stochastic objects in standard dynamical $\phi^4_3$ with operator $\lL_0$. 

\begin{prop} [\cite{homo_stochastic}, Theorem~3.1]
\label{prop:phi43_stochastic}
Let
\begin{equation*}
    \Upsilon_0 := \big( \<1>_0\;, \Pi_0 \<2>_0\;, \lL_0^{-1} \Pi_0^\perp \<2>_0\;, \<30>_0\;, \<31p>_0\;, \<22p>_0\;, \<32p>_0 \big)
\end{equation*}
be the collection of stochastic objects from the standard dynamical $\phi^4_3$ model with heat kernel $(\d_t - \lL_0)^{-1}$. Then for the choices
\begin{equation*}
    C_{\eps,1}^{(\delta)} = \E \big[ \<1>_\eps^{(\delta)} \big]^2\;, \qquad C_{\eps,2}^{(\delta)} = \E \big[ \<20>_\eps^{(\delta)} \cdot \<2>_\eps^{(\delta)} \big]\;,
\end{equation*}
for every $T>0$, the following convergences hold as $(\eps, \delta) \rightarrow (0,0)$:
\begin{enumerate}
\item $\Upsilon_\eps^{(\delta)} \rightarrow \Upsilon_0$ in $\yY_T$; 
\item $\fF_\eps^{(\delta)} \rightarrow \fF_0 := \big(\bar{A} \nabla \<1>_0, \, \bar{A} \nabla \<20>_0, \, \bar{A} \nabla \<30>_0 \big)$ in $\yY_T^{\flux}$, where $\<20>_0 = \I_0 (\<2>_0)$. 
\end{enumerate}
\end{prop}

% We make the following assumption on the stochastic objects $\Upsilon_\eps^{(\delta)}$ and $\fF_\eps^{(\delta)}$. 

% \begin{asp} \label{asp:phi4_stochastic}
% There exists a subset $\bB \subset [0,1] \times [0,1]$ containing $(0,0)$ as an accumulation point such that the followings are true for every realisation of the noise $\xi$. 
% \begin{enumerate}
%     \item There exists
%     \begin{equation*}
%     \Upsilon_0 := \big( \<1>_0\;, \Pi_0 \<2>_0\;, \lL_0^{-1} \Pi_0^\perp \<2>_0\;, \<30>_0\;, \<31p>_0\;, \<22p>_0\;, \<32p>_0 \big)
% \end{equation*}
%     such that for every $T>0$, $\Upsilon_\eps^{(\delta)} \rightarrow \Upsilon_0$ in $\yY_T$ as $(\eps,\delta) \rightarrow (0,0)$ when restricted to $\bB$. 

%     \item For every $T>0$, we have the convergence
%     \begin{equation*}
%     \fF_\eps^{(\delta)} \rightarrow \fF_0 := \big(\bar{A} \nabla \<1>_0, \, \bar{A} \nabla \<20>_0, \, \bar{A} \nabla \<30>_0 \big)
%     \end{equation*}
%     in $\yY_T^\flux$ as $(\eps,\delta) \rightarrow (0,0)$ when restricted to $\bB$, where $\<1>_0$ and $\<30>_0$ are the same objects from $\Upsilon_0$ above, and $\<20>_0 = \I_0 (\<2>_0)$ for the same $\<2>_0$ from above. 
% \end{enumerate}
% \end{asp}

We have the following theorem, whose proof will be provided at the end of this section. 

\begin{thm} \label{th:phi4_overall}
Let $C_\eps^{(\delta)}$ be such that Proposition~\ref{prop:phi43_stochastic} holds. Let $\phi_\eps(0) \in \cC^{-\frac{1}{2}-\eta}$ be such that
\begin{equation*}
    \|\phi_\eps(0) - \phi_0(0)\|_{-\sigma} \rightarrow 0
\end{equation*}
for some $\phi_0(0) \in \cC^{-\frac{1+\eta}{2}}$ as $\eps \rightarrow 0$. Let $\phi_\eps^{(\delta)}$ be the solution to \eqref{e:phi4_hom_full} with initial data $\phi_\eps(0)$ and enhanced stochastic objects $\Upsilon_\eps^{(\delta)}$, and $\phi_0$ be the solution to \eqref{e:phi4_limiting_equation} with initial data $\phi_0(0)$ and enhanced stochastic objects $\Upsilon_0$. Then for every $T>0$, $\phi_\eps^{(\delta)}$ converges in $L_T^\infty \cC_x^{-\sigma}$ to $\phi_0$ with the bound
\begin{equation*}
    \|\phi_\eps^{(\delta)} - \phi_0\|_{L_T^\infty \cC_x^{-\sigma}} \lesssim \eps^\theta + \|\phi_\eps(0) - \phi_0(0)\|_{-\sigma} + \|\Upsilon_\eps^{(\delta)} - \Upsilon_0\|_{\yY_T}
\end{equation*}
for some $\theta>0$. Furthermore, the flux also converges with the quantitative bound
\begin{equation*}
    \begin{split}
    &\phantom{11}\|A_\eps \nabla \phi_\eps^{(\delta)} - \bar{A} \nabla \phi_0\|_{\cC^{ - \kappa }_{ T } \cC^{ - \sigma - 1  - \kappa }_{ x }}\\
    &\lesssim \eps^\theta + \|\phi_\eps(0) - \phi_0(0)\|_{-\sigma} + \|\Upsilon_\eps^{(\delta)} - \Upsilon_0\|_{\yY_T} + \|\fF_\eps^{(\delta)} - \fF_0 \|_{\yY_{T}^{\flux}}
    \end{split}
\end{equation*}
for some $\theta>0$. 
\end{thm}

The above theorem is completely deterministic. For any choice of $C_{\eps,1}^{(\delta)}$ and $C_{\eps,2}^{(\delta)}$ such that the enhanced object $\Upsilon_\eps^{(\delta)}$ converges along a particular subsequence, then the solution also converges along that subsequence to a limiting solution, which solves the standard dynamical $\phi^4_3$ with constant coefficient $\bar{A}$ and enhanced stochastic object $\Upsilon_0$. In particular, it allows the situation that $\phi_\eps^{(\delta)}$ may converge to different limits along different relations $(\eps,\delta) \rightarrow (0,0)$. 

For readers' convenience, we summarise the stochastic objects in $\Upsilon_\eps^{(\delta)}$ and their regularities in the following table. 
\begin{center}
\begin{tabular}{p{1.5cm} p{1.2cm} p{1.2cm} p{1.8cm} p{2.4cm} p{1.2cm} p{1.2cm} p{1.2cm} p{1.2cm} }
\hline
Process & $\<1>_{\eps}^{(\delta)}$ & $\Pi_0 \<2>_{\eps}^{(\delta)}$ & $\lL_\eps^{-1} \Pi_0^\perp \<2>_\eps^{(\delta)}$ & $\<30>_{\eps}^{(\delta)}$ & $\<31p>_{\eps}^{(\delta)}$ & $\<22p>_{\eps}^{(\delta)}$ & $\<32p>_{\eps}^{(\delta)}$\\
\hline
Limit & $\<1>_0$ &$\Pi_0 \<2>_0$ & $\lL_0^{-1} \Pi_0^\perp \<2>_0$ & $\<30>_0$ & $\<31p>_0$&$\<22p>_0$&$\<32p>_0$\\
\hline
Reg. & $-\frac{1}{2}- $ & $L_T^\infty$ & $1-$ & $L_T^\infty \cC_x^{\frac{1}{2}-} \cap \cC_{\s}^{\eta}$ & $0-$&$0-$ &$-\frac{1}{2}-$\\
\hline
\end{tabular}
\end{center}
Here, $\cC_x^{\alpha-}$ denotes $\cC_x^{\alpha-\kappa}$, and a single ``$\alpha-$" denotes $L_T^\infty \cC_x^{\alpha-\kappa}$. We have the following theorem regarding the convergence of the stochastic objects.

% We use the space-time norm $\cC_{\mathfrak{s}}^{\alpha}$ in \cite{rs_theory} to measure $A_{\eps} \nabla \<1>_{\eps}$.
% Define the scaling index
%  \begin{equation*}
%      \mathfrak{s} := ( 2, 1, 1, 1)\;.
%  \end{equation*}
% Define a ball $B_{\mathfrak{s},(t,x)}$ in the sense of parabolic distance. More precisely,
%  \begin{equation*}
%       B_{\mathfrak{s},(t,x)}:= \{(\tau,y)\in \R\times \R^3; \sqrt{|\tau-t|}+|x-y|\leq 1\}\;.
%  \end{equation*}
% And for $r\in \N^{+}$,
% \begin{equation*}
%     \mathfrak{B}^{r}_{\mathfrak{s},(t,x)}:= \{g\in C_{0}^{r}(B_{\mathfrak{s},(t,x)})\;;\, \|g\|_{C^{r}}\leq 1\,\}\;.
% \end{equation*}
% As given in \cite[Definition~3.7]{rs_theory}, for $\alpha<0$, let $r = - {\lfloor{\alpha}}\rfloor$ and $\mathfrak{K}=[-1,1]\times [-\frac{1}{2},\frac{1}{2})^3$ and define for $f$, a space-time distribution defined on $[-2,2]\times \T^{3}$,
% \begin{equation*}
%     \|f\|_{\cC^{\alpha}_{\mathfrak{s}}([-1,1]\times \T^d)}:= \sup_{(t,x)\in \mathfrak{K}}\sup_{g\in \mathfrak{B}^{r}_{\mathfrak{s},(t,x)}}\sup_{\delta \in ( 0, 1/2)}\delta^{-\alpha}I_{(t,x)}(f,g;\delta)\;,
% \end{equation*}
% where
% \begin{equation*}
%     I_{(t,x)}(f,g;\delta):=\delta^{-5}\int_{-1}^{1}{\rm d}\tau \int_{[-\frac{1}{2},\frac{1}{2})^{3}}f(\tau-t,y-x)g(\delta^{-2}\tau,\delta^{-1}y){\rm d}y\;.
% \end{equation*}

\subsection{Derivation of the fixed point system}

In this subsection, we derive the fixed point system that encodes the decomposition of the solution $\phi_\eps^{(\delta)}$ to \eqref{e:phi4_hom_full} and its continuous dependence on the enhanced stochastic objects $\Upsilon_\eps^{(\delta)}$. 

Since the dependence of the solution $\phi_\eps^{(\delta)}$ on $\delta$ is via the stochastic object $\Upsilon_\eps^{(\delta)}$ only, we drop $\delta$ from the notations for simplicity. More precisely, with an abuse of notation, we write $\Upsilon_\eps := \Upsilon_\eps^{(\delta)}$ and $\phi_\eps := \phi_\eps^{(\delta)}$. We emphasise that these are not objects at the $\delta=0$ limit, but just with $\delta$ omitted in the notation. The choices $C_{\eps,j} = C_{\eps,j}^{(\delta)}$ for $j=1,2$ are encoded in the stochastic objects $\Upsilon_\eps = \Upsilon_\eps^{(\delta)}$.  

Recall that $\phi_\eps$ solves \eqref{e:phi4_hom_full} and the collection of stochastic objects $\Upsilon_\eps$ from \eqref{e:phi4_stochastic_collection}. Let
\begin{equation} \label{e:phi4_remainder_1}
    u_\eps := \phi_\eps - \<1>_\eps + \<30>_{\eps}\;, \qquad  \Lambda_\eps u := u - \<30>_\eps\;.
\end{equation}
Then $u_\eps$ satisfies the equation
\begin{equation*}
    \d_t u_\eps = (\lL_\eps - 1) u_\eps -  3 \Lambda_\eps u_\eps \,\cdot \, \<2>_\eps - 3  (\Lambda_\eps u_\eps)^2 \,\cdot\, \<1>_\eps - (\Lambda_{\eps} u_\eps)^3 - 9 C_{\eps,2} \,\phi_\eps\;.
\end{equation*}
This equation is still not ideal since the products $\Lambda_{\eps} u_\eps \cdot \<2>_\eps$ and $(\Lambda_{\eps}u_\eps)^2 \<1>_\eps$ are not uniformly bounded in $(\eps,\delta)$ given the regularities of $\<2>_\eps$ and $\<1>_\eps$. To be precise, the problematic parts of the above two products are
\begin{equation*}
    \Lambda_{\eps} u_\eps \circ_\eps \<2>_\eps \quad \text{and} \quad (\Lambda_{\eps}u_\eps \prec_\eps \Lambda_{\eps}u_\eps) \circ_\eps \<1>_\eps
\end{equation*}
respectively. We make the ansatz
\begin{equation*}
    u_\eps = e^{t\lL_\eps} u_\eps(0) -3 \iI_\eps \big( \Lambda_{\eps}u_\eps \prec_\eps \<2>_\eps \big) + u^{\#}_\eps\;,
\end{equation*}
and plug it back into the a priori ill-defined para-products $\Lambda_{\eps}u_\eps \circ_\eps \<2>_\eps$ and $(\Lambda_{\eps}u_\eps \prec_\eps \Lambda_{\eps} u_\eps) \circ_\eps \<1>_\eps$. We have the system of equations for $(u_\eps, u_\eps^{\#})$ as
\begin{equation} \label{e:phi4_system}
    \begin{split}
    u_\eps &= e^{t(\lL_\eps -1 )} u_\eps(0) - 3 \iI_\eps \big( \Lambda_\eps u_\eps \prec_\eps \<2>_\eps \big) - \iI_\eps \big( \kK_\eps (u_\eps, u_\eps^{\#}) \big)\;,\\
    u_\eps^{\#} &= e^{t (\lL_\eps -1)} u_\eps^{\#}(0) - \iI_\eps \big( \kK_\eps (u_\eps, u_\eps^{\#}) \big)\;,
    \end{split}
\end{equation}
where
\begin{equation} \label{e:phi4_K}
    \kK_\eps (u, u^{\#}) := (\Lambda_\eps u)^3 + 3 (\Lambda_\eps u) \succ_\eps \<2>_\eps + 3 \big( \sS_\eps^{(1)}(u) + \sS_\eps^{(2)}(u, u^{\#}) \big)\;,
\end{equation}
and the functionals $\sS_\eps^{(1)}$ and $\sS_\eps^{(2)}$ are respectively given by
\begin{equation} \label{e:phi4_S}
    \begin{split}
    \sS_\eps^{(1)}(u) = &\big( (\Lambda_\eps u) \circ_\eps (\Lambda_\eps u) \big) \<1>_\eps + 2 \big( \Lambda_\eps u \prec_\eps \Lambda_\eps u \big) \prec_\eps \<1>_\eps + 2 \big( \Lambda_\eps u \prec_\eps \Lambda_\eps u \big) \succ_\eps \<1>_\eps\\
    &+ 2 \, \Com_\eps \big( \Lambda_\eps u; \Lambda_\eps u; \<1>_\eps \big) +  2 (\Lambda_\eps u) \cdot \big( u \circ_\eps \<1>_\eps - \<31p>_\eps \big)\;,\\
    \sS_\eps^{(2)}(u,u^{\#}) = &\Big( e^{t (\lL_\eps-1)} \big(u_\eps(0) - u_\eps^{\#}(0) \big) + u^{\#} - 3 \,[ \, \iI_\eps, (\Lambda_\eps u) \prec_\eps \, ] \, \<2>_\eps \Big) \circ_\eps \<2>_\eps - \<32p>_\eps\\
    &- 3 \, \Com_\eps \big( \Lambda_\eps u; \iI_\eps(\<2>_\eps); \<2>_\eps \big) - 3 (\Lambda_\eps u) \, \Big( \<22p>_\eps - \big( e^{t(\lL_\eps - 1)} \<20>_\eps(0) \big) \circ_\eps \<2>_\eps \Big)\;.
    \end{split}
\end{equation}
They respectively represent the a priori ill-defined quantities
\begin{equation*}
    \sS_\eps^{(1)}(u_\eps) ``=" \Lambda_\eps^2 (u_\eps) \cdot \<1>_\eps\;, \quad \sS_\eps^{(2)}(u_\eps, u_\eps^\#) ``=" u_\eps \circ_\eps \<2>_\eps + 3 C_{\eps,2} \Lambda_\eps (u_\eps) - \<32p>_\eps\;.
\end{equation*}
Note that the heat kernel operator here corresponds to $\lL_\eps - 1$ instead of $\lL_\eps$. This guarantees that the stochastic objects are all stationary in time (since this allows integration in time from $-\infty$). One can check that if $\xi^{(\delta)}$ is smooth, then $\phi_\eps$ solves \eqref{e:phi4_hom_full} if and only if there exists $u_\eps^{\#}$ such that for $u_\eps$ given by \eqref{e:phi4_remainder_1}, the pair $(u_\eps, u_\eps^{\#})$ solves the system \eqref{e:phi4_system}. 

\begin{rmk} \label{rm:kernel_damped}
    Note that the damped heat kernel $\iI_\eps$ or its fixed time version $e^{t(\lL_\eps - 1)}$ still satisfy all the estimates in Section~\ref{sec:para_general}. The bounds for fixed time follow immediately since the damped kernel is the undamped one multiplied by a constant less than one. 

    As for continuity-in-time estimates, we have
    \begin{equation*}
        e^{t(\lL_\eps - 1)} - e^{s(\lL_\eps - 1)} = e^{-t} \big( e^{t \lL_\eps} - e^{s \lL_\eps} \big) + (e^{-t} - e^{-s}) \, e^{s \lL_\eps}\;.
    \end{equation*}
    The former has the same form as usual with an additional damping, and hence satisfies the same estimates. For the latter, continuity-in-time can be obtained from
    \begin{equation*}
        |e^{-t} - e^{-s}| \lesssim |t-s|, 
    \end{equation*}
    and the operator $e^{s \lL_\eps}$ has at least the same smoothing property as $e^{t \lL_\eps} - e^{s \lL_\eps}$. 
\end{rmk}

Hence, throughout this section, we will use the bounds in Section~\ref{sec:para_general} for heat kernels based on $\lL_\eps - 1$. For every $T>0$ and $\eps\in \N^{-1}$, define the spaces $\xX_{T}^{(1)}$ and $\xX_{T,\eps}^{(2)}$ by
\begin{equation} \label{e:phi4_space_components}
    \xX_{T}^{(1)} = L_{T}^\infty \cC_x^{-\sigma} \cap \cC_{2 (\sigma + 3\eta), T}^{\eta} \cC_{x}^{\sigma+2\eta}\;, \quad \xX_{T,\eps}^{(2)} = L_{T}^\infty \cC_x^{-\sigma} \cap \cC_{3(\sigma+\eta), T}^{\eta/2} \cC_{x}^{2\sigma; \, \eps}\;,
\end{equation}
and set
\begin{equation*}
    \xX_{T,\eps} := \xX_{T}^{(1)} \times \xX_{T,\eps}^{(2)}\;.
\end{equation*}
The space $\xX_{T,\eps}^{(1)}$ does not depend on $\eps$ since all regularity exponents are in the interval $(-1,1)$, and is in the range of equivalences as in Lemma~\ref{lem:Holder_equivalence}. We have the following bounds on the functional $\Lambda_\eps$ given in \eqref{e:phi4_remainder_1}.

\begin{lem} \label{le:phi4_lambda}
We have
\begin{equation*}
    \begin{split}
    \|\Lambda_\eps u_\eps\|_{L_{2(\sigma+3\eta),T}^{\infty} \cC_x^{\frac{1}{2}-\kappa}} &\lesssim \|u_\eps\|_{\xX_T^{(1)}} + \|\<30>_\eps\|_{L_T^\infty \cC_x^{\frac{1}{2}-\kappa}}\;,\\
    \|\Lambda_\eps u_\eps\|_{\cC_{\sigma + 3\eta, T}^{\eta/2} \cC_x^\eta} &\lesssim \|u_\eps\|_{\xX_T^{(1)}} + \|\<30>_\eps\|_{\cC_T^{\eta/2} \cC_x^\eta}\;,
    \end{split}
\end{equation*}
and
\begin{equation*}
    \begin{split}
    \|\Lambda_\eps u_\eps - \Lambda_0 u_0\|_{L_{2(\sigma+3\eta),T}^{\infty} \cC_x^{\frac{1}{2}-\kappa}} &\lesssim \|u_\eps - u_0\|_{\xX_T^{(1)}} + \|\<30>_\eps - \<30>_0\|_{L_T^\infty \cC_x^{\frac{1}{2}-\kappa}}\;,\\
    \|\Lambda_\eps u_\eps - \Lambda_0 u_0\|_{\cC_{\sigma + 3\eta, T}^{\eta/2} \cC_x^\eta} &\lesssim \|u_\eps - u_0\|_{\xX_T^{(1)}} + \|\<30>_\eps - \<30>_0\|_{\cC_T^{\eta/2} \cC_x^\eta}\;,
    \end{split}
\end{equation*}
\end{lem}
\begin{proof}
We consider the two bounds for $\Lambda_\eps u_\eps$. The $L_{2(\sigma+3\eta),T}^{\infty} \cC_x^{\frac{1}{2}-\kappa}$ bound follows directly from the definition of $\Lambda_\eps$. The $\cC_{\sigma+3\eta,T}^{\eta/2} \cC_x^\eta$ bound follows from the interpolation
\begin{equation*}
    \|u\|_{\cC_{\sigma+3\eta,T}^{\eta/2} \cC_x^\eta} \lesssim \|u\|_{L_T^\infty \cC_x^{-\sigma}}^{\frac{1}{2}} \|u\|_{\cC_{2(\sigma+3\eta)}^{\eta} \cC_x^{\sigma+2\eta}}^{\frac{1}{2}} \lesssim \|u\|_{\xX_T^{(1)}}\;.
\end{equation*}
The bounds for the difference $\Lambda_\eps u_\eps - \Lambda_0 u_0$ then follows from linearity. 
\end{proof}

\subsection{Uniform boundedness}

Recall $\eta \in (0,\frac{1}{1000})$ is fixed, and $\sigma := \frac{1}{2} + \eta$. Let
\begin{equation} \label{e:phi4_M_eps}
    \mM_\eps :=  \|u_\eps(0)\|_{\cC^{-\frac{1}{2}-\eta}} + \|u_\eps^{\#}(0)\|_{\cC^{-\frac{1}{2}-\eta}} + \|\Upsilon_\eps\|_{\yY_{T}}\;.
\end{equation}
We have the following theorem regarding the existence of the solution $(u_\eps, u_\eps^\#)$ to \eqref{e:phi4_system}. 

\begin{thm} \label{thm:phi4_fixed_pt}
    Fix $T>0$. For every $\eps\in \N^{-1}$, every realisation $\Upsilon_\eps \in \yY_{T,\eps}$ and every pair of initial data $\vec{u}_\eps(0) := \big(u_\eps(0), u_\eps^{\#}(0) \big) \in \cC_{x}^{-\sigma} \times \cC_x^{-\sigma}$, there exists $\tau_\eps \in (0,T)$ such that the system \eqref{e:phi4_system} has a unique solution $(u_\eps, u_\eps^{\#}) \in \xX_{T,\eps}$. 

    Furthermore, the local existence time $\tau_\eps$ and the norm $\|\vec{u}_\eps\|_{\xX_{\tau_\eps,\eps}}$ depends on $(\eps, \delta)$ via $\mM_\eps$ only (where $\delta$ is hidden in the notation). In particular, if $\mM_\eps$ is uniformly bounded in $\eps\in \N^{-1}$, then there exists $\tau>0$ independent of $\eps$ such that
    \begin{equation} \label{e:phi4_fixed_pt_bound}
        \|\vec{u}_\eps\|_{\xX_{\tau,\eps}} \lesssim 1 + \mM_\eps\;,
    \end{equation}
    where the proportionality constant is independent of $\eps$. 
\end{thm}
\begin{proof}
For $\tau \in (0,T)$ and $\eps\in \N^{-1}$, define the mild solution map $\Gamma_{\tau,\eps} = (\Gamma_{\tau,\eps}^{(1)}, \Gamma_{\tau,\eps}^{(2)})$ by
\begin{equation*}
    \begin{split}
    \Gamma_{\tau,\eps}^{(1)}(\vec{u}) &= e^{t(\lL_\eps-1)} u_\eps(0) - 3 \iI_\eps \big( \Lambda_\eps u \prec_\eps \<2>_\eps \big) - \iI_\eps \big( \kK_\eps (\vec{u}) \big)\;,\\
    \Gamma_{\tau,\eps}^{(2)}(\vec{u}) &= e^{t(\lL_\eps -1)} u_\eps^{\#}(0) - \iI_\eps \big( \kK_\eps (\vec{u}) \big)\;,
    \end{split}
\end{equation*}
where we write $\vec{u} = (u, u^{\#})$, and $\kK_\eps$ is defined in \eqref{e:phi4_K}, with some terms from $\sS_\eps^{(1)}$ and $\sS_\eps^{(2)}$ given in \eqref{e:phi4_S}. We first show that for sufficiently small $\tau>0$, $\Gamma_{\tau, \eps}$ is a contraction map from a bounded ball of $\xX_{\tau, \eps}$ into itself.

\begin{flushleft}
	\textit{Step 1.}
\end{flushleft}

We first check that for some $\rR$ and $\tau$ to be specified later, $\Gamma_{\tau,\eps}$ maps the ball in $\xX_{\tau, \eps}$ of radius $\rR$ centered at the origin into itself. For the initial data terms, by Lemma~\ref{le:heat_smoothing_space}, we have
\begin{equation} \label{e:phi4_fixed_pt_initial}
    \|e^{t (\lL_\eps-1)} u_\eps(0)\|_{\xX_\tau^{(1)}} \lesssim \|u_\eps(0)\|_{\cC^{-\sigma}}\;, \quad \|e^{t (\lL_\eps-1)} u_\eps^{\#}(0)\|_{\xX_{\tau,\eps}^{(2)}} \lesssim \|u_\eps^{\#}(0)\|_{\cC^{-\sigma}}\;.
\end{equation}
This gives the desired bounds for the two terms from the initial data. Also, by Proposition~\ref{pr:paraproducts_gBesov} and Lemma~\ref{le:phi4_lambda}, we have
\begin{equation*}
    \|(\Lambda_\eps u)(t) \prec_\eps \<2>_\eps(t) \|_{-1-\kappa;\eps} \lesssim \|(\Lambda_\eps u)(t)\|_{\cC_x^\eta} \cdot \|\<2>_\eps\|_{-1-\kappa;\eps} \lesssim \big( t^{-\frac{\sigma+3\eta}{2}} \|u\|_{\xX_{\tau}^{(1)}} + \mM_\eps \big) \mM_\eps
\end{equation*}
uniformly over $t \leq \tau$. It then follows from Lemma~\ref{le:heat_smoothing_spacetime} that there exists $\theta>0$ such that
\begin{equation} \label{e:phi4_fixed_pt_low_reg}
    \left\| \iI_\eps \big( \Lambda_\eps u \prec_\eps \<1>_\eps \big) \right\|_{\xX_{\tau}^{(1)}} \lesssim \tau^{\theta} \mM_\eps \big( \|u\|_{\xX_{\tau}^{(1)}} + \mM_\eps \big)\;.
\end{equation}
We now turn to the term $\kK_\eps$. If we can show that
\begin{equation} \label{e:phi4_fixed_pt_remainder_0}
    \| (\kK_\eps \vec{u})(r) \|_{-\sigma} \lesssim r^{-\frac{3}{2}(\sigma + 3\eta)} \big(1 + \mM_\eps + \|\vec{u}\|_{\xX_{\tau,\eps}} \big)^3
\end{equation}
uniformly over $r \leq \tau$, then it will again follow from Lemma~\ref{le:heat_smoothing_spacetime} that
\begin{equation} \label{e:phi4_fixed_pt_remainder}
    \left\| \iI_\eps (\kK_\eps \vec{u}) \right\|_{\xX_\tau^{(1)}} + \left\| \iI_\eps (\kK_\eps \vec{u}) \right\|_{\xX_{\tau, \eps}^{(2)}} \lesssim \tau^{\theta} \big( 1 + \mM_\eps + \|\vec{u}\|_{\xX_{\tau,\eps}} \big)^3\;.
\end{equation}
Combining \eqref{e:phi4_fixed_pt_initial}, \eqref{e:phi4_fixed_pt_low_reg} and \eqref{e:phi4_fixed_pt_remainder}, we see there exists $C_1, C_2>0$ such that for all $\vec{u}$ with $\|\vec{u}\|_{\xX_{\tau,\eps}} \leq \rR$, we have
\begin{equation} \label{e:phi4_fixed_pt_map_bound}
    \|\Gamma_{\tau, \eps} \vec{u}\|_{\xX_{\tau,\eps}} \leq C_1 \mM_\eps + C_2 \tau^{\theta} \big(1 + \mM_\eps + \rR \big)^3\;.
\end{equation}
Hence, if we take $\rR = 2 (1 + C_1 \mM_\eps)$ and $\tau>0$ be such that
\begin{equation} \label{e:phi4_fixed_pt_map_smalltime}
    C_2 \tau^\theta (1 + \mM_\eps + \rR)^3 < \frac{\rR}{4}\;,
\end{equation}
we see $\Gamma_{\tau,\eps}$ maps the ball of radius $\rR$ in $\xX_{\tau,\eps}$ into itself. Furthermore, we notice that the choice of $\rR$ and $\tau$ depends on $\eps$ via $\mM_\eps$ only. 

We now turn back to proving the bound \eqref{e:phi4_fixed_pt_remainder_0}. There are various terms from $\kK_\eps$; we give for example the term with $u^{\#}$, the two commutator terms $\Com_\eps$, and the commutator term involving $[\iI_\eps, \Lambda_\eps u \prec_\eps] \<2>_\eps$. 

For the ``high regularity" term $u^{\#}$, by Proposition~\ref{pr:paraproducts_gBesov}, we have
\begin{equation*}
    \|u^{\#}(r) \circ_\eps \<2>_\eps(r)\|_{\eta-\kappa;\eps} \lesssim \|u^{\#}(r)\|_{1+\eta;\eps} \|\<2>_\eps(r)\|_{-1-\kappa;\eps} \lesssim r^{-\frac{3}{2}(\sigma+\eta)} \|u^{\#}\|_{\xX_{\tau,\eps}^{(2)}} \mM_\eps\;,
\end{equation*}
and the bound follows from equivalence of the norms $\|\cdot\|_{\eta-\kappa;\eps} \sim \|\cdot\|_{\eta-\kappa}$ since $\eta - \kappa \in (0,1)$. For the two commutator terms $\Com_\eps$, by Corollary~\ref{cor:com_prec_circ}, we have
\begin{equation*}
    \begin{split}
    &\phantom{111}\left\| \Com_\eps\big( (\Lambda_\eps u)(r); (\Lambda_\eps u)(r); \<1>_\eps \big) \right\|_{\eta-2\kappa}\\
    &\lesssim \|(\Lambda_\eps u)(r)\|_{\frac{1}{2}-\kappa} \|(\Lambda_\eps u)(r)\|_{\eta} \|\<1>_\eps(r)\|_{-\frac{1}{2}-\kappa} \lesssim r^{-\frac{3}{2}(\sigma+3\eta)} (\mM_\eps + \rR)^3\;,
    \end{split}
\end{equation*}
and
\begin{equation*}
    \begin{split}
    &\phantom{111}\left\| \Com_\eps \big( (\Lambda_\eps u)(r); (\iI_\eps \<2>_\eps)(r); \<2>_\eps(r) \big) \right\|_{\eta-3\kappa}\\
    &\lesssim \|(\Lambda_\eps u)(r)\|_{\eta} \|(\iI_\eps \<2>_\eps)(r)\|_{1-2\kappa} \|\<2>_\eps(r)\|_{-1-\kappa;\eps} \lesssim r^{-\frac{\sigma+3\eta}{2}} (\mM_\eps + \rR)^3\;.
    \end{split}
\end{equation*}
Finally, control of the commutator term $[\iI_\eps, (\Lambda_\eps u) \prec_\eps] \<2>_\eps$ requires time continuity of $\Lambda_\eps u$ and hence for $u$, which is precisely given in Lemma~\ref{le:phi4_lambda}. Hence, by Lemma~\ref{le:heat_comm_spacetime}, we get
\begin{equation*}
    \left\| \big( [\iI_\eps, \Lambda_\eps u \prec_\eps] \<2>_\eps \big)(r) \right\|_{1+\eta-2\kappa; \eps} \lesssim r^{-\frac{\sigma+3\eta}{2}} \|\Lambda_\eps u\|_{\cC_{\sigma+3\eta}^{\eta/2} \cC_x^\eta} \|\<2>_\eps(r)\|_{-1-\kappa;\eps} \lesssim r^{-\frac{\sigma+3\eta}{2}} \mM_\eps^2\;,
\end{equation*}
which gives the desired bound. Other terms in $\kK_\eps$ can be treated similarly. This verifies \eqref{e:phi4_fixed_pt_low_reg}. Hence, one concludes that $\Gamma_{\tau,\eps}$ maps a ball of radius $\rR = 2 (1 + C_1 \mM_\eps)$ into itself provided $\tau$ sufficiently small (depending on $\mM_\eps$ only).

\begin{flushleft}
\textit{Step 2. }
\end{flushleft}

We now show that this map $\Gamma_{\tau,\eps}$ is a contraction in $\xX_{\tau,\eps}$ for sufficiently small $\tau$. For this, we need to get a bound for $\| \Gamma_{\tau,\eps}(\vec{u}) - \Gamma_{\tau,\eps}(\vec{v})\|_{\xX_{\tau,\eps}}$ for $\vec{u} = (u, u^\#)$ and $\vec{v} = (v, v^{\#})$ in $\xX_{\tau,\eps}$. The argument is essentially the same as in Step 1. Most of the terms are either ``constants" (the initial data and pure stochastic terms) or linear in $u$ or $u^{\#}$, and hence for the difference, they either cancel out or one just replaces $\|\vec{u}\|_{\xX_{\tau,\eps}}$ by $\|\vec{u} - \vec{v}\|_{\xX_{\tau,\eps}}$. 

The nonlinear terms are either quadratic or cubic in $u$. In these two situations, one replaces one factor of $\rR$ by $\|\vec{u}-\vec{v}\|_{\xX_{\tau,\eps}}$ and the other factors are kept changed. Hence, we get the bound
\begin{equation} \label{e:phi4_fixed_pt_map_contraction}
    \|\Gamma_{\tau,\eps}(\vec{u}) - \Gamma_{\tau,\eps}(\vec{v})\|_{\xX_{\tau,\eps}} \lesssim \tau^{\theta} \|\vec{u} - \vec{v}\|_{\xX_{\tau,\eps}} (1 + \mM_\eps + \rR)^{2}
\end{equation}
for some $\theta>0$. This shows that $\Gamma_{\tau,\eps}$ is a contraction in $\xX_{\tau,\eps}$ for sufficiently small $\tau = \tau_\eps$, and hence \eqref{e:phi4_system} admits a local solution $\vec{u}_\eps \in \xX_{\tau_\eps, \eps}$. 

Furthermore, it is clear from \eqref{e:phi4_fixed_pt_map_smalltime} and \eqref{e:phi4_fixed_pt_map_contraction} and the choice $\rR = 2 (1 + C_1 \mM_\eps)$ that the local existence time $\tau = \tau_\eps$ can be made proportional to $(1 + \mM_\eps)^{-2}$. In particular, it is independent of $(\eps,\delta)$ if $\mM_\eps = \mM_\eps^{(\delta)}$ is uniformly bounded. The bound \eqref{e:phi4_fixed_pt_bound} then follows from \eqref{e:phi4_fixed_pt_map_bound} and the choice of $\tau$ in \eqref{e:phi4_fixed_pt_map_smalltime}. 
\end{proof}

\subsection{Convergences of the solution and flux}

For $\eps\in \N^{-1}$, let $\vec{u}_\eps = (u_\eps, u_\eps^{\#})$ be the solution to the system \eqref{e:phi4_system} in $\xX_{\tau_\eps, \eps}$ with initial data $\big(u_\eps(0), u_\eps^{\#}(0) \big)$. We need to compare the solution $\vec{u}_\eps$ and $\vec{u}_0$ in one fixed space. Note that $u_\eps^\# \in \cC^{\alpha;\eps}$ for some $\alpha \in (1,2)$ which is outside the regularity exponent for equivalence of spaces, so we need to transform it in a way that it falls into the regularity space between $(-1,1)$. For this purpose, we decompose
\begin{equation*}
    u_\eps^\# = \Pi_0 u_\eps^\# + \lL_\eps^{-1} \lL_\eps u_\eps^\#\;,
\end{equation*}
where $\lL_\eps u_\eps^\# \in \cC^{\alpha-2}$ with $\alpha-2 \in (-1,0)$, and can be compared with $\lL_0 u_0^\#$. 

Hence, we define
\begin{equation*}
    \vec{U}_{ \eps } := \big( u_\eps, u_\eps^{\#}, \lL_\eps u_\eps^{\#} \big)\;,
\end{equation*}
and
\begin{equation*}
    \widetilde{\xX}_T := \xX_{T}^{(1)} \times L_T^\infty \cC_x^{-\sigma} \times \cC_{3(\sigma+\eta)}^{\eta/2} \cC_x^{-1+2\eta}\;,
\end{equation*}
where the space $\xX_T^{(1)}$ is defined in \eqref{e:phi4_space_components}. Let $\widetilde{\xX}_T^{+}$ be the same as $\widetilde{\xX}_T$ except that one replaces the two $L_T^\infty \cC_x^{-\sigma}$ components (one of them is hidden in the definition of $\xX_T^{(1)}$) by the slightly stronger space $L_T^\infty \cC_x^{-\frac{1+\eta}{2}}$. 

For $\eps = 0$, the operator $\lL_0 = \div(\bar{A} \nabla)$ is a constant coefficient elliptic operator. It is well known in this case that if $\big(u_0(0), u_0^{\#}(0) \big) \in \cC_x^{-\frac{1+\eta}{2}} \times \cC_x^{-\frac{1+\eta}{2}}$, then the solution $\vec{U}_0 = (u_0, u_0^{\#}, \lL_0 u_0^{\#})$ exists in $\widetilde{\xX}_T^{+}$ for arbitrary $T>0$. We have the following theorem regarding the convergence of $\vec{ v }_{ \eps }$ to $\vec{v}_0$.

\begin{thm} \label{thm:phi4_convergence_homo}
Fix $T>0$. Suppose $\big(u_0(0), u_0^{\#}(0)\big) \in \cC^{-\frac{1+\eta}{2}} \times \cC^{-\frac{1+\eta}{2}}$, $\big( u_\eps(0), u_\eps^{\#}(0) \big) \in \cC^{-\sigma} \times \cC^{-\sigma}$, and
\begin{equation*}
    \mM := 1 + \|\vec{U}_0\|_{\widetilde{\xX}_T^{+}} + \sup_{\eps\in \N^{-1}} \big( \|u_\eps(0)\|_{-\sigma} + \|u_\eps^{\#}(0)\|_{-\sigma} + \|\Upsilon_\eps\|_{\yY_{T}} \big) < +\infty\;.
\end{equation*}
Then for every $\eps\in \N^{-1}$, the solution $(u_\eps, u_\eps^{\#})$ exists in $\xX_{T,\eps}$ up to the same time $T$. Furthermore, there exists $\theta>0$ such that we have the bound
\begin{equation} \label{e:phi4_convergence_sol}
    \|\vec{U}_{\eps} - \vec{U}_0\|_{\widetilde{\xX}_T} \lesssim \eps^\theta + \|u_\eps(0) - u_0(0)\|_{-\sigma} + \|u_\eps^{\#}(0) - u_0^{\#}(0)\|_{-\sigma} + \| \Upsilon_\eps - \Upsilon_0 \|_{ { \yY }_{ T }}\;.
\end{equation}
The proportionality constant above depends on $\mM$ and $T$ only. 
\end{thm}
\begin{proof}
We first prove \eqref{e:phi4_convergence_sol} for a sufficiently small positive time (independent of $\eps$), and then iterate it to the whole time interval $[0,T]$ where the limiting solution $\vec{v}_0$ is defined. 

By Theorem~\ref{thm:phi4_fixed_pt}, there exists a universal $C_* > 1$ and $\tau_0>0$ independent of $\eps\in \N^{-1}$ such that
\begin{equation*}
    \|\vec{U}_{\eps}\|_{\widetilde{\xX}_{\tau}} < C_* \mM\;.
\end{equation*}
The existence of such a $\tau_0$ is guaranteed by the definition of $\mM$ and the bound \eqref{e:phi4_fixed_pt_bound}. 

The system of equations for $\vec{U}_{\eps} = (u_\eps,  u_\eps^{\#}, \lL_\eps u_\eps^{\#})$ can be derived immediately from \eqref{e:phi4_system} with the additional equation for $\lL_\eps u_\eps^{\#}$ arising from applying $\lL_\eps$ to the equation for $u_\eps^{\#}$. 

We need to control all the three components of $\vec{U}_{\eps} - \vec{U}_0$, up to a time $\tau < \tau_0$ specified later. For the terms with the initial data, it follows directly from Lemma~\ref{le:convergence_heat_normal} that
\begin{equation} \label{e:phi4_convergence_initial_1}
    \begin{split}
    \big\| e^{t(\lL_\eps-1)} u_\eps(0) - e^{t (\lL_0-1)} u_0(0) \big\|_{\xX_{\tau}^{(1)}} &\lesssim \eps^{\theta} \|u_0(0)\|_{-\frac{1+\eta}{2}} + \|u_\eps(0) - u_0(0)\|_{-\sigma}\;,\\
    \big\| e^{t(\lL_\eps-1)} u_\eps^{\#}(0) - e^{t (\lL_0-1)} u_0^{\#}(0) \big\|_{L_\tau^\infty \cC_x^{-\sigma}} &\lesssim \eps^\theta \|u_0^{\#}(0)\|_{-\frac{1+\eta}{2}} + \|u_\eps^{\#}(0) - u_0^{\#}(0)\|_{-\sigma},
    \end{split}
\end{equation}
and
\begin{equation} \label{e:phi4_convergence_initial_2}
    \begin{split}
    &\phantom{111}\big\| \lL_\eps e^{t(\lL_\eps-1)} u_\eps(0) - \lL_0 e^{t(\lL_0-1)} u_0(0) \big\|_{\cC_{3\sigma+\eta,\tau}^{\eta/2} \cC_{x}^{-1+2\eta}}\\
    &\lesssim \eps^{\theta} \|u_0(0)\|_{-\frac{1+\eta}{2}} + \|u_\eps(0) - u_0(0)\|_{-\sigma}\;.
    \end{split}
\end{equation}
This is where we need the slightly stronger space $L_T^\infty \cC_x^{-\frac{1+\eta}{2}}$ for $u_0$ and $u_0^{\#}$. 

For the second term in the equation for $u_\eps$, by Lemma~\ref{le:convergence_prec_heat_below}, we have the desired bound
\begin{equation} \label{e:phi4_convergence_rough}
    \begin{split}
    &\phantom{111}\big\| \iI_\eps (\Lambda_\eps u_\eps \prec_\eps \<2>_\eps) - \iI_0 (\Lambda_0 u_0 \prec \<2>_0) \big\|_{\xX_{\tau}^{(1)}}\\
    &\lesssim \tau^{\theta} \Big( \|\Lambda_\eps u_\eps - \Lambda_0 u_0\|_{L_{\sigma+3\eta,\tau}^{\infty} \cC_x^\eta} \|\Pi_0^\perp \<2>_\eps\|_{L_\tau^\infty \cC_{x}^{-1-\kappa;\eps}}\\
    &\phantom{111}+ \big( \eps^\theta \|\Pi_0^\perp \<2>_0\|_{L_\tau^\infty \cC_x^{-1-\kappa}} + \|\lL_{\eps}^{-1}\Pi_{0}^{\perp} \<2>_\eps-\lL_{0}^{-1}\Pi_{0}^{\perp} \<2>_0 \|_{L_\tau^\infty \cC_x^{1-\kappa}} \big) \, \|\Lambda_0 u_0\|_{L_{\sigma+3\eta,\tau}^\infty \cC_x^\eta} \Big)\\
    &\lesssim \tau^{\theta} \mM \Big( \eps^\theta \mM + \|u_\eps - u_0\|_{\xX_\tau^{(1)}} +\| \Upsilon_\eps- \Upsilon_0 \|_{{\yY}_{\tau}} \Big)\;.
     \end{split}
\end{equation}
As for the terms with $\kK_\eps$ and $\kK_0$, by Lemma~\ref{le:convergence_heat_spacetime_normal}, we have
\begin{equation} \label{e:phi4_convergence_IK}
    \begin{split}
    &\phantom{111}\|\iI_\eps (\kK_\eps) - \iI_0 (\kK_0)\|_{\xX_{\tau}^{(1)}} + \|\lL_\eps \iI_\eps (\kK_\eps) - \lL_0 \iI_0 (\kK_0)\|_{\cC_{3(\sigma+\eta),\tau}^{\eta/2} \cC_{x}^{-1+2\eta}}\\
    &\lesssim \tau^{\theta} \Big( \eps^\theta \|\kK_0\|_{L_{3(\sigma+3\eta),\tau}^\infty \cC_x^{-\sigma}} + \|\kK_\eps - \kK_0\|_{L_{3(\sigma+3\eta),\tau}^\infty \cC_x^{-\sigma}} \Big)\;.
    \end{split}
\end{equation}
It is standard that
\begin{equation} \label{e:phi4_convergence_K0}
    \|\kK_0\|_{L_{3(\sigma+3\eta),\tau}^\infty \cC_x^{-\sigma}} \lesssim  \mM^3\;.
\end{equation}
Hence, it remains to show that
\begin{equation} \label{e:phi4_convergence_K_difference}
    \|\kK_\eps - \kK_0\|_{L_{3(\sigma+3\eta),\tau}^\infty \cC_x^{-\sigma}} \lesssim \mM^2 \big( \eps^\theta \mM + \| \Upsilon_\eps - \Upsilon_0 \|_{{\yY}_{\tau}} + \|\vec{ v }_{ \eps } - \vec{v}_0\|_{\widetilde{\xX}_\tau} \big)\;.
\end{equation}
There are various terms in $\kK_\eps$. We give sketch for the following three, while all other terms can be treated in similar ways but only simpler: 
\begin{equation*}
    u^{\#}_{\eps}  \circ_\eps \<2>_\eps\;,\; \big([ \, \iI_\eps, \Lambda_\eps u_{\eps} \prec_\eps \, ] \, \<2>_\eps \big) \circ_\eps \<2>_\eps \;,\; \Com_\eps \big( u_\eps ; \iI_\eps(\<2>_\eps); \<2>_\eps \big)\;.
\end{equation*}
For the first one, by Lemma~\ref{lem:homo_diff_circ_estimate}, for fixed time, we have
\begin{equation*}
    \begin{split}
    \| u_\eps^{\#}(r) \circ_\eps  \<2>_\eps &(r) - u_0^{\#}(r) \circ \<2>_0(r) \|_{\cC_x^\eta} \lesssim \|u_0^{\#}(r)\|_{1+2\eta} |\!|\!| \<2>_\eps(r); \<2>_0(r)|\!|\!|_{-1-\kappa;\eps}\\
    &+ \Big( \eps^\theta \|u_\eps^{\#}(r)\|_{1+2\eta;\eps} + |\!|\!| u_\eps^{\#}(r); u_0^{\#}(r) |\!|\!|_{1+2\eta;\eps} \Big) \|\<20>_\eps(r)\|_{-1-\kappa;\eps}\;.
    \end{split}
\end{equation*}
Note that
\begin{equation*}
    \begin{split}
    |\!|\!| u_\eps^{\#}(r); u_0^{\#}(r) |\!|\!|_{1+2\eta; \eps} &\lesssim \|u_\eps^{\#}(r) - u_\eps^{\#}(r)\|_{\cC_x^{-\sigma}} + \|\lL_\eps u_\eps^{\#}(r) - \lL_0 u_0^{\#}(r)\|_{-1+2\eta}\;,\\
    |\!|\!| \<2>_\eps; \<2>_0 |\!|\!|_{-1-\kappa;\eps} &\lesssim \|\Upsilon_\eps - \Upsilon_0\|_{\yY_T}\;,
    \end{split}
\end{equation*}
we get the bound
\begin{equation*}
    \| u_\eps^{\#} \circ_\eps  \<2>_\eps - u_0^{\#} \circ \<2>_0 \|_{L_{3(\sigma+3\eta),\tau}^\infty \cC_x^\eta} \lesssim \mM \big( \eps^\theta \mM + \|\vec{ v }_{ \eps } - \vec{v}_0\|_{\widetilde{\xX}_\tau} + \|\Upsilon_\eps - \Upsilon_0\|_{\yY_T} \big)\;.
\end{equation*}
As for the second one, by Lemma~\ref{le:convergence_heat_comm_below}, we have
\begin{equation} \label{e:phi4_convergence_comm}
    \begin{split}
    &\phantom{111}\big\| \lL_\eps [\iI_\eps, \Lambda_\eps u_\eps \prec_\eps] \<2>_\eps - \lL_0 [\iI_0, \Lambda_0 u_0 \prec] \<2>_0 \big\|_{L_{3(\sigma+3\eta),\tau}^\infty \cC_x^{\eta-1-2\kappa}}\\
    &\lesssim \mM \big( \eps^\theta \mM + \|u_\eps - u_0\|_{\xX_\tau} + \|\Upsilon_\eps - \Upsilon_0\|_{\yY_T} \big)\;.
    \end{split}
\end{equation}
Another application of Lemma~\ref{lem:homo_diff_circ_estimate} yields the desired bound for
\begin{equation*}
    \big([ \, \iI_\eps, \Lambda_\eps u_{\eps} \prec_\eps \, ] \, \<2>_\eps \big) \circ_\eps \<2>_\eps - \big([ \, \iI_0, \Lambda_0 u_0 \prec \, ] \, \<2>_0 \big) \circ \<2>_0\;.
\end{equation*}
The bound for the difference
\begin{equation*}
    \Com_\eps \big( u_\eps ; \iI_\eps(\<2>_\eps); \<2>_\eps \big) - \Com_0 \big( u_0 ; \iI_0(\<2>_0); \<2>_0 \big)
\end{equation*}
follows from Lemma~\ref{le:convergence_comm_below} in essentially the same way. 

Combining the bounds \eqref{e:phi4_convergence_initial_1}, \eqref{e:phi4_convergence_initial_2}, \eqref{e:phi4_convergence_rough}, \eqref{e:phi4_convergence_IK}, \eqref{e:phi4_convergence_K0} and \eqref{e:phi4_convergence_K_difference}, we conclude there exists a $\theta>0$ and a universal $C_0 > 0$ such that
\begin{equation*}
    \|\vec{U}_{\eps} - \vec{U}_0\|_{\vec{\xX}_{\tau,\eps}} \leq C_0 \mM^2 \Big( \tau^{\theta} \|\vec{U}_{ \eps }- \vec{U}_0\|_{\vec{\xX}_{\tau,\eps}} + \eps^\theta \mM + \| \Upsilon_\eps - \Upsilon_0 \|_{{\yY}_{T}} + \|u_\eps(0) - u_0(0)\|_{-\sigma} \Big)\;.
\end{equation*}
We now choose $\tau \in (0, \tau_0)$ sufficiently small such that
\begin{equation*}
    C_0 \tau^\theta \mM^2 \leq \frac{1}{2}\;.
\end{equation*}
This then gives
\begin{equation*}
    \| \vec{U}_{\eps} - \vec{U}_0 \|_{\widetilde{\xX}_{\tau}} \lesssim \eps^\theta + \| \Upsilon_\eps-\Upsilon_0 \|_{{\yY}_{T}} + \|u_\eps(0) - u_0(0)\|_{-\sigma}\;, 
\end{equation*}
where the proportionality constant depends on $\mM$ only. For sufficiently small $\eps$, this further implies
\begin{equation*}
    \|u_\eps(\tau) - u_{0}(\tau)\|_{-\sigma} + \|u_\eps^{\#}(\tau) - u_0^{\#}(\tau)\|_{-\sigma} < \frac{1}{2}\;.
\end{equation*}
Since the smallness of such a time $\tau$ does not depend on $\eps$, one can iterate this procedure to extend the convergence to the whole time interval $[0,T]$ to obtain
\begin{equation*}
    \| \vec{U}_{\eps} - \vec{U}_0 \|_{\widetilde{\xX}_T} \lesssim \eps^\theta + \| \Upsilon_\eps-\Upsilon_0 \|_{{\yY}_{T}} + \|u_\eps(0) - u_0(0)\|_{-\sigma}\;.
\end{equation*}
This completes the proof of the theorem. 
\end{proof}

\begin{thm} \label{thm:phi4_convergence_flux}
For $\eps\in \N^{-1}$, let $(u_\eps, u_\eps^{\#})$ be the fixed point solution to \eqref{e:phi4_system}. For every $\kappa>0$, there exists $\theta>0$ such that
\begin{equation} \label{e:phi4_flux}
    \begin{split}
    \|A_{\eps} \nabla u_{\eps}-\bar{A}\nabla u_{0}\|_{ L^{\infty}_{3\sigma+\eta,T}\cC^{-\kappa}_{x}} \lesssim &\phantom{1}\eps^{\theta}+\|u_{\eps}(0)-u_{0}(0)\|_{-\frac{1}{2}-\eta}\\
    &+ \|\Upsilon_\eps - \Upsilon_0\|_{\yY_T} + \|A_\eps \nabla \<30>_\eps - \bar{A} \nabla \<30>_0\|_{L_T^\infty \cC_x^{-\frac{1}{2}-\kappa}}\;.
    \end{split}
\end{equation}
\end{thm}
\begin{proof}
Since $(u_\eps, u_\eps^{\#})$ solves the system \eqref{e:phi4_system} up to time $T$, we have
\begin{equation} \label{e:phi4_flux_decomposition}
    u_\eps =  e^{t \lL_\eps} u_\eps(0) - 3 \, (\Lambda_\eps u_\eps) \prec_\eps \iI_\eps (\<2>_\eps) - 3 [\iI_\eps, \Lambda_\eps u_\eps \prec_\eps] \<2>_\eps + u_\eps^\#\;.
\end{equation}
The convergence of the initial data term $A_\eps \nabla e^{t \lL_\eps} u_\eps(0)$ to $\bar{A} \nabla e^{t \lL_0} u_0(0)$ is immediate. The last two terms on the right hand side of \eqref{e:phi4_flux_decomposition} has regularity strictly above $1$. For general $G_\eps \in \cC_x^{1+\kappa;\eps}$ and $G_0 \in \cC^{1+\kappa}$, we have
\begin{equation*}
    \begin{split}
    \|A_\eps \nabla G_\eps - \bar{A} \nabla G_0\|_{-\kappa} &\leq \|(A_\eps \nabla \lL_\eps^{-1} - \bar{A} \nabla \lL_0^{-1}) \lL_\eps G_\eps\|_{-\kappa} + \|\bar{A} \nabla \lL_0^{-1} (\lL_\eps G_\eps - \lL_0 G_0)\|_{-\kappa}\\
    &\lesssim \eps^\theta \|\lL_\eps G_\eps\|_{-1+\kappa} + \|\lL_\eps G_\eps - \lL_0 G_0\|_{-1+\kappa}\;,
    \end{split}
\end{equation*}
where the bounds in the last line follow from flux estimate for standard elliptic periodic homogenisation. 

Taking $G_\eps = u_\eps^{\#}$ and $[\iI_\eps, \Lambda_\eps u_\eps \prec_\eps] \<2>_\eps$ and employing the bound for $\lL_\eps u_\eps^{\#} - \lL_0 u_0^{\#}$ in Theorem~\ref{thm:phi4_convergence_homo} for the former, and the bound \eqref{e:phi4_convergence_comm} for the latter, we obtain the desired bounds for
\begin{equation*}
    A_\eps \nabla u_\eps^{\#} - \bar{A} \nabla u_0^{\#} \quad \text{and} \quad A_\eps \nabla \big( [\iI_\eps, \Lambda_\eps u_\eps \prec_\eps] \, \<2>_\eps \big) - \bar{A} \nabla \big( [\iI_0, \Lambda_0 u_0 \prec] \, \<2>_0\big)\;.
\end{equation*}
It remains to control
\begin{equation*}
    A_\eps \nabla \big( \Lambda_\eps u_\eps \prec_\eps \iI_\eps (\<2>_\eps) \big) - \bar{A} \nabla \big( \Lambda_0 u_0 \prec \iI_0 (\<2>_0) \big)\;.
\end{equation*}
Note that we can write
\begin{equation*}
    \begin{split}
    A_\eps \nabla \big( &\Lambda_\eps u_\eps \prec_\eps \iI_\eps (\<2>_\eps) \big) = \Lambda_\eps (u_\eps) \cdot A_\eps \nabla \iI_\eps (\<2>_\eps)\\
    &+ \sum_{j \geq -1} \Big[ A_\eps \nabla S_{j,\eps} (\Lambda_\eps u_\eps) \cdot \Delta_{j,\eps} \iI_\eps (\<2>_\eps) + S_{j,\eps}^{\perp} (\Lambda_\eps u_\eps) \cdot A_\eps \nabla \Delta_{j,\eps} \iI_\eps (\<2>_\eps) \Big]\;.
    \end{split}
\end{equation*}
For the first term on the right hand side, it is standard that
\begin{equation*}
    \begin{split}
    &\phantom{111}\|\Lambda_\eps (u_\eps) \cdot A_\eps \nabla \iI_\eps (\<2>_\eps) - \Lambda_0 (u_0) \cdot \bar{A} \nabla \iI_0 (\<2>_0)\|_{ L^{\infty}_{3\sigma+\eta,T}\cC^{-\kappa}_{x}}\\
    &\lesssim \eps^{\theta} + \|u_{\eps} - u_{0}\|_{\xX_T^{(1)}} + \|\Upsilon_\eps - \Upsilon_0\|_{\yY_T} + \|A_\eps \nabla \<30>_\eps - \bar{A} \nabla \<30>_0\|_{L_T^\infty \cC_x^{-\frac{1}{2}-\kappa}}\;.
    \end{split}
\end{equation*}
We now turn to the terms with summation over $j \geq -1$. For each term in the sum, the $L_x^\infty$-norm can be controlled by
\begin{equation*}
    \begin{split}
    \| A_\eps \nabla S_{j,\eps} (\Lambda_\eps u_\eps) \cdot \Delta_{j,\eps} \iI_\eps (\<2>_\eps) \|_{L_x^\infty} &\lesssim \|\nabla S_{j,\eps} (\Lambda_\eps u_\eps)\|_{L^\infty} \|\Delta_{j,\eps} \iI_\eps (\<2>_\eps)\|_{L^\infty}\\
    &\lesssim 2^{(1-\alpha-\beta) j} \|\Lambda_\eps u_\eps\|_{\alpha} \|\iI_\eps (\<2>_\eps)\|_{\beta}\;,
    \end{split}
\end{equation*}
and
\begin{equation*}
    \begin{split}
    \| S_{j,\eps}^{\perp} (\Lambda_\eps u_\eps) \cdot A_\eps \nabla \Delta_{j,\eps} \iI_\eps (\<2>_\eps) \|_{L^\infty} &\lesssim \|S_{j,\eps}^{\perp} (\Lambda_\eps u_\eps)\|_{L^\infty} \cdot \|\nabla \Delta_{j,\eps} \iI_\eps (\<2>_\eps)\|_{L^\infty}\\
    &\lesssim 2^{(1-\alpha-\beta) j} \|\Lambda_\eps u_\eps\|_\alpha \|\iI_\eps (\<2>_\eps)\|_{\beta}\;,
    \end{split}
\end{equation*}
where we take $\alpha = \eta$ and $\beta = 1-\kappa$, and they are summable in $j$. 

On the other hand, one can show that the differences satisfy the bound
\begin{equation*}
    \begin{split}
    &\phantom{111}\| A_\eps \nabla S_{j,\eps} (\Lambda_\eps u_\eps) \cdot \Delta_{j,\eps} \iI_\eps (\<2>_\eps) - \bar{A} \nabla S_{j,0} (\Lambda_0 u_0) \cdot \Delta_{j,0} \iI_0 (\<2>_0) \|_{\cC_{x}^{-\kappa}}\\
    &\phantom{111}+ \| S_{j,\eps}^{\perp} (\Lambda_\eps u_\eps) \cdot A_\eps \nabla \Delta_{j,\eps} \iI_\eps (\<2>_\eps) - S_{j,0}^{\perp} (\Lambda_0 u_0) \cdot \bar{A} \nabla \Delta_{j,0} \iI_0 (\<2>_0) \|_{\cC_x^{-\kappa}}\\
    &\lesssim 2^{Mj} \Big( \eps^{\theta'}  \|\Lambda_\eps u_\eps\|_\alpha \|\iI_\eps (\<2>_\eps)\|_{\beta} + \|\Lambda_\eps u_\eps - \Lambda_0 u_0\|_\alpha \| \iI_0 (\<2>_0)\|_\beta\\
    &\phantom{111}+ \|\Lambda_0 u_0\|_\alpha \|\iI_\eps (\<2>_\eps) - \iI_0 (\<2>_0)\|_\beta \Big)
    \end{split}
\end{equation*}
for some $M>0$ and $\theta'>0$. Interpolating these two bounds for each $j$ and then summing over $j \geq -1$ then gives the desired bound for the summation on the right hand side. This completes the proof for the convergence of the flux. 
\end{proof}

\subsection{Proof of Theorem~\ref{th:phi4_overall}}

We now have all the ingredients to complete the proof of Theorem~\ref{th:phi4_overall}. 

\begin{proof} [Proof of Theorem~\ref{th:phi4_overall}]
Let $u_\eps^{(\delta)}$ be the first component of the solution $\vec{u}_\eps^{(\delta)}$ to \eqref{e:phi4_system} with initial data
\begin{equation*}
    u_\eps^{(\delta)}(0) = \phi_\eps^{(\delta)}(0) - \<1>_\eps^{(\delta)}(0) + \<30>_\eps^{(\delta)}(0)\;.
\end{equation*}
Then, it follows directly from the derivation of the system \eqref{e:phi4_system} that
\begin{equation*}
    \phi_\eps^{(\delta)} := \<1>_\eps^{(\delta)} - \<30>_\eps^{(\delta)} + u_\eps^{(\delta)}
\end{equation*}
solves the original equation \eqref{e:phi4_hom_full} with initial data $\phi_\eps(0)$. The desired convergence of the solution and fluxes then follow immediately from Theorems~\ref{thm:phi4_convergence_homo} and~\ref{thm:phi4_convergence_flux}. 
\end{proof}

\section{$2$D gPAM}
\label{sec:gpam}

In this section, we study the two dimensional generalised parabolic Anderson model given by
\begin{equation} \label{e:PAM_homo_full}
    \d_t u_\eps^{(\delta)} = \lL_\eps u_\eps^{(\delta)} + \Tr \Big( F_\eps^\sT \big( (\nabla u_\eps^{(\delta)})^{\otimes 2} - g^2(u_\eps^{(\delta)}) \widetilde{\fC}_\eps^{(\delta)} \big) \Big) + g(u_\eps^{(\delta)}) \, \big( \xi^{(\delta)} - C_\eps^{(\delta)} g'(u_\eps^{(\delta)}) \big)\;.
\end{equation}
Here, $F_\eps = \{f_{ij}(\cdot/\eps)\}$ for bounded periodic functions $f_{ij}$ with $f_{ij} = f_{ji}$, and $C_\eps^{(\delta)}$ and $\widetilde{\fC}_\eps^{(\delta)}$ are suitable renormalisation functions and matrices specified in Proposition~\ref{prop:PAM_stochastic} below. The candidate limiting equation for \eqref{e:PAM_homo_full} is
\begin{equation} \label{e:PAM_limiting_equation}
    \d_t u_0 = \div (\bar{A} \nabla u_0) + \Tr \Big( \bar{F}^\sT \big( (\nabla u_0)^{\otimes 2} - g^2(u_0) \, \widetilde{\boldsymbol{\infty}} \big) \Big) + g(u_0) \big( \xi - \infty g'(u_0) \big)\;,
\end{equation}
where $\bar{F}$ is the constant symmetric $2 \times 2$ matrix given by
\begin{equation} \label{e:PAM_limiting_matrix}
    \bar{F} := \Pi_0 (\Phi^\sT F \Phi) = \int_{\T^2} \big( \id + (\nabla \chi)(x) \big)^{\sT} \cdot F(x) \cdot \big( \id + (\nabla \chi)(x) \big) {\rm d}x\;.
\end{equation}
In the rest of this section, we fix two small parameters $\eta \in (0, \frac{1}{1000})$ and $\kappa \in (0, \frac{\eta}{1000})$. For two $d \times d$ matrices $P$ and $Q$, we define a new $d \times d$ matrix $P \odot Q$ by entry-wise multiplication as
\begin{equation*}
    (P \odot Q)_{i,j} = P_{ij} Q_{ij}\;.
\end{equation*}

\subsection{Definition of stochastic objects and main result}

Recall $\xi^{(\delta)}$ is a smooth approximation to $\xi$. For $\eps \in \N^{-1}$ and $\delta \in (0,1)$, define the two scalar-valued stochastic objects $\<1p>_\eps^{(\delta)}$ and $\<11p>_{\,\eps}^{(\delta)}$ by
\begin{equation*}
    \<1p>_\eps^{(\delta)} := - \lL_\eps^{-1} \Pi_0^\perp \xi^{(\delta)}\;, \qquad \<11p>_{\,\eps}^{(\delta)} := \<1p>_\eps^{(\delta)} \circ_\eps \xi^{(\delta)} - C_\eps^{(\delta)}\;.
\end{equation*}
Define the three vector-valued objects $\<1pd>_\eps^{(\delta)}$, $\<1pdF>_\eps^{(\delta)}$ and $\<1pdFP>_\eps^{(\delta)}$ by
\begin{equation*}
\<1pd>_{\eps}^{(\delta)} := \nabla \, \<1p>^{(\delta)}_{\eps}\;,\qquad \<1pdF>_{\eps}^{(\delta)} := F_\eps \nabla \, \<1p>^{(\delta)}_{\eps} = F_{\eps} \<1pd>_{\eps}^{(\delta)}\;, \qquad \<1pdFP>_{\eps}^{(\delta)} := \Phi_\eps^\sT F_\eps \nabla \, \<1p>^{(\delta)}_{\eps} = \Phi_{\eps}^{\sT} F_\eps \<1pd>_{\eps}^{(\delta)}\;.
\end{equation*}
Finally, define the matrix-valued objects
\begin{equation*}
    \<1pdF1pd0>_\eps^{(\delta)} := \nabla \<1p>_\eps^{(\delta)} \otimes \nabla \<1p>_\eps^{(\delta)} - \widetilde{\fC}_\eps^{(\delta)}\;, \quad F_\eps \odot \<1pdF1pd0>_{\eps}^{(\delta)} :=  F_\eps \odot \big( \nabla \<1p>_\eps^{(\delta)} \otimes \nabla \<1p>_\eps^{(\delta)} - \widetilde{\fC}_\eps^{(\delta)} \big)\;.
\end{equation*}
For every $\eps \in \N^{-1}$ and $\delta \in (0,1)$, define
\begin{equation*}
    \Upsilon_\eps^{(\delta)} := \big( \Pi_0\xi^{(\delta)}, \; \<1p>_\eps^{(\delta)}, \; \<11p>_\eps^{(\delta)}, \; \eps^{2\kappa} \, \<1pdF>_\eps^{(\delta)} ,\; \<1pdFP>_\eps^{(\delta)}, \; F_\eps \odot \<1pdF1pd0>_\eps^{(\delta)} \big)\;.
\end{equation*}
Let
\begin{equation*}
    \Upsilon_0 = \big( \Pi_0 \xi, \; \<1p>_0, \; \<11p>_0, \; \<1pdF1pd0>_0 \big)
\end{equation*}
be a choice of enhanced stochastic objects from standard $2$D g-PAM with the operator $\lL_0 = \div (\bar{A} \nabla)$. For a $2 \times 2$ matrix $P$, define the action
\begin{equation*}
    P \bullet \Upsilon_0 := \big( \Pi_0 \xi, \; \<1pd>_0, \; \<11p>_0, \; 0, \; P \, \<1pd>_0, \; P \odot \<1pdF1pd0>_0 \big)\;.
\end{equation*}
All the objects above are independent of time. Let
\begin{equation*}
    \yY := \R \times \cC^{1-\kappa}(\T^2; \R) \times \cC^{-\kappa} (\T^2; \R) \times \cC^{-\kappa} (\T^2; \R^2) \times \cC^{-\kappa}(\T^2; \R^2) \times \cC^{-\kappa}(\T^2; \R^{2 \times 2})\;.
\end{equation*}
The following proposition from \cite{homo_stochastic} gives convergences of the stochastic objects. 

\begin{prop} [\cite{homo_stochastic}, Theorem~4.1]
\label{prop:PAM_stochastic}
Recall $\bar{F}$ from \eqref{e:PAM_limiting_matrix}. For the choices
\begin{equation*}
    C_\eps^{(\delta)} = \E \big( \, \<1p>_\eps^{(\delta)} \cdot \xi^{(\delta)} \big)\;, \qquad (\widetilde{\fC}_\eps^{(\delta)})_{i,j} = \E \big( \d_i \<1p>_\eps^{(\delta)} \cdot \d_j \<1p>_\eps^{(\delta)} \big)\;,
\end{equation*}
we have the convergences $\Upsilon_\eps^{(\delta)} \rightarrow \bar{F} \bullet \Upsilon_0$ in $\yY$ and $A_\eps \nabla \<1p>_\eps^{(\delta)} \rightarrow \bar{A} \nabla \<1p>_0$ in $\cC^{-\kappa}(\T^2; \R^2)$ as $(\eps, \delta) \rightarrow (0,0)$. 
\end{prop}

For readers' convenience, we summarise the stochastic objects in $\Upsilon_\eps^{(\delta)}$ and their regularities in the following table. 
\begin{center}
{\renewcommand{\arraystretch}{1.2}
\begin{tabular}{p{1.5cm} p{1.2cm} p{1.0cm} p{1.2cm} p{1.2cm} p{1.2cm} p{1.7cm}}
\hline
Process & $\Pi_0 \xi^{(\delta)}$ & $\<1p>_\eps^{(\delta)}$ & $\<11p>_\eps^{(\delta)}$ & $\<1pdFP>_{\eps}^{(\delta)}$ & $ \eps^{2\kappa} \<1pdF>_\eps^{(\delta)}$ & $F_\eps \odot \<1pdF1pd0>_\eps^{(\delta)}$\\
\hline
Limit & $\Pi_0 \xi$ & $\<1p>_0$ & $\<11p>_0$ & $\bar{F} \, \<1pd>_0$ & $0$ & $\bar{F} \odot \<1pdF1pd0>_0$\\
\hline
Reg. & \text{const} & $1-$ & $0-$ & $0-$ & $0-$&$0-$\\
\hline
\end{tabular}
}
\end{center}
Here ``$\alpha-$" denotes $\cC_x^{\alpha-\kappa}$. Later we have the noise $\xi^{(\delta)}$ (including $\delta=0$) appearing in the fixed point system for both $\eps \in \N^{-1}$ and $\eps=0$, but only $\Pi_0 \xi^{(\delta)}$ is included in $\Upsilon_\eps^{(\delta)}$. This is because $\xi^{(\delta)}$ can be completely recovered from $\Pi_0 \xi^{(\delta)}$ and $\<1p>_\eps^{(\delta)}$ by
\begin{equation*}
    \xi^{(\delta)} = \Pi_0 \xi^{(\delta)} - \lL_\eps \<1p>_\eps^{(\delta)} = \Pi_0 \xi^{(\delta)} - \lL_0 \<1p>_0^{(\delta)}\;,
\end{equation*}
We now state our main theorem for $2$D gPAM. 

\begin{thm} \label{thm:PAM_overall}
Let $C_\eps^{(\delta)}$ and $\widetilde{\fC}_\eps^{(\delta)}$ be chosen such that Proposition~\ref{prop:PAM_stochastic} holds. Let $\{u_\eps(0)\} \in \cC^{1-\eta}$ such that
\begin{equation*}
    \|u_\eps(0) - u_0(0)\|_{\cC^{1-\eta}} \rightarrow 0
\end{equation*}
for some $u_0(0) \in \cC^{1-\eta'}$ for some $\eta' < \eta$ as $\eps \rightarrow 0$. Let $u_\eps$ be the solution to \eqref{e:PAM_homo_full} with initial data $u_\eps(0)$ and enhanced stochastic objects $\Upsilon_\eps^{(\delta)}$. Let $u_0$ be the solution to \eqref{e:PAM_limiting_equation} with initial data $u_0(0)$ and enhanced stochastic objects $\Upsilon_0$. 

Suppose $u_0$ exists up to (a random) time $T_0$. Then for every $T<T_0$ and all sufficiently small $\eps$ and $\delta$ (depending on the realisation of $\Upsilon_\eps^{(\delta)}$), $u_\eps^{(\delta)}$ exists up to time $T$, and there exists $\theta>0$ such that
\begin{equation*}
    \|u_\eps^{(\delta)} - u_0\|_{\cC_\fs^{1-\eta}([0,T] \times \T^2)} \lesssim \eps^\theta + \|u_\eps(0) - u_0(0)\|_{1-\eta} + \|\Upsilon_\eps^{(\delta)} - \bar{F} \bullet \Upsilon_0\|_{\cC^{-\kappa}}\;.
\end{equation*}
Furthermore, the flux also converges with the quantitative bound
\begin{equation*}
    \big\| A_\eps \nabla u_\eps^{(\delta)} - \bar{A} \nabla u_0 \big\|_{L_T^\infty \cC_x^{-2\eta}} \lesssim \eps^\theta + \|u_\eps(0) - u_0(0)\|_{1-\eta} + \|\Upsilon_\eps^{(\delta)} - \bar{F} \bullet \Upsilon_0\|_{\yY} + \|A_\eps \nabla \<1p>_\eps^{(\delta)} - \bar{A} \nabla \<1p>_0 \|_{\cC_x^{-\eta}}\;.
\end{equation*}
\end{thm}
\begin{proof}
It suffices to note that $u_\eps^{(\delta)}$ solves \eqref{e:PAM_homo_full} if and only if there exists $u_\eps^{\#; \delta}$ such that the pair $\vec{u}_\eps^{(\delta)} = \big(u_\eps^{(\delta)}, u_\eps^{\#; \delta}\big)$ solves the fixed point system \eqref{eq:PAM_system} (note that we omitted the $\delta$ in notation in the following). Hence, the claims of the theorem follows directly from Theorems~\ref{thm:PAM_convergence_homo} and~\ref{thm:PAM_convergence_flux}. 
\end{proof}

\subsection{Derivation of the fixed point system}

As in Section~\ref{sec:phi4}, we omit the superscript $\delta$ from the notations. So here, $\Upsilon_\eps$ and $u_\eps$ still denote the stochastic objects $\Upsilon_\eps^{(\delta)}$ and solution $u_\eps^{(\delta)}$ (possibly for both $\delta>0$ and $\delta=0$). 

We first decompose the product $g(u_\eps) \, \xi$ by
\begin{equation*}
    g(u_\eps ) \, \xi = g(u_\eps) \prec_\eps \xi + g(u_\eps) \succ_\eps \xi + g(u_\eps) \circ_\eps \xi\;.
\end{equation*}
For the resonance product $g(u_\eps) \circ_\eps \xi$, we use the para-linearisation operation $\PL_\eps^g$ from \eqref{e:para_linear_defn} and the commutator $\Com_\eps$ from \eqref{e:commutator_defn} to write it as
\begin{equation*}
    g(u_\eps) \circ_\eps \xi = g'(u_\eps) \, \big( u_\eps \circ_\eps \xi \big) + \Com_\eps \big( g'(u_\eps); \, u_\eps; \, \xi \big) + \PL_\eps^g (u_\eps) \circ_\eps \xi\;.
\end{equation*}
We can thus re-write the equation \eqref{e:PAM_homo_full} as
\begin{equation} \label{eq:PAM_rewrite}
    \begin{split}
    \d_t u_\eps &= \lL_\eps u_\eps + g(u_\eps) \prec_\eps \xi + \underbrace{g(u) \succ_\eps \xi + \Com_\eps \big( g'(u); \, u; \, \xi \big) + \PL_\eps^g (u) \circ_\eps \xi}_{\kK_\eps^{(0)}}\\
    &+ g'(u_\eps) \big( \underbrace{u_\eps \circ_\eps \xi - C_\eps g(u_\eps)}_{\kK_\eps^{(1)}} \big) + \big( \underbrace{\nabla^\sT u_\eps \cdot F_\eps \, \nabla u_\eps - \widetilde{C}_\eps g^2 (u_\eps)}_{\kK_\eps^{(2)}} \big)\;.
    \end{split}
\end{equation}
If $u \in \cC^{1-}$, then terms in $\kK_\eps^{(0)}$ above are all well behaved (the worst term having regularity $0-$). Hence, 
\begin{equation} \label{eq:PAM_K0}
    \kK_\eps^{(0)}(u) := g(u) \succ_\eps \xi + \Com_\eps \big( g'(u); \, u; \, \xi \big) + \PL_\eps^g (u) \circ_\eps \xi
\end{equation}
is a faithful representation. On the other hand, the resonance product $u_\eps \circ_\eps \xi$ and the normal product $\nabla^\sT u_{\eps} \cdot F_\eps \nabla u_\eps$ are both a priori ill-defined, and we need to explore more structures to find well-posed expressions for $\kK_\eps^{(1)}$ and $\kK_\eps^{(2)}$. 

We start with $\kK_\eps^{(1)}$. We make the ansatz
\begin{equation} \label{eq:PAM_ansatz}
    \begin{split}
    u_\eps (t) &= e^{t \lL_\eps} u_\eps(0) + \iI_\eps \big( g(u_\eps) \prec_\eps \xi \big) + \big( u_\eps^{\#}(t) - e^{t \lL_\eps} u_\eps^{\#}(0) \big)\\
    &= g(u_\eps) \prec_\eps \iI_\eps (\xi) + \gG_\eps (\vec{u}_\eps)\;,
    \end{split}
\end{equation}
where
\begin{equation} \label{eq:PAM_gG}
    \gG_{\eps}(\vec{u}):=e^{t \lL_\eps} \big( u_{\eps}(0) - u_{\eps}^{\#}(0) \big) + [\iI_\eps, g(u) \prec_\eps] \xi + u^{\#}\;,
\end{equation}
and we used the notation $\vec{u} = (u, u^{\#})$. We expect $u_\eps \in \cC^{1-}$ and $u_\eps^{\#} \in \cC^{1+; \eps}$ uniformly in $(\eps, \delta)$ (recall that we omitted the parameter $\delta$ for notational simplicity), so that $\gG_\eps (\vec{u}_\eps)$ has uniform $\cC^{1+;\eps}$ regularity. With the above ansatz, we have
\begin{equation*}
    \begin{split}
    u_\eps \circ_\eps \xi - C_\eps g(u_\eps) = \underbrace{\big( g(u_\eps) \prec_\eps \Pi_0^\perp \iI_\eps (\xi) \big) \circ_\eps \xi - C_\eps g(u_\eps)}_{g(u_\eps) \cdot ( \Pi_0^\perp \iI_\eps (\xi) \circ_\eps \xi - C_\eps ) \, + \, \Com_\eps ( g(u_\eps); \, \Pi_0^\perp \iI_\eps(\xi); \, \xi )} + \, \gG_{\eps}(\vec{u}_{\eps}) \circ_\eps \xi\;.
    \end{split}
\end{equation*}
Also note that by definition of $\<1p>_\eps$ and $\<11p>_\eps$, we have
\begin{equation*}
    \Pi_0^\perp \iI_\eps (\xi) \circ_\eps \xi - C_\eps = \big( \<1p>_\eps - e^{t \lL_\eps} \<1p>_\eps \big)  \circ_\eps \xi - C_\eps = \<11p>_\eps - (e^{t \lL_\eps} \<1p>_\eps) \circ_\eps \xi\;. 
\end{equation*}
Hence, adding back the initial data, this suggests us to define
\begin{equation} \label{eq:PAM_K1}
    \begin{split}
    \kK_\eps^{(1)}(\vec{u}) :=  g(u) \, \big( \<11p>_\eps - e^{t \lL_\eps} \<1p>_\eps \circ_\eps \xi \big) +\Com_\eps \big( g(u); \Pi_0^{\perp}\iI_\eps(\xi); \xi \big)+ \gG_{\eps}(\vec{u}) \circ_\eps \xi\;.
    \end{split}
\end{equation}
To proceed with $\kK_\eps^{(2)}$, we need to investigate finer structure of $\nabla u_\eps$. Applying $\nabla$ to both sides of \eqref{eq:PAM_ansatz} and re-writing $\nabla = \nabla \lL_\eps^{-1} \lL_\eps$, we have
\begin{equation} \label{eq:PAM_gradient_1st_expansion}
    \begin{split}
    \nabla u_\eps &=  \nabla \big( g(u_\eps) \prec_\eps \<1p>_\eps \big) + \underbrace{\nabla \lL_\eps^{-1} \lL_\eps \Big(  \gG_{\eps}( \vec{ u }_{ \eps } )- g(u_\eps) \prec_\eps \big( e^{t \lL_\eps} \<1p>_\eps \big) \Big)}_{\Phi_\eps \cdot \cC^{0+} \; + \text{small} \; \cC^{0+}}\;.
    \end{split}
\end{equation}
The quantity inside the big parenthesis above has regularity $\cC{1+; \,\eps}$. Hence, by properties of the operators $\lL_\eps$ and $\nabla \lL_\eps^{-1}$, we expect the second term on the right hand side above to be well behaved (having the structure indicated above). To deal with the first term, we need the following lemma. 

\begin{lem} \label{lem:comm_nabla_prec}
Let $\alpha, \beta \in (0,1)$ with $\alpha + \beta > 1$. Then, there exist $\theta, \nu \in (0,1)$ and bilinear maps $\Err_{\eps}^{(1)}, \Err_{\eps}^{(2)}:  \cC^\alpha \times \cC^\beta \rightarrow \cC^\nu$ with
\begin{equation} \label{eq:gradient_prec_commutator_error}
    \big\| \Err_\eps^{(1)}(f;h) \big\|_{\nu} + \|\Err_\eps^{(2)}(f;h)\|_{\nu} \lesssim \eps^\theta \|f\|_\alpha \|h\|_\beta\;,
\end{equation}
and that the following identity holds:
\begin{equation} \label{eq:gradient_prec_commutator}
    \nabla \big( f \prec_\eps h \big) = f \, \nabla h + \Phi_\eps \Big( \nabla f \prec h - f \succeq \nabla h + \Err_\eps^{(1)}(f;h) \Big) + \Err_{\eps}^{(2)}(f;h)\;.
\end{equation}
\end{lem}
\begin{proof}
Let
\begin{equation*}
    \begin{split}
    \Err_\eps^{(1)}(f;h) &:= \sum_j \bigg[  \nabla S_{j,0} f \cdot (\Delta_{j,\eps} - \Delta_{j,0}) h - \big( S_{j,\eps}^\perp - S_{j,0}^\perp \big) f \cdot \nabla \Delta_{j,0} h \bigg]\;,\\
    \Err_{\eps}^{(2)}(f;h) &:= \sum_j \bigg[ \big( \nabla S_{j,\eps} - \Phi_\eps \nabla S_{j,0} \big) f \cdot \Delta_{j,\eps} h- S_{j,\eps}^\perp f \cdot \big( \nabla \Delta_{j,\eps} - \Phi_{\eps}\nabla\Delta_{j,0} \big)  h \bigg]\;.
    \end{split}
\end{equation*}
One can check directly that the identity \eqref{eq:gradient_prec_commutator} holds with the above functionals $\Err_\eps^{(1)}$ and $\Err_\eps^{(2)}$. It then remains to show the bounds in \eqref{eq:gradient_prec_commutator_error}. By Corollary \ref{cor:convergence_single_block}, for every $\kappa>0$, $\nu\in(0,1)$, there exists $\theta>0$ such that for every $j\geq -1$, $\eps\in \N^{-1}$,
\begin{equation*}
    \|(\Delta_{j,\eps}-\Delta_{j,0})h\|_{\nu}\lesssim \eps^{\theta}2^{-(\beta-\nu-\kappa)j}\|h\|_{\beta}\;,\quad \|\big( S_{j,\eps}^\perp - S_{j,0}^\perp \big) f\|_{\nu}\lesssim \eps^{\theta}2^{-(\alpha-\nu-\kappa)j}\|f\|_{\alpha}\;.
\end{equation*}
By Lemma \ref{le:convergence_gradient_block}, for every $\kappa>0$, there exist $\theta,\nu>0$ such that for every $j\geq -1$, $\eps\in \N^{-1}$,
\begin{equation*}\begin{split}
    \|\big( \nabla S_{j,\eps} - \Phi_\eps \nabla S_{j,0} \big) f\|_{\nu}&\lesssim \eps^{\theta}2^{(1-\alpha+\kappa)j}\|f\|_{\alpha}\;,\\
    \|\big( \nabla \Delta_{j,\eps} - \Phi_{\eps}\nabla\Delta_{j,0} \big)  h\|_{\nu}&\lesssim \eps^{\theta}2^{(1-\beta+\kappa)j}\|h\|_{\beta}\;,
    \end{split}
\end{equation*}
Choosing $0<\nu<\frac{\alpha+\beta-1}{10}$ and sufficiently small $\kappa>0$, the above estimates yield the desired estimate for $\Err^{(1)}(f;h)$ and $\Err^{(2)}(f;h)$.
\end{proof}

Applying Lemma~\ref{lem:comm_nabla_prec} to the quantity $\nabla \big( g(u_\eps) \prec_\eps \<1p>_\eps \big)$, we see it has the expression
\begin{equation*}
    g(u_\eps) \nabla \, \<1p>_\eps + \Phi_\eps \Big( \underbrace{\nabla \big( g(u_\eps) \big) \prec \<1p>_\eps - g(u_\eps) \succeq \nabla \, \<1p>_\eps}_{\Lambda_\eps (u_\eps)} + \Err_\eps^{(1)} \big( g(u_\eps); \, \<1p>_\eps \big) \Big) + \Err_\eps^{(2)} \big( g(u_\eps); \, \<1p>_\eps \big)\;.
\end{equation*}
This suggests us to define the functional $\Lambda_\eps$ by
\begin{equation} \label{eq:PAM_Lambda}
    \Lambda_\eps (u) := \nabla \big( g(u) \big) \prec \<1p>_\eps - g(u) \succeq \nabla \, \<1p>_\eps\;.
\end{equation}
Hence, plugging it back to the right hand side of \eqref{eq:PAM_gradient_1st_expansion}, we get
\begin{equation}\label{eq:nabla_u}
    \nabla u_\eps = g(u_\eps) \nabla \<1p>_\eps + \Phi_\eps \v_\eps + \w_\eps\;,
\end{equation}
where $\v_\eps = \v_\eps (\vec{u}_\eps)$ and $\w_\eps = \w_\eps(\vec{u}_\eps)$ are respectively given by
\begin{equation} \label{eq:PAM_v}
    \begin{split}
    \v_\eps(\vec{u}):= \nabla \lL_0^{-1} \lL_\eps \Big( \gG_{\eps}(\vec{u}) - g(u) \prec_\eps e^{t \lL_\eps} \<1p>_\eps \Big)+ \Lambda_\eps (u) + \Err_\eps^{(1)} \big( g(u); \<1p>_\eps \big)\;,
    \end{split}
\end{equation}
and
\begin{equation} \label{eq:PAM_w}
\begin{split}
    \w_\eps(\vec{u}) :=& -\rR_\eps^\e \lL_\eps \Big( \gG_{\eps}(\vec{u}) - g(u) \prec_\eps e^{t \lL_\eps} \<1p>_\eps \Big)+\Err_\eps^{(2)} \big( g(u); \<1p>_\eps \big)\;.
    \end{split}
\end{equation}
Here, we recall from \eqref{eq:Green_remainder_elliptic} that $\rR_\eps^\e = - \big( \nabla \lL_\eps^{-1} - \Phi_\eps \nabla \lL_0^{-1} \big)$. We use boldface letters $\v_\eps$ and $\w_\eps$ to indicate that they are vector-valued. With the above decomposition of $\nabla u_\eps$, we define
\begin{equation} \label{eq:PAM_K2}
    \begin{split}
    \kK_\eps^{(2)}(\vec{u}) &:= g^2(u) \, \sum_{i,j} \big(F_\eps \odot \<1pdF1pd0>_{\eps})_{ij}  + 2 g(u) \, \big( \v_\eps^\sT(\vec{u}) \cdot \<1pdFP>_\eps + \w_\eps^\sT(\vec{u}) F_\eps \cdot \<1pd>_\eps \big)\\
    &\phantom{11111}+ \big( \v_\eps^\sT(\vec{u}) \Phi_\eps^\sT + \w_\eps^\sT(\vec{u}) \big) F_\eps \big( \Phi_\eps \v_\eps(\vec{u}) + \w_\eps(\vec{u}) \big)\;.
    \end{split}
\end{equation}
The forms of $\v_0$, $\w_0$ and $\kK_0^{(2)}$ are not just the same by putting $\eps = 0$. In fact, we have
\begin{equation*}
\v_0 (\vec{u}) := \nabla \big( \gG_0 (\vec{u}) - g(u) \prec e^{t \lL_0} \<1p>_0 \big) + \Lambda_0 (u)\;, \quad \w_0 \equiv 0\;,
\end{equation*}
and
\begin{equation} \label{eq:PAM_K2_limit}
\kK_0^{(2)}(\vec{u}) = g^2(u) \, \sum_{i,j} \big( \bar{F} \odot \<1pdF1pd0>_0 \big)_{i,j} + 2 \, g(u) \, \v_0^\sT(\vec{u}) \, \bar{F} \, \<1pd>_0 + \v_0^\sT(\vec{u}) \, \bar{F} \, \v_0(\vec{u})\;.
\end{equation}
Now, with $\kK_\eps^{(j)}$ for $j=0,1,2$ and $\eps \in [0,1]$ given in \eqref{eq:PAM_K0}, \eqref{eq:PAM_K1}, \eqref{eq:PAM_K2} and \eqref{eq:PAM_K2_limit}, we define
\begin{equation} \label{eq:PAM_K}
    \kK_\eps (\vec{u}) := \kK_\eps^{(0)}(u) + g'(u) \cdot \kK_\eps^{(1)}(\vec{u}) + \kK_\eps^{(2)}(\vec{u})\;.
\end{equation}
Thus, we are led to the fixed point system
\begin{equation} \label{eq:PAM_system}
    \begin{split}
    u_\eps &= e^{t \lL_\eps} u_\eps(0) + \iI_\eps \big( g(u_\eps) \prec_\eps \xi \big) + \iI_\eps \big( \kK_\eps (\vec{u}_\eps) \big)\;,\\
    u_\eps^{\#} &= e^{t \lL_\eps} u_\eps^{\#}(0) + \iI_\eps \big( \kK_\eps(\vec{u}_\eps) \big)
    \end{split}
\end{equation}
for the pair $\vec{u}_\eps = (u_\eps, u_\eps^{\#})$, where $\kK_\eps^{(0)}$, $\kK_\eps^{(1)}$ and $\kK_\eps^{(2)}$ are given by \eqref{eq:PAM_K0}, \eqref{eq:PAM_K1} and \eqref{eq:PAM_K2}. Also the functionals $\Lambda_\eps$ and $\gG_\eps$ appearing in the expressions for $\kK_\eps^{(j)}$ are given by \eqref{eq:PAM_Lambda} and \eqref{eq:PAM_gG} respectively. 

Let $\alpha = 1-\eta$. Define the spaces $\xX_T^{(1)}$ and $\xX_{T,\eps}^{(2)}$ by
\begin{equation} \label{eq:PAM_space}
    \xX_T^{(1)} = \cC_\fs^\alpha \big( [0,T] \times \T^2 \big)\;, \qquad \xX_{T,\eps}^{(2)} = L_{T}^{\infty} \cC_x^{\alpha+2\eta; \eps}\;,
\end{equation}
and set
\begin{equation*}
    \xX_{T,\eps} := \xX_{T}^{(1)} \times \xX_{T,\eps}^{(2)}\;.
\end{equation*}
This is the space we will find the solution $\vec{u}_\eps = (u_\eps, u_\eps^\#)$ to our fixed point system \eqref{eq:PAM_system}.

\subsection{Preliminary bounds}

\begin{lem} \label{lem:PAM_noise_Self}
We have
\begin{equation*}
    |\!|\!| \xi; \xi |\!|\!|_{-1-\kappa; \eps} \lesssim \| \<1p>_\eps - \<1p>_0 \|_{1-\kappa}\;.
\end{equation*}
\end{lem}
\begin{proof}
By definition of $\<1p>_\eps$ and $\<1p>_0$, we have
\begin{equation*}
    \xi = \Pi_0 \xi - \lL_\eps \<1p>_\eps = \Pi_0 \xi - \lL_0 \<1p>_0\;.
\end{equation*}
The claim then follows directly from definition of the norm $|\!|\!| \cdot; \, \cdot |\!|\!|_{\beta}$ for $\beta \in (-2,-1)$. 
\end{proof}

\begin{lem} \label{lem:PAM_g_composite}
Let $\beta \in (0,1)$. We have the bounds
\begin{equation*}
    \begin{split}
    \|g(u_\eps)\|_{\cC_\fs^\beta} &\lesssim 1 + \|u_\eps\|_{\cC_\fs^\beta}\;,\\
    \|g(u_\eps) - g(u_0)\|_{\cC_\fs^\beta} &\lesssim \big( 1 + \|u_0\|_{\cC_\fs^\beta} \big) \, \|u_\eps - u_0\|_{\cC_\fs^\beta}\;.
    \end{split}
\end{equation*}
\end{lem}
\begin{proof}
This is identical to Lemma~\ref{lem:Holder_composition}, or \cite[Lemma~3.5]{homo_gPAM}. We omit the details. 
\end{proof}

\begin{lem} \label{lem:PAM_Lambda}
Recall the definition of $\Lambda_\eps$ in \eqref{eq:PAM_Lambda}. We have
\begin{equation*}
    \|\Lambda_\eps (u_\eps)\|_{L_T^\infty \cC_x^{\alpha-\kappa}} \lesssim \|u_\eps\|_{L_T^\infty \cC_x^\alpha} \|\<1p>_\eps\|_{1-\kappa}\;,
\end{equation*}
and
\begin{equation*}
    \begin{split}
    &\phantom{111}\|\Lambda_\eps (u_\eps) - \Lambda_0 (u_0)\|_{L_T^\infty \cC_x^{\alpha-\kappa}}\\
    &\lesssim \big( 1 + \|u_0\|_{L_T^\infty \cC_x^\alpha} \big) \Big( \|\<1p>_\eps\|_{1-\kappa} \|u_\eps - u_0\|_{L_T^\infty \cC_x^\alpha} + \|\<1p>_\eps - \<1p>_0\|_{1-\kappa} \Big)\;.
    \end{split}
\end{equation*}
\end{lem}
\begin{proof}
The claims are immediate from Lemma~\ref{lem:PAM_g_composite} (same true with $\cC_\fs^\beta$ norm replaced by $L_T^\infty \cC_x^\alpha$) and the estimate for the para-product $\prec$ from Proposition~\ref{pr:paraproducts_gBesov}. 
\end{proof}

\begin{lem} \label{lem:PAM_gG}
Recall the definition of $\gG_\eps$ from \eqref{eq:PAM_gG}. We have the bound
\begin{equation*}
    \|\gG_\eps (\vec{u}_\eps)\|_{L_{3\eta,T}^\infty \cC_x^{1+\eta;\eps}} \lesssim \|\vec{u}_\eps\|_{\xX_{T,\eps}} + \big( 1 + \|u_\eps\|_{\xX_T^{(1)}} \big) \|\<1p>_\eps\|_{1-\kappa}\;.
\end{equation*}
Furthermore, there exists $\theta>0$ such that
\begin{equation*}
    \begin{split}
    &\phantom{111}\| \lL_\eps \gG_\eps (\vec{u}_\eps) - \lL_0 \gG_0 (\vec{u}_0) \|_{L_{3\eta,T}^\infty \cC_x^{-1+\eta-2\kappa}}\\
    &\lesssim \eps^\theta \|\vec{u}_0\|_{\xX_{T,0}} + \|u_\eps - u_0\|_{\xX_T^{(1)}} + \|\lL_\eps u_\eps^\# - \lL_0 u_0^\#\|_{L_{3\eta,T}^\infty \cC_x^{-1+\eta-\kappa}}\\
    &\phantom{111}+ \big( 1 + \|u_0\|_{\xX_T^{(1)}} \big) \big( \eps^\theta \|\<1p>_0\|_{1-\kappa} + \|\<1p>_\eps - \<1p>_0\|_{1-\kappa} + \|\<1p>_\eps\|_{1-\kappa} \|u_\eps - u_0\|_{\xX_T^{(1)}} \big)\;.
    \end{split}
\end{equation*}
\end{lem}
\begin{proof}
The bounds follow from Lemmas~\ref{le:convergence_heat_normal}, ~\ref{le:convergence_heat_comm_below} and~\ref{lem:PAM_g_composite}. 
\end{proof}

\begin{lem} \label{lem:PAM_K0}
Recall the definition of $\kK_\eps^{(0)}$ from \eqref{eq:PAM_K0}. We have the bound
\begin{equation*}
    \|\kK_\eps^{(0)}(u_\eps)\|_{L_T^\infty \cC_x^{\alpha-1-\kappa}} \lesssim \big( 1 + \|u_\eps\|_{L_T^\infty \cC_x^\alpha}^2 \big) \,\big( |\Pi_0 \xi| + \|\<1p>_\eps\|_{1-\kappa} \big)\;.
\end{equation*}
Furthermore, there exists $\theta>0$ such that
\begin{equation*}
    \begin{split}
    \| &\kK_\eps^{(0)}(u_\eps) - \kK_0^{(0)} (u_0) \|_{L_T^\infty \cC_x^{\alpha-1-2\kappa}} \lesssim \big( 1 + \|u_0\|_{L_T^\infty \cC_x^\alpha}^4 + \|u_\eps\|_{L_T^\infty \cC_x^\alpha}^4 \big)\\
    &\cdot \big( 1 + |\Pi_0 \xi| + \|\<1p>_\eps\|_{1-\kappa} + \|\<1p>_0\|_{1-\kappa} \big) \cdot \big( \eps^\theta + \|\<1p>_\eps - \<1p>_0\|_{1-\kappa} + \|u_\eps - u_0\|_{L_T^\infty \cC_x^\alpha} \big)\;.
    \end{split}
\end{equation*}
\end{lem}
\begin{proof}
The uniform bound is straightforward. For the difference $\kK_\eps^{(0)} - \kK_0^{(0)}(u_0)$, the desired bound follows from Lemmas~\ref{le:convergence_prec_below}, ~\ref{le:convergence_comm_below}, \ref{lem:para_linearisation_diff} and~\ref{lem:homo_diff_circ_estimate}. 
\end{proof}

\begin{lem} \label{lem:PAM_K1}
Recall the definition of $\kK_\eps^{(1)}$ from \eqref{eq:PAM_K1}. We have the bound
\begin{equation*}
    \|\kK_\eps^{(1)}(\vec{u}_\eps)\|_{L_{3\eta,T}^\infty \cC_x^{-\kappa}} \lesssim \Big( \big(1 + \|u_\eps\|_{\xX_T^{(1)}} \big) \big( 1 + \|\Upsilon_\eps\|_{\yY} \big) + \|\vec{u}_\eps\|_{\xX_{T,\eps}} \Big) \, \|\Upsilon_\eps\|_{\yY}\;.
\end{equation*}
Furthermore, there exists $\theta>0$ such that
\begin{equation*}
    \begin{split}
    \|\kK_\eps^{(1)} (\vec{u}_\eps) - &\kK_0^{(1)}(\vec{u}_0)\|_{L_{3\eta,T}^\infty \cC_x^{-2\kappa}} \lesssim \big( 1 + \|\vec{u}_\eps\|_{\xX_{T,\eps}} + \|\vec{u}_0\|_{\xX_{T,0}} \big) \big( 1 + \|\Upsilon_0\|_{\yY} + \|\Upsilon_\eps\|_{\yY} \big)^4\\
    &\Big( \eps^\theta + \|\Upsilon_\eps - \bar{F} \bullet \Upsilon_0\|_{\yY} + \|u_\eps - u_0\|_{\xX_T^{(1)}} + \|\lL_\eps u_\eps^\# - \lL_0 u_0^\# \|_{L_{3\eta,T}^\infty \cC_x^{-1+\eta}} \Big)\;.
    \end{split}
\end{equation*}
\end{lem}
\begin{proof}
By the expression \eqref{eq:PAM_K1}, $\kK_\eps^{(1)}(\vec{u}_\eps) - \kK_0^{(1)}(\vec{u}_0)$ is the sum of the differences of three pairs. The bounds for the first two follow directly from Lemmas~\ref{lem:homo_diff_circ_estimate}, ~\ref{le:convergence_comm_below} and~\ref{le:convergence_heat_normal}. For the last one, again by Lemma~\ref{le:convergence_comm_below}, we have
\begin{equation*}
    \begin{split}
    \big\| \gG_\eps (\vec{u}_\eps) &\circ_\eps \xi - \gG_0 (\vec{u}_0) \circ \xi \big\|_{L_{3\eta}^\infty \cC_x^{\eta-2\kappa}} \lesssim \|\gG_0 (\vec{u}_0)\|_{L_{3\eta,T}^\infty \cC_x^{1+\eta}} \, |\!|\!| \xi; \, \xi |\!|\!|_{-1-\kappa;\eps}\\
    &+\Big( \eps^\theta \|\gG_\eps(\vec{u}_\eps)\|_{L_{3\eta,T}^\infty \cC_x^{1+\eta; \eps}} + \big\| \lL_ \eps \gG_\eps (\vec{u}_\eps) - \lL_0 \gG_0 (\vec{u}_0) \big\|_{L_{3\eta,T}^\infty \cC_x^{1+\eta;\eps}} \Big) \|\xi\|_{-1-\kappa;\eps}
    \end{split}
\end{equation*}
The desired bound then follows from Lemmas~\ref{lem:PAM_noise_Self} and~\ref{lem:PAM_gG}. 
\end{proof}

\begin{lem} \label{lem:PAM_v_w}
There exists $\nu, \theta > 10 \kappa$ such that
\begin{equation*}
    \|\v_\eps (\vec{u}_\eps)\|_{L_{3\eta,T}^\infty \cC_x^\nu} \lesssim \|\vec{u}_\eps\|_{\xX_T^{(\eps)}} + \big( 1 + \|u_\eps\|_{\xX_T^{(1)}} \big) \|\<1p>_\eps\|_{\cC^{-\kappa}}\;,
\end{equation*}
and
\begin{equation*}
    \begin{split}
    \|&\v_\eps (\vec{u}_\eps) - \v (\vec{u}_0) \|_{L_{3\eta,T}^\infty \cC_x^\nu} \lesssim \eps^\theta \|\vec{u}_0\|_{\xX_T^{(0)}} + \|u_\eps - u_0\|_{\xX_T^{(1)}} + \|\lL_\eps u_\eps^\# - \lL_0 u_0^\#\|_{L_T^\infty \cC_x^{\eta-1}}\\
    &+ \big( 1 + \|u_0\|_{\xX_T^{(1)}} + \|u_\eps\|_{\xX_T^{(1)}} \big) \big( 1 + \|\<1p>_0\|_{1-\kappa} + \|\<1p>_\eps\|_{1-\kappa} \big) \Big( \eps^\theta + \|\<1p>_\eps - \<1p>_0\|_{1-\kappa} + \|u_\eps - u_0\|_{\xX_T^{(1)}} \Big)\;,
    \end{split}
\end{equation*}
and
\begin{equation*}
    \| \w_\eps (\vec{u}_\eps) \|_{L_{3\eta,T}^\infty \cC_x^\nu} \lesssim \eps^\theta \Big( \|\vec{u}_\eps\|_{\xX_{T,\eps}} + \big( 1 + \|u_\eps\|_{\xX_T^{(1)}} \big) \, \|\<1p>_\eps\|_{\cC^{1-\kappa}} \Big)\;.
\end{equation*}
\end{lem}
\begin{proof}
The uniform-in-$\eps$ bound for $\v_\eps$ follows from the uniform boundedness parts of Lemmas~\ref{lem:PAM_g_composite}, ~\ref{lem:PAM_gG} and heat kernel estimates. The bound for $\w_\eps$ follows from the same ingredients with an additional application of Lemma~\ref{lem:elliptic_homo_Green_remainder}. 

For the difference $\v_\eps (\vec{u}_\eps) - \v_0 (\vec{u}_0)$, the desired bound follows from the corresponding parts in Lemmas~\ref{lem:PAM_g_composite}, ~\ref{lem:PAM_gG} and~\ref{lem:prec_normal_high}. All these are immediate applications of the mentioned lemmas. So we omit the details. 
\end{proof}

\begin{lem} \label{lem:PAM_K2}
Recall the definition of $\kK_\eps^{(2)}$ and $\kK_0^{(2)}$ from \eqref{eq:PAM_K2} and \eqref{eq:PAM_K2_limit}. We have the bound
\begin{equation*}
    \|\kK_\eps^{(2)}(\vec{u}_\eps)\|_{L_{3\eta,T}^\infty \cC_x^{-\kappa}} \lesssim \big(1 + \|u_\eps\|_{\xX_{T}^{(1)}}^2\big) (1 + \|\vec{u}_\eps\|_{\xX_{T,\eps}}) \big( 1 + \|\Upsilon_\eps\|_{\yY} \big)^2\;.
\end{equation*}
Furthermore, for every $\kappa>0$, there exists $\theta>0$ such that
\begin{equation*}
    \begin{split}
    \|\kK_\eps^{(2)}(\vec{u}_\eps) - &\kK_0^{(2)}(\vec{u}_0)\|_{L_{3\eta,T}^\infty \cC_x^{-3\kappa}} \lesssim \big(1 + \|\vec{u}_0\|_{\xX_{T,0}} + \|\vec{u}_\eps\|_{\xX_{T,\eps}} \big)^{5} \big( 1 + \|\Upsilon_0\|_{\yY} + \|\Upsilon_\eps\|_{\yY} \big)^5\\
    &\times \Big( \eps^\theta + \|\Upsilon_\eps - \bar{F} \bullet \Upsilon_0\|_{\yY} + \|u_\eps - u_0\|_{\xX_T^{(1)}} + \|\lL_\eps u_\eps^\# - \lL_0 u_0^\#\|_{L_{3\eta,T}^\infty \cC_x^{-1+\eta}} \Big)\;.
    \end{split}
\end{equation*}
\end{lem}
\begin{proof}
The bounds follow directly from the expressions of $\kK_\eps^{(2)}$ and $\kK_0^{(2)}$ in \eqref{eq:PAM_K2} and \eqref{eq:PAM_K2_limit} as well as the bounds for $\v_\eps$ and $\w_\eps$ from Lemma~\ref{lem:PAM_v_w}. We omit the details. 
\end{proof}

\subsection{Uniform boundedness}

Let
\begin{equation}
    \mM_\eps := \|u_\eps(0)\|_{1-\eta} + \|u_\eps^\#(0)\|_{1+\eta;\eps} + \|\Upsilon_\eps\|_{{\cC^{-\kappa}}}\;.
\end{equation}
We have the following theorem regarding the existence of the solution $\vec{u}_\eps = (u_\eps, u_\eps^\#)$ to \eqref{eq:PAM_system}. 

\begin{thm} \label{thm:PAM_fixed_point}
For every $\eps \in \N^{-1}$, every realisation of $\Upsilon_\eps \in {\cC^{-\kappa}}$ and every pair of initial data $\vec{u}_\eps(0) = \big( u_\eps(0), u_\eps^{\#}(0) \big) \in \cC^{1-\eta} \times \cC^{1+\eta;\eps}$, there exists $\tau_\eps >0$ such that the system \eqref{eq:PAM_system} has a unique solution $\vec{u}_\eps \in \xX_{\tau_\eps, \eps}$. 

Furthermore, the local existence time $\tau_\eps$ and the norm $\|\vec{u}_\eps\|_{\xX_{\tau_\eps,\eps}}$ depends on $(\eps,\delta)$ via $\mM_\eps$ only (recall that $\delta$ is hidden in the notation). In particular, if $\mM_\eps$ is uniformly bounded, then there exists $\tau>0$ independent of $\eps$ such that
\begin{equation} \label{eq:PAM_sol_uniform_bound}
    \|\vec{u}_\eps\|_{\xX_{\tau,\eps}} \lesssim 1 + \sup_\eps \mM_\eps\;,
\end{equation}
where the proportionality constant is independent of $\eps$ (and the hidden $\delta$). 
\end{thm}
\begin{proof}
For $\tau>0$ and $\eps \in \N^{-1}$, define the mild solution map $\Gamma_{\tau,\eps} := \big(\Gamma_{\tau,\eps}^{(1)}, \Gamma_{\tau,\eps}^{(2)} \big)$ by
\begin{equation*}
    \begin{split}
    \Gamma_{\tau,\eps}^{(1)}(\vec{u}) &= e^{t \lL_\eps} u_\eps(0) + \iI_\eps \big( g(u) \prec_\eps \xi \big) + \iI_\eps \big( \kK_\eps (\vec{u}) \big)\;,\\
    \Gamma_{\tau,\eps}^{(2)}(\vec{u}) &= e^{t \lL_\eps} u_\eps^{\#}(0) + \iI_\eps \big( \kK_\eps (\vec{u}) \big)\;,
    \end{split}
\end{equation*}
where $\kK_\eps$ is given in \eqref{eq:PAM_K}. We show that for some $R, \tau > 0$ to be specified later, the map $\Gamma_{\tau,\eps}$ is a contraction from the centered ball in $\xX_{\tau,\eps}$ with radius $R$ into itself, hence yielding a fixed point. 

Let $\vec{u} \in \xX_{\tau,\eps}$. Also recall the spaces $\xX_{\tau}^{(1)}$ and $\xX_{\tau,\eps}^{(2)}$ from \eqref{eq:PAM_space}. First, by Lemma~\ref{le:heat_smoothing_space}, we have
\begin{equation*}
    \|e^{t\lL_\eps} u(0)\|_{\xX_{\tau}^{(1)}} \lesssim \|u_\eps(0)\|_{1-\eta}\;, \qquad \|e^{t \lL_\eps} u_\eps^\#(0)\|_{\xX_{\tau,\eps}^{(2)}} \lesssim \|u_\eps^\#(0)\|_{1+\eta;\eps}\;.
\end{equation*}
By \eqref{e:bony_generalised_low_high} in Proposition~\ref{pr:paraproducts_gBesov} and Lemma~\ref{le:heat_smoothing_spacetime}, we have
\begin{equation*}
    \big\| \iI_\eps \big( g(u) \prec_\eps \xi \big) \big\|_{\xX_\tau^{(1)}} \lesssim \tau^\theta (1+\rR) \mM_\eps\;.
\end{equation*}
Also, by Lemmas~\ref{lem:PAM_K0}, \ref{lem:PAM_K1} and~\ref{lem:PAM_K2}, we have
\begin{equation*}
    \big\| \iI_\eps \big( \kK_\eps (\vec{u}) \big)  \big\|_{\xX_{\tau}^{(1)} \cap \xX_{\tau,\eps}^{(2)}} \lesssim \tau^\theta (1 + \rR)^5 (1 + \mM_\eps)^3\;.
\end{equation*}
Combining the above, we see there exists a universal constant $C>0$ such that
\begin{equation} \label{eq:PAM_sol_map_bound}
    \|\Gamma_{\tau,\eps}(\vec{u})\|_{\xX_{\tau,\eps}} \leq C \big( \mM_\eps + \tau^\theta (1+\rR)^5 (1 + \mM_\eps)^3 \big)\;.
\end{equation}
Hence, as long as one chooses
\begin{equation*}
    \rR > 2 C \mM_\eps \qquad \text{and} \qquad \tau^\theta (1 + \rR)^5 (1+ \mM_\eps)^3 < \frac{C \mM_\eps}{2}\;,
\end{equation*}
the map $\Gamma_{\tau,\eps}$ maps the ball of radius $\rR$ (and centered at the origin) in $\xX_{\tau,\eps}$ into itself. To show it is also a contraction, we repeat the above procedures to get
\begin{equation*}
    \| \Gamma_{\tau,\eps}(\vec{u}_{1}) - \Gamma_{\tau,\eps}(\vec{u}_{2}) \|_{\xX_{\tau,\eps}} \lesssim \tau^\theta \big( 1 + \rR \big)^5 \big( 1 + \mM_\eps \big)^5 \|\vec{u}_{1} - \vec{u}_{2}\|_{\xX_{\tau,\eps}}\;.
\end{equation*}
Fix the above choice of $\rR$, one can further choose a smaller $\tau = \tau_\eps$ to make $\Gamma_{\tau,\eps}$ a contraction. Hence, one gets a unique fixed point solution $\vec{u}_\eps \in \xX_{\tau_\eps}$. 

Furthermore, the above argument also shows that the existence time $\tau_\eps$ depends on $\eps$ via $\mM_\eps$ only. Hence, using \eqref{eq:PAM_sol_map_bound} and that $\vec{u}_\eps$ is a fixed point, we obtain \eqref{eq:PAM_sol_uniform_bound}. This completes the proof. 
\end{proof}

\subsection{Convergences of the solution and flux}

For $\eps \in \N^{-1}$, let $\vec{u}_\eps \in \xX_{\tau_\eps,\eps}$ be the solution to \eqref{eq:PAM_system} from Theorem~\ref{thm:PAM_fixed_point}. We define
\begin{equation*}
    \vec{ U }_\eps := \big( u_\eps, \Pi_0 u_\eps^\#, \lL_\eps u_\eps^\# \big)\;,
\end{equation*}
and the spaces
\begin{equation*}
    \widetilde{\xX}_T := \xX_T^{(1)} \times L_T^\infty \R \times L_{T}^\infty \cC_x^{-1+\eta}\; \quad \widetilde{\xX}_T^{+} := \cC_{\fs}^{1-\frac{\eta}{2}}([0,T] \times \T^2) \times L_T^\infty \R \times L_{T}^\infty \cC_x^{-1+\frac{3\eta}{2}}\;.
\end{equation*}
For $\eps=0$, the homogenised operator $\lL_0 = \div (\bar{A} \nabla)$ is a constant coefficient elliptic operator. It is well known that if $\big(u_0(0), u_0^\#(0) \big) \in \cC^{1-\frac{\eta}{2}} \times \cC^{1+\frac{3\eta}{2}}$, then there exists a random $T>0$ such that
\begin{equation*}
    \vec{U}_0 = \big( u_0, \Pi_0 u^\#(0), \lL_0 u^\#(0) \big) \in \widetilde{\xX}_T^+
\end{equation*}
in the slightly stronger space $\widetilde{\xX}_T^+$. We have the following theorem regarding convergence of the solutions.

\begin{thm} \label{thm:PAM_convergence_homo}
Fix $T>0$ as above (so that $\|\vec{U}_0\|_{\xX_T^+} < +\infty$). Suppose
\begin{equation*}
    \mM := 1 + \|\vec{ U }_{ 0 }\|_{\widetilde{\xX}_T^+} + \sup_{\eps \in \N^{-1}} \big( \|u_\eps(0)\|_{1-\eta} + \|u_\eps^\#(0)\|_{1+\eta;\eps} + \|\Upsilon_\eps\|_{\yY} \big) < +\infty\;.
\end{equation*}
Then for all sufficiently small $\eps \in \N^{-1}$, the solution $\vec{u}_\eps = (u_\eps, u_\eps^\#)$ exists in $\xX_{T,\eps}$ for the same $T$. Furthermore, there exists $\theta>0$ such that
\begin{equation} \label{eq:PAM_convergence_homo}
    \|\vec{ U }_{ \eps } - \vec{ U }_{ 0 }\|_{\widetilde{\xX}_T} \lesssim \eps^\theta + \|u_\eps(0) - u_0(0)\|_{1-\eta} + |\!|\!| u_\eps^\#(0); u_0^\#(0) |\!|\!|_{1+\eta;\eps} + \|\Upsilon_\eps - \bar{F} \bullet \Upsilon_0\|_{\yY}\;.
\end{equation}
The proportionality constant depends on $\mM$ and $T$ only. 
\end{thm}
\begin{proof}
As in the case of dynamical $\phi^4_3$, we first prove \eqref{eq:PAM_convergence_homo} for sufficiently small positive time $\tau>0$ (independent of $\eps$), and then iterate it to the whole time interval $[0,T]$ where the limiting solution $\vec{ U }_{ 0 }$ exists. 

By Theorem~\ref{thm:PAM_fixed_point}, there exists $\tau>0$ independent of $\eps$ such that
\begin{equation*}
    \|\vec{ U }_{ \eps }\|_{\widetilde{\xX}_{\tau}} \lesssim \mM\;.
\end{equation*}
To control difference $\vec{ U }_{ \eps } - \vec{ U }_{ 0 }$, we first note that by Lemmas~\ref{le:convergence_heat_normal}, we have
\begin{equation*}
    \begin{split}
    \|e^{t \lL_\eps} u_\eps(0) - e^{t \lL_0} u_0(0)\|_{\cC_\fs^{1-\eta}} &\leq \|e^{t \lL_\eps} (u_\eps(0) - u_0(0))\|_{\cC_\fs^{1-\eta}} + \|\big(e^{t \lL_\eps} - e^{t \lL_0} \big) u_0(0) \|_{\cC_\fs^{1-\eta}}\\
    &\lesssim \|u_\eps(0) - u_0(0)\|_{1-\eta} + \eps^\theta \|u_0\|_{1-\frac{\eta}{2}}\;.
    \end{split}
\end{equation*}
Here, we get a positive power $\eps^\theta$ for $\eps$ with the norm $\cC_\fs^{1-\eta}$ since we assumed $u_0(0)$ is in the higher regularity space $\cC^{1-\frac{\eta}{2}}$. Similarly, we also have
\begin{equation*}
    \begin{split}
    &\phantom{111}\|\lL_\eps e^{t \lL_\eps} u_\eps^{\#}(0) - \lL_0 e^{t \lL_0} u_0^\#(0)\|_{L_\tau^\infty \cC_x^{-1+\eta}}\\
    &\lesssim \|\lL_\eps u_\eps^\#(0) - \lL_0 u_0^\#(0)\|_{-1+\eta} + \eps^\theta \|\lL_0 u_0^\#(0)\|_{-1+\frac{3\eta}{2}}\;.
    \end{split}
\end{equation*}
The reason we get a factor $\eps^\theta$ for the regularity space $ \cC^{ - 1 + \eta } $ is that we assume $\lL_0 u_0^\#(0)$ lives in the higher regularity space $\cC^{-1+\frac{3\eta}{2}}$. 

Also, by Lemma~\ref{le:convergence_prec_heat_below}, we have (noting that $\Pi_0^\perp \xi = \lL_\eps \<1p>_\eps = \lL_0 \<1p>_0$)
\begin{equation*}
    \begin{split}
    &\phantom{111}\big\| \iI_\eps \big( g(u_\eps) \prec_\eps \xi \big) - \iI_0 \big( g(u_0) \prec \xi \big) \big\|_{\cC_\fs^{1-\eta}([0,\tau] \times \T^2)}\\
    &\lesssim \tau^\theta \mM^3 \big( \eps^\theta + \|u_\eps - u_0\|_{L_\tau^\infty \cC_x^{1-\eta}} + \|\<1p>_\eps - \<1p>_0\|_{1-\kappa} \big)
    \end{split}
\end{equation*}
for some $\theta>0$. Finally, combining Lemma~\ref{le:convergence_heat_spacetime_normal} together with Lemmas~\ref{lem:PAM_K0}, \ref{lem:PAM_K1} and~\ref{lem:PAM_K2}, we have
\begin{equation*}
    \begin{split}
    &\phantom{111}\big\| \iI_\eps \big( \kK_\eps (\vec{u}_\eps) \big) - \iI_0 \big( \kK_0 (\vec{u}_0) \big) \big\|_{\cC_\fs^{1-\eta}([0,\tau] \times \T^2)} + \big\| \lL_\eps \iI_\eps \big( \kK_\eps (\vec{u}_\eps) \big) - \lL_0 \iI_0 \big( \kK_0 (\vec{u}_0) \big) \big\|_{L_{3\eta,\tau}^\infty \cC_x^{-1+\eta}}\\
    &\lesssim \tau^\theta \Big( \eps^\theta \|\kK_0 (\vec{u}_0)\|_{L_{3\eta,\tau}^\infty \cC_x^{-3\kappa}} + \|\kK_\eps (\vec{u}_\eps) - \kK_0 (\vec{u}_0)\|_{L_{3\eta,\tau}^\infty \cC_x^{-3\kappa}} \Big)\\
    &\lesssim \tau^\theta \mM^{10} \big( \eps^\theta + \|\vec{ U }_{ \eps } - \vec{U}_0\|_{\widetilde{\xX}_\tau} + \|\Upsilon_\eps - \bar{F} \bullet \Upsilon_0\|_{\yY} \big)\;,
    \end{split}
\end{equation*}
The relevant bound for $|\Pi_0 (u_\eps^\# - u_0^\#)|$ is straightforward and hence we omit the details. Combining all the above, we conclude there exists $C > 0$ and $\theta>0$ such that
\begin{equation*}
    \begin{split}
    \|\vec{ U }_{ \eps } - \vec{U}_0\|_{\widetilde{\xX}_\tau} \leq C \mM^{10} \big(&\tau^\theta \|\vec{ U }_{ \eps } - \vec{U}_0\|_{\widetilde{\xX}_\tau} + \eps^\theta + \|u_\eps(0) - u_0(0)\|_{1-\eta}\\
    &+ |\!|\!| u_\eps^\#(0); u_0^\#(0) |\!|\!|_{1+\eta;\eps} + \|\Upsilon_\eps - \bar{F} \bullet \Upsilon_0\|_{\yY}  \big)\;,
    \end{split}
\end{equation*}
which implies
\begin{equation*}
    \|\vec{ U }_{ \eps } - \vec{U}_0\|_{\widetilde{\xX}_T} \lesssim \eps^\theta + \|u_\eps(0) - u_0(0)\|_{1-\eta} + |\!|\!| u_\eps^\#(0); u_0^\#(0) |\!|\!|_{1+\eta;\eps} + \|\Upsilon_\eps - \bar{F} \bullet \Upsilon_0\|_{\yY}\;.
\end{equation*}
if $\tau$ is sufficiently small (independent of $\eps$). For sufficiently small $\eps$, the assumption on the convergence of initial data further implies
\begin{equation*}
    \|u_\eps(\tau) - u_0(\tau)\|_{1-\eta} + \big| \Pi_0 \big( u_\eps^\# (\tau) - u_0^\#(\tau) \big) \big| + \| \lL_\eps u_\eps^\#(\tau) - \lL_0 u_0^\#(\tau) \|_{-1+\eta} < \frac{1}{2}\;.
\end{equation*}
Since the time $\tau>0$ does not depend on $\eps$, one can iterate this procedure to the whole time interval $[0,T]$ and conclude \eqref{eq:PAM_convergence_homo}. 
\end{proof}

We next deal with convergence of the flux. 

\begin{thm} \label{thm:PAM_convergence_flux}
Let $(u_\eps, u_\eps^\#)$ be the solution to the fixed point system \eqref{eq:PAM_system}. Under the above assumptions, we have
\begin{equation*}
    \big\| A_\eps \nabla u_\eps - \bar{A} \nabla u_0  \big\|_{L_T^\infty \cC_x^{-2\eta}} \lesssim \eps^\theta + \|u_\eps(0) - u_0(0)\|_{1-\eta} + \|\Upsilon_\eps - \bar{F} \bullet \Upsilon_0\|_{\yY} + \| A_{ \eps } \nabla \<1p>_{ \eps } - \bar{ A } \nabla \<1p>_{ 0 } \|_{ - \eta } \; .
\end{equation*}
The proportionality constant depends on $\mM$ and $T$ only. 
\end{thm}
\begin{proof}
By \eqref{eq:nabla_u}, we have
\begin{equation*}
    A_{ \eps } \nabla u_{ \eps } = A_{ \eps } \nabla \<1p>_{ \eps } \, \cdot \, g ( u_{ \eps } ) + A_{ \eps } \Phi_{ \eps } \, \cdot \, v_{ \eps } ( \vec{ u }_{ \eps } ) + A_{ \eps } \, \cdot \, w_{ \eps } ( \vec{ u }_{ \eps } ) 
\end{equation*}
where $ v_{ \eps } $, $ w_{ \eps } $  are given in \eqref{eq:PAM_v}, \eqref{eq:PAM_w}. Lemma \ref{lem:PAM_K2} and Theorem \ref{thm:PAM_convergence_homo} give quantitative convergences of $  v_{ \eps } ( \vec{ u }_{ \eps } ) $ and $ w_{ \eps } ( \vec{ u }_{ \eps } )  $ to $  v_{ 0 } ( \vec{ u }_{ 0 } ) $ and $ 0 $ in $ \cC^{ 0 + } $ respectively. 
Using that $ \| A_{ \eps } \Phi_{ \eps } - \bar{ A } \|_{ - \kappa } \lesssim \eps^{ \theta } $ and the convergence results above, we can deduce the convergence of the three terms in the expression for $ A_{ \eps } \nabla u_{ \eps } - \bar{ A } \nabla u_{ 0 }$. This completes the proof of the theorem. 
\end{proof}

\section{The KPZ equation}
\label{sec:kpz}

In this section, we consider the periodic homogenisation problem for KPZ equation on the one dimensional torus with oscillatory coupling coefficient $g_\eps = g(\cdot / \eps)$. For notational simplicity, we write $\b := \frac{g}{A}$ and $\b_\eps := \b (\cdot / \eps)$. Hence, the equation we consider is
\begin{equation} \label{e:kpz_homo_full}
    \d_t h_\eps^{(\delta)} = \lL_\eps h_\eps^{(\delta)} + A_\eps \b_\eps \big( \d_x h_\eps^{(\delta)} \big)^2 + \xi^{(\delta)} - C_\eps^{(\delta)} -  4 \widetilde{C}_\eps^{(\delta)} \cdot A_{\eps} \d_{x}h_{\eps}^{(\delta)}\;.
\end{equation}
The candidate for the limiting equation is
\begin{equation} \label{e:kpz_limiting_equation}
    \d_t h_0 = \bar{A} \d_x^2 h_0 + \lambda (\d_x h_0)^2 + \xi - \infty\;,
\end{equation}
where $\bar{A} = \big( \Pi_0 (A^{-1}) \big)^{-1}$, and the coupling constant $\lambda$ is given by
\begin{equation} \label{e:kpz_coupling_const}
    \lambda = \bar{A}^2 \int_\T \frac{g(x)}{A^2(x)} {\rm d}x = \bar{A}^2 \int_\T \frac{\textbf{b}(x)}{A(x)} {\rm d}x\;.
\end{equation}
To prove uniform boundedness in $\eps$ and convergence, the main new analytic tools, in addition to the general framework developed in Sections~\ref{sec:para_general} and~\ref{sec:para_convergence}, are the bounds related to derivatives specific to $d=1$ given in Section~\ref{sec:bounds_derivatives_1d}. 

For $\delta \in (0,1)$, let $\xi^{(\delta)}$ be a smooth approximation of $\xi$ at scale $\delta$. In the rest of the section, we fix the parameters
\begin{equation} \label{eq:kpz_parameters}
    \eta \in \Big(0, \frac{1}{1000} \Big)\;, \qquad \kappa \in \Big(0, \frac{\eta}{1000} \Big)\;, \qquad \sigma := \frac{1}{2} + 10 \eta\;.
\end{equation}

\subsection{Definition of stochastic objects and main result}
\label{sec:kpz_stochastic}

Let $\iI_\eps$ and $\I_\eps$ denote
\begin{equation*}
    (\iI_\eps f)(t) := \int_{0}^{t} e^{(t-r) \lL_\eps} f(r) {\rm d}r\;, \qquad (\I_\eps f)(t) := \int_{-\infty}^{t} e^{(t-r) \lL_\eps} \Pi_0^\perp f(r) {\rm d}r
\end{equation*}
respectively. As usual, we introduce graphic notations to represent stochastic objects, but this time more complicated than the standard KPZ situation since they involve not only spatial derivatives but also multiplications by $A_\eps$ and $\b_\eps$. We use a green edge to denote the operation $\I_\eps$, a blue edge to denote the operation ${\b}_{\eps} \cdot \d_x \I_\eps$, and a red edge to denote the operation $A_\eps \cdot \d_x \I_\eps$. 

Let $\xi^{(\delta)}$ be a smooth approximation to $\xi$ at scale $\delta$. With the above graphic notations, the ``free fields" $\<1gk>^{(\delta)}_\eps$, $\<1bk>^{(\delta)}_\eps$ and $\<1rk>^{(\delta)}_\eps$ are respectively given by
\begin{equation} \label{e:kpz_stochastic_ff}
    \<1gk>_\eps^{(\delta)}(t) := \int_{-\infty}^{t} e^{(t-r) \lL_\eps} \Pi_0^\perp \xi^{(\delta)}(r) {\rm d}r\;, \qquad \<1bk>_\eps^{(\delta)} := \b_\eps \d_x \, \<1gk>_\eps^{(\delta)}\;, \qquad \<1rk>_\eps^{(\delta)} :=  A_\eps \d_x \, \<1gk>_\eps^{(\delta)}\;.
\end{equation}
Next, define the objects $\<2k>_\eps^{(\delta)}$, $\<20gk>_\eps^{(\delta)}$, $\<20bk>_\eps^{(\delta)}$ and $\<20rk>_\eps^{(\delta)}$ by
\begin{equation} \label{e:kpz_stochastic_square}
    \begin{split}
    \<2k>_\eps^{(\delta)} := \<1rk>_\eps^{(\delta)} \cdot \<1bk>_\eps^{(\delta)} - C_{\eps,1}^{(\delta)}\;, \quad \<20gk>_\eps^{(\delta)} (t) := (\I_\eps \<2k>_\eps^{(\delta)})(t) = \int_{-\infty}^{t} e^{(t-r) \lL_\eps} \Pi_0^\perp \<2k>_\eps^{(\delta)} (r) {\rm d}r\;
    \end{split}
\end{equation}
    and
    \begin{equation*}
    \begin{split}
    \<20bk>_\eps^{(\delta)} := \b_\eps \d_x \, \<20gk>_\eps^{(\delta)}\;, \qquad \<20rk>_\eps^{(\delta)} := A_\eps \d_x \, \<20gk>_\eps^{(\delta)}\;. 
    \end{split}
\end{equation*}
We also define two other terms $\<10rk>_{\eps}^{(\delta)}$ and $\<101k>_{\eps}^{(\delta)}$ by
\begin{equation} \label{e:kpz_stochastic_12}
   \<10rk>_{\eps}^{(\delta)} := A_{\eps} \d_{x} \I_{\eps}(\<1bk>_{\eps}^{(\delta)})\;, \qquad \<101k>_{\eps}^{(\delta)} := \<10rk>_{\eps}^{(\delta)} \circ_{\eps} \<1bk>_{\eps}^{(\delta)} - \widetilde{C}_\eps^{(\delta)}\;.
\end{equation}
These are all the objects with one or two occurrences of the white noise. There are two terms with three occurrences, $\<21k>_\eps^{(\delta)}$ and $\<210rk>_\eps^{(\delta)}$, respectively given by
\begin{equation*}
    \<21k>_\eps^{(\delta)} := \<20rk>_\eps^{(\delta)} \cdot \<1bk>_\eps^{(\delta)} - 2 \widetilde{C}_\eps^{(\delta)} \<1rk>_\eps^{(\delta)}\;, \quad \<210rk>_\eps^{(\delta)}(t) := A_{\eps}\d_{x}\I_\eps (\<21k>_\eps^{(\delta)})(t)
\end{equation*}
where $\widetilde{C}_\eps^{(\delta)}$ is the same as in \eqref{e:kpz_stochastic_12}. We introduce two final objects $\<211k>_\eps^{(\delta)}$ and $\<22k>_\eps^{(\delta)}$ by
\begin{equation*}
    \<211k>_\eps^{(\delta)} := \<210rk>_\eps^{(\delta)} \circ_\eps \<1bk>_\eps^{(\delta)} - \widetilde{C}_\eps^{(\delta)} \cdot \<20rk>_{\eps}^{(\delta)} - C_{\eps,2}^{(\delta)}\;, \qquad \<22k>_\eps^{(\delta)} := \<20rk>_\eps^{(\delta)} \cdot \<20bk>_\eps^{(\delta)} -  C_{\eps,3}^{(\delta)}\;.
\end{equation*}
The renormalisation function $\widetilde{C}_\eps^{(\delta)}$ appearing in \eqref{e:kpz_homo_full} is precisely the one in \eqref{e:kpz_stochastic_12}, while $C_\eps^{(\delta)}$ is
\begin{equation*}
    C_{\eps}^{(\delta)} := C_{\eps,1}^{(\delta)} + 4 C_{\eps,2}^{(\delta)} + C_{\eps,3}^{(\delta)}\;.
\end{equation*}
These renormalisation functions depend on $x \in \T^1$ but are independent of time. Now we re-organise the objects as follows. Define the collection of random objects
\begin{equation} \label{e:kpz_stochastic_collection}
    \Upsilon_\eps^{(\delta)} := \big( \; \<1bk>_\eps^{(\delta)}, \; \<20bk>_\eps^{(\delta)}, \; \<101k>_\eps^{(\delta)}, \; \<211k>_\eps^{(\delta)}, \; \<22k>_\eps^{(\delta)}, \; \<10rk>_\eps^{(\delta)}, \; \<210rk>_\eps^{(\delta)}  \, \big)\;.
\end{equation}
Define two additional collection of objects $\hH_\eps^{(\delta)}$ and $\fF_\eps^{(\delta)}$ by
\begin{equation*}
    \hH_\eps^{(\delta)} := \Big( \; \<1gk>_\eps^{(\delta)}, \<20gk>_\eps^{(\delta)},\; \int_0^{\cdot} \Pi_0 \big( \xi^{(\delta)} + \<2k>_\eps^{(\delta)} + \<21k>_\eps^{(\delta)} \big)(r) {\rm d}r \Big)\;, \qquad \fF_\eps^{(\delta)} := \big( \, \<1rk>_\eps^{(\delta)}, \, \<20rk>_\eps^{(\delta)} \, \big)\;.
\end{equation*}
We define the spaces $\yY_T$, $\yY_T^{\text{h}}$ and $\yY_T^{\flux}$ by
\begin{equation*}
    \begin{split}
    \yY_T &:= \cC_T^{-\kappa} \cC_x^{-\frac{1}{2}-\kappa} \times L_T^\infty \cC_x^{-\kappa} \times \cC_T^{-\kappa} \cC_x^{-\kappa} \times \cC_T^{-\kappa} \cC_x^{-\kappa} \times L_T^\infty \cC_x^{-\kappa} \times \cC_T^{2\kappa} \cC_x^{\frac{1}{2}-5\kappa} \times \cC_{T}^{2\eta} \cC_x^{\frac{1}{2}-5\eta}\;,\\
    \yY_T^{\text{h}} &:= L_T^\infty \cC_x^{\frac{1}{2}-\kappa} \times L_T^\infty \cC_x^{1-\kappa} \times L_T^\infty \R\;, \qquad \yY_T^{\flux} := \cC_T^{-\kappa} \cC_x^{-\frac{1}{2}-\kappa} \times L_T^\infty \cC_x^{-\kappa}\;.
    \end{split}
\end{equation*}
Let
\begin{equation}
    \Upsilon_0 := \big( \; \<1bks>_0, \; \<20bks>_0,\; \<101ks>_0, \; \<211ks>_0, \; \<22ks>_0,  \; \<10bks>_0, \; , \<210bks>_0\big)
\end{equation}
be a collection of stochastic objects for the standard KPZ equation with operator $\lL_0 = \bar{A} \d_x^2$. Let $\hH_0$ and $\fF_0$ be
\begin{equation*}
    \hH_0 := \Big( \; \I_0 (\xi), \; \I_0 (\<2ks>_0), \;   \int_0^{\cdot} \Pi_0 \big( \xi + \<2ks>_0 + \<21ks>_0 \big)(r) {\rm d}r \Big)\;, \quad \fF_0 := \big( \, \<1bks>_0, \; \, \<20bks>_0 \big)\;.
\end{equation*}
For $c \in \R$, define the actions $c \bullet \Upsilon_0$, $c \bullet \hH_0$ and $c \bullet \fF_0$ by
\begin{equation*}
    c \bullet \Upsilon_0 := \Big( \; \frac{c}{\bar{A}} \cdot \<1bks>_0, \; \frac{c^2}{\bar{A}} \cdot \<20bks>_0,\; \frac{c^2}{\bar{A}} \cdot \<101ks>_0, \; c^3 \, \<211ks>_0, \; c^3 \, \<22ks>_0,  \; c \, \<10bks>_0, \;  c^2 \, \bar{A} \, \<210bks>_0 \; \Big)
\end{equation*}
and
\begin{equation*}
    c \bullet \hH_0 := \Big( \; \I_0 (\xi), \; c \, \I_0 (\<2ks>_0), \; \int_0^{\cdot} \Pi_0 \big( \xi + c \, \<2ks>_0 + c^2 \, \<21ks>_0 \big)(r) {\rm d}r  \Big)\;, \quad c \bullet \fF_0 := \big( \bar{A} \, \<1bks>_0, \; c \bar{A} \, \<20bks>_0 \big)\;.
\end{equation*}
The following proposition from \cite{homo_stochastic} gives convergence of the stochastic objects $\Upsilon_\eps^{(\delta)}$, $\hH_\eps^{(\delta)}$ and $\fF_\eps^{(\delta)}$.

\begin{prop} [\cite{homo_stochastic}, Theorem~5.1]
\label{prop:kpz_stochastic}
For the choices
\begin{equation*}
    \begin{split}
    C_{\eps,1}^{(\delta)} &= \E \big( \<1rk>_\eps^{(\delta)} \cdot \<1bk>_\eps^{(\delta)} \big)\;, \quad \widetilde{C}_\eps^{(\delta)} = \E \big( \<10rk>_\eps^{(\delta)} \cdot \<1bk>_\eps^{(\delta)} \big)\;,\\
    C_{\eps,2}^{(\delta)} &= \E \big( \<210rk>_\eps^{(\delta)} \cdot \<1bk>_\eps^{(\delta)} \big)\;, \quad  C_{\eps,3}^{(\delta)} = \E \big( \<20rk>_\eps^{(\delta)} \cdot \<20bk>_\eps^{(\delta)} \big)\;,
    \end{split}
\end{equation*}
we have
\begin{equation*}
    \|\Upsilon_\eps^{(\delta)} - \lambda \bullet \Upsilon_0\|_{\yY_T} + \|\hH_\eps^{(\delta)} - \lambda \bullet \hH_0\|_{\yY_T^{\text{h}}} + \|\fF_\eps^{(\delta)} - \lambda \bullet \fF_0\|_{\yY_T^\flux} \rightarrow 0
\end{equation*}
as $(\eps, \delta) \rightarrow (0,0)$. Here, $\lambda$ is the constant given by \eqref{e:kpz_coupling_const}. Furthermore, the limit $\Upsilon_0$ satisfies the stronger bound
\begin{equation*}
    \|\Upsilon_0\|_{\yY_T^+} < +\infty\;,
\end{equation*}
where $\yY_T^+$ is the same as $\yY_T$ except that one replaces all appearances of $\cC_T^{-\kappa}$ by $L_T^\infty$. 
\end{prop}

The reason we separate the stochastic objects into collections $\Upsilon_\eps^{(\delta)}$, $\hH_\eps^{(\delta)}$ and $\fF_\eps^{(\delta)}$ is the following. Later in \eqref{eq:kpz_dpd}, we will decompose the solution $h_\eps^{(\delta)}$ to
\begin{equation*}
    h_\eps^{(\delta)} = \<1gk>_\eps^{(\delta)} + \<20gk>_\eps^{(\delta)} + 2 \<210gk>_\eps^{(\delta)} + \int_0^t \Pi_0 \big( \xi^{(\delta)} + \<2k>_\eps^{(\delta)} + \<21k>_\eps^{(\delta)} \big)(r) {\rm d}r + u_\eps^{(\delta)}\;.
\end{equation*}
The objects in $\Upsilon_\eps^{(\delta)}$ are precisely those governing the convergence of the remainder $u_\eps^{(\delta)}$ and its flux. The objects from $\hH_\eps^{(\delta)}$ and $\fF_\eps^{(\delta)}$ give additional information encoded by $h_\eps^{(\delta)}$ and its flux. The object $\<210gk>_\eps^{(\delta)}$ is in the decomposition of $h_\eps^{(\delta)}$ but does not appear in the collection $\hH_\eps^{(\delta)}$ since it can be determined from $\<210rk>_\eps^{(\delta)}$ (in $\Upsilon_\eps^{(\delta)}$) by
\begin{equation*}
    \<210gk>_\eps^{(\delta)} := \lL_\eps^{-1} \d_x \<210rk>_\eps^{(\delta)}\;.
\end{equation*}
For readers' convenience, we list the objects from $\Upsilon_\eps^{(\delta)}$ together with their limits and regularities in the table below. 
{\renewcommand{\arraystretch}{1.2}
\begin{center}
\begin{tabular}{p{1.3cm}  p{1.4cm}  p{1.7cm} p{1.7cm} p{1.2cm} p{1.2cm} p{1.6cm} p{1.5cm}}
\hline
Process & $\<1bk>_\eps^{(\delta)}$ &  $ \<20bk>_\eps^{(\delta)}$ & $ \<101k>_\eps^{(\delta)}$ & $\<211k>_\eps^{(\delta)}$ & $\<22k>_\eps^{(\delta)}$ &  $\<10rk>_\eps^{(\delta)}$ &  $\<210rk>_\eps^{(\delta)}$\\
\hline
Limit &  $(\lambda / \bar{A}) \, \<1bks>_0$ & $(\lambda^2 / \bar{A}) \, \<20bks>_0$ & $(\lambda^2 / \bar{A}) \, \<101ks>_0$ & $\lambda^3 \, \<211ks>_0$ &$\lambda^2 \, \<22ks>_0$ &  $\lambda \, \<10bks>_0$ &  $\lambda^2 \bar{A} \, \<210bks>_0$\\
\hline
Reg. & $\cC_T^{0-}\cC_x^{-\frac{1}{2}-}$ & $0-$ & $\cC_T^{0-}\cC_x^{0-}$ & $\cC_T^{0-}\cC_x^{0-}$ & $0-$ &  $\cC_T^{2\kappa} \cC_x^{\frac{1}{2}- 5\kappa} $ &  $\cC_T^{2\eta} \cC_x^{\frac{1}{2}-5\eta}$\\
\hline
\end{tabular}
\end{center}
Here, $\alpha-$ denotes the space $L_T^\infty \cC_x^{\alpha-\kappa}$, and $\cC_T^{\alpha-} \cC_x^{\beta-}$ denotes $\cC_T^{\alpha-\kappa} \cC_x^{\beta-\kappa}$. The following theorem is the main result regarding the KPZ equation. 
}

\begin{thm} \label{thm:kpz_overall}
Let $C_\eps^{(\delta)}$ and $\widetilde{C}_\eps^{(\delta)}$ be chosen such that Proposition~\ref{prop:kpz_stochastic} holds. Let $\{h_\eps(0)\}_{\eps \in \N^{-1}}$ be a sequence of functions in $\cC^\sigma(\T)$ such that
\begin{equation*}
    \|h_\eps(0) - h_0(0)\|_{\frac{1}{2}-\eta} \rightarrow 0
\end{equation*}
as $\eps \rightarrow 0$ for some $h_0(0) \in \cC^{\frac{1-\eta}{2}}$. Let $h_\eps^{(\delta)}$ be the solution to \eqref{e:kpz_homo_full} with initial data $h_\eps(0)$ and enhanced stochastic objects $\hH_\eps^{(\delta)}$ and $\Upsilon_\eps^{(\delta)}$. Let $h_0$ be the solution to \eqref{e:kpz_limiting_equation} with coupling constant $\lambda$ in \eqref{e:kpz_coupling_const}, initial data $h_0(0)$ and enhanced stochastic objects $\hH_0$ and $\Upsilon_0$. Then for every $T>0$, $h_\eps$ converges to $h_0$ in $L_T^\infty \cC_x^{\frac{1}{2}-\eta}$ with
\begin{equation*}
    \begin{split}
    \|h_\eps^{(\delta)} - h_0\|_{L_T^\infty \cC_x^{\frac{1}{2}-\eta}} &\lesssim \eps^\theta + \|h_\eps(0) - h_0(0)\|_{\frac{1}{2}-\eta}\\
    &\phantom{11}+ \|\Upsilon_\eps^{(\delta)} - \lambda \bullet \Upsilon_0\|_{\yY_T} + \|\hH_\eps^{(\delta)} - \lambda \bullet \hH_0\|_{\yY_T^{\text{h}}} 
    \end{split}
\end{equation*}
for some $\theta>0$. Furthermore, the flux also converges with the quantitative bound
\begin{equation*}
    \begin{split}
    \|A_\eps \d_x h_\eps^{(\delta)} - \bar{A} \d_x h_0 \|_{\cC^{-\kappa}_{\sigma,T} \cC_{x}^{-\frac{1}{2}-2\kappa}} &\lesssim \eps^\theta + \|h_\eps(0) - h_0(0)\|_{\frac{1}{2}-\eta}\\
    &\phantom{11}+\|\Upsilon_\eps^{(\delta)} - \lambda \bullet \Upsilon_0\|_{\yY_T} + \|\fF_\eps^{(\delta)} - \lambda \bullet \fF_0\|_{\yY_{T}^{\flux}}\;.
    \end{split}
\end{equation*}
Here, the proportionality constants are polynomials of uniform bounds of the corresponding norms of initial data and stochastic objects $\Upsilon_{\eps}^{(\delta)}$, $\hH_{\eps}^{(\delta)}$ and  $\fF_{\eps}^{(\delta)}$.
\end{thm}

The above theorem is purely deterministic if one assumes Proposition~\ref{prop:kpz_stochastic}, and one allows flexibilities in choosing the renormalisation functions and different limiting objects along different ways of $(\eps, \delta) \rightarrow (0,0)$.

\subsection{Derivation of the fixed point system}

In this subsection, we derive the fixed point system that encodes the decomposition of the solution $h_\eps^{(\delta)}$ to \eqref{e:kpz_homo_full} and its continuous dependence on the enhanced stochastic objects $\Upsilon_\eps^{(\delta)}$. Similar as before, we drop $\delta$ from the notations for simplicity. More precisely, with an abuse of notation, we write $\Upsilon_\eps := \Upsilon_\eps^{(\delta)}$ and $h_\eps := h_\eps^{(\delta)}$. Again, these are not objects at the $\delta = 0$ limit, but just with $\delta$ omitted in the notation. The choices $C_{\eps,j} = C_{\eps,j}^{(\delta)}$ for $j=1,2,3$ and $\widetilde{C}_{\eps}^{(\delta)}$ are encoded in the stochastic objects $\Upsilon_\eps = \Upsilon_\eps^{(\delta)}$.  

With this notation, the KPZ equation is given by
\begin{equation*}
    \d_t h_\eps = \lL_\eps h_\eps + A_\eps \b_\eps (\d_x h_\eps)^2 + \xi - C_\eps - 4 \widetilde{C}_{\eps} \cdot  A_{\eps}\d_{x}h_{\eps}\;, 
\end{equation*}
Define the functional $\Lambda_\eps$ (including $\eps=0$) by
\begin{equation} \label{e:KPZ_functional_Lambda}
    \Lambda_\eps u := 2 A_\eps \d_x u + 4 \<210rk>_\eps\;, \quad \Lambda_0 u := 2 \d_x u + 4 \lambda^2 \, \<210bks>_0\;.
\end{equation}
Then
\begin{equation} \label{eq:kpz_dpd}
    \begin{split}
    u_\eps &:= h_\eps - \I_\eps \big( \xi + \<2k>_\eps + 2 \<21k>_\eps \big) - \int_{0}^{\cdot} \Pi_0 \big( \xi + \<2k>_\eps + 2 \<21k>_\eps \big)(r) {\rm d}r\\
    &= h_\eps - \<1gk>_\eps - \<20gk>_\eps - 2 \<210gk>_\eps - \int_{0}^{\cdot} \Pi_0 \big( \xi + \<2k>_\eps + 2 \<21k>_\eps \big)(r) {\rm d}r
    \end{split}
\end{equation}
satisfies the equation
\begin{equation} \label{e:kpz_u_formal}
    \begin{split}
    \d_t u_\eps = &\phantom{1}\lL_\eps u_\eps + (\Lambda_\eps u_\eps) \prec_\eps \<1bk>_\eps + (\Lambda_\eps u) \succ_\eps \<1bk>_\eps + 4 \, \<211k>_\eps + \<22k>_\eps + \<20bk>_{\eps} \cdot (\Lambda_\eps u_\eps)\\
    &+ \frac{\b_\eps}{4 A_\eps} \cdot ( \Lambda_\eps u_\eps )^2 + 2 \big( (A_\eps \d_x u_\eps) \circ_\eps \<1bk>_\eps - \widetilde{C}_\eps \cdot \Lambda_\eps u_\eps \big)\;.
    \end{split}
\end{equation}
Note that $\<2k>_\eps$ makes sense as a spacetime distribution, but not as a distribution in space for any fixed time. Hence, it always comes with integration in time. The above equation for $u_\eps$ is still formal since the resonance product $( A_\eps \d_x u_\eps) \circ_\eps \<1bk>_\eps$ is a priori ill defined for any fixed $\eps$ (uniformly in $\delta$), and one needs to subtract the divergent part $\widetilde{C}_\eps \, \Lambda_\eps u_\eps$ from it. 

To be more precise, we make the ansatz
\begin{equation} \label{e:kpz_ansatz}
    u_\eps = e^{t \lL_\eps} \big(u_\eps(0) - u_\eps^{\#}(0) \big) + \iI_\eps \big( (\Lambda_\eps u_\eps) \prec_\eps \<1bk>_\eps \big) + u_\eps^{\#}\;,
\end{equation}
where one expects $A_\eps \d_x u_\eps^{\#} \in \cC^{\frac{1}{2}+; \, \eps}$ and its resonance product $\circ_\eps$ with $\<1bk>_\eps$ be well defined. With this ansatz, we have the expression
\begin{equation*}
    A_{ \eps } \d_{ x } u_{ \eps } = \Lambda_{ \eps } (u_{\eps}) \prec_{ \eps } \<10rk>_{ \eps } + F_{ \eps } ( u_{ \eps } , u_{ \eps }^{ \# } )
\end{equation*}
where
\begin{equation} \label{eq:kpz_F_rhs}
    \begin{split}
    F_\eps(u, u^\#) = &\phantom{1}A_\eps \d_x \Big[ e^{t \lL_\eps} \big(u_\eps(0) - u_\eps^{\#}(0) \big) + u_\eps^{\#} \Big] + [A_\eps \iI_\eps', \Lambda_\eps(u) \prec_\eps] \, \<1bk>_\eps\\
    &+ (\Lambda_\eps u) \prec_\eps \big(A_\eps \iI_\eps' (\<1bk>_\eps) - \<10rk>_\eps \big)\;.
    \end{split}
\end{equation}
We expect $ F_{ \eps } ( u_{ \eps } , u^{ \# }_{ \eps } ) \in \cC^{\frac{1}{2}+} $. Using commutators, we have
% \begin{equation*}
%    ( \Lambda_{ \eps } u_{ \eps } \prec_{ \eps } \<1bk>_{ \eps }) \circ_{ \eps } \<10rk>_{ \eps } - \tilde{ C }_{ \eps } \Lambda_{ \eps } u_{ \eps } \; .
% \end{equation*}
% \begin{equation*}
%     \Lambda_\eps (u_{\eps}) \cdot \<101k>_\eps + \Com_\eps \big( \Lambda_\eps u_{\eps}; \, \<10rk>_\eps; \, \<1bk>_\eps \big) \; .
% \end{equation*}
\begin{equation} \label{eq:kpz_resonance_expression}
    A_\eps \d_x u_\eps \circ_\eps \<1bk>_\eps - \widetilde{C}_\eps \Lambda_\eps u_\eps = \Lambda_\eps (u_{\eps}) \cdot \<101k>_\eps + \Com_\eps \big( \Lambda_\eps u_{\eps}; \, \<10rk>_\eps; \, \<1bk>_\eps \big) + \<1bk>_\eps \circ_\eps F_\eps (u_\eps, u_\eps^{\#})\;.
\end{equation}
Plugging the expression \eqref{eq:kpz_resonance_expression} back into the right hand side of \eqref{e:kpz_u_formal}, and combining with \eqref{e:kpz_ansatz}, we deduce the system of equations
\begin{equation} \label{eq:kpz_system}
    \begin{split}
    u_\eps &= e^{t \lL_\eps} \big( u_\eps(0) - u_\eps^\#(0) \big) + \iI_\eps \big( (\Lambda_\eps u_\eps) \prec_\eps \<1bk>_\eps \big) + u_\eps^{\#}\;;\\
    u_\eps^{\#} &= e^{t \lL_\eps} u_\eps^\#(0) + \iI_\eps \Big( (\Lambda_\eps u_\eps) \succ_\eps \<1bk>_\eps + \frac{\b_\eps}{4 A_\eps} (\Lambda_\eps u_\eps)^2 + \kK_\eps^{(1)}(u_\eps) + 2 \, \kK_\eps^{(2)}(\vec{u}_\eps) \Big)\;,
    \end{split}
\end{equation}
where
\begin{equation} \label{eq:kpz_K}
    \begin{split}
    \kK_\eps^{(1)}(u) &= 4 \, \<211k>_\eps + \<22k>_\eps +  \Lambda_\eps (u) \, \big( \<20bk>_{\eps} + 2 \, \<101k>_\eps \big)\;,\\
    \kK_\eps^{(2)}(\vec{u}) &= \Com_\eps \big( \Lambda_\eps u; \, \<10rk>_\eps, \, \<1bk>_\eps \big) + F_\eps (\vec{u}) \circ_\eps \<1bk>_\eps\;,
    \end{split}
\end{equation}
and $F_\eps(\vec{u}) = F_\eps(u, u^\#)$ is given by \eqref{eq:kpz_F_rhs}. 

Recall the definition of $\Lambda_0$ from \eqref{e:KPZ_functional_Lambda}. For $\eps=0$, the system $(u_0, u_0^\#)$ is given by
\begin{equation} \label{eq:kpz_system_0}
    \begin{split}
    u_0 &= e^{t \lL_0} \big( u_0(0) - u_0^\#(0) \big) + \lambda \iI_0 \big( (\Lambda_0 u_0) \prec \<1bks>_0 \big) + u_0^\#\;;\\
    u_0^{\#} &= e^{t \lL_0} u_0^\#(0) + \iI_0 \Big( \lambda(\Lambda_0 u_0) \succ \<1bks>_0 + \frac{\lambda}{4} (\Lambda_0 u_0)^2 + \kK_0^{(1)}(u_0) + 2 \, \kK_0^{(2)}(\vec{u}_0) \Big)\;,
    \end{split}
\end{equation}
where
\begin{equation} \label{eq:kpz_K_0}
    \begin{split}
    \kK_0^{(1)}(u) &= 4 \lambda^3 \, \<211ks>_0 + \lambda^3 \, \<22ks>_0 +  \lambda^2 \Lambda_0 (u) \, \big( \<20bks>_0 + 2 \, \<101ks>_0 \big)\;,\\
    \kK_0^{(2)}(\vec{u}) &= \lambda^2 \, \Com_0 \big( \Lambda_0 u; \, \<10bks>_0; \, \<1bks>_0 \big) + F_0 (\vec{u}) \circ \<1bks>_0\;,
    \end{split}
\end{equation}
and
\begin{equation} \label{eq:kpz_F_rhs_0}
    \begin{split}
    F_0(u, u^\#) = &\phantom{1}\lambda \, \d_x \Big[ e^{t \lL_0} \big(u_0(0) - u_0^{\#}(0) \big) + u_0^{\#} \Big] + \lambda^2 \, [\iI_0', \,  (\Lambda_0 u) \prec] \, \<1bks>_0\\
    &+ \lambda^2 \, (\Lambda_0 u) \prec \big(\iI_0' (\<1bks>_0) - \<10bks>_0 \big)\;.
    \end{split}
\end{equation}

\subsection{Preliminary bounds}
\label{sec:kpz_preliminary}

Recall from \eqref{eq:kpz_parameters} that $\eta \in (0, \frac{1}{1000})$ and $\sigma := \frac{1}{2} + 7\eta$. For $\eps\in \N^{-1}$, define the spaces
\begin{equation} \label{eq:kpz_space_uniform}
    \xX_{T,\eps}^{(1)} := L_{T}^\infty \cC_{x}^{\frac{1}{2}-\eta} \cap \cC_{\sigma, T}^{2\eta} \cC_{x}^{1+4\eta; \, \eps}\;, \qquad \xX_{T,\eps}^{(2)} := \xX_{T, \eps}^{(1)} \cap \cC_{1+4\eta, T}^{2\kappa} \cC_{x}^{\frac{3}{2} + \eta; \, \eps}\;,
\end{equation}
and
\begin{equation*}
    \xX_{T,\eps} := \xX_{T,\eps}^{(1)} \times \xX_{T,\eps}^{(2)}\;.
\end{equation*}
We seek solution $\vec{u}_\eps = (u_\eps, u_\eps^{\#})$ to the system \eqref{eq:kpz_system} in $\xX_{T,\eps}$. Since $\Lambda_\eps u_\eps$ appears many times on the right hand side of the equation, we first give a bound for that quantity.

% \begin{lem} \label{lem:kpz_stochasitc_210}
% We have the uniform bound
% \begin{equation} \label{e:kpz_stochasitc_210_uniform}
%     \|\<210gk>_\eps\|_{\cC_T^{\eta/2} \cC_x^{\frac{3}{2}-2\eta; \eps}} + \|\<210rk>_\eps\|_{\cC_T^{\eta/2} \cC_x^{\frac{1}{2}-2\eta}} \lesssim \|\lL_{\eps}\<210gk>_\eps\|_{\cC_T^{\eta/2} \cC_x^{-\frac{1}{2}-2\eta}}\;,
% \end{equation}
% For the difference, there exists $\theta>0$ such that
% \begin{equation} \label{e:kpz_stochastic_210_difference}
%     \|\<210rk>_\eps - \<210rk>_0\|_{\cC_T^{\eta/2} \cC_x^{\frac{1}{2}-3\eta}} \lesssim \eps^{\theta} \|\lL_0 \<210gk>_0\|_{\cC_T^{\eta/2} \cC_x^{-\frac{1}{2}-2\eta}} + \| \lL_\eps \<210gk>_\eps - \lL_0 \<210gk>_0 \|_{\cC_T^{\eta/2} \cC_x^{-\frac{1}{2}-2\eta}}\;.
% \end{equation}
% \end{lem}
% \begin{proof}
% The first bound in \eqref{e:kpz_stochasitc_210_uniform} directly follows from definition. The second one follows from Lemma~\ref{le:derivative_a_1d}. The bound \eqref{e:kpz_stochastic_210_difference} follows from Lemma~\ref{le:convergence_derivative_a_1d}. 
% \end{proof}

\begin{lem} \label{le:kpz_lambda}
The functional $\Lambda_\eps$ given in \eqref{e:KPZ_functional_Lambda} satisfies the bounds
\begin{equation*}
    \|\Lambda_\eps u_\eps\|_{\cC_{\sigma, T}^{2\eta} \cC_{x}^{4\eta}} \lesssim \|u_\eps\|_{\cC_{\sigma, T}^{2\eta} \cC_{x}^{1+4\eta; \eps}} + \|\<210rk>_\eps\|_{\cC_{T}^{2\eta} \cC_{x}^{\frac{1}{2}-5\eta}}\;,
\end{equation*}
and
\begin{equation*}
    \begin{split}
    \|\Lambda_\eps u_\eps - \bar{A} \Lambda_0 u_0\|_{\cC_{\sigma, T}^{2\eta} \cC_{x}^{4\eta}} &\lesssim \eps^\theta \|\lL_0 u_0\|_{\cC_{\sigma, T}^{2\eta} \cC_{x}^{4\eta-1}} + \|\lL_\eps u_\eps - \lL_0 u_0\|_{\cC_{\sigma, T}^{2\eta} \cC_{x}^{4\eta-1}}\\
    &\phantom{11}+ \|\<210rk>_\eps - \bar{A} \lambda^2 \, \<210bks>_0\|_{\cC_{T}^{2\eta} \cC_{x}^{\frac{1}{2}-5\eta}}\;. 
    \end{split}
\end{equation*}
\end{lem}
\begin{proof}
The bounds follow immediately from the definition of $\Lambda_\eps$ and Lemmas~\ref{le:derivative_a_1d} and~\ref{le:convergence_derivative_a_1d}. 
\end{proof}

In what follows, we collect some preliminary bounds for quantities that appear on the right hand side of \eqref{eq:kpz_system}. 

\begin{lem} \label{lem:kpz_prec_succ}
We have
\begin{equation*}
    \begin{split}
    \|\Lambda_\eps (u_\eps) \prec_\eps \<1bk>_\eps\|_{\cC_{\sigma,T}^{-\kappa} \cC_x^{-\frac{1}{2}-\kappa}} &\lesssim \big( \|\Upsilon_\eps\|_{\yY_T} + \|u_\eps\|_{\xX_{T,\eps}^{(1)}} \big) \, \| \, \<1bk>_\eps \|_{\cC_T^{-\kappa} \cC_x^{-\frac{1}{2}-\kappa}}\;,\\
    \|\Lambda_\eps (u_\eps) \succ_\eps \<1bk>_\eps\|_{\cC_{\sigma,T}^{-\kappa} \cC_x^{3\eta-\frac{1}{2}}} &\lesssim \big( \|\Upsilon_\eps\|_{\yY_T} + \|u_\eps\|_{\xX_{T,\eps}^{(1)}} \big) \, \| \, \<1bk>_\eps \|_{\cC_T^{-\kappa} \cC_x^{-\frac{1}{2}-\kappa}}\;.
    \end{split}
\end{equation*}
Furthermore, there exists $\theta>0$ such that
\begin{equation*}
    \begin{split}
    &\phantom{111}\big\| \Lambda_\eps (u_\eps) \prec_\eps \<1bk>_\eps - \lambda\Lambda_0(u_0) \prec \<1bks>_0 \big\|_{\cC_{\sigma,T}^{-\kappa} \cC_x^{\frac{1}{2}-3\kappa}} + \big\| \Lambda_\eps (u_\eps) \succ_\eps \<1bk>_\eps - \lambda\Lambda_0(u_0) \succ \<1bks>_0 \big\|_{\cC_{\sigma,T}^{-\kappa} \cC_x^{3\eta-\frac{1}{2}}}\\
    &\lesssim \|\Lambda_\eps (u_\eps) - \bar{A}\Lambda_0 (u_0)\|_{\cC_{\sigma,T}^{2\eta} \cC_x^{4\eta}} \|\Upsilon_\eps\|_{\yY_T} + \|\Lambda_0 (u_0)\|_{\cC_{\sigma,T}^{2\eta} \cC_x^{4\eta}} \big( \eps^\theta \|\Upsilon_0\|_{\yY_T^+} + \|\Upsilon_\eps - \lambda \bullet  \Upsilon_0\|_{\yY} \big)\;,
    \end{split}
\end{equation*}
where $\yY_T^+$ is the strenghening of $\yY_T$ by replacing all appearances of $\cC_T^{-\kappa}$ by $L_T^\infty$. 
\end{lem}
\begin{proof}
This is a direct consequence of Proposition~\ref{prop:Young_multilinear} and the para-product estimates in Proposition~\ref{pr:paraproducts_gBesov}. So we omit the details. 
\end{proof}

\begin{lem} \label{lem:kpz_quadratic}
We have
\begin{equation*}
    \Big\| (\b_\eps / A_\eps) \cdot (\Lambda_\eps u_\eps)^2  \Big\|_{L_{2\sigma,T}^\infty L_x^\infty} \lesssim \big( \|\Upsilon_{\eps}\|_{\yY_T} + \|u_\eps\|_{\xX_{T,\eps}^{(1)}} \big)^2\;.
\end{equation*}
Furthermore, there exists $\theta>0$ such that
\begin{equation*}
    \begin{split}
    &\phantom{111} \Big\| (\b_\eps / A_\eps) \cdot (\Lambda_\eps u_\eps)^2 - \lambda \cdot (\Lambda_0 u_0)^2 \Big\|_{L_{2\sigma,T}^\infty \cC_x^{-\kappa}}\\
    &\lesssim \eps^\theta \|\Lambda_\eps u_\eps\|_{\cC_{\sigma,T}^{2\eta} \cC_x^{4\eta}}^2 + \Big( \|\Lambda_\eps u_\eps + \bar{A}\Lambda_0 u_0\|_{\cC_{\sigma,T}^{2\eta} \cC_x^{4\eta}} \Big) \, \Big( \|\Lambda_\eps u_\eps - \bar{A}\Lambda_0 u_0\|_{\cC_{\sigma,T}^{2\eta} \cC_x^{4\eta}} \Big)
    \end{split}
\end{equation*}
\end{lem}
\begin{proof}
The uniform-in-$\eps$ bound is straightforward. For the difference, it suffices to have the decomposition
\begin{equation*}
    \begin{split}
    &\phantom{111} (\b_\eps / A_\eps) \cdot (\Lambda_\eps u_\eps)^2 - \Pi_0 (\b / A) \cdot (\bar{A}\Lambda_0 u_0)^2\\
    &= \Big( (\b_\eps / A_\eps) - \Pi_0 (\b / A) \Big) \cdot (\Lambda_\eps u_\eps)^2 + \Pi_0 (\b/A) \, \big( \Lambda_\eps u_\eps - \bar{A}\Lambda_0 u_0 \big) \, \big( \Lambda_\eps u_\eps + \bar{A}\Lambda_0 u_0 \big)\;,
    \end{split}
\end{equation*}
and note that $\b_\eps / A_\eps$ converges to $\lambda / \bar{A}^2$ in $L_T^\infty \cC_x^{-\kappa}$. The claim then follows. 
\end{proof}

\begin{lem} \label{lem:kpz_K1}
Recall the expressions $\kK_\eps^{(1)}$ and $\kK_0^{(1)}$ from \eqref{eq:kpz_K} and \eqref{eq:kpz_K_0}. We have the uniform bound
\begin{equation*}
    \big\| \kK_\eps^{(1)} (u_\eps) \big\|_{\cC_{\sigma, T}^{-\kappa} \cC_x^{-2\kappa}} \lesssim \|\Upsilon_\eps\|_{\yY_T} \big( 1 + \|\Upsilon_\eps\|_{\yY_T} + \|u_\eps\|_{\xX_{T,\eps}^{(1)}} \big)\;.
\end{equation*}
Furthermore, there exists $\theta>0$ such that
\begin{equation*}
    \begin{split}
    \big\| \kK_\eps^{(1)}(u_\eps) - \kK_0^{(1)}( u_0 )  \big\|_{\cC_{\sigma, T}^{-\kappa} \cC_x^{-2\kappa}} &\lesssim \big( 1 + \|\Upsilon_0\|_{\yY_T} + \|u_0\|_{\xX_{T,0}^{(1)}} \big) \|\Upsilon_\eps - \lambda \bullet \Upsilon_0\|_{\yY_T}\\
    &\phantom{11}+ \|\Upsilon_\eps\|_{\yY_T} \big\|\Lambda_\eps(u_\eps) - \bar{A}\Lambda_0(u_0) \big\|_{\cC_{\sigma,T}^{2\eta} \cC_x^{4\eta}}\;. 
    \end{split}
\end{equation*}
\end{lem}
\begin{proof}
The first claim is straightforward from Lemma~\ref{le:kpz_lambda}. For the second one, we have
\begin{equation*}
    \begin{split}
    \kK_\eps^{(1)}(u_\eps) - \kK_0^{(1)}(u_0) &= 4 \, \big( \<211k>_\eps -  \lambda^3 \, \<211ks>_0 \big) + \big( \<22k>_\eps - \lambda^3 \, \<22ks>_0 \big) +  \big( \Lambda_\eps (u_\eps) - \bar{A}\Lambda_0(u_0) \big) \, \big( \<20bk>_{\eps} + 2 \, \<101k>_\eps \big)\\
    &\phantom{111}+ \bar{A}\Lambda_0 (u_0) \big[ \left(\<20bk>_{\eps} - (\lambda^2 / \bar{A}) \cdot \<20bks>_0\right) + 2 \, \left(\<101k>_\eps -   \, (\lambda^2 / \bar{A}) \cdot\<101ks>_0\right) \big]\;.
    \end{split}
\end{equation*}
The claim then again follows from Lemma~\ref{le:kpz_lambda}. This completes the proof. 
\end{proof}

\begin{lem} \label{lem:kpz_F_bound}
Recall the expressions $F_\eps$ and $F_0$ from \eqref{eq:kpz_F_rhs} and \eqref{eq:kpz_F_rhs_0}. We have the bound
\begin{equation} \label{eq:kpz_F_uniform}
    \| F_\eps (\vec{u}_\eps) \|_{\cC_{1+4\eta,T}^{2\kappa} \cC_x^{\frac{1}{2}+2\kappa}} \lesssim \|\Upsilon_\eps\|_{\yY_T} \big( \|\Upsilon_\eps\|_{\yY_T} + \|\vec{u}_\eps\|_{\xX_{T,\eps}} \big) + \|u_\eps(0) - u_\eps^{\#}(0)\|_{\cC_x^{\frac{1}{2}-\eta}}\;.
\end{equation}
Furthermore, there exists $\theta>0$ such that
\begin{equation*}
    \begin{split}
    \|F_\eps &(\vec{u}_\eps) - (\bar{A} / \lambda) F_0 (\vec{u}_0)\|_{L_{1+4\eta,T}^\infty \cC_x^{\frac{1}{2}+2\kappa}} \lesssim \eps^\theta \big( \|u_0(0) - u_0^{\#}(0)\|_{\frac{1-\eta}{2}} + \|u_0^\#\|_{\xX_{T,0}^{(2)}} \big)\\
    &+ \big\| \big( u_\eps(0) - u_\eps^{\#}(0) \big) - \big( u_0(0) - u_0^{\#}(0) \big) \big\|_{\frac{1}{2}-\eta} + \|\lL_\eps u_\eps^\# - \lL_0 u_0^\#\|_{L_{1+4\eta,T}^\infty \cC_x^{\eta-\frac{1}{2}}}\\
    &+ \|\Upsilon_\eps - \lambda \bullet \Upsilon_0\|_{\yY_T}  \|\Lambda_\eps (u_\eps)\|_{\cC_{\sigma,T}^{2\eta} \cC_x^{4\eta}} + \big( \eps^\theta \|\Lambda_0 (u_0)\|_{\cC_{\sigma,T}^{2\eta} \cC_x^{4\eta}}\\
    &+ \|\Lambda_\eps (u_\eps) - \bar{A}\Lambda_0 (u_0)\|_{\cC_{\sigma,T}^{2\eta}\cC_x^{4\eta}} \big) \cdot \|\Upsilon_0\|_{\yY_T^+}\;.\\
    \end{split}
\end{equation*}
\end{lem}
\begin{proof}
We first deal with the uniform bound \eqref{eq:kpz_F_uniform}. We treat the right hand side in \eqref{eq:kpz_F_rhs} term by term. First, by definition of $\xX_{T,\eps}^{(2)}$ and Lemmas~\ref{le:heat_smoothing_space} and~\ref{le:derivative_a_1d}, we have
\begin{equation} \label{eq:kpz_F_a_derivative}
    \Big\| A_\eps \d_x \Big( e^{t \lL_\eps} \big(u_\eps(0) - u_\eps^{\#}(0) \big) + u_\eps^{\#} \Big) \Big\|_{\cC_{1+4\eta,T}^{2\kappa} \cC_x^{\frac{1+\eta}{2}}} \lesssim \|u_\eps(0) - u_\eps^{\#}(0)\|_{\cC_x^{\frac{1}{2}-\eta}} + \|u_\eps^{\#}\|_{\xX_{T,\eps}^{(2)}}\;.
\end{equation}
For the commutator term, by Proposition~\ref{prop:commutator_derivative_heat_negative_time} (with $\gamma=-\frac{1}{2}-\kappa$) and Lemma~\ref{le:kpz_lambda}, we have the bound
\begin{equation} \label{eq:kpz_F_comm_heat_derivative_prec}
    \big\| \big[ A_\eps \iI_\eps', \; \Lambda_\eps(u_\eps) \prec_\eps \big] \, \<1bk>_\eps \big\|_{\cC_{\sigma,T}^{2\kappa} \cC_x^{\frac{1}{2}+\eta}} \lesssim \|\Upsilon_\eps\|_{\yY_T} \big( \|u_\eps\|_{\xX_{T,\eps}^{(1)}} + \|\Upsilon_\eps\|_{\yY_T} \big)\;.
\end{equation}
Finally, note that
\begin{equation*}
    \big\| A_\eps \d_x \iI_\eps (\<1bk>_\eps) - \<10rk>_\eps \big\|_{\cC_{3\eta,T}^{\eta/2} \cC_x^{\frac{1}{2}+\eta}} = \Big\| A_\eps \d_x e^{t \lL_\eps} \int_{-\infty}^{0} e^{-r \lL_\eps} \<1bk>_\eps(r) {\rm d}r \Big\|_{\cC_{3\eta,T}^{\eta/2} \cC_x^{\frac{3}{2}+\eta; \eps}} \lesssim \|\<1bk>_\eps\|_{\cC_T^{-\kappa} \cC_x^{-\frac{1}{2}-\kappa}}\;,
\end{equation*}
we have
\begin{equation} \label{eq:kpz_F_remainder}
    \big\| (\Lambda_\eps u_\eps) \prec_\eps \big(A_\eps \d_x \iI_\eps (\<1bk>_\eps) - \<10rk>_\eps \big) \big\|_{\cC_{\sigma+3\eta,T}^{\eta/2} \cC_x^{\frac{1}{2}+\eta}} \lesssim \|\Upsilon_\eps\|_{\yY_T} \big( \| \Upsilon_{\eps} \|_{\yY_T} + \|u_\eps\|_{\xX_{T,\eps}^{(1)}} \big)\;.
\end{equation}
Combining \eqref{eq:kpz_F_a_derivative}, \eqref{eq:kpz_F_comm_heat_derivative_prec} and \eqref{eq:kpz_F_remainder} and noting that $1 + 4 \eta > \sigma + 3 \eta$, we get the uniform bound \eqref{eq:kpz_F_uniform}. 

As for the difference $F_\eps (\vec{u}_\eps) - (\bar{A} / \lambda) \, F_0(\vec{u}_0)$, the desired bounds for most of the terms follow by the same arguments together with corresponding convergence lemmas in Sections~\ref{sec:para_convergence} and~\ref{sec:bounds_derivatives_1d}. We give details for the term involving the commutator $\big[ A_\eps \iI_\eps', \; \Lambda_\eps(u_\eps) \prec_\eps \big] \, \<1bk>_\eps$. We decompose the difference into
\begin{equation} \label{eq:kpz_F_commutator_difference}
    \begin{split}
    \big[ A_\eps \iI_\eps', \; \Lambda_\eps(u_\eps) \prec_\eps \big] \, \<1bk>_\eps - &\lambda \bar{A} \, \big[ \iI_0', \,  (\Lambda_0 u) \prec \big] \, \<1bks>_0 = \big[ A_\eps \iI_\eps', \; \Lambda_\eps(u_\eps) \prec_\eps \big] \, \big(\<1bk>_\eps - (\lambda / \bar{A}) \, \<1bks>_0 \big)\\
    &+ \frac{\lambda}{\bar{A}} \, \big( \big[ A_\eps \iI_\eps', \; \Lambda_\eps(u_\eps) \prec_\eps \big] \, \<1bks>_0 - \bar{A}^2 \, \big[ \iI_0', \,  (\Lambda_0 u_0) \prec \big] \, \<1bks>_0 \big)\;.
    \end{split}
\end{equation}
The advantage of decomposing in this way is that although $\<1bk>_\eps - (\lambda / \bar{A}) \, \<1bks>_0$ still converges in $\cC_T^{-\kappa} \cC_x^{-\frac{1}{2}-\kappa}$, the object $\<1bks>_0$ actually lives in $L_T^\infty \cC_x^{-\frac{1}{2}-\kappa}$. For the first term above, it follows directly from Proposition~\ref{prop:commutator_derivative_heat_negative_time} (just replacing $\<1bk>_\eps$ with $\<1bk>_\eps - (\lambda / \bar{A}) \, \<1bks>_0$) that
\begin{equation*}
    \big\| [ A_\eps \iI_\eps', \; \Lambda_\eps(u_\eps) \prec_\eps ] \, (\<1bk>_\eps - (\lambda / \bar{A}) \, \<1bks>_0 ) \big\|_{L_{\sigma,T}^{\infty} \cC_x^{\frac{1}{2}+\eta}} \lesssim  \|\Lambda_\eps (u_\eps)\|_{\cC_{\sigma,T}^{2\eta} \cC_x^{4\eta}} \| \, \<1bk>_\eps - (\lambda / \bar{A}) \, \<1bks>_0\|_{\cC_T^{-\kappa} \cC_x^{-\frac{1}{2}-\kappa}}\;.
\end{equation*}
As for the second one, by Proposition~\ref{prop:derivative_heat_prec_comm} (with $\alpha = 2\eta$, $\beta = -\frac{1}{2}-\kappa$), we have
\begin{equation*}
    \begin{split}
    &\phantom{111}\big\| [ A_\eps \iI_\eps', \; \Lambda_\eps(u_\eps) \prec_\eps ] \, \<1bks>_0 - \bar{A}^2 \, [ \iI_0', \,  (\Lambda_0 u_0) \prec ] \, \<1bks>_0 \big) \big\|_{L_{\sigma,T}^\infty \cC_x^{\frac{1}{2}+\eta}}\\
    &\lesssim \Big( \eps^\theta \|\Lambda_0 (u_0) \|_{\cC_{\sigma,T}^{2\eta} \cC_x^{4\eta}} + \big\| \Lambda_\eps (u_\eps) - \bar{A} \Lambda_0 (u_0) \big\|_{\cC_{\sigma,T}^{2\eta} \cC_x^{4\eta}} \Big) \cdot \|\<1bks>_0\|_{L_T^\infty \cC_x^{\frac{1}{2}-\kappa}}\;,
    \end{split}
\end{equation*}
where we also used that $\eta > 1000\kappa$. The above two bounds together give the desired bound for the difference of the commutators --- the $L_{\sigma,T}^\infty \cC_x^{\frac{1}{2}+\eta}$-norm of the right hand side of \eqref{eq:kpz_F_commutator_difference}. The bounds for all other terms can be obtained in similar (but simpler) ways. This completes the proof of the lemma. 
\end{proof}

\begin{lem} \label{lem:kpz_K2}
Recall the expressions for $\kK_\eps^{(2)}$ and $\kK_0^{(2)}$ from \eqref{eq:kpz_K} and \eqref{eq:kpz_K_0}. We have
\begin{equation*}
    \big\| \kK_\eps^{(2)}(\vec{u}_\eps) \big\|_{\cC_{1+4\eta,T}^{-\kappa} \cC_x^{\kappa}} \lesssim \|\Lambda_\eps (u_\eps)\|_{\cC_{\sigma,T}^{2\eta} \cC_x^{4\eta}} \|\Upsilon_\eps\|_{\yY_T}^2 + \|F_\eps (\vec{u}_\eps)\|_{\cC_{1+4\eta,T}^{2\kappa} \cC_x^{\frac{1}{2}+4\kappa}} \|\Upsilon_\eps\|_{\yY_T}\;.
\end{equation*}
Furthermore, there exists $\theta>0$ such that
\begin{equation*}
    \begin{split}
    &\phantom{111}\big\| \kK_\eps^{(2)}(\vec{u}_\eps) - \kK_0^{(2)}(\vec{u}_0) \big\|_{\cC_{1+4\eta,T}^{-\kappa} \cC_x^\kappa}\\
    &\lesssim \|\Lambda_\eps (u_\eps)\|_{\cC_{\sigma,T}^{2\eta} \cC_x^{4\eta}} \big( \|\Upsilon_\eps\|_{\yY_T} + \|\Upsilon_0\|_{\yY_T^+} \big) \, \|\Upsilon_\eps - \lambda \bullet \Upsilon_0\|_{\yY_T}\\
    &\phantom{111}+  \Big( \eps^\theta \|\Lambda_0 (u_0)\|_{\cC_{\sigma,T}^{2\eta} \cC_x^{4\eta}} + \|\Lambda_\eps(u_\eps) - \bar{A}\Lambda_0 (u_0)\|_{\cC_{\sigma,T}^{2\eta} \cC_x^{4\eta}} \Big) \|\Upsilon_0\|_{\yY_T^+}^2\\
    &\phantom{111}+ \Big( \eps^\theta \|F_0(\vec{u}_0)\|_{L_{1+4\eta,T}^\infty \cC_x^{\frac{1}{2}+2\kappa}} + \left\|  F_\eps(\vec{u}_\eps) -  (\bar{A}/\lambda) \, F_0(\vec{u}_0)\right\|_{L_{1+4\eta,T}^{\infty} \cC_x^{\frac{1}{2}+2\kappa}} \Big) \, \| \Upsilon_0 \|_{\yY_T^+}\\
    &\phantom{111}+ \|F_\eps (\vec{u}_\eps)\|_{\cC_{1+4\eta,T}^{2\kappa} \cC_x^{\frac{1}{2}+2\kappa}} \|\Upsilon_\eps - \lambda \bullet \Upsilon_0\|_{\yY_T}\;.
    \end{split}
\end{equation*}
\end{lem}
\begin{proof}
The uniform bound for $\kK_\eps^{(2)}$ follows directly from Proposition~\ref{prop:Young_multilinear} with $\T_m = \Com_\eps$, $f_1 = \Lambda_\eps (u_\eps)$, $f_2 = \<10rk>_\eps$ and $g = \<1bk>_\eps$. For the difference, we write
\begin{equation*}
    \begin{split}
    \Com_\eps \big( &\Lambda_\eps u_\eps; \, \<10rk>_\eps, \, \<1bk>_\eps \big) - \lambda^2 \, \Com_0 \big( \Lambda_0 u_0; \, \<10bks>_0; \, \<1bks>_0 \big) = \Com_\eps \big( \Lambda_\eps u_\eps; \, \<10rk>_\eps, \, \<1bk>_\eps - (\lambda / \bar{A}) \, \<1bks>_0 \big)\\
    &+ \frac{\lambda}{\bar{A}} \Big[ \Com_\eps \big( \Lambda_\eps u_\eps; \, \<10rk>_\eps, \, \<1bks>_0 \big) - \Com_0 \big( \bar{A} \, \Lambda_0 u_0; \, \lambda \, \<10bks>_0; \, \<1bks>_0 \big) \Big]\;.
    \end{split}
\end{equation*}
The bound for the first term on the right hand side above follows from Proposition~\ref{prop:Young_multilinear} in the same way except that for $g = \<1bk>_\eps - (\lambda / \bar{A}) \, \<1bks>_0$. As for the second term (the difference between two commutators), since the third input $\<1bks>_0 \in L_T^\infty \cC_x^{-\frac{1}{2}-\kappa}$ is bounded in time now, the desired bound then follows from Lemma~\ref{le:convergence_comm_good}. 

Similarly, we also decompose
\begin{equation*}
    \begin{split}
    F_\eps (\vec{u}_\eps) &\circ_\eps \<1bk>_\eps - F_0 (\vec{u}_0) \circ \<1bks>_0 = F_\eps (\vec{u}_\eps) \circ_\eps \big( \<1bk>_\eps - (\lambda / \bar{A}) \, \<1bks>_0 \big)\\
    &+ \lambda / \bar{A} \Big( F_\eps (\vec{u}_\eps) \circ_\eps \<1bks>_0 - (\bar{A} / \lambda) F_0 (\vec{u}_0) \circ \<1bks>_0 \Big)\;,
    \end{split}
\end{equation*}
and the desired bound then follows in the same way as above (using Proposition~\ref{eq:Young_multilinear} for the first term, and Lemma~\ref{le:convergence_circ_good} for the second term). 
\end{proof}

\subsection{Uniform boundedness}
\label{sec:kpz_uniform}

Fix $T>0$. Similar as before, let
\begin{equation*}
    \mM_\eps := \|u_\eps(0)\|_{\cC^{\frac{1}{2}-\eta}} + \|u_\eps^{\#}(0)\|_{\cC^{\frac{1}{2}-\eta}} + \|\Upsilon_{\eps}\|_{\yY_T}\;.
\end{equation*}
We have the following theorem. 

\begin{thm}\label{thm:kpz_fixed_pt}
    Let $\eta \in (0, \frac{1}{1000})$ and the space $\xX_{T,\eps} = \xX_{T,\eps}^{(1)} \times \xX_{T,\eps}^{(2)}$ be as given in \eqref{eq:kpz_space_uniform}. For every $\eps \in \N^{-1}$ and every pair of initial data $\big( u_\eps(0), u_\eps^\#(0) \big) \in \cC^{\frac{1}{2}-\eta} \times \cC^{\frac{1}{2}-\eta}$, there exists $\tau_\eps \in (0,T)$ such that the system \eqref{eq:kpz_system} admits a unique solution $\vec{u}_\eps = (u_\eps, u_\eps^{\#}) \in \xX_{\tau_\eps, \eps}$. 

    Furthermore, the existence time $\tau_\eps$ and the norm $\|\vec{u}_\eps\|_{\xX_{\tau_\eps, \eps}}$ depends on $\eps$ via $\mM_\eps$ only. In particular, if $\mM_\eps$ is uniformly bounded in $\eps$, then there exists $\tau>0$ independent of $\eps$ such that
    \begin{equation} \label{e:kpz_sol_uniform_bound}
        \|\vec{u}_\eps\|_{\xX_{\tau,\eps}} \lesssim 1 + \mM_\eps\;,
    \end{equation}
    where the proportionality constant is independent of $\eps$. 
\end{thm}
\begin{proof}
For $\tau \in (0,T)$ and $\eps\in \N^{-1}$, define the mild solution map $\Gamma_{\tau,\eps} = (\Gamma_{\tau,\eps}^{(1)}, \Gamma_{\tau,\eps}^{(2)})$ by
\begin{equation*}
    \begin{split}
    \Gamma_{\tau,\eps}^{(1)}(\vec{u}) &= e^{t \lL_\eps} u_\eps(0) + \iI_\eps \big( (\Lambda_\eps u) \prec_\eps \<1bk>_\eps \big)\\
    &\phantom{11111}+ \iI_\eps \Big( (\Lambda_\eps u) \succ_\eps \<1bk>_\eps + \frac{\b_\eps}{4 A_\eps} (\Lambda_\eps u)^2 + \kK_\eps^{(1)}(u) + 2 \, \kK_\eps^{(2)}(\vec{u}) \Big)\;,\\
    \Gamma_{\tau,\eps}^{(2)}(\vec{u}) &= e^{t \lL_\eps} u^{\#}_\eps(0) + \iI_\eps \Big( (\Lambda_\eps u) \succ_\eps \<1bk>_\eps + \frac{\b_\eps}{4 A_\eps} (\Lambda_\eps u)^2 + \kK_\eps^{(1)}(u) + 2 \, \kK_\eps^{(2)}(\vec{u}) \Big)\;,
    \end{split}
\end{equation*}
where we write $\vec{u} = (u, u^{\#})$, and $\kK_{\eps}$ is defined in \eqref{eq:kpz_K}. We first show that for sufficiently small $\tau>0$, $\Gamma_{\tau, \eps}$ is a contraction map from a bounded ball of $\xX_{\tau, \eps}$ into itself. 

\begin{flushleft}
	\textit{Step 1.}
\end{flushleft}

We first check that $\Gamma_{\tau,\eps}$ maps the ball in $\xX_{\tau, \eps}$ of radius $\rR$ centered at the origin into itself for some suitable $\rR$ and $\tau$ to be specified later. For the initial data term, by Lemma~\ref{le:heat_smoothing_space}, we have
\begin{equation*}
    \|e^{t \lL_\eps} u_\eps(0)\|_{\xX_{\tau,\eps}^{(1)}} \lesssim \|u_\eps(0)\|_{\frac{1}{2}-\eta}\;, \quad \|e^{t \lL_\eps} u^{\#}_\eps(0)\|_{\xX_{\tau,\eps}^{(2)}} \lesssim \|u^{\#}_\eps(0)\|_{\frac{1}{2}-\eta}\;.
\end{equation*}
This gives the desired bounds for terms from the initial data. Also, by  Proposition~\ref{prop:Young_heat_convolution} and Lemma~\ref{lem:kpz_prec_succ}, we see there exists $\theta>0$ such that
\begin{equation*}
    \| \iI_\eps \big( \Lambda_\eps u \prec_\eps \<1bk>_\eps \big) \|_{\xX_{\tau,\eps}^{(1)}} \lesssim \tau^\theta \|\Lambda_\eps u \prec_\eps \<1bk>_\eps\|_{\cC_{\sigma,\tau}^{-\kappa} \cC_x^{-\frac{1}{2}-\kappa}} \lesssim \tau^{\theta} \big( \|u\|_{\xX_{\tau,\eps}^{(1)}} + \mM_\eps \big) \mM_\eps\;.
\end{equation*}
We now turn to the terms that contributing to $u_\eps^\#$; that is, those on the right hand side of $\Gamma_{\tau,\eps}^{(2)}$. Again by Proposition~\ref{prop:Young_heat_convolution} and Lemma~\ref{lem:kpz_prec_succ}, we have
\begin{equation*}
    \begin{split}
    \big\| \iI_\eps \big( \Lambda_\eps u \succ_\eps \<1bk>_\eps \big) \big\|_{\xX_{\tau,\eps}^{(1)} \cap \xX_{\tau,\eps}^{(2)}} \lesssim \tau^\theta \|\Lambda_\eps u \succ_\eps \<1bk>_\eps\|_{\cC_{\sigma,\tau}^{-\kappa} \cC_x^{3\eta - \frac{1}{2}}} \lesssim \tau^\theta \big( \|u\|_{\xX_{\tau,\eps}^{(1)}} + \mM_\eps \big) \mM_\eps\;.
    \end{split}
\end{equation*}
Also, the term with $(\Lambda_\eps u)^2$ satisfies the uniform bound
\begin{equation} \label{eq:kpz_fixed_pt_quadratic}
    \big\| \iI_\eps \big( (\b_\eps / 4 A_\eps) \cdot (\Lambda_\eps u)^2 \big)  \big\|_{\xX_{\tau,\eps}^{(1)} \cap \xX_{\tau,\eps}^{(2)}} \lesssim \tau^\theta \big\| (\Lambda_\eps u)^2 \big\|_{L_{2\sigma,T}^{\infty} L_x^{\infty}} \lesssim \tau^\theta \big( \|u\|_{\xX_{\tau,\eps}^{(1)}} + \mM_\eps \big)^2\;.
\end{equation}
As for $\kK_\eps^{(1)}$ and $\kK_\eps^{(2)}$, by Proposition~\ref{prop:Young_heat_convolution} and Lemmas~\ref{lem:kpz_K1} and~\ref{lem:kpz_K2}, we have
\begin{equation*}
    \big\| \iI_\eps \big( \kK_\eps^{(1)}(u_\eps) + \kK_\eps^{(2)}(\vec{u}_\eps) \big) \big\|_{\xX_{\tau,\eps}^{(1)} \cap \xX_{\tau,\eps}^{(2)}} \lesssim \tau^\theta \mM_\eps^2 \big( \|\vec{u}\|_{\xX_{\tau,\eps}} + \mM_\eps + 1 \big)\;.
\end{equation*}
Combining all of the above, we deduce that
\begin{equation*}
    \| \Gamma_{\tau,\eps} \vec{u} \|_{\xX_{\tau,\eps}} \lesssim \mM_\eps + \tau^\theta \big( 1 + \mM_\eps + \|\vec{u}\|_{\xX_{\tau,\eps}} \big)^3\;.
\end{equation*}
Hence, we can choose $C_0$ sufficiently large and $\rR = 1 + 2 C_0 \mM_\eps$ such that
\begin{equation} \label{e:kpz_fixed_pt_map_space}
    \|\Gamma_{\tau, \eps} \vec{u}\|_{\xX_{\tau,\eps}} \leq \frac{\rR}{2} + C \tau^{\theta} \rR^3
\end{equation}
for some $C>0$ and all $\vec{u}$ with $\|\vec{u}\|_{\xX_{\tau,\eps}} \leq \rR$. Hence, $\Gamma_{\tau,\eps}$ maps a ball of radius $\rR$ in $\xX_{\tau,\eps}$ into itself provided $\tau$ is sufficiently small (depending on $\mM_\eps$ only). This shows that $\Gamma_{\tau,\eps} = \Gamma^{(\delta)}_{\tau,\eps}$ maps a ball of radius $\rR$ into itself provided $\tau$ sufficiently small (depending on $\mM_\eps = \mM_\eps^{(\delta)}$ only).

\begin{flushleft}
\textit{Step 2. }
\end{flushleft}

We now show that this map $\Gamma_{\tau,\eps}$ is a contraction in $\xX_{\tau,\eps}$ for sufficiently small $\tau$. For this, we need to get a bound for $\| \Gamma_{\tau,\eps}(\vec{u}) - \Gamma_{\tau,\eps}(\vec{v})\|_{\xX_{\tau,\eps}}$ for $\vec{u} = (u, u^\#)$ and $\vec{v} = (v, v^{\#})$ in $\xX_{\tau,\eps}$. The argument is essentially the same as in Step 1. Most of the terms are either ``constants" (the initial data and pure stochastic terms) or linear in $u$ or $u^{\#}$. Hence for the difference, they either cancel out or one replaces $\|\vec{u}\|_{\xX_{\tau,\eps}}$ by $\|\vec{u} - \vec{v}\|_{\xX_{\tau,\eps}}$. 

The only nonlinear term is $\frac{\b_\eps}{A_\eps} (\Lambda_\eps u)^2$ which is quadratic in $u$. For this term, one simply replaces one factor of $\rR$ by $\|u-v\|_{\xX_{\tau,\eps}^{(1)}}$. Hence, we get the bound
\begin{equation} \label{e:kpz_fixed_pt_contraction}
    \|\Gamma_{\tau,\eps}(\vec{u}) - \Gamma_{\tau,\eps}(\vec{v})\|_{\xX_{\tau,\eps}} \lesssim \tau^{\theta} \|\vec{u} - \vec{v}\|_{\xX_{\tau,\eps}} (1 + \mM_\eps + \rR)^2
\end{equation}
for some $\theta>0$. Combining \eqref{e:kpz_fixed_pt_map_space} and \eqref{e:kpz_fixed_pt_contraction}, we see if one chooses $\tau>0$ such that
\begin{equation} \label{e:kpz_local_existence_time}
    \tau^\theta \rR^2 < c
\end{equation}
for some sufficiently small (but universal) $c$, $\Gamma_{\tau,\eps}$ is a contraction in $\xX_{\tau,\eps}$ for sufficiently small $\tau = \tau^{(\delta)}_\eps$. Furthermore, this existence time $\tau_\eps$ depends on $\mM_\eps$ and $\rR$ only (and hence $\mM_\eps$ only), and is independent of $(\eps, \delta)$ if $\mM_\eps = \mM^{(\delta)}_\eps$ is uniformly bounded. 

\begin{flushleft}
\textit{Step 3. }
\end{flushleft}
Now, suppose $\mM_\eps = \mM_\eps^{(\delta)}$ is uniformly bounded in $(\eps, \delta)$ when restricted to $\bB$. Then the local existence time $\tau$ from \eqref{e:kpz_local_existence_time} is also independent of $(\eps, \delta)$. Furthermore, for this $\tau$, we have the bound
\begin{equation*}
    \|\vec{u}_\eps\|_{\xX_{\tau,\eps}} \lesssim 1 + \mM_\eps\;.
\end{equation*}
This completes the proof. 
\end{proof}

\subsection{Convergences of the solution and flux}

Let $\vec{u}_\eps = (u_\eps, u_\eps^{\#}) \in \xX_{\tau_\eps,\eps}$ be the solution to the system \eqref{eq:kpz_system}. Define
\begin{equation*}
    \vec{U}_{\eps} := (u_{\eps}, \, \lL_{\eps}u_{\eps}, \, u_\eps^{\#}, \, \lL_{\eps}u_{\eps}^{\#})\;,
\end{equation*}
and
\begin{equation*}
    \widetilde{\xX}_{T} := L_{T}^{\infty}{\cC^{\frac{1}{2}-\eta}_{x}}  \times  \cC_{\sigma, T}^{2\eta} \cC_{x}^{-1+4\eta}  \times  L_T^\infty \cC_x^{\frac{1}{2}-\eta}  \times  \cC_{1+4\eta, T}^{2\kappa} \cC_{x}^{-\frac{1}{2} + \eta}\;.
\end{equation*}
We further let $\widetilde{\xX}_T^{+}$ be the space same as $\widetilde{\xX}_T$, except that one replaces the two $L_T^\infty \cC_x^{\frac{1}{2}-\eta}$ components by the slightly stronger spaces $L_T^\infty \cC_x^{\frac{1-\eta}{2}}$. 

For $\eps = 0$, it is well known that $\vec{U}_0 \in \widetilde{\xX}_{T}^{+}$ can be defined up to arbitrary time $T>0$ provided that the initial data are in $\cC_x^{\frac{1-\eta}{2}}$. We have the following theorem regarding the convergence.

\begin{thm} \label{thm:kpz_convergence_homo_flux}
Fix $T>0$. Suppose
\begin{equation*}
    \big( u_\eps(0), u_\eps^{\#}(0) \big) \in \cC^{\frac{1}{2}-\eta} \times \cC^{\frac{1}{2}-\eta}\;, \quad \big( u_0(0), u_0^{\#}(0) \big) \in \cC^{\frac{1-\eta}{2}} \times \cC^{\frac{1-\eta}{2}}\;,
\end{equation*}
and
\begin{equation*}
    \mM := 1 + \|\vec{U}_0\|_{\widetilde{\xX}_{T}^{+}} + \|\<1bk>_0\|_{L_T^\infty \cC_x^{-\frac{1}{2}-\kappa}} + \sup_{\eps\in \N^{-1}} \big( \|u_\eps(0)\|_{\frac{1}{2}-\eta} + \|u_\eps^{\#}(0)\|_{\frac{1}{2}-\eta} + \|\Upsilon_\eps\|_{\yY_{T}} \big) < +\infty\;.
\end{equation*}
Then for every $\eps\in \N^{-1}$, the solution $(u_\eps, u_\eps^{\#})$ to \eqref{eq:kpz_system} exists in $\xX_{T,\eps}$ up to the same time $T$, and there exists $\theta>0$ such that
\begin{equation} \label{e:kpz_convergence_sol}
    \|\vec{U}_\eps - \vec{U}_0\|_{\widetilde{\xX}_T} \lesssim \eps^\theta + \|u_\eps(0) - u_0(0)\|_{\frac{1}{2}-\eta} + \|u_\eps^{\#}(0) - u_0^{\#}(0)\|_{\frac{1}{2}-\eta} + \|\Upsilon_\eps - \lambda \bullet \Upsilon_0\|_{\yY_{T}}\;.
\end{equation}
Furthermore, the fluxes also satisfy the bound
\begin{equation} \label{e:kpz_convergence_flux}
    \begin{split}
    &\phantom{111}\|A_{\eps}\d_{x}u_{\eps}-\bar{A}\d_{x}u_{0}\|_{\cC^{2\eta}_{\sigma,T} \cC_x^{3\eta}} + \|A_{\eps}\d_{x}u^{\#}_{\eps}-\bar{A}\d_{x}u_{0}^{\#}\|_{\cC_{1+4\eta,T}^{2\kappa} \cC_x^{\frac{1+\eta}{2}}}\\
    &\lesssim \eps^{\theta} + \|u_{\eps}(0)-u_{0}(0)\|_{\frac{1}{2}-\eta} + \|u_{\eps}^{\#}(0)-u_{0}^{\#}(0)\|_{\frac{1}{2}-\eta} + \|\Upsilon_{\eps}-\lambda \bullet \Upsilon_{0}\|_{\yY_T}\;.
    \end{split}
\end{equation}
All proportionality constants above depend on $\mM$ and $T$ only. 
\end{thm}
\begin{proof}
We first note that by Lemma~\ref{le:convergence_derivative_a_1d}, there exists $\theta>0$ such that
\begin{equation*}
    \begin{split}
    \| A_\eps \d_x u_\eps - \bar{A} \d_x u_0 \|_{\cC_{\sigma,T}^{2\eta} \cC_x^{3\eta}} &\lesssim \eps^\theta \|\lL_0 u_0\|_{\cC_{\sigma,T}^{2\eta} \cC_x^{-1+4\eta}} + \|\lL_\eps u_\eps - \lL_0 u_0\|_{\cC_{\sigma,T}^{2\eta} \cC_x^{-1+4\eta}}\\
    \| A_\eps \d_x u_\eps^{\#} - \bar{A} \d_x u_0^{\#} \|_{\cC_{1+4\eta,T}^{2\kappa} \cC_x^{\frac{1+\eta}{2}}} &\lesssim \eps^\theta \|\lL_0 u_0^{\#}\|_{\cC_{1+4\eta,T}^{2\kappa} \cC_{x}^{-\frac{1}{2}+\eta}} + \|\lL_\eps u_\eps^{\#} - \lL_0 u_0^{\#}\|_{\cC_{1+4\eta, T}^{2\kappa} \cC_{x}^{-\frac{1}{2}+\eta}}\;,
    \end{split}
\end{equation*}
where the proportionality constants are independent of $\eps$. Hence, the convergence of fluxes \eqref{e:kpz_convergence_flux} follows from the bound \eqref{e:kpz_convergence_sol}. It then remains to prove the convergence of the solutions. 

By Theorem~\ref{thm:kpz_fixed_pt} and the bound \eqref{e:kpz_sol_uniform_bound}, there exists a universal $C_*>1$ and $\tau_0>0$ independent of $\eps\in \N^{-1}$ such that $\vec{U}_\eps$ exists up to time $\tau_0$ with the bound
\begin{equation} \label{e:kpz_convergence_local_uniform}
    \|\vec{U}_\eps\|_{\widetilde{\xX}_{\tau_0}} < C_* \mM\;.
\end{equation}
We first prove the bound \eqref{e:kpz_convergence_sol} for sufficiently small $\tau \in (0, \tau_0)$ specified later, and then iterate it to $T>0$ where the solution $\vec{U}_0$ is defined. Finally, convergence of the flux follows automatically from the bounds on $\|\vec{U}_\eps - \vec{U}_0\|_{\widetilde{\xX}_T}$. Further smallness of $\tau$ will be needed and specified later in the proof (but still depending on $\mM$ only). 

We need to compare $\vec{U}_\eps$ and $\vec{U}_0$ for all their four components. For terms with the initial data, it follows from Lemma~\ref{le:convergence_heat_normal} that
\begin{equation} \label{e:kpz_convergence_initial_1}
    \begin{split}
    \| e^{t \lL_\eps} u_\eps(0) - e^{t \lL_0} u_0(0) \|_{L_\tau^\infty \cC_x^{\frac{1}{2}-\eta}} &\lesssim \eps^\theta \mM + \|u_\eps(0) - u_0(0)\|_{\frac{1}{2}-\eta}\;,\\
    \| \lL_\eps e^{t \lL_\eps} u_\eps(0) - \lL_0 e^{t \lL_0} u_0(0) \|_{\cC_{\sigma, \tau}^{2\eta} \cC_{x}^{-1+4\eta}} &\lesssim \eps^\theta \mM + \|u_\eps(0) - u_0(0)\|_{\frac{1}{2}-\eta}\;,
    \end{split}
\end{equation}
and
\begin{equation} \label{e:kpz_convergence_initial_2}
    \begin{split}
    \| e^{t \lL_\eps} u_\eps^{\#}(0) - e^{t \lL_0} u_0^{\#}(0) \|_{L_\tau^\infty \cC_x^{\frac{1}{2}-\eta}} &\lesssim \eps^\theta \mM + \|u_\eps^{\#}(0) - u_0^{\#}(0)\|_{\frac{1}{2}-\eta}\;,\\
    \| \lL_\eps e^{t \lL_\eps} u_\eps^{\#}(0) - \lL_0 e^{t \lL_0} u_0^{\#}(0) \|_{\cC_{1+4\eta, \tau}^{2\kappa} \cC_{x}^{-\frac{1}{2}+\eta}} &\lesssim \eps^\theta \mM + \|u_\eps^{\#}(0) - u_0^{\#}(0)\|_{\frac{1}{2}-\eta}\;,
    \end{split}
\end{equation}
where we have replaced $\|u_0(0)\|_{\frac{1-\eta}{2}}$ and $\|u_0^{\#}(0)\|_{\frac{1-\eta}{2}}$ by $\mM$. This is the place where we used that $\big( u_0(0), u_0^{\#}(0) \big)$ have slightly better regularity than $\frac{1}{2}-\eta$. 

All other terms on the right hand sides of \eqref{eq:kpz_system} are of the form $\iI_\eps (h_\eps)$ for $h_\eps$ being one of the followings: 
\begin{equation} \label{eq:kpz_rhs_list}
    \Lambda_\eps (u_\eps) \prec_\eps \<1bk>_\eps\;, \quad \Lambda_\eps (u_\eps) \succ_\eps \<1bk>_\eps\;, \quad \frac{\b_\eps}{A_\eps} \cdot (\Lambda_\eps u_\eps)^2\;, \quad \kK_\eps^{(1)}(u_\eps)\;, \quad \kK_\eps^{(2)}(\vec{u}_\eps)\;.
\end{equation}
Only the quadratic term $\frac{\b_\eps}{A_\eps} \cdot (\Lambda_\eps u_\eps)^2$ is a continuous-in-time evolution in some distribution space, while all others have weighted $\cC^{-\kappa}$ regularity in time. Hence we deal with the quadratic term first. By Lemmas~\ref{lem:kpz_quadratic}, \ref{le:convergence_heat_below} and~\ref{le:convergence_heat_spacetime_normal}, we have
\begin{equation*}
    \begin{split}
    &\phantom{111}\left\| \iI_\eps \big( (\b_\eps / A_\eps) \cdot (\Lambda_\eps u_\eps)^2 \big) - \iI_0 \big( \Pi_0 (\b / A) \cdot (\bar{A}\Lambda_0 u_0)^2 \big) \right\|_{L_{\tau}^\infty \cC_x^{\frac{1}{2}-\eta}}\\
    &\lesssim \tau^\theta \mM \Big( \eps^\theta \mM + \|\vec{U}_\eps - \vec{U}_0\|_{\widetilde{\xX}_{\tau}} + \|\Upsilon_\eps - \lambda \bullet \Upsilon_0\|_{\yY_T} \Big)\;,
    \end{split}
\end{equation*}
and
\begin{equation*}
    \begin{split}
    &\phantom{111}\left\| \lL_\eps \iI_\eps \big( (\b_\eps / A_\eps) \cdot (\Lambda_\eps u_\eps)^2 \big) - \lL_0 \iI_0 \big( \Pi_0 (\b / A) \cdot (\bar{A}\Lambda_0 u_0)^2 \big) \right\|_{\cC_{\sigma,\tau}^{2\eta} \cC_x^{4\eta-1}}\\
    &\lesssim \tau^\theta \mM \Big( \eps^\theta \mM + \|\vec{U}_\eps - \vec{U}_0\|_{\widetilde{\xX}_{\tau}} + \|\Upsilon_\eps - \lambda \bullet \Upsilon_0\|_{\yY_T} \Big)\;.
    \end{split}
\end{equation*}
All other terms $h_\eps$ in the list \eqref{eq:kpz_rhs_list} are spacetime distributions. To control the differences, we write
\begin{equation*}
    \begin{split}
    \iI_\eps (h_\eps) - \iI_0 (h_0) &= \iI_\eps (h_\eps - h_0) + (\iI_\eps - \iI_0) h_0\;,\\
    \lL_\eps \iI_\eps (h_\eps) - \lL_0 \iI_0 (h_0) &= \lL_\eps \iI_\eps (h_\eps - h_0) + (\lL_\eps \iI_\eps - \lL_0 \iI_0) h_0\;.
    \end{split}
\end{equation*}
The $\cC_{\sigma',T}^{-\kappa} \cC_x^{\gamma}$-norm of the difference $h_\eps - h_0$ can be controlled by corresponding lemmas in Section~\ref{sec:kpz_preliminary}, where $\gamma \in (-1,1)$ depending on the actual term. Hence, the quantities
\begin{equation*}
    \iI_\eps (h_\eps - h_0) \quad \text{and} \quad \lL_\eps \iI_\eps (h_\eps - h_0)
\end{equation*}
can be controlled using Proposition~\ref{prop:Young_heat_convolution} as well as the difference bound for $h_\eps - h_0$. 

On the other hand, we note that although $h_\eps$ is uniformly bound in $\eps$ only as spacetime distributions, $h_0$ is actually a continuous-in-time evolution in corresponding space of distributions. In particular, $h_0$ satisfies the bounds in the preliminary lemmas in Section~\ref{sec:kpz_preliminary} with $\cC_{\sigma',\tau}^{-\kappa}$ strengthened to $L_{\sigma',\tau}^\infty$. This is because the stochastic object $\<1bk>_0 \in L_T^\infty \cC_{x}^{-\frac{1}{2}-\kappa}$. Hence, the quantities
\begin{equation*}
    (\iI_\eps - \iI_0) h_0 \quad \text{and} \quad (\lL_\eps \iI_\eps - \lL_0 \iI_0) h_0
\end{equation*}
can be controlled with Lemmas~\ref{le:convergence_heat_below} and~\ref{le:convergence_heat_spacetime_normal}.

Combining all of the above together, we see there exist $C_1, C_2 > 0$ and $\theta > 0$ such that
\begin{equation*}
    \begin{split}
    \|\vec{U}_\eps - \vec{U}_0\|_{\widetilde{\xX}_\tau} \leq &C_1 \Big( \eps^\theta \mM + \|u_\eps(0) - u_0(0)\|_{\frac{1}{2}-\eta} + \|u_\eps^{\#}(0) - u_0^{\#}(0)\|_{\frac{1}{2}-\eta} \Big)\\
    &+ C_2 \tau^{\theta} \mM^2 \big( \eps^\theta \mM + \|\Upsilon_\eps - \lambda \bullet \Upsilon_0\|_{\yY_T} + \|\vec{U}_\eps - \vec{U}_0\|_{\widetilde{\xX}_\tau} \big)\;.
    \end{split}
\end{equation*}
Taking $\tau$ sufficiently small such that $C_2 \tau^\theta \mM^2 < \frac{1}{2}$ and that \eqref{e:kpz_convergence_local_uniform} is still satisfied, we can absorb $\|\vec{U}_\eps - \vec{U}_0\|_{\widetilde{\xX}_\tau}$ into the left hand side, and deduce that
\begin{equation} \label{e:kpz_convergence_local}
    \begin{split}
    \|\vec{U}_\eps - \vec{U}_0\|_{\widetilde{\xX}_\tau} \lesssim \mM^2 \Big( \eps^\theta \mM &+ \|u_\eps(0) - u_0(0)\|_{\frac{1}{2}-\eta}\\
    &+ \|u_\eps^{\#}(0) - u_0^{\#}(0)\|_{\frac{1}{2}-\eta} + \|\Upsilon_\eps - \lambda \bullet \Upsilon_0\|_{\yY_T} \Big)\;.
    \end{split}
\end{equation}
Since the local existence time $\tau>0$ above is uniform in $\eps$, one can re-start the procedure and iterate the bound \eqref{e:kpz_convergence_local} finitely many times to reach time $T$. This proves the convergence \eqref{e:kpz_convergence_sol} and hence the theorem as well. 
\end{proof}

\subsection{Proof of Theorem~\ref{thm:kpz_overall}}

We are now in position to combine all ingredients above to prove Theorem~\ref{thm:kpz_overall}. We include the mollification parameter $\delta$ in this final argument. 

\begin{proof} [Proof of Theorem~\ref{thm:kpz_overall}]
Let $h_\eps(0) \in \cC^{\frac{1}{2}-\eta}$ and $h_0(0) \in \cC^{\frac{1-\eta}{2}}$ be given. Let $(u_\eps^{(\delta)}, u_\eps^{(\delta); \#}$ be the solution to the system \eqref{eq:kpz_system} with initial data
\begin{equation*}
    u_\eps^{(\delta)}(0) = h_\eps(0) - \<1gk>_\eps^{(\delta)}(0) - \<20gk>_\eps^{(\delta)}(0) - 2 \, \<210gk>_\eps^{(\delta)}(0)
\end{equation*}
and $u_\eps^{(\delta); \#}(0) = 0$. Then, it follows from the derivation of the system \eqref{eq:kpz_system} that
\begin{equation} \label{e:kpz_sol_decomposition}
    h_\eps^{(\delta)} = \<1gk>_\eps^{(\delta)} + \<20gk>_\eps^{(\delta)} + 2 \, \<210gk>_\eps^{(\delta)} + u_\eps^{(\delta)} + \int_{0}^{t} \Pi_0 \big( \xi + \<2k>_\eps^{(\delta)} +2 \, \<21k>_\eps^{(\delta)} \big)(r) {\rm d}r
\end{equation}
solves the original equation \eqref{e:kpz_homo_full} with initial data $h_\eps(0)$. Then we have
\begin{equation*}
    \| h_{ \eps }^{ ( \delta ) } - h_0 \|_{ L_T^\infty \cC_x^{\frac{1}{2}-\eta}} \lesssim \|u_\eps^{(\delta)} - u_0\|_{L_T^\infty \cC_x^{\frac{1}{2}-\eta}} + \|\hH_\eps^{(\delta)} - \lambda \bullet \hH_0\|_{\yY_T^h} + \|\Upsilon^{(\delta)}_\eps - \lambda \bullet \Upsilon_0\|_{\yY_T}\;.
\end{equation*}
There is the additional term on the difference between $\Upsilon_\eps^{(\delta)}$ and $\lambda \bullet \Upsilon_0$ here since the term $\<210gk>_\eps^{(\delta)}$ is not included in $\hH_\eps^{(\delta)}$ but can be completely determined from $\<210rk>_\eps^{(\delta)} \in \Upsilon_\eps^{(\delta)}$. The desired bound for $h_{\eps}^{(\delta)} - h_0$ then follows from Theorem~\ref{thm:kpz_convergence_homo_flux}. 

As for the flux, from \eqref{e:kpz_sol_decomposition}, we have the expression
\begin{equation*}
    A_\eps \d_x h_\eps^{(\delta)} = \<1rk>_\eps^{(\delta)} + \<20rk>_\eps^{(\delta)} + 2 \, \<210rk>_\eps^{(\delta)} + A_\eps \d_x u_\eps^{(\delta)}\;,
\end{equation*}
which gives
\begin{equation*}
    \begin{split}
    &\phantom{111}\big\| A_\eps \d_x h_\eps^{(\delta)} - \bar{A} \d_x h_0 \big\|_{\cC_{\sigma,T}^{-\kappa} \cC_x^{-\frac{1}{2}-2\kappa}}\\
    &\lesssim \|\fF_\eps^{(\delta)} - \lambda \bullet \fF_0\|_{\yY_T^{\flux}} + \|\Upsilon_\eps^{(\delta)} - \lambda \bullet \Upsilon_0\|_{\yY_T} + \|A_\eps \d_x u_\eps^{(\delta)} - \bar{A} \d_x u_0\|_{\cC_{\sigma,T}^{2\eta} \cC_x^{3\eta}}\;.
    \end{split}
\end{equation*}
Note that $\fF_\eps^{(\delta)}$ only includes $\<1rk>_\eps^{(\delta)}$, and the convergence of $\<20rk>_\eps^{(\delta)}$ is encoded in $\Upsilon_\eps^{(\delta)}$. The desired claim then follows from directly from the bound for $A_\eps \d_x u_\eps^{(\delta)} - \bar{A} \d_x u_0$ in Theorem~\ref{thm:kpz_convergence_homo_flux}. This completes the proof of the theorem. 
\end{proof}

%\section{Periodic homogenisation for the stochastic heat equation}\label{sec: Periodic homogenisation for SHE}
%\include{SHE}    

\appendix

\section{An oscillation lemma}

\begin{lem} \label{le:convergence_oscillation}
For every $\alpha>0$, we have
\begin{equation*}
    \| f( \cdot / \eps ) \|_{\cC^{-\alpha}} \lesssim  \eps^\alpha \|f\|_{L^{\infty}}
\end{equation*}
uniformly over $f \in L^\infty(\T^d)$ with $\Pi_0 f = 0$ and $\eps \in \N^{-1}$. 
\end{lem}
\begin{proof}
Write $f_\eps = f(\cdot / \eps)$ for simplicity. Fix $\eps = \frac{1}{N}$. 
Since $\Pi_0 f = 0$, the Fourier series of $f_\eps$, which we denote by $\widehat{f_\eps}$, is supported on the set $\{Nk\}_{k \in \Z^d \setminus \{0\}}$. 

Let $\{\Delta_j\}_{j \geq -1}$ be the Littlewood-Paley block operators associated to the Laplacian $\Delta$. Hence, $\Delta_j f_\eps = 0$ whenever $2^j < \frac{N}{16}$. For $j$ such that $2^j \geq \frac{N}{16}$, we have
\begin{equation*}
    \|\Delta_j f_\eps\|_{L^\infty} \lesssim \|f_\eps\|_{L^\infty} \lesssim N^{-\alpha} 2^{\alpha j} \|f\|_{L^\infty}\;,
\end{equation*}
which in turn implies
\begin{equation*}
    \|f_\eps\|_{\cC^{-\alpha}} = \sup_{j \geq -1} \Big( 2^{-\alpha j} \|\Delta_j f_\eps\|_{L^\infty} \Big) \lesssim \eps^\alpha \|f\|_{L^\infty}\;.
\end{equation*}
This completes the proof. 
\end{proof}

\section{Results for elliptic periodic homogenisation}
\label{sec:elliptic_homogenisation}

In this section, we fix a coefficient matrix $A$ satisfying Assumption~\ref{as:A}. Let $A_{\eps} = A(\cdot/\eps)$ for $\eps\in \N^{-1}$, and $\bar{A}$ be the homogenised matrix given in \eqref{def:bar_a_homogenised_matrix}. Define $\lL_{\eps} := \div(A_{\eps}\nabla)$ and $\lL_0 = \div (\bar{A} \nabla)$. The following result is standard (see \cite[Theorem~3.4.5]{Shen2018} for example).

\begin{lem} \label{lem:elliptic_homo_convergence}
We have
\begin{equation*}
    \|(\lL_\eps^{-1} - \lL_0^{-1}) f\|_{L^2(\T^d)} \lesssim \eps \|f\|_{L^2(\T^d)}
\end{equation*}
uniformly over $f \in L^2(\T^d)$ with $\Pi_0 f = 0$. 
\end{lem}

We have the following corollary. 

\begin{cor} \label{cor:resolvent_convergence}
We have the bound
    \begin{equation*}
        \| (\id - \lL_\eps)^{-1} - (\id - \lL_0)^{-1} \|_{L^{2}(\T^{d})\to L^{2}(\T^{d})} \lesssim \eps
    \end{equation*}
uniformly in $\eps \in \N^{-1} \cup \{0\}$. 
\end{cor}
\begin{proof}
Since $(\id - \lL_\eps)^{-1} \1 = (\id - \lL_0)^{-1} \1 = \1$, it suffices to consider the action of these operators on the orthogonal complement of constant functions in $L^2(\T^d)$. We denote this orthogonal complement by $\eE$. 

On $\eE$, the operator $(\id - \lL_\eps)^{-1}$ can be represented by
\begin{equation*}
    (\id - \lL_\eps)^{-1} = F(\lL_\eps^{-1})
\end{equation*}
with $F(x) = \frac{x}{x-1}$, which is analytic on $\mathbf{C} \setminus \{1\}$. Since the spectrum of $\lL_\eps^{-1}$ on $\eE$ is completely contained in $[-c, 0]$ for some $c>0$ independent of $\eps$ (see proof of Lemma~\ref{lem: pointwise estimate for lL G eps} below), the desired bound follows from contour integration. 
\end{proof}

The following result is from \cite[Lemma~4.4.6]{Shen2018}.

\begin{lem} \label{lem:elliptic_homo_convergence_gradient}
We have
\begin{equation*}
    \| \big( \nabla \lL_\eps^{-1} - (\id + \nabla\chi(\cdot/\eps)) \nabla \lL_0^{-1} \big) f\|_{L^{\infty}(\T^{d})} \lesssim \sqrt{\eps} \|f\|_{\cC^{\frac{1}{2}}}
\end{equation*}
uniformly over $f \in \cC^{\frac{1}{2}}(\T^d)$ with $\Pi_0 f = 0$. 
\end{lem}

Let $G_{\eps}$ be the Green's function of $-\lL_{\eps}$ on $\T^{d}$. The following results are from \cite[Theorem~13]{FLinCompactness} and \cite[Theorem~5.6.1]{Shen2018}.

\begin{lem} \label{lem:elliptic_homo_Green_pointwise}
We have the pointwise bound
\begin{equation*} 
|G_\eps(x,y)| \lesssim \left \{
\begin{array}{rl}
&\log \big(2 + |x-y|^{-1}\big)\,, \quad d=2\\
&|x-y|^{-d+2}\,, \qquad \qquad d \geq 3
\end{array} \right.,
\end{equation*}
and the gradient bounds
    \begin{equation*}
        |\nabla_{x} G_{\eps}(x,y)| + |\nabla_{y}G_{\eps}(x,y)| \lesssim |x-y|^{-d+1}\;,\quad d\geq 2\;.
    \end{equation*}
Finally, the cross derivative also satisfies the pointwise bound
    \begin{equation*}
        |\nabla_{x}\nabla_{y} G_{\eps}(x,y)|\lesssim |x-y|^{-d}\;.
    \end{equation*}
All proportionality constants are uniform in $\eps \in \N^{-1}$ and $x,y \in \T^d$. 
\end{lem}

\begin{lem} \label{lem:elliptic_homo_Green_pointwise_expansion}
Fix $d \geq 2$. We have the bounds
    \begin{equation*}
        |G_{\eps}(x,y)-G_{0}(x,y)| \lesssim \frac{\eps}{|x-y|^{d-1}}\;, \quad |\nabla_x G_{\eps}(x,y) - \Phi_\eps (x) \nabla_x G_{0}(x,y)| \lesssim \frac{\eps}{|x-y|^d}\;.
    \end{equation*}
Both proportionality constants are uniform in $\eps \in \N^{-1}$ and $x,y \in \T^d$. 
\end{lem}

\begin{lem} \label{lem:elliptic_homo_Green_remainder}
For every $\nu' \in (0,1)$, there exist $\nu, \theta > 0$ such that
\begin{equation*}
    \| \big( \nabla \lL_\eps^{-1} - \Phi_\eps \nabla \lL_0^{-1} \big) f \|_{\cC^{\nu}(\T^d)} \lesssim \eps^\theta \|f\|_{\cC^{\nu'-1}(\T^d)}\;.
\end{equation*}
The proportionality constant is uniform in $\eps \in \N^{-1}$ and $x,y \in \T^d$. 
\end{lem}

\section{Bounds on parabolic Green's functions}
\label{sec: parabolic_Green_function}

Same as before, we fix a coefficient matrix $A$ satisfying Assumption~\ref{as:A}, and let $A_\eps = A(\cdot/\eps)$. Let $\bar{A}$ be the homogenised matrix. Let $\lL_\eps = \div (A_\eps \nabla)$ and $\lL_0 = \div (\bar{A} \nabla)$. 

Let $Q_\eps$ be the Green's function of $\partial_t-\lL_\eps$ in $\R^+ \times \T^d$. In this section, we aim at pointwise estimates for Green's functions $Q_\eps$. Let $\Gamma_\eps$ be the fundamental solution of $\partial_t-\lL_\eps$ in $\R^+ \times \R^d$. For $x,y\in[0,1)^{d}$, we have the relation
\begin{equation} \label{e:periodise}
    Q_{\eps}( t ; x , y )=\sum_{ k\in\Z^d } \Gamma_{\eps}(t; {x}, {y} + k )\;.
\end{equation}
The following proposition for $\Gamma_\eps$ is well known. 

\begin{prop} \label{prop:fundamental_Gaussian_bounds}
There exists $c>0$ such that the fundamental solution $\Gamma_\eps$, its gradients and cross derivatives satisfy the bounds
\begin{equation*}
    \begin{split}
    |\Gamma_\eps(t;x,y)| &\lesssim t^{-\frac{d}{2}} e^{-\frac{c|x-y|^2}{t}}\;,\\
    |\nabla_x \Gamma_\eps(t;x,y)| + |\nabla_y \Gamma_\eps(t;x,y)| &\lesssim t^{-\frac{d+1}{2}} e^{-\frac{c|x-y|^2}{t}}\;,\\
    |\nabla_x \nabla_y \Gamma_\eps(t;x,y)| &\lesssim t^{-\frac{d+2}{2}} e^{-\frac{c|x-y|^2}{t}}\;.
    \end{split}
\end{equation*}
All proportionality constants are independent of $\eps \in \N^{-1}$, $x,y \in \R^d$ and $t \in \R^+$. 
\end{prop}
\begin{proof}
This is the content of \cite[(1-6)]{GengShen2020} and \cite[Theorem 2.7]{GengShen2020}. 
\end{proof}

\begin{lem} \label{lem:summation_Gaussian_bound}
For every $c>0$, we have
	\begin{equation*}
		\sum_{n\in\Z^{d}} e^{-\frac{c|z+n|^2}{t}} \lesssim (1 + t^{\frac{d}{2}}) e^{-\frac{c |z|^2}{t}}
	\end{equation*}
uniformly over $z \in [-\frac{1}{2}, \frac{1}{2}]^d$ and $t \in \R^+$. The proportionality constant depends on $c$. 
\end{lem}
\begin{proof}
Note that $|z| \leq |z+n|$ for all $z \in [-\frac{1}{2}, \frac{1}{2}]^d$ and $n \in \N$, we immediately have
	\begin{equation*}
		\sum_{|n|< 10d} e^{-\frac{c|z+n|^2}{t}} \lesssim e^{-\frac{c|z|^2}{t}}\;.
	\end{equation*}
For $|n|\geq 10d$, we have $|z+n| \geq \frac{|n|}{2}|n|+|z|$, and hence
	\begin{equation*}
		\sum_{|n| \geq 10d} e^{-\frac{c|z+n|^2}{t}} \lesssim e^{-\frac{c|z|^2}{t}} \sum_{|n|\geq 10d} e^{-\frac{c|n|^{2}}{4t}} \lesssim (1 + t^{\frac{d}{2}}) e^{-\frac{c|z|^2}{t}}\;.
	\end{equation*}
This completes the proof. 
\end{proof}

We are now ready to derive the following Gaussian bounds for the periodic Green's function $Q_\eps$. 

\begin{prop} \label{prop:Green_point}
The Green's function $Q_\eps$ satisfies the pointwise bound
    \begin{equation}\label{eq: estimate 1 for Green_point}
        |Q_\eps ( t ; x , y )| \lesssim (t^{-\frac{d}{2}} +1 ) e^{ - \frac{c |x-y|^2 }{ t } }\;.
    \end{equation}
Its gradients satisfy the bound
    \begin{equation}\label{eq: estimate 2 for Green_point}
        |\nabla_x Q_\eps ( t ; x , y )| +|\nabla_{y} Q_{\eps}( t ; x , y )|\lesssim t^{-\frac{1}{2}} ( t^{-\frac{d}{2}} + 1 )e^{-\frac{c |x-y|^2}{t}}\;.
    \end{equation}
Finally, its cross derivative also satisfies the pointwise bound
    \begin{equation}\label{eq: estimate 3 for Green_point}
        |\nabla_{x} \nabla_{y} Q_{\eps}( t ; x , y )| \lesssim (t^{-\frac{d}{2}}+ 1) e^{-\frac{c |x-y|^2}{t}}\;.
    \end{equation}
All proportionality constants are independent of $\eps \in \N^{-1}$, $x,y \in \T^d$ and $t \in \R^+$. 
\end{prop}
\begin{proof}
We choose representations of $x,y \in \T^d$ such that $x-y \in [-\frac{1}{2}, \frac{1}{2}]^d$. The three bounds then follow immediately from the representation \eqref{e:periodise}, the corresponding bounds for $M_\eps$ in Proposition~\ref{prop:fundamental_Gaussian_bounds} and Lemma~\ref{lem:summation_Gaussian_bound}. 
\end{proof}

We have the following immediate corollary. 

\begin{cor} \label{cor:coefficient_uniform}
For every coefficient matrix $A$ satisfying Assumption~\ref{as:A}, there exist $\Lambda, M>0$ such that $A_\eps \in \sS_d(\Lambda, M)$ for all $\eps \in \N^{-1}$. Furthermore, the value $M$ depends on $A$ via $\|A\|_{\cC^\alpha}$ only. 
\end{cor}

\begin{lem} \label{lem:convolution}
Let $\alpha \in (0,d)$ and $\lambda>0$. We have
\begin{equation} \label{e:convolution_1}
	\int_{\T^d} \frac{2^{dn} |y-z|^{-\alpha}}{(1+2^{n}|x-y|)^{d+\lambda}} {\rm d}y \lesssim  |z-x|^{-\alpha}\;,
\end{equation}
and
\begin{equation} \label{e:convolution_2}
	\int_{\T^d} \frac{2^{dn}}{(1+2^{n}|x-y|)^{d+\lambda}} \cdot \frac{2^{dn}}{(1+2^{n}|z-y|)^{d+\lambda}}{\rm d}y \lesssim \frac{2^{dn}}{(1+2^{n}|x-z|)^{d+\lambda}}\;.
\end{equation}
Both proportionality constants are independent of $n \in \Z$ and $x,z \in \T^d$. 
\end{lem}
\begin{proof}
We start with \eqref{e:convolution_1}. We split the domain of integration into $\{y: |y-z| \geq \frac{|x-z|}{2}\}$ and $\{y: |y-z| < \frac{|x-z|}{2}\}$. We respectively have
\begin{equation*}
    \int_{|y-z| \geq \frac{|x-z|}{2}} \frac{2^{dn} |y-z|^{-\alpha}}{(1+2^{n}|x-y|)^{d+\lambda}} {\rm d}y \lesssim |x-z|^{-\alpha} \int_{\T^d} \frac{2^{dn}}{(1 + 2^n |x-y|)^{d+\lambda}} {\rm d}y \lesssim |x-z|^{-\alpha}\;,
\end{equation*}
and
\begin{equation*}
    \begin{split}
    \int_{|y-z| < \frac{|x-z|}{2}} \frac{2^{dn} |y-z|^{-\alpha}}{(1+2^{n}|x-y|)^{d+\lambda}} {\rm d}y &\lesssim \frac{2^{dn}}{(1 + 2^n |x-z|)^{d+\lambda}} \int_{|y-z| < \frac{|x-z|}{2}} |y-z|^{-\alpha} {\rm d}y\\
    &\lesssim |x-z|^{-\alpha}\;.
    \end{split}
\end{equation*}
They together yield the bound \eqref{e:convolution_1}. We now turn to \eqref{e:convolution_2}. On the domain $\{y: |y-z| < \frac{|x-z|}{2}\}$, we have $|x-y| \sim |x-z|$, and hence the integral on the left hand side of \eqref{e:convolution_2} is bounded by
\begin{equation*}
    \frac{2^{dn}}{(1 + 2^n |x-z|)^{d+\lambda}} \int_{\T^d} \frac{2^{dn}}{(1 + 2^n |z-y|)^{d+\lambda}} {\rm d}y \lesssim \frac{2^{dn}}{(1 + 2^n |x-z|)^{d+\lambda}}\;.
\end{equation*}
The same bound holds on the domain $\{|y-x| < \frac{|x-z|}{2}\}$ just with $x$ and $z$ swapped. Finally, on the domain $\{y: |y-z| \geq \frac{|x-z|}{2}\} \cup \{|y-x| \geq \frac{|x-z|}{2}\}$, the integral is controlled by
\begin{equation*}
    \frac{2^{2dn}}{(1 + 2^n |x-z|)^{2(d+\lambda)}} \lesssim \frac{2^{dn}}{(1 + 2^n |x-z|)^{d+\lambda}}\;.
\end{equation*}
Combining the above gives the bound \eqref{e:convolution_2} and also completes the proof of the theorem. 
\end{proof}

We now introduce the homogeneous Littlewood-Paley blocks $\widetilde{\Delta}_{n,\eps}$ for $n \in \Z$ by
\begin{equation*}
    \widetilde{\Delta}_{n,\eps} := \varphi_0 \big( 2^{-n} \sqrt{-\lL_\eps} \big)\;, 
\end{equation*}
where $\varphi_0$ given in \eqref{e:dyadic_partition} is supported on an annulus. Note that $\widetilde{\Delta}_{n,\eps} = \Delta_{n,\eps}$ for $n \geq 0$, but $\widetilde{\Delta}_{n,\eps}$ are also defined for all negative integers. 

\begin{lem}\label{lem:semigroup_operator_norm}
    There exists $C>0$ such that for every $t>0$, $\alpha\in \R$, $\gamma\in [0,2]$, $ \eps \in \N^{ - 1 }$,
        \begin{equation*}
           \| e^{ t \lL_{ \eps } } - \id \|_{ \cC^{ \alpha ; \eps } \to \cC^{ \alpha - \gamma ; \eps } } \leq \, C \, t^{ \frac{ \gamma }{ 2 } } \; . 
        \end{equation*}
    \end{lem}
    \begin{proof}
        By Green function estimates, we have
        \begin{equation*}
            \| \Delta_{ k , \eps }( e^{ t \lL_{ \eps } } - \id )f \|_{ L^{ \infty } } \lesssim 2^{ - \alpha k }\| f \|_{ {\alpha ; \eps } } \; . 
        \end{equation*}
        On the other hand, we have
        \begin{equation*}
            \Delta_{ k , \eps } ( e^{ t \lL_{ \eps } } - \id ) f \, = \, \int_{ 0 }^{ t } e^{ s \lL_{ \eps } } \lL_{ \eps } \Delta_{ k , \eps } ( f ) { \rm d } s \; .
        \end{equation*}
        Then we have
        \begin{equation*}
            \| \Delta_{ k , \eps }( e^{ t \lL_{ \eps } } - \id ) f \|_{ L^{ \infty } } \lesssim t \, 2^{ ( 2 - \alpha ) k } \, \| f \|_{ \alpha ; \eps } \; .
        \end{equation*}
        By interpolation, we have
        \begin{equation*}
            \| \Delta_{ k , \eps } ( e^{ t \lL_{ \eps } }-\id ) f\|_{ L^{ \infty } } \lesssim t^{ \frac{ \gamma }{ 2 } }\, 2^{ ( \gamma - \alpha ) k } \, \| f \|_{ \alpha ; \eps }\;.
        \end{equation*}
        This completes the proof.
    \end{proof}
    
\begin{lem} \label{lem:heat_convolution_dyadic_bound}
The homogeneous Littlewood-Paley block operators satisfy the bound
    \begin{equation*}
        \left|  \big( e^{t\lL_\eps} \widetilde{\Delta}_{n,\eps} \big)(x,y) \right| \lesssim e^{-ct 2^{2n}} \cdot \frac{2^{dn}}{(1 + 2^n |x-y|)^{d+\lambda}}\;.
    \end{equation*}
The proportionality constant depends on $d$ and $\lambda$ but independent of $\eps \in \N^{-1}$, $n \in \Z$, $x,y \in \T^d$ and $t \in \R^+$. 
\end{lem}
\begin{proof}
By definition, we have
	\begin{equation*}
		\big( e^{t\lL_\eps} \widetilde{\Delta}_{n,\eps} \big)(x,y) = \big( F_{n,t}(2^{-n} \sqrt{-\lL_\eps}) \big)(x,y)\;,
	\end{equation*}
where
	\begin{equation*}
		F_{n,t}(\theta) = e^{-t 2^{2n} \theta^{2}} \varphi_0(\theta)\;.
	\end{equation*}
	The claim then follows from Proposition~\ref{pr:kernel_pointwise_bound} and that $\varphi_0$ is supported on an annulus. 
\end{proof}	 

For $f: \T^d \times \T^d \rightarrow \R$, we write
\begin{equation*}
	(\lL_{\eps,x}f)( x , y ): = \div_{x}( A_{\eps}(x) \nabla_{ x }f)(x,y)\;, \quad (\lL_{\eps,y}f)(x,y) := \div_{y} (A_{\eps}(y) \nabla_{y}f)(x,y)\;.
\end{equation*}
For the Green's function $Q_{\eps}$, by definition, we have
\begin{equation*}
	\lL_{\eps,x}Q_{\eps}(t; x, y) = -\sum_{n\in \N^{+}} e^{-t\lambda_{n,\eps}}\lambda^{2}_{n,\eps}\psi_{n,\eps}(x)\psi_{n,\eps}(y) = (\lL_{\eps,y} Q_\eps)(t;x,y)\;.
\end{equation*}
Hence, we simply write $\lL_\eps Q_\eps$.

\begin{lem} \label{lem: pointwise estimate for lL G eps}
There exists $c>0$ such that one has the pointwise bound
\begin{equation*}
	|(\lL_{\eps}Q_{\eps})( t; x, y )| \lesssim  e^{-ct} (\sqrt{t} + |x-y|)^{-2-d}
\end{equation*}	
uniformly over $\eps \in \N^{-1}$, $x,y \in \T^d$ and $t \in \R^+$. 
\end{lem}
\begin{proof}
Integrating by parts, we see there exists $\Lambda > 0$ independent of $\eps$ such that
\begin{equation} \label{e:reverse_Poincare}
    \lambda_{1,\eps}^2 \|\psi_{1,\eps}\|^{2}_{L^{2}(\T^d)} = \scal{-\lL_{\eps}\psi_{1,\eps},\psi_{1,\eps}} \geq \Lambda^{-1} \|\nabla \psi_{1,\eps}\|^2_{L^{2}(\T^{d})}
\end{equation}
where $\lambda_{1,\eps}>0$ is the smallest \textit{positive} eigenvalue of $\sqrt{-\lL_\eps}$ and $\psi_{1,\eps}$ is its associated $L^2$-normalised eigenfunction. Since the eigenfunction associated to $\lambda_{0,\eps} = 0$ is the constant $1$, and that eigen-functions associated to different eigenvalues are orthogonal, we have $\Pi_0 \psi_{1,\eps} = 0$. Thus, Poincar\'e's inequality implies
\begin{equation} \label{e:Poincare}
    \|\nabla \psi_{1,\eps}\|^2_{L^{2}(\T^{d})}\geq 4\pi^2 \|\psi_{1,\eps}\|^{2}_{L^{2}(\T^{d})}\;.
\end{equation}
Hence, \eqref{e:reverse_Poincare} and \eqref{e:Poincare} together give
\begin{equation*}
    \lambda_{1,\eps} \geq \frac{2 \pi}{\sqrt{\Lambda}}\;,
\end{equation*}
which in turn implies that there exists $K \in \N$ such that for all $\eps \in \N^{-1}$, we have $\widetilde{\Delta}_{n,\eps} \lL_\eps = 0$ for all $n < -K$. Hence, we have
\begin{equation*}
	(\lL_{\eps}Q_{\eps})( t; x, y ) = \sum_{n=-K}^{\infty} \Big(\sum_{\ell=n-1}^{n+1} \lL_{\eps} \widetilde{\Delta}_{\ell,\eps} \Big) (e^{t\lL_{\eps}} \widetilde{\Delta}_{n,\eps})( x, y )\;.
\end{equation*}
By Lemmas~\ref{lem:blocks_derivatives_pointwise}, ~\ref{lem:convolution} and~\ref{lem:heat_convolution_dyadic_bound}, there exist $c_1, c_2 > 0$ dependent on $K$ such that
\begin{equation*}
    \begin{split}
	|\lL_{\eps}Q_{\eps}( t; x, y )| &\lesssim \sum_{n \geq -K} \int_{ \T^{ d } }\frac{2^{ (d+2) n }}{(1+2^n |z-x|)^{d+3}}\cdot \frac{2^{ d n }e^{- c t 2^{2n}}}{(1+2^n |z-y|)^{d+3}}{\rm d}z\\
    &\lesssim e^{-c_1 t} \cdot \sum_{n \geq -K} \frac{2^{(d+2)n} e^{-c_2 t 2^{2n}}}{(1 + 2^n |x-y|)^{d+3}}\;.
    \end{split}
\end{equation*}
We fix $t \in \R^+$ and $x,y \in \T$, and seek bounds for the above sum that are independent of $t$, $x$ and $y$. First, for $x \neq y$, let $k \in \N$ be the unique integer such that $2^{-k-1} \leq |x-y| < 2^{-k}$. Then 
\begin{equation*}
	\frac{ 2^{ (d+2) n } e^{-c_2 t 2^{2n}} }{ ( 1+2^{n} |x-y| )^{d+3} } \leq 2^{(d+2)(k+1)} \cdot \frac{2^{(d+2)(n-k-1)}}{(1 + 2^{n-k-1})^{d+3}}\;.
\end{equation*}
Summing the above bound over $n \in \Z$ gives
\begin{equation} \label{e:G_sg_space}
    |(\lL_\eps Q_\eps)(t;x,y)| \lesssim e^{-c_1 t} \cdot 2^{(d+2)(k+1)} \lesssim e^{-c_1 t} \cdot |x-y|^{-(d+2)}\;.
\end{equation}
On the other hand, for every $t>0$, there exists $\ell \in \Z$, such that $2^{\ell-1} < \sqrt{t} \leq 2^{\ell}$. We then have
\begin{equation*}
    \frac{ 2^{ (d+2) n } e^{-c_2 t 2^{2n}} }{ ( 1+2^{n} |x-y| )^{d+3}} \leq 2^{ (d+2)n} e^{ - c_2 2^{ 2 (n+ \ell -1) } } \leq 2^{-(d+2)(\ell-1)} \cdot 2^{(d+2)(n+\ell-1)} e^{-c_2 2^{2(n+\ell-1)}}\;.
\end{equation*}
Summing the above bound over $\ell \in \Z$ gives
\begin{equation} \label{e:G_sg_time}
    |(\lL_\eps Q_\eps)(t;x,y)| \lesssim e^{-c_1 t} \cdot 2^{-(d+2)\ell} \lesssim e^{- c_1 t} \cdot t^{-\frac{d+2}{2}}\;.
\end{equation}
Combining \eqref{e:G_sg_space} and \eqref{e:G_sg_time} gives the desired claim and thus completes the proof. 
\end{proof}

Similarly, we also have the following bounds. Its proof is very similar as above, and hence we omit the details. 
	
\begin{lem} \label{eq: global point wise estimate for Green function}
There exists $c>0$ such that
\begin{equation*}
	|Q_{\eps}( t ; x , y )-1| \lesssim e^{- c t }(\sqrt{t}+| x - y | )^{-d}\;.
\end{equation*}
Furthermore, the gradients and cross derivatives also satisfy the bounds
\begin{equation*}
	|\nabla_{ x } Q_{\eps}( t ; x , y )| + |\nabla_{ y } Q_{\eps}( t ; x , y )| \lesssim e^{- c t }(\sqrt{t}+| x - y | )^{-d-1}\;,
\end{equation*}
and
\begin{equation*}\begin{split}
	|\nabla_{ x } \nabla_{y} Q_{\eps}( t ; x , y )| \lesssim e^{- c t }(\sqrt{t}+| x - y | )^{-d-2}\;,\\
	|\nabla_{ x } \lL_{\eps,y} Q_{\eps}( t ; x , y )| \lesssim  e^{- c t }(\sqrt{t}+| x - y | )^{-d-3}\;.
	\end{split}
\end{equation*}
All proportionality constants above are independent of $\eps \in \N^{-1}$, $x,y \in \T^d$ and $t \in \R^+$. 
\end{lem}

\section{Young integrals}\label{sec:Young_integral}

\subsection{Definition and multi-linear operation}

For $\beta \in (0,1)$ and $\sigma \in (0,2)$, recall from \eqref{eq:spacetime_Holder_norm_weight} the weighted H\"older norm
\begin{equation*}
    \|f\|_{\cC_{\sigma,T}^\beta \zZ} := \|f\|_{L_{\sigma,T}^\infty \zZ} + \sup_{0<s<t<T} \left( s^{\frac{\sigma}{2}} \cdot \frac{\|f(t) - f(s)\|_{\zZ}}{(t-s)^\beta} \right)\;.
\end{equation*}
We define the $\cC_{\sigma,T}^{-\beta} \zZ$-norm for $f \in L_{\sigma,T}^\infty \zZ$ by
\begin{equation} \label{eq:spacetime_negative_Holder_norm_weight}
    \|f\|_{\cC_{\sigma,T}^{-\beta} \zZ} := \sup_{t \in [0,T]} \left( t^{\frac{\sigma}{2}} \left\| \int_0^t f(r) {\rm d}r \right\|_{\zZ} \right) + \sup_{0<s<t<T} \bigg( s^{\frac{\sigma}{2}} \cdot \frac{\big\| \int_s^t f(r) {\rm d}r \big\|_{\zZ}}{(t-s)^{1-\beta}} \bigg)\;.
\end{equation}
The above time integrals of $f$ are well defined since $f \in L_{\sigma,T}^\infty \zZ$. The space $\cC_{\sigma,T}^{-\beta} \zZ$ is the completion of $L_{\sigma,T}^\infty \zZ$ under the above norm.

\begin{prop} \label{prop:Young_multilinear}
Let $\mathsf{T}_m : \xX_1 \times \xX_2 \times \dots \times \xX_{m-1} \times \yY \rightarrow \zZ$ be a multi-linear form. Let $\nu_1, \dots, \nu_{m-1}, \nu \in (0,1)$ such that
\begin{equation*}
    \min \{\nu_1, \dots, \nu_{m-1}\} > \nu\;.
\end{equation*}
Also let $\sigma \in [0,2)$. Then we have the bound
\begin{equation} \label{eq:multilinear_pointwise_bound}
    \begin{split}
    &\phantom{111}\big\| \mathsf{T}_m \big( f_1, \dots, f_{m-1}, g \big) \big\|_{\cC_{\sigma,T}^{-\nu} \zZ}\\
    &\lesssim \|f_1\|_{\cC_{\sigma,T}^{\nu_1} \xX_1} \bigg( \prod_{j=2}^{m-1} \|f_j\|_{\cC_T^{\nu_j} \xX_j} \bigg) \|g\|_{\cC_T^{-\nu} \yY} \cdot \|\mathsf{T}_m\|_{\xX_1 \times \cdots \times \xX_{m-1} \times \yY \rightarrow \zZ}\;.
    \end{split}
\end{equation}
\end{prop}
\begin{proof}
It suffices to show that the integral
\begin{equation} \label{eq:Young_multilinear}
    \int_s^t \mathsf{T}_m \big( f_1, \dots, f_{m-1}, g \big)(r) {\rm d}r := \int_s^t \mathsf{T}_m \big( f_1(r), \dots, f_{m-1}(r), {\rm d}G(r) \big)
\end{equation}
is well defined as a Young integral, and satisfies desired bounds. The only difference with the classical Young integral is that there is a weight $\sigma>0$ in the input $f_1$, which gives an integrable singularity near time $0$. For sake of completeness, we sketch an argument here. 

Fix $[a,b] \subset (0,T]$, let $\{\pP_n\}_n$ be the dyadic partition of $[a,b]$ with intervals of equal length $\delta_n = \frac{b-a}{2^n}$. For every $n$, let
\begin{equation*}
    \I^n_{[a,b]} := \sum_{j} \mathsf{T}_{m} \Big( f_1(\tau_j), \dots, f_{m-1}(\tau_j), \; G(\tau_{j+1}) - G(\tau_j) \Big)\;,
\end{equation*}
where the sum is taken over all points in the partition $\pP_n$. Then
\begin{equation*}
    \begin{split}
    \I^{n+1}_{[a,b]} - &\I^{n}_{[a,b]} = \sum_j \Big( \mathsf{T}_m \big( f_1 (\tau_j + \delta_{n+1})\;, \dots, f_{m-1}(\tau_j + \delta_{n+1}), G(\tau_{j+1}) - G(\eta_j) \big)\\
    &- \mathsf{T}_m  \big( f_1 (\tau_j)\;, \dots, f_{m-1}(\tau_j), G(\tau_{j+1}) - G(\eta_j) \big) \Big)\;,
    \end{split}
\end{equation*}
where the sum is taken over all $j$ such that $\tau_j \in \pP_n$, and $\eta_j := \frac{1}{2}(\tau_j + \tau_{j+1})$. By assumption, we have
\begin{equation*}
    \| \I^{n+1}_{[a,b]} - \I^{n}_{[a,b]} \|_{\zZ} \lesssim \sum_j \tau_{j}^{-\frac{\sigma}{2}} (\tau_{j+1} - \tau_j)^{1+\nu_{\min} - \nu} \lesssim 2^{-n(\nu_{\min} - \nu)} \int_a^b r^{-\frac{\sigma}{2}} {\rm d}r\;,
\end{equation*}
where $\nu_{\min} = \min\{\nu_1, \dots, \nu_{m-1}\}$, and the proportionality constant is proportional to the right hand side of \eqref{eq:multilinear_pointwise_bound}. 

Hence, the Young integral in \eqref{eq:Young_multilinear} is well defined as the $n \rightarrow +\infty$ of $\I_{[s,t]}^n$, and satisfies the bound
\begin{equation*}
    \begin{split}
    \Big\| \int_s^t \mathsf{T}_m &\big( f_1(r), \dots, f_{m-1}(r), {\rm d}G(r) \big) \Big\|_{\zZ} \lesssim \big\| \mathsf{T}_m \big( f_1(s), \dots, f_{m-1}(s), G(t) - G(s) \big) \big\|_\zZ\\
    &+ \Big( \int_s^t r^{-\frac{\sigma}{2}} {\rm d}r \Big) \cdot \|f_1\|_{\cC_{\sigma,T}^{\nu_1} \xX_1} \bigg( \prod_{j=2}^{m-1} \|f_j\|_{\cC_T^{\nu_j} \xX_j} \bigg) \|g\|_{\cC_T^{-\nu} \yY} \cdot \|\mathsf{T}_m\|_{\xX_1 \times \cdots \times \xX_{m-1} \times \yY \rightarrow \zZ}\\
    &\lesssim s^{-\frac{\sigma}{2}} (t-s)^{1-\nu} \|f_1\|_{\cC_{\sigma,T}^{\nu_1} \xX_1} \bigg( \prod_{j=2}^{m-1} \|f_j\|_{\cC_T^{\nu_j} \xX_j} \bigg) \|g\|_{\cC_T^{-\nu} \yY} \cdot \|\mathsf{T}_m\|_{\xX_1 \times \cdots \times \xX_{m-1} \times \yY \rightarrow \zZ}\;.
    \end{split}
\end{equation*}
This completes the proof. 
\end{proof}

\subsection{Spacetime integration with heat kernel}

This section aims to deal with the spacetime convolution of heat kernel with a spacetime distribution $h$, where $h$ is one of the objects in \eqref{eq:kpz_rhs_list} (except the quadratic one) from the KPZ equation. All of those objects belong to $\cC_{\sigma_1,T}^{-\kappa} \cC_x^{\gamma;\eps}$ for some $\sigma_1 \in (0,2)$, $\gamma \in (-1,1)$ and $\kappa > 0$ arbitrarily small. Hence, we keep the choice of $\kappa$ sufficiently small to avoid further assumptions on various exponents. 

\begin{prop} \label{prop:Young_heat_convolution}
Let $\alpha, \beta \in [0,2)$ and $\gamma \in \R$ with $\alpha + \beta - \gamma \in (0,2)$. Let $\sigma_1 \in (0,2)$. Then for all sufficiently small $\kappa$ and $h \in \cC_{\sigma_1,T}^{-\kappa} \cC_{x}^{\gamma;\eps}$, the integral
\begin{equation*}
    (\iI_\eps h)(t) := \int_0^t e^{(t-r) \lL_\eps} {\rm d} H(r)
\end{equation*}
is well defined as a weighted Young integral in $\cC_x^{\beta;\eps}$, where $H$ is the time integral of $h$. Furthermore, there exists $\theta>0$ such that
\begin{equation*}
    \|\iI_\eps h\|_{\cC_{\sigma_2,T}^{\alpha/2} \cC_x^{\beta;\eps}} \lesssim T^\theta \|h\|_{\cC_{\sigma_1,T}^{-\kappa} \cC_x^{\gamma;\eps}}\;,
\end{equation*}
where $\sigma_2 = (\sigma_1 + \alpha + \beta - \gamma - 2 + 7 \kappa) \vee 0$. 
\end{prop}
\begin{proof}
The argument is essentially the same as that in Proposition~\ref{prop:Young_multilinear} except that in the current case, there is also a weight in time. So we only give a sketch argument. 

We first give bounds on $\|\cdot\|_{\beta+\alpha;\eps}$-norm for fixed time $t$, and then extend it to $\cC_{\sigma_2,T}^{\alpha/2}\cC_x^{\beta;\eps}$-norm. For any $[a,b] \subset (0,t) \subset [0,T]$, we have
\begin{equation*}
    \I_{[a,b]}^{n+1} - \I_{[a,b]}^{n} = - \sum_{j} \big( e^{(\eta_j-\tau_j) \lL_\eps} - \id \big) e^{(t-\eta_j) \lL_\eps} \big( H(\tau_{j+1}) - H(\eta_j) \big)\;,
\end{equation*}
where the sum is taken over all $j$ such that $\tau_j \in \pP_n$, the dyadic partition of $[a,b]$, and $\eta_j = \frac{1}{2} (\tau_j + \tau_{j+1})$. By standard heat kernel estimates, we have
\begin{equation*}
    \begin{split}
    \| \I_{[a,b]}^{n+1} - \I_{[a,b]}^n \|_{\alpha+\beta; \eps} &\lesssim \sum_{j} (\eta_j - \tau_j)^{2\kappa} (t - \eta_j)^{-\frac{\alpha+\beta-\gamma}{2}-2\kappa} \|H(\tau_{j+1}) - H(\eta_j)\|_{\gamma;\eps}\\
    &\lesssim 2^{-\kappa n} \bigg( \int_a^b r^{-\frac{\sigma_1}{2}} (t-r)^{-\frac{\alpha+\beta-\gamma}{2} - 2\kappa} {\rm d}r \bigg) \cdot \|h\|_{\cC_{\sigma_1,T}^{-\kappa} \cC_x^{\gamma;\eps}}\\
    &\lesssim 2^{-\kappa n} a^{-\frac{\sigma_1}{2}} (t-a)^{-\frac{\alpha+\beta-\gamma}{2}-2\kappa} \, (b-a)^{1-\kappa} \cdot \|h\|_{\cC_{\sigma_1,T}^{-\kappa} \cC_x^{\gamma;\eps}}\;.
    \end{split}
\end{equation*}
This implies that the Young integral is well defined on $[a,b]$, and satisfies the bound
\begin{equation*}
    \begin{split}
    \left\| \int_a^b e^{(t-r) \lL_\eps} {\rm d} H(r) \right\|_{\alpha+\beta;\eps} &\leq \sum_{n \geq 0} \| \I_{[a,b]}^{n+1} - \I_{[a,b]}^n \|_{\alpha+\beta; \eps} + \big\| e^{(t-a) \lL_\eps} \big( H(b) - H(a) \big) \big\|_{\alpha+\beta; \eps}\\
    &\lesssim a^{-\frac{\sigma_1}{2}} (t-a)^{-\frac{\alpha+\beta-\gamma}{2}-2\kappa} \, (b-a)^{1-\kappa} \cdot \|h\|_{\cC_{\sigma_1,T}^{-\kappa} \cC_x^{\gamma;\eps}}\;.
    \end{split}
\end{equation*}
To give well-posedness and desired bound for $(\iI_\eps h)(t)$, we bootstrap the above bound. Fix $ \delta > 0 $ small, and let $\tau_0 = \delta$. For $1\leq j \leq m-1$, define $\tau_{j} := 2 \tau_{j-1}$, where $m$ is the unique integer such that $\tau_{m-1} \in [\frac{t}{4}, \frac{t}{2})$. Then define $\tau_m := t-\delta$. Since, 
\begin{equation*}
    \int_\delta^{t-\delta} e^{(t-r) \lL_\eps} {\rm d}H(r) = \sum_j \int_{\tau_j}^{\tau_{j+1}} e^{(t-r) \lL_\eps} {\rm d}H(r)\;,
\end{equation*}
we have
\begin{equation*}
    \begin{split}
    \left\| \int_\delta^{t-\delta} e^{(t-r) \lL_\eps} {\rm d}H(r) \right\|_{\alpha+\beta; \eps} &\lesssim  \bigg( \sum_j \tau_j^{-\frac{\sigma_1}{2}} (t-\tau_{j})^{-\frac{\alpha+\beta-\gamma}{2}-2\kappa} (\tau_{j+1} - \tau_j)^{1-\kappa} \bigg) \, \|h\|_{\cC_{\sigma_1,T}^{-\kappa} \cC_x^{\gamma;\eps}}\\
    &\lesssim t^{1-\frac{\sigma_1 + \alpha + \beta - \gamma}{2}-3\kappa} \|h\|_{\cC_{\sigma_1,T}^{-\kappa} \cC_x^{\gamma;\eps}} \lesssim t^{\theta - \frac{\sigma_2}{2}} \|h\|_{\cC_{\sigma_1,T}^{-\kappa} \cC_x^{\gamma;\eps}}\;.
    \end{split}
\end{equation*}
where the second inequality follows from our choice of $\tau_j$ so that $\tau_{j+1} - \tau_j \sim \tau_j$. Since the bound does not depend on $\delta$, sending $\delta \rightarrow 0$ yields
\begin{equation*}
    \|(\iI_\eps h)(t)\|_{\alpha+\beta;\eps} \lesssim t^{\theta - \frac{\sigma_2}{2}} \|h\|_{\cC_{\sigma_1,T}^{-\kappa} \cC_x^{\gamma;\eps}}\;.
\end{equation*}
We now to turn to the estimate for time-continuity. If $s \leq \frac{t}{2}$ (so that $s \lesssim t-s$), we can control the difference directly by
\begin{equation*}
    \|(\iI_\eps h)(t) - (\iI_\eps h)(s)\|_{\beta;\eps} \leq \|(\iI_\eps h)(t)\|_{\beta;\eps} + \|(\iI_\eps h)(s)\|_{\beta;\eps} \lesssim t^{\theta} \, s^{-\frac{\sigma_2}{2}} \|h\|_{\cC_{\sigma_1,T}^{-\kappa} \cC_x^{\gamma;\eps}}\;.
\end{equation*}
For $s > \frac{t}{2}$ (so that $t-s \lesssim s$), we write the difference as
\begin{equation*}
    (\iI_\eps h)(t) - (\iI_\eps h)(s) = \big( e^{(t-s) \lL_\eps} - \id \big) (\iI_\eps h)(s) + \int_s^t e^{(t-r) \lL_\eps} {\rm d}H(r)\;.
\end{equation*}
For the first one, it follows from the bound for fixed time that
\begin{equation*}
    \begin{split}
    \Big\| \big( e^{(t-s) \lL_\eps} - \id \big) (\iI_\eps h)(s)  \Big\|_{\beta;\eps} &\lesssim (t-s)^{\alpha/2} \|(\iI_\eps h)(s)\|_{\alpha+\beta;\eps}\\
    &\lesssim (t-s)^{\alpha/2} t^{\theta - \frac{\sigma_2}{2}} \|h\|_{\cC_{\sigma_1,T}^{-\kappa} \cC_x^{\gamma;\eps}}\;.
    \end{split}
\end{equation*}
For the second one, we control it by
\begin{equation*}
    \begin{split}
    \left\| \int_s^t e^{(t-r) \lL_\eps} {\rm d}H(r) \right\|_{\beta;\eps} &\lesssim s^{-\frac{\sigma_1}{2}} (t-s)^{1-\frac{\beta-\gamma}{2}-2\kappa} \, \|h\|_{\cC_{\sigma_1,T}^{-\kappa} \cC_x^{\gamma;\eps}} + \left\| e^{(t-s) \lL_\eps} \big( H(t) - H(s) \big) \right\|_{\beta;\eps}\\
    &\lesssim  t^\theta (t-s)^{\alpha/2} s^{-\frac{\sigma_2}{2}} \|h\|_{\cC_{\sigma_1,T}^{-\kappa} \cC_x^{\gamma;\eps}}\;.
    \end{split}
\end{equation*}
This completes the proof. 
\end{proof}

\subsection{Commutator with the heat operator}

The aim of this subsection is to establish the bounds for $[A_\eps \iI_\eps ', f \prec_\eps] g$ with $f = \Lambda_\eps^{(\delta)} (u_\eps^{(\delta)})$ and $g = \<1bk>_\eps^{(\delta)}$ or $g = \<1bk>_\eps^{(\delta)} - \frac{\lambda}{\bar{A}} \<1bks>_0$ in the KPZ equation in mind. Note that
\begin{equation*}
    [A_\eps \iI_\eps', f \prec_\eps] g = A_\eps \d_x \, \big( [\iI_\eps, f \prec] g \big) + [A_\eps \d_x, f \prec_\eps] \, \iI_\eps(g)\;.
\end{equation*}
We focus on the estimate for $[\iI_\eps, f \prec_\eps] g$ since this is where the negative-in-time regularity of $g$ brings technical differences. For this term, we have
\begin{equation*}
    \big( [\iI_\eps, f \prec_\eps] g \big)(t) = \int_0^t \hH_\eps(t;r) {\rm d} G(r)\;,
\end{equation*}
where
\begin{equation*}
    \hH_\eps (t;r) g = e^{(t-r) \lL_\eps} \big( f(r) \prec_\eps g \big) - f(t) \prec_\eps \big( e^{(t-r) \lL_\eps} g \big)\;.
\end{equation*}
Since we have this specific object in mind, to avoid complexity, we use the same exponents as those in the KPZ equation and do not try to generalise. Hence, in what follows, we always assume
\begin{equation*}
    \eta \in \left(0, \frac{1}{1000}\right)\;, \quad \kappa \in \left(0, \frac{\eta}{1000}\right)\;, \quad \sigma := \frac{1}{2} + 7\eta\;.
\end{equation*}
Also, $\beta, \gamma \in \R$ such that $\beta - \gamma \in (1, 2 - 10 \kappa)$. We have the following estimate on time difference of $\hH_\eps$. 

\begin{lem} \label{lem:heat_comm_negative_fixed_t}
We have the pointwise bound
\begin{equation*}
    \begin{split}
    &\phantom{111}\big\| \big( \hH_\eps(t;r) - \hH_\eps(t;r') \big) \, g  \big\|_{\beta + 3\eta - 6\kappa;\eps}\\
    &\lesssim (r-r')^{3\kappa} (t-r)^{-\frac{\beta-\gamma}{2}} |r'|^{-\frac{\sigma}{2}} \, \|f\|_{\cC_{\sigma,T}^{2\eta} L_x^\infty \cap L_{\sigma,T}^\infty \cC_x^{4\eta}} \, \|g\|_{\gamma;\eps}\;.
    \end{split}
\end{equation*}
\end{lem}
\begin{proof}
Write $\hH_\eps = \hH_\eps^{(1)} + \hH_\eps^{(2)}$, where
\begin{equation*}
    \hH_\eps^{(1)}(t;r) g = e^{(t-r) \lL_\eps} \big( \big( f(r) - f(t) \big) \prec_\eps g \big)\;, \quad \hH_\eps^{(2)}(t;r) g = \big[ e^{(t-r) \lL_\eps}, f(t) \prec_\eps \big] \, g\;.
\end{equation*}
For $\hH_\eps^{(1)}$, we have
\begin{equation*}
    \begin{split}
    \|\hH_\eps^{(1)}(t;r) g\|_{\beta + 4\eta;\eps} &\lesssim (t-r)^{-\frac{\beta + 4\eta -\gamma}{2}} \big\| \big( f(r) - f(t) \big) \prec_\eps g \big\|_{\gamma;\eps}\\
    &\lesssim (t-r)^{-\frac{\beta-\gamma}{2}} r^{-\frac{\sigma}{2}} \|f\|_{\cC_{\sigma,t}^{2\eta} L_x^\infty} \|g\|_{\gamma;\eps}\;,
    \end{split}
\end{equation*}
and hence
\begin{equation*}
    \big\| \Delta_{k,\eps} \big( \hH_\eps^{(1)}(t;r) g \big) \big\|_{L^\infty} \lesssim 2^{-(\beta + 4\eta) k} \, (t-r)^{-\frac{\beta-\gamma}{2}} \, r^{-\frac{\sigma}{2}} \, \|f\|_{\cC_{\sigma,t}^{2\eta} L_x^\infty} \, \|g\|_{\gamma;\eps}\;.
\end{equation*}
Also, by Lemma~\ref{le:heat_comm_space}, we have
\begin{equation*}
    \big\| \Delta_{k,\eps} \big( \hH_\eps^{(2)}(t;r) g \big) \big\|_{L^\infty} \lesssim 2^{- (\beta + 4\eta) k} \, (t-r)^{-\frac{\beta-\gamma}{2}} \, t^{-\frac{\sigma}{2}} \, \|f\|_{L_{\sigma,t}^{\infty} \cC_x^{4\eta;\eps}} \, \|g\|_{\gamma;\eps}\;.
\end{equation*}
Combining the above, we get
\begin{equation} \label{eq:comm_heat_block_fixed_time}
    \|\Delta_{k,\eps} \big( \hH_\eps(t;r) g \big)\|_{L^\infty} \lesssim 2^{-(\beta+ 4\eta) k} (t-r)^{-\frac{\beta-\gamma}{2}} \, r^{-\frac{\sigma}{2}} \, \|f\|_{\cC_{\sigma,T}^{2\eta} L_x^\infty \cap L_{\sigma,t}^{\infty} \cC_x^{4\eta;\eps}} \, \|g\|_{\gamma;\eps}\;.
\end{equation}
This directly implies the pointwise bound
\begin{equation} \label{eq:comm_heat_block_opt_reg}
    \begin{split}
    &\phantom{111}\big\| \Delta_{k,\eps} \big( \big( \hH_\eps(t;r) - \hH_\eps(t;r') \big) g \big) \big\|_{L^\infty}\\
    &\lesssim 2^{- (\beta + 4\eta) k} (t-r)^{-\frac{\beta-\gamma}{2}} \, |r'|^{-\frac{\sigma}{2}} \, \|f\|_{\cC_{\sigma,T}^{2\eta} L_x^\infty \cap L_{\sigma,T}^\infty \cC_x^{4\eta;\eps}} \, \|g\|_{\gamma;\eps}\;.
    \end{split}
\end{equation}
To get a factor of a positive power of $(r-r')$, we decompose
\begin{equation*}
    \begin{split}
    \big( \hH_\eps &(t;r) - \hH_\eps(t;r') \big) g = \big( \id - e^{(r-r') \lL_\eps} \big) \, e^{(t-r) \lL_\eps} \, \big( f(r) \prec_\eps g \big)\\
    &+ e^{(t-r') \lL_\eps} \Big( \big( f(r) - f(r') \big) \prec_\eps g \Big) + f(r) \prec_\eps \Big( \big( \id - e^{(r-r')\lL_\eps} \big) e^{(t-r) \lL_\eps} g \Big)\;.
    \end{split}
\end{equation*}
It then follows from the above decomposition that
\begin{equation*}
    \big\| \big( \hH_\eps(t;r) - \hH_\eps (t;r') \big) g \big\|_{\beta-24\kappa;\eps} \lesssim (r-r')^{12 \kappa} (t-r)^{-\frac{\beta-\gamma}{2}} |r'|^{-\frac{\sigma}{2}} \, \|f\|_{\cC_{\sigma,t}^{2\eta} L_x^\infty} \, \|g\|_{\gamma;\eps}\;,
\end{equation*}
and hence
\begin{equation} \label{eq:comm_heat_block_time_diff}
    \begin{split}
    &\phantom{111}\big\| \Delta_{k,\eps} \big( \big( \hH_\eps (t;r) - \hH_\eps (t;r') \big) g \big) \big\|_{L^\infty}\\
    &\lesssim 2^{-(\beta-24\kappa)k} \, (r-r')^{12 \kappa} (t-r)^{-\frac{\beta-\gamma}{2}} |r'|^{-\frac{\sigma}{2}} \, \|f\|_{\cC_{\sigma,t}^{2\eta} L_x^\infty} \, \|g\|_{\gamma;\eps}\;.
    \end{split}
\end{equation}
The claim of the lemma then follows by interpolating \eqref{eq:comm_heat_block_opt_reg} and \eqref{eq:comm_heat_block_time_diff}. 
\end{proof}

\begin{lem} \label{lem:heat_comm_negative_time_diff}
We have
\begin{equation*}
    \begin{split}
    &\phantom{111}\big\| \big( \hH_\eps(t;r) - \hH_\eps(s;r) \big) g - \big( \hH_\eps(t;r') - \hH_\eps(s;r') \big) g \big\|_{\beta + 2\eta - 8\kappa; \eps}\\
    &\lesssim (t-s)^{2\kappa} (r-r')^{2\kappa} (s-r)^{-\frac{\beta-\gamma}{2}} |r'|^{-\frac{\sigma}{2}} \, \|f\|_{\cC_{\sigma,T}^{2\eta} L_x^\infty \cap L_{\sigma,T}^\infty \cC_x^{4\eta}} \, \|g\|_{\gamma;\eps}
    \end{split}
\end{equation*}
uniformly over all $0<s<t<T$ and all $\eps \in [0,1]$. 
\end{lem}
\begin{proof}
First, by Lemma~\ref{lem:heat_comm_negative_fixed_t}, we have
\begin{equation} \label{eq:heat_comm_time_diff_opt_reg}
    \begin{split}
    &\phantom{111}\big\| \Delta_{k,\eps} \big( \big( \hH_\eps(t;r) - \hH_\eps(s;r) \big) g - \big( \hH_\eps(t;r') - \hH_\eps(s;r') \big) g \big) \big\|_{L^\infty}\\
    &\lesssim 2^{-(\beta + 3\eta - 6 \kappa) k} \, (r-r')^{3\kappa} (s-r)^{-\frac{\beta-\gamma}{2}} |r'|^{-\frac{\sigma}{2}} \, \|f\|_{\cC_{\sigma,T}^{2\eta} L_x^\infty \cap L_{\sigma,T}^\infty \cC_x^{4\eta}} \, \|g\|_{\gamma;\eps}\;.
    \end{split}
\end{equation}
To get a factor with positive power of $t-s$, we decompose
\begin{equation*}
    \begin{split}
    \big( \hH_\eps(t;r) - &\hH_\eps(s;r) \big) g = \big( e^{(t-s) \lL_\eps} - \id \big) e^{(s-r) \lL_\eps} \big(f(r) \prec_\eps g \big)\\
    &- \big( f(t) - f(s) \big) \prec_\eps e^{(t-r) \lL_\eps} g - f(t) \prec_\eps \big( e^{(t-s) \lL_\eps} - \id \big) e^{(s-r) \lL_\eps} g\;.
    \end{split}
\end{equation*}
It then follows from this decomposition that
\begin{equation*}
    \big\| \big( \big( \hH_\eps (t;r) - \hH_\eps(s;r) \big) g \big) \big)\|_{\beta - 12 \kappa; \eps} \lesssim (t-s)^{6\kappa} (s-r)^{-\frac{\beta-\gamma}{2}} r^{-\frac{\sigma}{2}} \|f\|_{\cC_{\sigma,T}^{2\eta} L_x^\infty} \, \|g\|_{\gamma;\eps}\;,
\end{equation*}
and hence
\begin{equation} \label{eq:heat_comm_time_diff}
    \begin{split}
    &\phantom{111}\big\| \Delta_{k,\eps} \big( \big( \hH_\eps(t;r) - \hH_\eps(s;r) \big) g - \big( \hH_\eps(t;r') - \hH_\eps(s;r') \big) g \big) \big\|_{L^\infty}\\
    &\lesssim 2^{-(\beta- 12\kappa) k} \, (t-s)^{6\kappa} (s-r)^{-\frac{\beta-\gamma}{2}} |r'|^{-\frac{\sigma}{2}} \, \|f\|_{\cC_{\sigma,T}^{2\eta} L_x^\infty} \, \|g\|_{\gamma;\eps}\;.
    \end{split}
\end{equation}
The claim of the lemma then follows by interpolating \eqref{eq:heat_comm_time_diff_opt_reg} and \eqref{eq:heat_comm_time_diff}. 
\end{proof}

\begin{prop} \label{prop:Young_heat_comm}
The quantity
\begin{equation*}
    [\iI_\eps, f \prec_\eps] g := \int_{0}^{t} \hH_\eps(t;r) {\rm d} G(r)
\end{equation*}
is well defined as a Young integral, and satisfies the bound
\begin{equation*}
    \big\| \, [\iI_\eps, f \prec_\eps] g \, \big\|_{\cC_{\sigma', T}^{2\kappa} \cC_x^{\beta + 2\eta - 8\kappa; \eps}} \lesssim \|f\|_{\cC_{\sigma,T}^{2\eta} L_x^\infty \cap L_{\sigma,T}^\infty \cC_x^{4\eta}} \, \|g\|_{\cC_T^{-\kappa} \cC_x^{\gamma;\eps}}\;,
\end{equation*}
where $\sigma' = \beta-\gamma+\sigma-2$. 
\end{prop}
\begin{proof}
Given the estimate of the integrand in Lemma~\ref{lem:heat_comm_negative_time_diff}, the proof follows in essentially the same way as that of Proposition~\ref{prop:Young_heat_convolution}, but only simpler since $g$ does not have a weight at $0$. We omit the detailed argument. 
\end{proof}

The above proposition holds for general spatial dimension $d$. On the other hand, as long as $A_\eps \d_x$ is involved, we restrict ourselves to $d=1$. We have the following proposition. 

\begin{prop} \label{prop:commutator_derivative_heat_negative_time}
Let $d=1$, $\gamma \in (-1,0)$, and $\sigma \in (0,2)$. Then we have
\begin{equation*}
    \big\| [A_\eps \iI_\eps', \, f \prec_\eps] \, g \big\|_{\cC_{\sigma,T}^{2\kappa} \cC_x^{\gamma+1+2\eta - 10\kappa}} \lesssim \|f\|_{\cC_{\sigma,T}^{2\eta} L_x^\infty \cap L_{\sigma,T}^\infty \cC_x^{4\eta}} \, \|g\|_{\cC_T^{-\kappa} \cC_x^{\gamma;\eps}}\;.
\end{equation*}
\end{prop}
\begin{proof}
By Lemma~\ref{le:derivative_a_1d} and Proposition~\ref{prop:Young_heat_comm}, we have
\begin{equation*}
    \big\| A_\eps \d_x \big( [\iI_\eps, f \prec_\eps] \, g \big) \, \big\|_{\cC_{\sigma',T}^{2\kappa} \cC_x^{\gamma+1+2\eta - 10\kappa}} \lesssim \|f\|_{\cC_{\sigma,T}^{2\eta} L_x^\infty \cap L_{\sigma,T}^\infty \cC_x^{4\eta}} \, \|g\|_{\cC_T^{-\kappa} \cC_x^{\gamma;\eps}}\;.
\end{equation*}
Furthermore, by Proposition~\ref{pr:a_derivative_prec_comm}, we have
\begin{equation*}
    \big\| [A_\eps \d_x, f \prec_\eps] \, \iI_\eps(g) \big\|_{\cC_{\sigma,T}^{2\kappa} \cC_x^{2\eta+\gamma+1-8\kappa}} \lesssim \|f\|_{\cC_{\sigma,T}^{2\eta} L_x^\infty \cap L_{\sigma,T}^\infty \cC_x^{4\eta}} \, \|g\|_{\cC_T^{-\kappa} \cC_x^{\gamma;\eps}}\;.
\end{equation*}
The conclusion then follows by combining the above two bounds. 
\end{proof}

\section{Comparison with the previous version}
\label{sec:comparison}

The previous version of this manuscript was submitted to and appeared on arXiv in August 2025. In that version, the statement for the KPZ equation of oscillatory coefficient was restricted to the case $g_\eps = A_\eps$. To cover the more general situation (with general $g$), one needs to deal with spacetime distributions instead of continuous-in-time evolutions in the space of distributions. Hence, we did not submit that preprint to journal while trying to find a clean way to cover the general case. Meanwhile, we have also found improvements / simplifications for other parts of the article. These all together constitute the four main changes of this version. We list them below. 

\begin{enumerate}

\item \textbf{Adding the $2$D gPAM model} (including a quadratic nonlinearity of the derivation of the solution). This is the content of Theorems~\ref{thm:gpam_informal} (informal) and~\ref{thm:PAM_overall} (complete statement). 

$2$D gPAM was not treated in our 2025 arXiv version. Adding it in this version was inspired by \cite{homo_spde_rs}, but the treatment is quite different. Without the quadratic nonlinearity of the derivative (that is, $F_\eps \equiv 0$ in \eqref{e:gPAM_formal_intro}), the ansatz (together with uniform bounds and convergence) fits neatly with the generalised para-controlled distribution framework here. To treat the additional product $\nabla^\sT u_\eps^{\#} \, F_\eps \, \nabla u_\eps^{\#}$, we employed the integration by parts technique as well as the ``completion of product" trick from \cite{homo_gPAM}. 

\item \textbf{Extension to general oscillatory coupling coefficient $g_\eps = g(\cdot/\eps)$ in KPZ} (in the previous version, we only treated the special case $g_\eps = A_\eps$). 

The extension requires treating spacetime distributions (not only distribution in space for fixed times), and hence some serious modifications of the previous content. After comparing several different possible treatments, we ended up formulating this part as Young integration in time. This is the content of Appendix~\ref{sec:Young_integral}. The main content of the article is not affected. 

\item \textbf{Simplification of the proof of Proposition~\ref{pr:decoupling_spectrum}} -- uniform decay estimate at high spectrum. 

In the previous version, the proof of Proposition~\ref{pr:decoupling_spectrum} goes via a spectral expansion of $\Delta_{i,a} f \cdot \Delta_{j,a} g$. In the current version, we utilise the pointwise bound on the kernel from Proposition~\ref{e:kernel_pointwise_bound} as well as duality to give a direct and shorter proof. 

\item \textbf{Significant simplification of the proof of Theorem~\ref{th:convergence_single_block}} -- convergence of the block operators. 

The proof of Theorem~\ref{th:convergence_single_block} in the previous version involves a combination of the following ingredients: quantitative convergence of individual eigenvalues from \cite{KLS2013}, Weyl's asymptotic formula together with a clustering argument. The proof of that single theorem took almost ten pages. 

In March--April 2026, Siqi Fu gave a series of lectures on the spectral theory for complex Laplacian at Westlake University. We learned the Helffer-
Sj\"ostrand’s representation for self-adjoint operators. This representation enables us to give a much simpler proof of the convergence of block operations (Theorem~\ref{th:convergence_single_block}). 

\end{enumerate}

\bigskip
\bigskip
\bigskip
\bigskip
\bigskip

\newpage

\section{Index of symbols}

In this appendix, we give a list of frequently used symbols in this article as well as their meanings. 

\begin{longtable}{cl}
    \hline
        Symbol & Meaning\\
        \hline
         $\Pi_0$  & Average of functions on $\T^d$\\
         $\Pi_0^{\perp}$  & $\id-\Pi_0$\\
         $\sS_d(\Lambda)$  & Set of symmetric matrix coefficients on $\T^d$ with ellipticity between $\frac{1}{\Lambda}$ and $\Lambda$\\
         $\sS_d(\Lambda,M)$  & Subset of $\sS_d(\Lambda)$ specified in Definition~\ref{def:operator_class}\\
         % $\Z$& Set of integers \\
         % $\N$& Set of positive integers\\
         $\N^{-1}$  & Set of the inverse of positive integers\\
         $\cC^{\alpha}_{\fs}$& Space-time parabolic H\"older spaces\\
         $[\cdot,\cdot]$  & Commutator $[A,B]f = AB(f)-BA(f)$\\
         $a$ & Generic element in $\sS_d(\Lambda, M)$\\
         $Q_a$ & Green function of $\d_t-\lL_{a}$\\
        $\lL_{a}$ & Elliptic operator $\div(a\nabla)$ on $\T^d$\\
        $\Delta_{j,a}$ & Littlewood-Paley blocks in \eqref{eq: paraproduct_a} based on $a\in \sS_d(\Lambda,M)$\\
        $\bB_{p,q}^{\alpha;a}$& Besov spaces based on $\lL_{a}$ in \eqref{eq: gBesov_defn}\\
        $\cC^{\alpha;a}$& H\"older-Besov spaces based on $\lL_a$ in \eqref{eq: Holder_Besov_defn}\\
        $\prec_a$, $\circ_a$, $\succ_a$ & Para-products based on $a\in \sS_d(\Lambda,M)$, defined in \eqref{eq: paraproduct_a} \\
        $\iI_{a}$& Integration operator $(\d_t-\lL_{a})^{-1}$\\
        % $\mathbf{\Pi}_a$& Sums of Littlewood Paley blocks based on $a$, defined in \eqref{eq: mathbf Pi}\\
        $\P_a$ & Sums of gradients of Littlewood Paley blocks based on $a$, defined in \eqref{e:operator_Pa}\\
        $\Com_a$ & Commutator operation based on $a$, defined in \eqref{e:commutator_defn}\\
        $\|\cdot\|_{\mathfrak{L}^{\alpha;a}_{\sigma,T}}$& Weighted-in-time space-time Besov spaces based on $a$, defined in \eqref{eq: weighted_mathfrakL}\\
        $A$   & Coefficient matrix satisfying Assumption~\ref{as:A}\\
        $A_\eps$ & Rescaled coefficient $A(\cdot / \eps)$\\
        $\bar{A}$ & Homogenised matrix for $A$, given in \eqref{def:bar_a_homogenised_matrix}\\
        $\chi$ & Homogenisation corrector associated with $A$, defined in \eqref{def:chi_homo_corrector}\\
        $\cC^{\alpha;\eps}$& H\"older-Besov spaces based on $\lL_\eps$\\
        $|\!|\!|\cdot;\cdot|\!|\!|_{\alpha;\eps}$ & Modified norms for comparison of objects in $\cC^{\alpha;\eps}$ and $\cC^{\alpha;0}$\\
        \phantom{1} & with $|\alpha| \in (1,2)$, defined in \eqref{e:norm_comparison_high_reg} and \eqref{e:norm_comparison_low_reg}\\
        $\lL_{\eps}$ & Elliptic operator $\div(A_{\eps}\nabla )$ on $\T^d$\\
        $\lL_{0}$ & Elliptic operator $\div(\bar{A}\nabla )$ on $\T^d$\\
        $\Delta_{j,\eps}$ & Littlewood-Paley blocks in \eqref{eq: paraproduct_a} based on $A_{\eps}$ where $A$ satisfies Assumption \ref{as:A}\\
        $\prec_\eps$, $\circ_\eps$, $\succ_\eps$ & Paraproducts based on $A_{\eps}$ for the matrix $A$ satisfying Assumption \ref{as:A}\\
        $\iI_{\eps}$& Integration operator $(\d_t-\lL_{\eps})^{-1}$\\
        % $\mathbf{\Pi}_\eps$& Sums of Littlewood Paley blocks based on $A_{\eps}$, defined in \eqref{eq: mathbf Pi}\\
        $\P_\eps$ & Sums of gradients of Littlewood-Paley blocks based on $A_{\eps}$, defined in \eqref{e:operator_Pa}\\
        $\Com_{\eps}$ & Commutator based on $A_{\eps}$ in \eqref{e:commutator_defn}\\
        $Q_\eps$ & Green function of $\d_t-\lL_{\eps}$\\
        $G_\eps$ & Green function of $-\lL_{\eps}$\\
        % & \\
        %  & \\
        %  & \\
        %  & \\
        %  & \\
        %  & \\
        %  & \\
        %  & \\
        %  & \\
        %  & \\
        %  & \\
        %  & \\
        %  & \\
        %  & \\
        %  & \\
        %  & \\
        %  & \\
        %  & \\
        %  & \\
        $(\lambda_{n,\eps})_{n \geq 0}$& Eigenvalues of $\sqrt{-\lL_{\eps}}$ in increasing order, with multiplicities repeated\\
        $\psi_{n,\eps}$& Eigenfunction of $\sqrt{-\lL_{\eps}}$ with eigenvalue $\lambda_{n,\eps}$\\
        % $P_{n,\eps}$& Projection onto the one dimensional space spanned by $\psi_{n,\eps}$ \\
        $(\tilde{\lambda}_{n,\eps})_{n \geq 0}$ & Eigenvalues of $\sqrt{-\lL_{\eps}}$ in strictly increasing order\\
        % $\tilde{P}_{n,\eps}$ & Projection to the eigenspace of $\lambda_{n,\eps}$\\
         % & \\
         % & \\
         % & \\
         % & \\
         $\Upsilon_\eps^{(\delta)}$& Enhanced stochastic objects from SPDEs\\
         $\fF_{\eps}^{(\delta)}$& Flux related stochastic objects from SPDEs\\
         % & \\
         % & \\
         % & \\
         % & \\
         \hline
    %\caption{Caption}
\end{longtable}

\textsc{Beijing International Center for Mathematical Research, Peking University, 5 Yiheyuan Road, Haidian District, Beijing, 100871, China}. 
\\
Email: yilin\underline{\phantom{1}}chen@pku.edu.cn
\\
\\
\textsc{Institute for Theoretical Sciences, Westlake University, 600 Dunyu Road, Xihu District, Hangzhou, Zhejiang, 310030, China}. 
\\
Email: xuweijun@westlake.edu.cn
\end{document}